\documentclass[11pt]{article}
\usepackage[margin=1in]{geometry}
\usepackage{comment,slashed,tensor}
\usepackage[T1]{fontenc}
\usepackage{ stmaryrd }
\usepackage{cancel}
\usepackage[symbol]{footmisc}
\usepackage[title]{appendix}

\usepackage[bookmarksnumbered, bookmarksopen, colorlinks, citecolor=red, linkcolor=black]{hyperref}
 \usepackage{amsfonts}
 \usepackage{comment}
  \usepackage{cancel}
\usepackage{graphicx}
 \usepackage{amsmath, amssymb}
 \usepackage{amsthm}
 \usepackage{mathtools}

 \numberwithin{equation}{section}
\newtheorem{thm}{Theorem}[section]
\newtheorem{lem}[thm]{Lemma}
\newtheorem{cor}[thm]{Corollary}

\newtheorem{rmk}[thm]{Remark}
\newtheorem{defn}[thm]{Definition}

\newtheorem{nota}[thm]{Notation}
\newtheorem{conv}[thm]{Convention}

\allowdisplaybreaks

\newcommand{\R}{\mathbb{R}}
\newcommand{\N}{\mathbb{N}}

\newcommand{\calH}{\mathcal{H}}
\newcommand{\calF}{\mathcal{F}}
\newcommand{\calM}{\mathcal{M}}

\newcommand{\calC}{\mathcal{C}}

\newcommand{\calL}{\mathcal{L}}
\newcommand{\calN}{\mathcal{N}}

\newcommand{\calR}{\mathcal{R}}
\newcommand{\calU}{\mathcal{U}}

\newcommand{\mfg}{\mathfrak{g}}

\newcommand{\mfC}{\mathfrak{C}}
\newcommand{\mfR}{\mathfrak{R}}
\newcommand{\mfA}{\mathfrak{A}}

\newcommand{\mfF}{\mathfrak{F}}

\newcommand{\mfc}{\mathfrak{c}}
\newcommand{\mfa}{\mathfrak{a}}
\newcommand{\mfL}{\mathfrak{L}}
\newcommand{\mfN}{\mathfrak{N}}
\newcommand{\mfG}{\mathfrak{G}}
\newcommand{\mfM}{\mathfrak{M}}

\makeatletter
\newcommand*{\rom}[1]{\expandafter\@slowromancap\romannumeral #1@}
\makeatother

\DeclareMathOperator{\tr}{tr}
\DeclareMathOperator{\supp}{supp}

\DeclareMathOperator{\Id}{Id}
\DeclareMathOperator{\adj}{adj}
\DeclareMathOperator{\ord}{ord}

\newcommand{\p}{\partial}

\newcommand{\dist}{\operatorname{dist}}

\newcommand{\abs}[1]{\left\lvert#1\right\rvert}
\newcommand{\norm}[1]{\left\lVert#1\right\rVert}
\newcommand{\const}[1]{\calC\left(#1\right)}

\newcommand{\sptdomain}{\mathcal{D}_T}

\usepackage{tensor}

\title{The Israel--Stewart equations with a physical vacuum boundary}
\author{Runzhang Zhong\thanks{Department of Mathematics, University of Missouri, Columbia, MO, 65211, United States. Email: \url{rz47r@missouri.edu}}}

\begin{document}

\maketitle
\begin{abstract}
    We establish local existence for the nonlinear Israel--Stewart equations describing a relativistic viscous fluid with a physical vacuum free boundary. Compared with ideal fluids with a physical vacuum, the relaxation equation for the bulk viscosity introduces further derivative loss. The problem is compounded by the degenerate nature of the equations of motion on the free boundary. In this work, we construct solutions in weighted Sobolev spaces via a Nash--Moser iteration. Since the domain moves, the iteration is carried out in Lagrangian coordinates, while the linearized operator at each iterate is inverted on the Eulerian domain with the local well-posedness of the linearized equations established in our earlier work.
\end{abstract}
\tableofcontents
\section{Introduction}
Relativistic viscous hydrodynamics describes dissipative fluids in regimes where relativistic effects cannot be neglected. Its main applications are the quark--gluon plasma produced in heavy-ion collisions and the matter in neutron star mergers, where a more accurate description requires the inclusion of viscosity, diffusion and heat conduction. The standard models are of Israel--Stewart type, in which the viscous contributions are not given as functions of the hydrodynamic variables but are new dynamical variables in their own right. For the relevant physical background and for the mathematical results, we refer to the survey \cite{disconzi2024recent}, the introduction of \cite{Zhong_2025}, and the references therein.

In this paper, we consider the Israel--Stewart model for relativistic viscous fluids in a Minkowski background $(\R^{1+3},g)$. We consider the case where dissipation arises only from bulk viscosity, and neglect shear viscosity and particle diffusion.

The fluid state is represented by the energy density $\rho\geq0$, the $4$-velocity $u$, which is assumed to be a forward time-like vector field, normalized by
\begin{equation}\label{eq:normalization of four velocity}
	u^\mu u_\mu = -1,
\end{equation}
and the bulk viscosity $\Pi$.

The evolution of the fluid is governed by the system
\begin{equation}\label{eq:general equation of motion}
\left\{
    \begin{aligned}
        &\partial_\mu T^{\mu\nu} = 0,\\
        &\partial_\mu J^\mu = 0,\\
        &\tau_\Pi u^\mu \partial_\mu \Pi + \Pi + \lambda \Pi^2 + \zeta\,\partial_\mu u^\mu = 0,
    \end{aligned}
    \right.
\end{equation}
where $T$ is the energy-momentum tensor for a viscous fluid, defined by
\begin{equation*}
	T_{\mu\nu}=\rho u_\mu u_\nu+(p+\Pi)\Delta_{\mu\nu},
\end{equation*}
and $J$ is the rest-mass current, defined by
\begin{equation*}
	J^\mu = n_b u^\mu.
\end{equation*}
Above, $n_b$ is the baryon density, the $2$-tensor $\Delta_{\mu\nu}:= g_{\mu\nu} + u_\mu u_\nu$ is the projection tensor onto the subspace orthogonal to $u$, and $p$ is the pressure of the fluid subject to an equation of state
\begin{equation*}
	p=p(\rho,n_b).
\end{equation*}
The three parameters in the evolution of the bulk viscosity $\Pi$ are the relaxation time $\tau_\Pi=\tau_\Pi(\rho)$, the transport coefficient $\lambda=\lambda(\rho)$, and the bulk viscosity coefficient $\zeta=\zeta(\rho)$.

Throughout the paper we use rectangular space-time coordinates $\{x^\alpha\}_{\alpha=0}^{3}$ on the Minkowski spacetime $\mathbb{R}^{1+3}$, identifying the coordinate $x^0=t$ as time and $x=(x^1,x^2,x^3)$ as the spatial components. The derivative $\p_\mu$ denotes the coordinate vector field $\frac{\p}{\p x^\mu}$. We adopt the Einstein summation convention. Greek indices range from $0$ to $3$, while Latin indices range from $1$ to $3$. Indices are raised and lowered using the Minkowski metric $g$.

\subsection{Physical vacuum boundary condition and Eulerian formulation}
It was shown in \cite{bemfica2019causality} that the system \eqref{eq:general equation of motion} is symmetric hyperbolic provided the causality condition holds:
\begin{equation}\label{eq:causality}
    \left[\frac{\zeta}{\tau_\Pi} + n_b\left(\frac{\partial p}{\partial n_b}\right)_{\rho}\right]\frac{1}{\rho + p + \Pi}\le1 - \left(\frac{\partial p}{\partial \rho}\right)_{n_b}.
\end{equation}
Here, we assume $\rho + p + \Pi > 0$ and $\tau_\Pi>0$. Under these conditions, the standard theory of first-order symmetric hyperbolic systems yields local well-posedness. A more complete analysis of nonlinear causality and strong hyperbolicity for a general class of Israel--Stewart-type theories, including shear viscosity and coupling to gravity, has recently been given in \cite{cordeiro2026nonlinearcausalitystronghyperbolicity}. However, none of these results addresses situations relevant for gaseous stars, in which a body of gas is surrounded by vacuum, whose boundary moves with the fluid. As we will see, the standard theory does not apply in this case due to the degeneracy of the equations of motion caused by the vacuum.

The gas occupies the domain
\begin{equation*}
    \Omega_t :=\{x\in\R^3:\rho(t,x)>0\},
\end{equation*}
at time $t$, whose boundary $\Gamma_t:=\partial\Omega_t$ is the vacuum boundary advected by the fluid's velocity $u$.

A distinguishing character of a gas, as opposed to a liquid, is that the density $\rho$ vanishes on the free boundary $\Gamma_t$,
\begin{equation*}
    \rho|_{\Gamma_t} = 0, \quad \text{for all $t>0$.}
\end{equation*}
The physical interpretation of this condition is that no matter is present on the boundary. Due to the absence of matter, it is natural to impose that the bulk viscosity $\Pi$ also vanishes on $\Gamma_t$,
\begin{equation*}
    \Pi|_{\Gamma_t} = 0,\quad \text{ for all $t>0$.}
\end{equation*}

We work with the power law equation of state,
\begin{equation}\label{eq:power law}
    p(\rho)=\rho^{\kappa+1}, \quad \kappa>0,
\end{equation}
under which the baryon density $n_b$ decouples from the system
\eqref{eq:general equation of motion}. Contracting
\eqref{eq:general equation of motion} with $u_\nu$ and with
$\Delta_{\nu\alpha}$ yields the energy and momentum equations,
respectively, and the full system can be reformulated as the following
\textit{initial boundary value problem}:
\begin{subequations}\label{eq:pre Eulerian formulation}
    \begin{align}
        u^\mu \partial_\mu \rho + (\rho^{\kappa+1} + \rho + \Pi)\partial_\mu u^\mu &= 0 && \text{in } \sptdomain, \\
        (\rho^{\kappa+1}+\rho+\Pi)u^\mu \partial_\mu u^\alpha + (\kappa + 1)\rho^{\kappa} \Delta^{\alpha\mu} \partial_\mu \rho + \Delta^{\alpha\mu}\partial_\mu \Pi &= 0 && \text{in } \sptdomain, \\
        \tau_\Pi u^\mu \partial_\mu \Pi + \Pi + \lambda \Pi^2 + \zeta\partial_\mu u^\mu &= 0 && \text{in } \sptdomain,\\
        u^\mu u_\mu &= -1 && \text{in } \sptdomain, \\
        \rho &= 0 && \text{on } \p\sptdomain, \\
        \Pi &=0 && \text{on } \p\sptdomain,\\
        u^\mu \partial_\mu|_{\Gamma} &\in T\Gamma,
        \\
        \sptdomain &= \bigcup_{0\leq t\leq T}\{t\}\times\Omega_t,\\
        \Gamma &= \bigcup_{0\leq t\leq T}\{t\}\times\Gamma_t,\\
        \rho(0,\cdot)=\rho_0(\cdot),\qquad u(0,\cdot)=u_0(\cdot),\qquad\Pi(0,\cdot)&=\Pi_0(\cdot).
\end{align}
\end{subequations}
As noted above, the system \eqref{eq:pre Eulerian formulation} is degenerate in the presence of vacuum.

\begin{rmk}
    Technically, the boundary of $\sptdomain$ consists of the ``cap'' $\Omega_T$ and the ``bottom'' $\Omega_0$, together with the lateral boundary
    \[ 
        \bigcup_{0<t<T}\{t\}\times\partial\Omega_t.
    \]
    Above, with a slight abuse of notation, we write $\p\sptdomain$ to denote only the lateral part.
\end{rmk}

It turns out that the decay rates of $\rho$ and $\Pi$ near the free boundary are important to ensure that the vacuum boundary evolves with a finite non-zero acceleration. We determine the admissible decay rates of the density $\rho$ and the bulk viscosity $\Pi$ near the free boundary $\Gamma_t$. By expressing the fluid $4$-acceleration in terms of $\rho$ and $\Pi$ and analyzing its behavior as the distance to $\Gamma_t$ tends to zero, we obtain the following decay rates that guarantee a finite and nonzero boundary acceleration (see \cite[Section~1.2]{Zhong_2025}):
\begin{equation}
    \rho\sim \mathrm{dist}(x,\Gamma_t)^\frac{1}{\kappa}, \quad \Pi\sim \mathrm{dist}(x,\Gamma_t)^{\frac{1}{\kappa}+2}.
\end{equation}
Here and throughout, for nonnegative functions $f$ and $g$, we write $f\sim g$ if there is a constant $C\geq 1$ such that $C^{-1}g\leq f\leq Cg$. These decay rates are referred to as the \textit{physical vacuum boundary condition}. In particular, we point out that, if $\Pi = 0$, this condition reduces to the physical vacuum condition for \textit{relativistic ideal fluids} in \cite{disconzi2022relativistic}.

To track the boundary behavior, it is convenient to introduce variables whose leading-order behavior is comparable to the distance function $\dist(\cdot,\Gamma_t)$. For this purpose, we define new variables
\begin{equation*}
    r:=\rho^\kappa,\quad \pi:=\Pi^\frac{\kappa}{2\kappa+1},
\end{equation*}
so that both $r$ and $\pi$ decay linearly near the free boundary. We can show that these decay rates are preserved by the flow.

We impose assumptions on the relaxation time $\tau_\Pi$, the bulk viscosity coefficient $\zeta$ and the transport coefficient $\lambda$. We assume the relaxation time $\tau_\Pi$ is a constant. For convenience, we take $\tau_\Pi=\frac{\kappa}{2\kappa+1}$. Define the functions $\hat{\zeta}$ and $\hat{\lambda}$
\[
\hat{\zeta}(r):=\frac{\zeta(r^{1/\kappa})}{r^{2+1/\kappa}},\qquad
\hat{\lambda}(r):=\lambda(r^{1/\kappa})\,r^{2+1/\kappa}.
\]
As in \cite{Zhong_2025}, we assume that $\hat{\zeta}$ and $\hat{\lambda}$ extend smoothly to $r=0$ and that $\hat\zeta$ is positive up to $r=0$:
\begin{equation}\label{eq:coefficient-assumptions}
    \hat{\zeta},\,\hat{\lambda}\in C^\infty([0,\infty)),\qquad \hat{\zeta}>0\quad\text{on }[0,\infty).
\end{equation}
In particular, $\hat\zeta(0)>0$ implies that $\zeta(r^{1/\kappa})\sim r^{2+1/\kappa}$ near the free boundary. Away from $r=0$, \eqref{eq:coefficient-assumptions} only requires $\zeta$ and $\lambda$ to be smooth and $\zeta$ to be positive. Since $r$ is bounded, $\hat\zeta$ is bounded from below by a positive constant on the range of $r$, and $\hat\zeta$, $\hat\lambda$ and their derivatives are bounded there. These bounds depend only on $\norm{r}_{L^\infty}$.

The smoothness of $\hat\zeta$ and $\hat\lambda$ up to $r=0$ prevents factors of $1/r$, which are singular at the boundary, from appearing when the coefficients are differentiated in the higher-order estimates. The lower bound on $\hat\zeta$ is used in the energy estimates of \cite{Zhong_2025}, where the $\tilde\pi$-equation is multiplied by a factor $\hat\zeta^{-1}$ to symmetrize the linearized system, and it enters the present paper through Theorem~\ref{thm:Eulerian-Linearized-Wellposedness}.

Rewriting the system \eqref{eq:pre Eulerian formulation} in the new variables $(r,u,\pi)$, we obtain the system in Eulerian coordinates:
\begin{subequations}\label{eq:Eulerian-Main-System}
    \begin{align}
        u^\mu \partial_\mu r + \kappa ra_1\p_\mu u^\mu &= 0 && \text{in } \sptdomain,\label{eq:r-equation} \\
        u^\mu \p_\mu u^0 + \left(1+\frac{1}{\kappa}\right)a_2\Delta^{0\mu}\partial_\mu r + \left(2+\frac{1}{\kappa}\right)ra_3\Delta^{0\mu}\partial_\mu\pi &=0 && \text{in }\sptdomain,\label{eq:u-zero-equation}\\
        u^\mu \p_\mu u^i + \left(1+\frac{1}{\kappa}\right)a_2\Delta^{i\mu}\partial_\mu r + \left(2+\frac{1}{\kappa}\right)ra_3\Delta^{i\mu}\partial_\mu\pi &=0 && \text{in }\sptdomain,\label{eq:vec-u-equation}\\
        u^\mu\partial_\mu\pi + \pi + r\hat{\lambda}a_0^{3+\frac{1}{\kappa}} + ra_4\partial_\mu u^\mu &= 0 && \text{in } \sptdomain,\label{eq:pi-equation}\\
        u^\mu u_\mu &= -1 && \text{in } \sptdomain, \label{eq:Normalization-Condition}\\
        r &= 0 && \text{on } \Gamma_t, \\
        \pi &=0 && \text{on } \Gamma_t,\\
        r,\pi &\sim \dist(\cdot,\Gamma_t) && \text{near } \Gamma_t, \\
        u^\mu \partial_\mu|_{\Gamma} &\in T\Gamma,
        \\
        \sptdomain &= \bigcup_{0\leq t \leq T}\{t\}\times\Omega_t,\\
        \Gamma &= \bigcup_{0\leq t\leq T}\{t\}\times\Gamma_t,\\
        r(0,\cdot)=r_0(\cdot),\qquad u(0,\cdot)&=u_0(\cdot),\qquad \pi(0,\cdot)=\pi_0(\cdot),
\end{align}
\end{subequations}
where the coefficients $a_0$, $a_1$, $a_2$, $a_3$, $a_4$ are given by:
\begin{equation}\label{eq:Coefficients}
    a_0 = \frac{\pi}{r}, \quad a_1 = r + 1 + r^2a_0^{2+\frac{1}{\kappa}},\quad a_2 = \frac{1}{a_1}, \quad a_3 = a_0^{1+\frac{1}{\kappa}}a_2,\quad a_4 = \frac{\hat{\zeta}}{a_0^{1+\frac{1}{\kappa}}}.
\end{equation}
Note that the normalization condition \eqref{eq:Normalization-Condition} can be shown to be preserved by the flow. Here, the equation \eqref{eq:u-zero-equation} for $u^0$ is written separately from the spatial part on purpose. At first sight, the system \eqref{eq:Eulerian-Main-System} appears overdetermined: the unknowns $(r,u,\pi)$ consist of six scalar functions, while the system contains seven equations including the normalization condition $u^\mu u_\mu =-1$. However, the system \eqref{eq:Eulerian-Main-System} is not inconsistent in that if one solves the system \eqref{eq:Eulerian-Main-System} excluding the equation \eqref{eq:u-zero-equation}, then $u^0$ can be defined through the normalization with the spatial part $\vec u$. A direct computation shows that $u^0$ solves the equation \eqref{eq:u-zero-equation}.

\subsection{Lagrangian formulation and main result}
We pass to Lagrangian coordinates in order to transform the free-boundary problem into one posed on a fixed domain. We denote by $y=(y_1,y_2,y_3)$ the Lagrangian coordinates on the fixed domain $\Omega_0$. Let $\eta(t,y):[0,T]\times\Omega_0\to\R^3$ be the flow map associated with the $3$-dimensional velocity field $\frac{\vec{u}}{u^0}$, defined by
\begin{equation}\label{eq:Flow-Map-Equation}
    \left\{
        \begin{aligned}
            \partial_t \eta(t,y)&=\frac{\vec{u}}{u^0}(t,\eta(t,y)),\\
            \eta(0,y)&=y,
        \end{aligned}
    \right.
\end{equation}
where $u^0$ is defined through \eqref{eq:Normalization-Condition} in terms of $\vec u$.

\begin{rmk}\label{rmk:three-velocity-flow}
The flow map is defined through $\vec u/u^0$ rather than the $4$-velocity $u$. Both have the fluid worldlines as integral curves, but \eqref{eq:Flow-Map-Equation} parametrizes them by the coordinate time $t$, so that Lagrangian and Eulerian time coincide and $\eta(t,\cdot)$ maps $\Omega_0$ onto the Eulerian domain $\Omega_t$, on which the weight $r(t,\cdot)$ and hence the weighted norms used in Sections~\ref{sec:linearized-eqs}--\ref{sec:tame-estimate-linearization-Eulerian} are defined. Parametrizing by proper time would make the level sets of the Lagrangian time solution-dependent hypersurfaces, and the slice-by-slice Eulerian estimates of Sections~\ref{sec:linearized-eqs}--\ref{sec:tame-estimate-linearization-Eulerian} could not be transferred to the Lagrangian domain. 
\end{rmk}

We further introduce the Lagrangian variables corresponding to the Eulerian variables $(r, u,\pi)$:
\begin{align*}
	R(t,y)&:=r(t,\eta(t,y)),\\
	\vec{U}(t,y)&:=\vec{u}(t,\eta(t,y)),\\
	U^0(t,y)&:= u^0(t,\eta(t,y)),\\
	P(t,y) &:=\pi(t,\eta(t,y)),\\
	F(t,y) &:= (\nabla_y\eta)(t,y).
\end{align*}

By the chain rule, the Eulerian derivative $\partial_i = \nabla_i$ can be written in the Lagrangian frame as
\begin{equation}
	\nabla_i = (F^{-1})^j_i \frac{\p}{\p y_j}.
\end{equation}

We define the metric
\begin{equation*}
	H^{ij}:=H^{ij}(\vec u) = \delta^{ij}-\frac{u^iu^j}{(u^0)^2}.
\end{equation*}
Then, in terms of the Lagrangian variables $(R,\vec{U},P,\eta)$, the system \eqref{eq:Eulerian-Main-System} becomes
\begin{subequations}\label{eq:Lagrangian-formulation}
    \begin{align}
        &U^0\p_t R + \kappa RA_1\frac{U^i}{U^0}\p_tU_i+\kappa RA_1(F^{-1})_i^l H^{ij}(\vec{U})\p_l U_j = 0\qquad \text{in } [0,T]\times\Omega_0, \label{eq:R-equation}\\
        &U^0\p_t U_i + \left(1+\frac{1}{\kappa}\right)A_2(U_iU^0\p_t R+(F^{-1})_i^k\p_kR)\notag\\
        &\qquad\qquad+ \left(2+\frac{1}{\kappa}\right)RA_3(U_iU^0\p_t P+(F^{-1})_i^k\p_k P) =0 \qquad \text{in }[0,T]\times\Omega_0,\\      
        &U^0\p_t P + P + R\hat{\lambda}(R)A_0^{3+\frac{1}{\kappa}}  + RA_4\frac{U^i}{U^0}\p_t U_i + RA_4(F^{-1})^l_iH^{ij}(\vec{U})\p_l U_j = 0 \qquad \text{in } [0,T]\times\Omega_0,\label{eq:P-eqaution}\\
        &\p_t \eta - \frac{\vec{U}}{U^0} = 0 \qquad \text{in } [0,T]\times\Omega_0,\label{eq:Lagrangian-flow-map-equation}\\
        &U^0=\sqrt{1+\abs{\vec{U}}^2} \qquad \text{in } [0,T]\times\Omega_0, \label{eq:Lagrangian-normalization}\\
        &R = 0 \qquad \text{on } \Gamma_0, \\
        &P =0 \qquad \text{on } \Gamma_0,\\
        &R,P \sim \dist(\cdot,\Gamma_0) \qquad \text{near } \Gamma_0, \\
        &\eta(t,\Gamma_0) =\Gamma_t,
        \\       
        &R(0,\cdot)=r_0(\cdot),\quad \vec{U}(0,\cdot)=\vec{u}_0(\cdot),\quad P(0,\cdot)=\pi_0(\cdot),\quad\eta(0,y)=y,
    \end{align}
\end{subequations}
where the coefficients $A_0,A_1,A_2,A_3,A_4$ are 
\begin{equation}\label{eq:Lagrangian-Coeffcients}
	A_0 = \frac{P}{R}, \quad A_1 = R + 1 + R^2A_0^{2+\frac{1}{\kappa}},\quad A_2 = \frac{1}{A_1}, \quad A_3 = A_0^{1+\frac{1}{\kappa}}A_2,\quad A_4 = \frac{\hat{\zeta}(R)}{A_0^{1+\frac{1}{\kappa}}}.
\end{equation}
\begin{rmk}
	We denote by $\p_i := \frac{\partial}{\partial x^i}$ the Eulerian spatial derivatives. Above and throughout, by an abuse of notation, we will sometimes use the same symbol to denote spatial Lagrangian derivatives. We also use the same notation $\partial_t$ for both Eulerian and Lagrangian time derivatives. The intended meaning will be clear from the argument of the function being differentiated.
\end{rmk}

Our main result is the local existence of solutions to the system \eqref{eq:Lagrangian-formulation} for sufficiently regular data satisfying
the physical vacuum condition. The functional framework is introduced in Section~\ref{sec:function-space} and the precise statement is Theorem~\ref{thm:main}. Informally:

\begin{thm}[Main result, informal]
Let $(r_0,\vec u_0,\pi_0)$ be initial data on a bounded domain $\Omega_0$ satisfying the physical vacuum condition $r_0,\pi_0\sim\dist(\cdot,\Gamma_0)$, the causality condition \eqref{eq:causality}, and lying in a weighted Sobolev space of sufficiently high order. Then there exist $T>0$ and a solution $(R,\vec U,P,\eta)$ of \eqref{eq:Lagrangian-formulation} on $[0,T]\times\Omega_0$, in the same scale of weighted Sobolev spaces as the data but with a loss of derivatives.
\end{thm}

\subsection{Difficulty and strategy}\label{sec:difficulty-and-strategy}
Compared with previous works on ideal fluids with a physical vacuum boundary \cite{disconzi2022relativistic,ifrim2023compressible}, the new difficulty is that the bulk viscosity in the Israel--Stewart system is a dynamical variable, governed by a transport equation sourced by the divergence of the velocity. A direct derivative count and application of transport estimates reveals a loss of derivatives. This loss of regularity is compounded by the need to control the fluid variables in the weighted Sobolev spaces required to handle the degeneracy of the equations due to the presence of vacuum.

In \cite{Zhong_2025} we showed that the linearized equations of the Eulerian formulation are nevertheless well-posed in the weighted spaces, with estimates that lose derivatives with respect to the background at which one linearizes. Heuristically, this is the situation in which Nash--Moser iteration applies: it produces a solution of the nonlinear equation from a linear theory whose estimates lose derivatives with respect to the background. The key ingredients are tame estimates for the linearized equations, which we establish in Sections~\ref{sec:linearized-eqs}--\ref{sec:Lagrangian-tame}.

However, implementing the Nash--Moser iteration directly on the Eulerian domain presents a further difficulty: since the domain itself moves, each iterate, constructed by solving the linearized equations, is defined on a different domain, making it difficult to formulate convergence in a fixed function space. Moreover, the smoothing operator in the Nash--Moser iteration scheme is, in general, difficult to construct in a way that preserves the physical vacuum boundary condition, which is essential for solving the linearized equations. For these reasons, we instead carry out the Nash--Moser iteration over the Lagrangian domain, defined through the flow of the spatial velocity $\frac{\vec u}{u^0}$ (see Remark~\ref{rmk:three-velocity-flow}).

Implementing the Nash--Moser iteration over the Lagrangian domain presents its own difficulty: closing the estimates for the linearized equations of the Lagrangian formulation is hard because of the inverse deformation tensor $F^{-1}$. In other words, linearization and the change of coordinates between the Eulerian and Lagrangian formulations do not commute. In our setting, however, this non-commutativity is, in a suitable sense, of lower order: the Lagrangian iterates can be constructed through the local well-posedness theory for the linearized equations of Eulerian formulation, and the estimates over the Eulerian domain can be carried over to the Lagrangian domain. 

Finally, the smoothing operators required by the scheme must act on weighted Sobolev spaces with a degenerate weight. We obtain them by adapting the regularization kernels of \cite{ifrim2023compressible} to the scales of the scheme.

With the tame estimates and the smoothing operators in hand, the Nash--Moser scheme converges, and we obtain the first nonlinear local existence result for the free-boundary Israel--Stewart system.

\subsection{Related work}
For the history of the free-boundary relativistic Euler equations, including results with symmetry and for liquids, we refer to \cite[Section~1.6]{disconzi2022relativistic} and \cite{disconzi2024recent}. Here, we focus on the three lines of work that are most relevant to the present paper.

\emph{Physical vacuum boundaries.} For a gas surrounded by vacuum, the relevant boundary condition is the \emph{physical vacuum boundary} condition, which goes back to the work of Liu and Liu--Yang on vacuum states in compressible gas dynamics \cite{liu1996compressible,liu1997compressible,liu2000compressible}. Local well-posedness for the compressible Euler equations with a physical vacuum boundary was established independently by Jang and Masmoudi \cite{jang2009well,Jang_Masmoudi_2014} and by Coutand and Shkoller \cite{coutand2011well,coutand2012well} (see also \cite{coutand2010priori} for a priori
estimates), in Lagrangian coordinates, with energies built from the wave-like structure of the equations and weighted by powers of the density. Ifrim and Tataru \cite{ifrim2023compressible} then developed an Eulerian approach with a sharper functional framework, giving local well-posedness at lower regularity and a continuation criterion. This result was recently extended to the non-isentropic case by Liu and Luo
\cite{liu2025physicalvacuumproblemscompressible}. For the relativistic Euler equations, a priori estimates were obtained by Jang, LeFloch and Masmoudi \cite{Jang-LeFloch-Masmoudi-2016} and by Had\v{z}i\'{c}, Shkoller and Speck \cite{Hadzic-Shkoller-Speck-2019}, and local well-posedness by Disconzi, Ifrim and Tataru \cite{disconzi2022relativistic} along the lines of \cite{ifrim2023compressible}. The ideal gas equation of state was treated at the linearized level recently in \cite{luczak2024priori}. In our work, we take from \cite{ifrim2023compressible} the scale of weighted Sobolev spaces as our functional framework and the good regularization kernels to build the smoothing operators for Nash--Moser iteration (Appendix~\ref{app:smoothing}).

\emph{Relativistic viscous fluids.} For the Israel--Stewart equations with bulk viscosity coupled to Einstein's equations, causality and local well-posedness in Sobolev spaces were established in \cite{bemfica2019causality}. This was recently extended in \cite{cordeiro2026nonlinearcausalitystronghyperbolicity} to a general class of
Israel--Stewart-type theories including shear viscosity. We also refer to the survey \cite{disconzi2024recent} for a comprehensive account of the mathematical theory of relativistic viscous fluids. None of these results involves a free boundary. In \cite{Zhong_2025}, we took the first step towards the free boundary problem: we determined the physical vacuum condition appropriate for the Israel--Stewart equations with bulk viscosity and proved that the equations linearized about a solution satisfying it are well posed in weighted Sobolev spaces. The estimates obtained there lose derivatives with respect to the background solution, as explained in Section~\ref{sec:difficulty-and-strategy}, so a Picard iteration does not close. This prompts the use of Nash--Moser in the present paper.

\emph{Nash--Moser iteration.} A loss of derivatives in the linearized equations is the typical situation in which Nash--Moser iteration is used and the method has been applied to a number of free boundary problems in fluid dynamics. Among these works, the closest to our setting is Makino's treatment of the gaseous star with a physical vacuum boundary \cite{Makino-2017}. In our work, the abstract Nash--Moser theorem we use is the one formulated in \cite{secchi2016nash}.

\subsection{Outline of the paper}

\emph{Function spaces, Section~\ref{sec:function-space}.}
The system degenerates on the free boundary, and the natural spaces are the weighted Sobolev spaces $H^{s,\sigma}_w(\Omega)$ with weight $w$, which is a non-degenerate defining function of $\Omega$. We recall their basic properties: embeddings and interpolation inequalities from \cite{ifrim2023compressible}. We also introduce the space-time versions, which serve as the scales in the Nash--Moser iteration.

\emph{Nash--Moser iteration scheme, Section~\ref{sec:Nash-Moser-iteration-scheme}.}
We set up the Nash--Moser iteration. From the Lagrangian formulation, we define the iteration map $\mfF$, which only acts on the spatial Lagrangian velocity $\vec U$, and the remaining variables, including $U^0$, the flow map $\eta$ and $(R,P)$, are determined from $\vec U$ by solving the constraint and the transport equations they satisfy. Then, a solution of the Lagrangian formulation is given by $\vec U$ solving $\mfF(\vec U)=0$ with prescribed initial data $\vec U(0,\cdot)=\vec u_0$. Notice that the initial data $(r_0,\pi_0)$ of the remaining variables are built into the definition of $\mfF$ as the data of the transport equations that determine $(R,P)$ from $\vec U$. Next, we introduce smoothing operators on the space-time scales, which are built from the spatial ones of Appendix~\ref{app:smoothing}. Since the unknown $\vec U$ at $t=0$ carries the prescribed initial data, the iteration map $\mfF$ does not map between Banach spaces. To fix this, we construct from the initial data an approximate solution which solves the equation to high order at $t=0$, and reduce to a modified iteration map $\tilde{\mfF}$. Finally, we compute the linearization of $\tilde\mfF$ and show that, after moving to the Eulerian coordinates, it is the linearized Eulerian system coupled to a transport equation coming from the variation of the flow map.

\emph{Invertibility of the linearization over the Eulerian domain,
Section~\ref{sec:linearized-eqs}.}
We construct the right inverse of the linearized iteration map required by the Nash--Moser theorem. We recall from \cite{Zhong_2025} the local well-posedness of the Eulerian linearized equations with a physical vacuum boundary. The linearized iteration map over the Eulerian domain differs from the equations in \cite{Zhong_2025} by the transport equation for the variation of the flow map, and thus a similar energy method applies. The resulting estimates have derivative losses with respect to the background states, which is the reason for the Nash--Moser iteration, but at this stage the losses are not of tame type yet.

\emph{First-order tame estimate over the Eulerian domain,
Section~\ref{sec:tame-estimate-linearization-Eulerian}.}
The energy estimate of Section~\ref{sec:linearized-eqs} is refined into a tame estimate by commuting the modified linearized equations with the weighted elliptic operators of Section~\ref{subsec:elliptic-estimates}. The commutation produces a large number of terms. We introduce a bookkeeping scheme that records the powers of the weight, the number of derivatives falling on the linearized variables and the total number of derivatives for each term, and we show that every commutator term is of lower order in this scheme. The tame estimates then follow.

\emph{Tame estimates in Lagrangian coordinates, Section~\ref{sec:Lagrangian-tame}.}
The estimates of Sections~\ref{sec:linearized-eqs}
and~\ref{sec:tame-estimate-linearization-Eulerian} are transferred from the Eulerian domain, where they are proved, to the Lagrangian domain, where the iteration is defined. This requires tame bounds for the change of coordinates, which we establish in Section~\ref{subsec:change-of-coordinate}. We then obtain the first-order and second-order tame estimates for the linearized iteration map.

\emph{Proof of the main theorem, Section~\ref{sec:proof-main-theorem}.}
We verify the hypotheses of the abstract Nash--Moser theorem of \cite{secchi2016nash} are satisfied in our case. Following the convergence argument in \cite{secchi2016nash}, we obtain a solution to $\mfF(\vec U)=0$ with $\vec U(0,\cdot)=\vec u_0$, and hence a solution of the Lagrangian formulation.

\emph{Appendices.}
Appendix~\ref{sec:solving-time-derivative} shows the time derivatives of the modified linearized system can be expressed in terms of spatial derivatives, which is needed to control the space-time norms by the spatial ones. Appendix~\ref{app:smoothing} constructs the spatial smoothing operators on the weighted Sobolev scale, using the regularization kernels in \cite{ifrim2023compressible}.

\noindent\textbf{List of notations}
\begin{itemize}
	\item $\const{\cdot}$ denotes constants depending non-decreasingly on the arguments.
	\item $\calF(\cdot)$ denotes functions smooth in a neighborhood of the range of its arguments.
\end{itemize}

\section{Function spaces}\label{sec:function-space}
The system \eqref{eq:Eulerian-Main-System} is a first-order hyperbolic system that degenerates on the free boundary $\Gamma_t$. This degeneracy makes standard Sobolev spaces insufficient for our analysis. We therefore employ weighted Sobolev spaces adapted to the degeneracy. In this section, we introduce weighted Sobolev norm and the time-dependent norms that will be useful later in the estimates and the construction of solutions.

\subsection{Weighted Sobolev spaces}\label{sec:weighted-sobolev}

Let $\Omega\subset\mathbb{R}^3$ be a bounded domain and $w=w(x)$ be a non-degenerate defining function of $\Omega$, i.e. the function $w$ satisfies
\begin{equation*}
	w\geq 0, \qquad \Omega=\{x\in\R^3: w(x)>0\},
\end{equation*}
and $w$ is Lipschitz on $\overline\Omega$ such that $\abs{\nabla w}>0 \ \text{on } \partial\Omega$. In particular, near the boundary $w$ is comparable to the distance: there exists $0<c\leq C$ such that one has 
\begin{equation*}
	c\,\dist(x,\p\Omega)\leq w(x)\leq C\,\dist(x,\p\Omega).
\end{equation*}

\begin{defn}
Let $N\geq0$ be an integer, $\sigma>-1/2$ and $w$ be a defining function for $\Omega$.
\begin{itemize}
\item For a scalar distribution $f$ on $\Omega$, the weighted Sobolev
space $H^{N,\sigma}_w$ consists of those $f$
with
\begin{equation*}
	\|f\|_{H^{N,\sigma}_w}^2
:=\sum_{|\alpha|\le N}
\|w^\sigma \partial^\alpha f\|_{L^2}^2
<\infty.
\end{equation*}
\item For a finite array $X$ of variables (vectors, matrices, or higher-order
tensor valued), the weighted Sobolev norm $\norm{\cdot}_{H^{N,\sigma}_w}$ of $X$ is the sum of the weighted
Sobolev norms of its components.
\end{itemize}
\end{defn}

\begin{nota}
For tensorial objects $X_1,\ldots,X_k$, the product $X_1\cdots X_k$ denotes any contraction of them. Each of its components is a finite sum of pointwise products of components of the $X_i$, so when estimating norms we treat $X_1\cdots X_k$ as a single such product and sum over the finitely many components at the end, at the cost of a combinatorial constant absorbed into $\calC(\cdot)$. Likewise $\p^kf$ denotes the collection $\{\p^\alpha f\mid|\alpha|=k\}$, treated as a single component in estimates. In particular the order of the factors in a schematic product does not matter, and we write them in whichever order is convenient.
\end{nota}

The following properties of the weighted Sobolev spaces can be found in \cite{ifrim2023compressible,disconzi2022relativistic}. We have a Hardy-type inequality:
\begin{lem}\label{lem:WeightedEmbedding}
Let $N \ge 0$, $\sigma>-\frac{1}{2}$ and $k\geq0$.
Then,
\[
H^{N+k,\sigma+k}_w \subset H^{N,\sigma}_w.
\]
\end{lem}

As an immediate consequence, the weighted spaces embed into standard Sobolev spaces.

\begin{cor}
Let $\sigma > 0$ and $\sigma \le N$. Then
\[
H^{N,\sigma}_w \subset H^{N-\sigma}.
\]
\end{cor}

Combining this with standard Sobolev embeddings yields Morrey-type bounds.

\begin{lem}\label{lem:Morrey}
For $N \ge 0$ and $\sigma > -\tfrac12$,
\[
H^{N,\sigma}_w \subset C^s,
\qquad
0 \le s \le N - \sigma - \frac{3}{2},
\]
where equality can hold only if $s$ is not an integer.
\end{lem}

We next record an interpolation inequality between weighted Sobolev norms.

\begin{lem}\label{lem:Weighted-Interpolation}
Let $\sigma_0,\sigma_m\in\mathbb{R}$ and $1\le p_0,p_m\le\infty$.  
For $0<j<m$, define
\[
\theta_j = \frac{j}{m},\qquad
\frac1{p_j} = \frac{1-\theta_j}{p_0} + \frac{\theta_j}{p_m},\qquad
\sigma_j = (1-\theta_j)\sigma_0 + \theta_j\sigma_m.
\]
Assume
\[
m-\sigma_m - 3\!\left(\frac1{p_m}-\frac1{p_0}\right) > -\sigma_0,
\qquad
\sigma_j > -\frac1{p_j}.
\]
Then for $0<j<m$,
\[
\| w^{\sigma_j}\partial^j f \|_{L^{p_j}}
\lesssim
\| w^{\sigma_0} f \|_{L^{p_0}}^{1-\theta_j}
\| w^{\sigma_m}\partial^m f \|_{L^{p_m}}^{\theta_j}.
\]
\end{lem}

A particular case that will be used repeatedly is when $p_0=\infty$, $\sigma_0=0$ and $p_m=2$. We record it here for later use.

\begin{cor}\label{cor:L2Linfty-interpolation}
Let $\sigma_m>-\tfrac12$ and assume
\[
m-\sigma_m-\frac{3}{2}>0.
\]
For $0<j<m$, define
\[
\theta_j=\frac{j}{m},\qquad
\frac1{p_j}=\frac{\theta_j}{2},\qquad
\sigma_j=\theta_j\sigma_m.
\]
Then
\[
\| w^{\sigma_j}\partial^j f \|_{L^{p_j}}
\lesssim
\|f\|_{L^\infty}^{1-\theta_j}
\| w^{\sigma_m}\partial^m f \|_{L^2}^{\theta_j}.
\]
\end{cor}

Lemma~\ref{lem:Weighted-Interpolation} with $p_0=p_m=2$ and $\sigma_0=\sigma_m=\sigma$ implies that the weighted Sobolev norm is equivalent to the sum of the two end-point terms:
\begin{cor}\label{cor:weighted-L2-end-point-equivalence}
	Let $N\geq 0$ be an integer and $\sigma>-1/2$. Then,
	\begin{equation}
		\norm{f}_{H^{N,\sigma}_w}\lesssim \norm{w^\sigma f}_{L^2} + \norm{w^\sigma\p^Nf}_{L^2}.
	\end{equation}
\end{cor}

In the energy estimates, we often need to bound multilinear expressions of
the form $r^{\sigma}\prod_{i,j_i}\partial^{\alpha_{j_i}}f_i$ in $L^2$,
where $r$ solves \eqref{eq:r-equation} along $\vec{u}$.

\begin{cor}\label{Cor:Multilinear-Interpolation}
Let $\sigma>-\frac{1}{2}$ and $k$ be an integer such that
\[
    k-\sigma >\frac{3}{2}.
\]
Let $f_1,\ldots,f_n$ be sufficiently regular functions, $J_1,\ldots,J_n\geq 1$ and $\alpha_{j_i}\geq 1$ be integers such that
\[
    \sum_{i=1}^{n} \sum_{j_i=1}^{J_i} \alpha_{j_i} = k.
\]
Then,
\[
    \norm{w^\sigma \prod_{i=1}^{n} \prod_{j_i=1}^{J_i} \partial^{\alpha_{j_i}}f_i}_{L^2}\lesssim \calC\sum_{i=1}^n \prod_{j \neq i}\left(\norm{f_j}_{L^\infty}\right)\norm{w^\sigma\p^kf_i}_{L^2},
\]
where $\calC$ depends on $\norm{f_i}_{L^\infty}$ for $i=1,\dots,n$.
\end{cor}

\begin{proof}
Choose
\[
  \theta_{j_i} = \frac{\alpha_{j_i}}{k}, \qquad
  \Lambda_i = \sum_{j_i} \theta_{j_i}, \qquad
  \nu_{j_i} = \sigma\, \theta_{j_i}.
\]
By construction, $\sum_{i=1}^{n}\Lambda_i = \sum_{i,j_i}\theta_{j_i} = 1$.
Since $\alpha_{j_i} \geq 1$, we have $\theta_{j_i} > 0$, and since
$\sigma > -\tfrac{1}{2}$ and $\theta_{j_i} \leq 1$, we have
$\nu_{j_i} = \sigma\theta_{j_i} > -\frac12\theta_{j_i}$ and $\sum_{i,j_i}\nu_{j_i}=\sigma$.

We decompose the weight $w^\sigma$ in the following way
\[
  w^\sigma = \prod_{i=1}^{n}\prod_{j_i=1}^{J_i} w^{\nu_{j_i}},
  \qquad \nu_{j_i} = \sigma\,\theta_{j_i}
\]
and distribute to each factor.

Choose the H\"{o}lder exponents $p_{j_i} = 1/\theta_{j_i}$, so that $\sum_{i,j_i}1/p_{j_i} = 1$. Squaring the $L^2$ norm and applying H\"older's inequality,

\[
\norm{w^\sigma\prod_{i,j_i}\p^{\alpha_{j_i}}f_i}_{L^2}^2 = \int \prod_{i=1}^n\prod_{j_i=1}^{J_i}\abs{w^{\nu_{j_i}}\p^{\alpha_{j_i}}f_i}^2dx\leq \prod_{i=1}^n\prod_{j_i=1}^{J_i}\norm{w^{\nu_{j_i}}\p^{\alpha_{j_i}}f_i}_{L^{2p_{j_i}}}^2.
\]

Fix $(i,j_i)$ and apply Lemma~\ref{lem:Weighted-Interpolation} with
\[\quad m=k, \quad j=\alpha_{j_i},\quad p_0=\infty,\quad p_m =2,\quad \sigma_0 = 0,\quad\sigma_m = \sigma.\]
We also verify that the parameters satisfy
\begin{align*}
  m-\sigma_m - 3\!\left(\frac1{p_m}-\frac1{p_0}\right) > -\sigma_0 
  &\Longleftarrow k-\sigma > \tfrac{3}{2},\\
  \sigma_j > -\frac1{p_j}&\Longleftarrow \sigma > -\tfrac{1}{2}.
\end{align*}

Thus, we obtain
\[
    \norm{w^{\nu_{j_i}}\p^{\alpha_{j_i}}f_i}_{L^{2p_{j_i}}}\lesssim \norm{f_i}_{L^\infty}^{1-\theta_{j_i}}\norm{w^\sigma\p^kf_i}_{L^2}^{\theta_{j_i}}.
\]
Substituting and collecting exponents by index $i$,
\begin{equation}\label{eq:collected}
  \norm{w^{\sigma}
    \prod_{i,j_i}\p^{\alpha_{j_i}} f_i}_{L^2} \lesssim
  \prod_{i=1}^{n}
  \norm{f_i}_{L^{\infty}}^{J_i-\Lambda_i}
  \,\norm{w^{\sigma}\p^k f_i}_{L^{2}}^{\Lambda_i}.
\end{equation}

Writing $J_i - \Lambda_i = (J_i-1)+(1-\Lambda_i)$, we have
\[
  \prod_{i=1}^{n}
  \norm{f_i}_{L^{\infty}}^{J_i-\Lambda_i}
  \norm{w^{\sigma}\p^{k} f_i}_{L^{2}}^{\Lambda_i}
  = \underbrace{\prod_{i=1}^{n}\norm{f_i}_{L^{\infty}}^{J_i-1}}_{=\calC}
  \;\cdot\;
  \prod_{i=1}^n
  \norm{f_i}_{L^{\infty}}^{1-\Lambda_i}
  \norm{w^{\sigma}\p^k f_i}_{L^{2}}^{\Lambda_i}.
\]
Using $1-\Lambda_i = \sum_{j\neq i}\Lambda_j$ and setting
$a_i = \prod_{j\neq i}\norm{f_j}_{L^\infty}
\cdot\norm{w^{\sigma}\p^k f_i}_{L^2}$, the remaining product equals
$\prod_i a_i^{\Lambda_i}$.  Since $\Lambda_i > 0$ and
$\sum_i \Lambda_i = 1$, Young's inequality gives
\[
  \prod_{i=1}^{n} a_i^{\Lambda_i}
  \leq
  \sum_{i=1}^{n}\Lambda_i a_i
  \lesssim
  \sum_{i=1}^{n}
  \left(\prod_{j\neq i}\norm{f_j}_{L^{\infty}}\right)
  \norm{w^{\sigma}\p^{k}f_i}_{L^{2}}.
\]
Combining with $\calC$ yields 
\[
    \norm{w^\sigma \prod_{i=1}^{n} \prod_{j_i=1}^{J_i} \partial^{\alpha_{j_i}}f_i}_{L^2}\lesssim \calC\sum_{i=1}^n \prod_{j \neq i}\left(\norm{f_j}_{L^\infty}\right)\norm{w^\sigma\p^kf_i}_{L^2}.
\]

\end{proof}
\begin{lem}\label{lem:coefficient-product-estimate}
Let $f_1,\ldots,f_n$ be state variables such that
$\mathrm{Range}(f_1,\ldots,f_n)\subset K\subset\subset\R^n$. Let
$F=F(f_1,\ldots,f_n)$ be smooth in a neighborhood of the range of its arguments, and let $\tilde f$ be a
function defined on $\Omega$.
 
\begin{enumerate}
\item[(a)] For $\sigma>-1/2$ and any integer $k\geq 0$,
\begin{equation}\label{eq:coeff-product-Wkinfty}
    \norm{w^\sigma\p^k\big(F(f_1,\ldots,f_n)\,\tilde f\big)}_{L^2}
    \lesssim \const{\norm{f_1}_{W^{k,\infty}},\ldots,\norm{f_n}_{W^{k,\infty}}}
    \norm{\tilde f}_{H^{k,\sigma}_w}.
\end{equation}
\item[(b)] If further $k-\sigma>\frac32$, then
\begin{equation}\label{eq:coeff-product-weighted}
\begin{aligned}
    \norm{w^\sigma\p^k\big(F(f_1,\ldots,f_n)\,\tilde f\big)}_{L^2}
    &\lesssim \const{\norm{f_1}_{L^\infty},\ldots,\norm{f_n}_{L^\infty}}\\
    &\quad\times\Big(\norm{w^\sigma\p^k\tilde f}_{L^2}
    + \norm{\tilde f}_{L^\infty}\sum_{j=1}^n\norm{w^\sigma\p^kf_j}_{L^2}\Big).
\end{aligned}
\end{equation}
\end{enumerate}
\end{lem}
 
\begin{proof}
By the Leibniz rule, $\p^k\big(F(f_1,\ldots,f_n)\tilde f\big)$ is a finite
sum of terms multilinear in derivatives of $f_1,\ldots,f_n,\tilde f$, with
coefficients smooth in $f_1,\ldots,f_n$ (by the chain rule applied to $F$)
and total derivative count $k$. Since
$\mathrm{Range}(f_1,\ldots,f_n)\subset K\subset\subset\R^n$ and $F$ is smooth on
a neighborhood of $K$, every such coefficient is bounded on $K$ by a constant depending on $\norm{f_1}_{L^\infty},\ldots,\norm{f_n}_{L^\infty}$.

Bounding all factors involving $f_1,\ldots,f_n$ directly in $L^\infty$ (using
$\norm{f_i}_{W^{k,\infty}}$) and the remaining factor involving $\tilde f$ in $H^{k,\sigma}_w$, each term is bounded by
\[
    \const{\norm{f_1}_{W^{k,\infty}},\ldots,\norm{f_n}_{W^{k,\infty}}}
    \norm{w^\sigma\p^k\tilde f}_{L^2}
    \le \const{\norm{f_1}_{W^{k,\infty}},\ldots,\norm{f_n}_{W^{k,\infty}}}
    \norm{\tilde f}_{H^{k,\sigma}_w},
\]
which gives \eqref{eq:coeff-product-Wkinfty} upon summing over the finitely
many terms in the expansion.

If $k-\sigma>\frac{3}{2}$, applying Corollary~\ref{Cor:Multilinear-Interpolation} to each term gives
\begin{equation*}
	\const{\norm{f_1}_{L^\infty},\ldots,\norm{f_n}_{L^\infty}}
    \Big(\norm{w^\sigma\p^k\tilde f}_{L^2}
    + \norm{\tilde f}_{L^\infty}\sum_{j=1}^n\norm{w^\sigma\p^kf_j}_{L^2}\Big)
\end{equation*}
for each term, and summing over the finitely many terms in the expansion
yields \eqref{eq:coeff-product-weighted}.
\end{proof}

\subsection{Time-dependent norms}\label{subsec:time-dependent-norms}
For each $t\in [0,T]$, let $\Omega_t\subset\R^3$ be a bounded domain with a non-degenerate
defining function $w(t)=w(t,\cdot)$ as in Section~\ref{sec:weighted-sobolev}. Below, the norms $X$, $Y$ and $Z$ are norms on $\Omega_t$ at a fixed time. Their counterparts $\tilde X$, $\tilde Y$, $\tilde Z$ are obtained by taking the sup over $t\in[0,T]$.

\begin{defn}\label{defn:time-dependent-norms}
Let $l\ge0$ (with $l\geq 1$ for the $X$-norms), $k\ge1$ be integers and $\alpha>-\frac12$. Let $f$ be defined on
$\bigcup_{t\in[0,T]}\{t\}\times\Omega_t$ and let $w(t)$ be a defining function of $\Omega_t$ for each $t$. Set
\begin{align*}
    \norm{f}_{X^{l,\alpha}_w}&:=\norm{f}_{H^{2l,l+\alpha}_w}+\norm{\p_tf}_{H^{2l-2,l+\alpha-1}_w},
    & \norm{f}_{\tilde X^{l,\alpha}_w}&:=\sup_{t\in[0,T]}\norm{f(t)}_{X^{l,\alpha}_{w(t)}},\\
    \norm{f}_{Y^{l,\alpha}_w}&:=\norm{f}_{H^{2l,l+\alpha}_w},
    & \norm{f}_{\tilde Y^{l,\alpha}_w}&:=\sup_{t\in[0,T]}\norm{f(t)}_{Y^{l,\alpha}_{w(t)}},\\
    \norm{f}_{Z^k_x}&:=\norm{f}_{W^{k,\infty}(\Omega_t)}+\norm{\p_tf}_{W^{k-1,\infty}(\Omega_t)},
    & \norm{f}_{\tilde Z^k_x}&:=\sup_{t\in[0,T]}\norm{f(t)}_{Z^k_x}.
\end{align*}
For $F$ defined on $[0,T]\times\Omega_0$ and $w(t)$ a defining function of $\Omega_0$ for each $t$, the
norms $X^{l,\alpha}_w$, $\tilde X^{l,\alpha}_w$, $Y^{l,\alpha}_w$, $\tilde Y^{l,\alpha}_w$ are defined by
the same formulas, and
\begin{align*}
    \norm{F}_{Z^k_y}&:=\norm{F}_{W^{k,\infty}(\Omega_0)}+\norm{\p_tF}_{W^{k-1,\infty}(\Omega_0)},
    & \norm{F}_{\tilde Z^k_y}&:=\sup_{t\in[0,T]}\norm{F(t)}_{Z^k_y}.
\end{align*}
In both cases $\p_t$ is the time derivative at fixed spatial variable. We call $X,\tilde X$ the solution
norms, $Y,\tilde Y$ the source norms and $Z,\tilde Z$ the control norms.
\end{defn}

The subscripts $x$ and $y$ in the $Z$-norms indicate that the norms are taken over the Eulerian domain $\Omega_t$, whose coordinate is $x$, and over the Lagrangian domain $\Omega_0$, whose coordinate is $y$. The norms $X$, $Y$ need no such subscript, since the domain is determined by the weight. In this paper the only defining function of $\Omega_t$ we use is the Eulerian density $r$, and the only defining functions of $\Omega_0$ we use are the Lagrangian density $R$ and the initial density $r_0$ (see Lemma~\ref{lem:comparability}). When the subscript is omitted from $H^{N,\sigma}$, $X$ or $Y$, the corresponding weight is $r$.

By Lemma~\ref{lem:WeightedEmbedding} and Lemma~\ref{lem:Morrey}, we have the following relations among the above norms:
\begin{lem}\label{lem:time-dependent-norms-properties}
Let $l\ge1$, $k\ge1$ be integers and $\alpha>-\frac12$.
\begin{enumerate}
\item[(i)] $\norm{f}_{X^{l,\alpha}_w}\lesssim\norm{f}_{X^{l+1,\alpha}_w}$ and $\norm{f}_{Y^{l,\alpha}_w}\lesssim\norm{f}_{Y^{l+1,\alpha}_w}$.
\item[(ii)] $\norm{\p f}_{H^{2l-2,l+\alpha-1}_w}\lesssim\norm{f}_{H^{2l,l+\alpha}_w}$.
\item[(iii)] If $l-\alpha-\frac32>k$, then $\norm{f}_{\tilde Z^k}\lesssim\norm{f}_{\tilde X^{l,\alpha}_w}$.
\end{enumerate}
The implicit constants depend on $\Omega_t$ and $w(t)$ only through those of Lemmas~\ref{lem:WeightedEmbedding}
and Lemma~\ref{lem:Morrey}.
\end{lem}

The solution norms $X$ measure the iterates of the Nash--Moser scheme and the solutions of the linearized equations. The source norms $Y$ measure sources: the energy estimates of Section~\ref{sec:linearized-eqs} need only $L^1_t$ of the source, and $\norm{g}_{L^1_tH^{2l,l+\alpha}_w}\le T\norm{g}_{\tilde Y^{l,\alpha}_w}\le\norm{g}_{\tilde Y^{l,\alpha}_w}$ since $T\le1$. The control norms $Z$ are the Lipschitz-type norms measuring the background states on which the constants in the estimates depend. See, for example, Theorem~\ref{thm:modified-linearized-with-xi-local-wellposedness}.

\begin{nota}\label{nota:l0}
Let $l_0$ be the smallest integer greater than $\frac1{2\kappa}+2$. By Lemma~\ref{lem:Morrey},
\[
    \norm{f}_{L^\infty}\lesssim\norm{f}_{H^{2l_0,\,l_0+\frac1{2\kappa}-\frac12}_w},\qquad
    \norm{f}_{\tilde Z^1}\lesssim\norm{f}_{\tilde X^{l_0+1,\frac1{2\kappa}}_w}.
\]
\end{nota}

\section{Nash--Moser iteration scheme}\label{sec:Nash-Moser-iteration-scheme}
In this section, we define a Nash--Moser iteration scheme for solving the nonlinear Israel--Stewart system \eqref{eq:Lagrangian-formulation}. 

Given a sufficiently regular vector field $\vec U$ defined on $[0,T]\times\Omega_0$, where $T$ will be chosen later on depending on the initial data, we define $U^0$ and $\eta$ through \eqref{eq:Lagrangian-normalization} and \eqref{eq:Lagrangian-flow-map-equation}, respectively, i.e.,
\begin{align}
	U^0&=\sqrt{1+\abs{\vec U}^2},\label{eq:U^0-defining-eq}\\
	\eta(t,\cdot) &= \Id_{\Omega_0} + \int_0^t \left(\frac{\vec U}{U^0}\right)(s,\cdot) ds.\label{eq:eta-defining-eq}
\end{align}

Next, we define the variables $(R,A_0)$ in terms of $\vec U$ as the solutions to the following sub-system:
\begin{equation}\label{eq:R-A0-defining}
	\begin{aligned}
		&U^0\p_t R + \kappa RA_1\frac{U^i}{U^0}\p_tU_i+\kappa RA_1(F^{-1})^{l}_i H^{ij}(\vec{U})\p_l U_j = 0\qquad \text{in } [0,T]\times\Omega_0\\
		&U^0\p_t A_0 + A_0 + \hat\lambda(R)A_0^{3+\frac{1}{\kappa}}\\
		&\qquad + (A_4-\kappa A_0A_1)\frac{U^i}{U^0}\p_t U_i + (A_4-\kappa A_0A_1)(F^{-1})^{l}_i H^{ij}(\vec{U})\p_l U_j=0 \qquad \text{in } [0,T]\times\Omega_0\\
		&R(0,\cdot)=r_0(\cdot),\qquad A_0(0,\cdot)=\left(\frac{\pi_0}{r_0}\right)(\cdot).
	\end{aligned}
\end{equation}
For fixed $\vec U$, this can be treated as a first-order ODE system and therefore admits a unique solution.

Setting $P:=RA_0$, a direct computation shows that $(R,P)$ solve the following sub-system formed by \eqref{eq:R-equation} and \eqref{eq:P-eqaution} with their corresponding initial data:
\begin{equation}\label{eq:R-P-defining}
    \begin{aligned}
        &U^0\p_t R + \kappa RA_1\frac{U^i}{U^0}\p_tU_i+\kappa RA_1(F^{-1})^{l}_i H^{ij}(\vec{U})\p_l U_j = 0\qquad \text{in } [0,T]\times\Omega_0, \\
        &U^0\p_t P + P + R\hat{\lambda}(R)A_0^{3+\frac{1}{\kappa}}  + RA_4\frac{U^i}{U^0}\p_t U_i + RA_4(F^{-1})^{l}_iH^{ij}(\vec{U})\p_l U_j = 0 \qquad \text{in } [0,T]\times\Omega_0,\\
        &R(0,\cdot)=r_0(\cdot),\quad P(0,\cdot)=\pi_0(\cdot).
    \end{aligned}
\end{equation}

\begin{rmk}
We define $(R,A_0)$ first and set $P=RA_0$, rather than defining $(R,P)$ directly by \eqref{eq:R-P-defining}, because \eqref{eq:R-P-defining} involves the ratio $P/R$, which is not defined on $\p\Omega_0$, where both $R$ and $P$ vanish. In contrast, the system \eqref{eq:R-A0-defining} for $(R,A_0)$ is non-degenerate up to the boundary since its data $\frac{\pi_0}{r_0}$ is well defined up to the boundary by Definition~\ref{defn:admissible-data}.
\end{rmk}

Now, we define the iteration map $\mfF=\mfF(\vec U)$ acting on the Lagrangian velocity field $\vec U$. Given $\vec U$, let $(R,U^0,P,\eta)$ be as defined above. We define
\begin{equation}\label{eq:iteration-map}
    \left(\mfF(\vec U)\right)_i:=U^0\p_t U_i + \left(1+\frac{1}{\kappa}\right)A_2(U_iU^0\p_t R+(F^{-1})^k_i\p_kR)
        + \left(2+\frac{1}{\kappa}\right)RA_3(U_iU^0\p_t P+(F^{-1})^k_i\p_k P).
\end{equation}

Then, a solution to the nonlinear Israel--Stewart system \eqref{eq:Lagrangian-formulation} is given by solving
\begin{equation}\label{eq:iteration-map-equation}
    \left\{\begin{aligned}
        \mfF(\vec U) &= 0\\
        \vec U(0,\cdot) &= \vec u_0.
    \end{aligned}
    \right.
\end{equation}
Note that the initial data of $R,P$ and $\eta$ do not appear in \eqref{eq:iteration-map-equation} because they are already prescribed in \eqref{eq:R-A0-defining} and \eqref{eq:eta-defining-eq}, which are part of the definition of $\mfF$.

\begin{rmk}\label{rmk:r-u-pi-defining}
    If the initial data $(r_0,\pi_0)$ satisfy the physical vacuum boundary conditions, then the functions $(R,P)$ defined through $\vec U$ preserve the physical vacuum boundary conditions for as long as the flow exists, which are essential for inverting the linearized operator of the iteration map using the local well-posedness of the linearized equations of the Eulerian formulation. Iterating only on the variable $\vec U$, rather than all fluid variables $(R,\vec U,P)$, allows us to avoid losing these conditions in the smoothing process of Nash--Moser iteration.

    Moreover, setting
    \[r=R\circ\eta^{-1},\qquad u^0 = U^0\circ\eta^{-1}, \qquad \vec u = \vec U\circ\eta^{-1},\qquad \pi = P\circ\eta^{-1},\] 
    and defining the moving domain $\Omega_t=\Omega_t(\vec U)$ at time $t\in[0,T]$ as
\begin{equation*}
	\Omega_t:=\eta(t,\Omega_0),
\end{equation*}
    one can verify that $(r,\pi)$ solve the corresponding Eulerian part
    \begin{equation}\label{eq:r-pi-defining}
    \begin{aligned}
        &u^\mu \partial_\mu r + \kappa ra_1\p_\mu u^\mu = 0 \qquad \text{in } \sptdomain, \\
        &u^\mu\partial_\mu\pi + \pi + r\hat{\lambda}a_0^{3+\frac{1}{\kappa}} + ra_4\partial_\mu u^\mu = 0  \qquad \text{in } \sptdomain,\\
        &r(0,\cdot)=r_0(\cdot),\qquad \pi(0,\cdot)=\pi_0(\cdot),
    \end{aligned}
\end{equation}
and thus $(r,\pi)$ satisfy the physical vacuum boundary condition, and $r=r(t,\cdot)$ is the defining function for $\Omega_t$.
\end{rmk}

\subsection{Smoothing operator}
The estimates are obtained in scales of weighted Sobolev spaces. Consequently, a smoothing operator adapted to weighted Sobolev spaces is required, which we construct in this subsection.

\begin{lem}[Smoothing operators]\label{lem:spatial-smoothing-operator}
Let $m_*\ge1$ be an integer and let $S_\theta$, $\theta\ge1$, be the operators of
Lemma~\ref{lem:spatial-smoothing-operator-appendix} for $(\Omega,w)=(\Omega_0,r_0)$, acting on space-time functions at each fixed time, $(S_\theta u)(t,\cdot):=S_\theta(u(t,\cdot))$, and componentwise on vector-valued functions. Then, the family of operators
\begin{equation*}
    S_\theta: \tilde Y^{0,\frac{1}{2\kappa}}_{r_0}\to\bigcap_{l\geq 0}\tilde Y^{l,\frac{1}{2\kappa}}_{r_0},\qquad \theta\geq 1,
\end{equation*}
have the following estimates: for all integers $0\leq\alpha,\beta\leq m_*$ and $u\in \tilde Y^{\alpha,\frac{1}{2\kappa}}_{r_0}$,
\begin{align}
    \norm{S_\theta u}_{\tilde Y^{\beta,\frac1{2\kappa}}_{r_0}}&\le C\theta^{(\beta-\alpha)_+}\norm{u}_{\tilde Y^{\alpha,\frac1{2\kappa}}_{r_0}},\label{eq:spacetime-smoothing-S1}\\
    \norm{S_\theta u-u}_{\tilde Y^{\beta,\frac1{2\kappa}}_{r_0}}&\le C\theta^{\beta-\alpha}\norm{u}_{\tilde Y^{\alpha,\frac1{2\kappa}}_{r_0}},\qquad\beta\le\alpha,\label{eq:spacetime-smoothing-S2}\\
    \norm{\tfrac{d}{d\theta}S_\theta u}_{\tilde Y^{\beta,\frac1{2\kappa}}_{r_0}}&\le C\theta^{\beta-\alpha-1}\norm{u}_{\tilde Y^{\alpha,\frac1{2\kappa}}_{r_0}},\label{eq:spacetime-smoothing-S3}
\end{align}
where $(\beta-\alpha)_+:=\max{(\beta-\alpha,0)}$ and $C=C(m_*,\kappa,\Omega_0,r_0)>0$ is independent of $\theta,\alpha,\beta$ and $u$. Moreover, the same estimates hold with the scale $\tilde Y^{l,\frac{1}{2\kappa}}_{r_0}$ replaced by $\tilde X^{l,\frac{1}{2\kappa}}_{r_0}$ for integers $1\le\alpha,\beta\le m_*$.
\end{lem}
\begin{proof}
The estimates on $\tilde Y^{\alpha,\frac1{2\kappa}}_{r_0}$ follow by applying \eqref{eq:S1}--\eqref{eq:S3} of Lemma~\ref{lem:spatial-smoothing-operator-appendix} at each $t\in[0,T]$ and taking the supremum in $t$.

For the estimates on $\tilde X^{\alpha,\frac1{2\kappa}}_{r_0}$, fix integers $1\leq\alpha,\beta\leq m_*$. Since $S_\theta$ is independent of $t$ and bounded on $H^{2(\beta-1),\beta-1+\frac{1}{2\kappa}}_{r_0}$ by \eqref{eq:S1}, for every $\varphi\in C_c^\infty((0,T))$, one has, as an identity in $H^{2(\beta-1),\beta-1+\frac{1}{2\kappa}}_{r_0}$ between Bochner integrals,
\begin{align*}
    \int_0^T S_\theta u(t)\,\varphi^\prime(t)\,dt
    =S_\theta\left(\int_0^T u(t)\,\varphi^\prime(t)\,dt\right)
    =-S_\theta\left(\int_0^T\p_tu(t)\,\varphi(t)\,dt\right)
    =-\int_0^T S_\theta\p_tu(t)\,\varphi(t)\,dt.
\end{align*}

Consequently,
\begin{equation}
\label{eq:time-derivative-smoothing-commutation}
    \p_t(S_\theta u)=S_\theta\p_tu
\end{equation}
in the sense of distributions on $(0,T)$. The same calculation, with $\frac{d}{d\theta}S_\theta$ in place of $S_\theta$, gives
\begin{equation}
\label{eq:time-derivative-parameter-derivative-commutation}
    \p_t\left(\frac{d}{d\theta}S_\theta u\right)=\frac{d}{d\theta}S_\theta\p_tu.
\end{equation}

Applying \eqref{eq:S1}--\eqref{eq:S3} at each $t$ with input index $\alpha-1$ and output index $\beta-1$, whose difference is again $\beta-\alpha$, and using \eqref{eq:time-derivative-smoothing-commutation} and \eqref{eq:time-derivative-parameter-derivative-commutation} gives
\begin{align*}
    \norm{\p_t S_\theta u}_{L^\infty_tH^{2(\beta-1),\beta-1+\frac{1}{2\kappa}}_{r_0}}
    &\leq C\theta^{(\beta-\alpha)_+}\norm{\p_t u}_{L^\infty_tH^{2(\alpha-1),\alpha-1+\frac{1}{2\kappa}}_{r_0}},
    \\
    \norm{\p_t (S_\theta u-u)}_{L^\infty_tH^{2(\beta-1),\beta-1+\frac{1}{2\kappa}}_{r_0}}
    &\leq C\theta^{\beta-\alpha}\norm{\p_t u}_{L^\infty_tH^{2(\alpha-1),\alpha-1+\frac{1}{2\kappa}}_{r_0}},
    \qquad \beta\leq\alpha,
    \\
    \norm{\p_t \tfrac{d}{d\theta}S_\theta u}_{L^\infty_tH^{2(\beta-1),\beta-1+\frac{1}{2\kappa}}_{r_0}}
    &\leq C\theta^{\beta-\alpha-1}\norm{\p_t u}_{L^\infty_tH^{2(\alpha-1),\alpha-1+\frac{1}{2\kappa}}_{r_0}}.
\end{align*}
Adding the above to the inequalities for the scale $\tilde Y^{l,\frac{1}{2\kappa}}_{r_0}$ proves the three estimates on $\tilde X^{\alpha,\frac{1}{2\kappa}}_{r_0}$.

\end{proof}

\begin{rmk}\label{rmk:smoothing-range}
Definition~2.3 of \cite{secchi2016nash} states the smoothing estimates for all regularity indices. The proof of Theorem~2.4 in that reference invokes only finitely many of them. More precisely, Lemma~3.2 of \cite{secchi2016nash} uses the $X$-scale through $\tilde\alpha+\max\{r^\prime,s^\prime\}$, while Lemma~3.7 of that reference uses the $Y$-scale through $\tilde\alpha+s$, where $r,r^\prime,s,s^\prime$ are the losses in Assumptions~2.1 and~2.2 of that reference. For the tame estimates stated in this paper, these losses (see Section~\ref{sec:proof-main-theorem}) are
\begin{equation*}
    r = 2,\qquad r'=3,\qquad s = 0,\qquad s'=4.
\end{equation*}
Thus, for the base iteration in Theorem~2.4(i) of \cite{secchi2016nash}, one takes $m_*=\tilde\alpha+4$. The kernels are
selected after $m_*$ has been fixed, with $N_*\geq m_*$.
\end{rmk}

\subsection{Admissible data}
Here, we specify the class of initial data. Since the physical vacuum boundary condition is stated in terms of the distance function $\dist(\cdot,\p\Omega_0)$, which is smooth only near $\p\Omega_0$, we first introduce a globally smooth substitute for it. We then list the conditions on $(r_0,\vec u_0,\pi_0)$. They ensure that the physical vacuum boundary condition holds at $t=0$ and propagates, that the coefficient $a_0=\pi/r$ stays bounded, and that the data are regular enough for the Nash--Moser iteration.

\begin{defn}\label{defn:reference-function}
Let $\Omega_0\subset\mathbb R^3$ be a bounded domain with smooth boundary. There exists $\varepsilon_0>0$ such that $\mathrm{dist}(\cdot,\p\Omega_0)$ is well-defined and smooth on $\Omega_0^{\varepsilon_0}:=\{x\in\Omega_0 \mid\mathrm{dist}(x,\p\Omega_0)<\varepsilon_0\}$. Fix bump functions $\chi_{\varepsilon_0}\equiv1$ on $\Omega_0\setminus\Omega_0^{\varepsilon_0}$, vanishing outside an open subset of $\Omega_0$ containing $\Omega_0\setminus\Omega_0^{\varepsilon_0}$, and $\chi_{\frac{\varepsilon_0}{2}}\equiv1$ on
$\Omega_0\setminus\Omega_0^{\frac{\varepsilon_0}{2}}$, vanishing outside an open subset of
$\Omega_0$ containing $\Omega_0\setminus\Omega_0^{\frac{\varepsilon_0}{2}}$. The
reference function $d$ is defined as
\[
    d := \mathrm{dist}(\cdot,\p\Omega_0)(1-\chi_{\varepsilon_0}) + \chi_{\frac{\varepsilon_0}{2}}.
\]
\end{defn}

\begin{rmk}
	Since $\mathrm{dist}(\cdot,\p\Omega_0)$ is defined and smooth only near $\p\Omega_0$, it cannot by itself serve as a weight over all of $\Omega_0$. The reference function $d$ is a globally smooth substitute that agrees with the true distance function near the boundary while being bounded away from $0$ in the interior. In particular, any function comparable to the reference function $d$ is positive, bounded away from $0$ away from the boundary and satisfies the physical vacuum boundary condition.
\end{rmk}

\begin{defn}\label{defn:admissible-data}
Let $k\ge1$ be an integer. We say that initial data $(r_0, u_0,\pi_0)$ on $\Omega_0$ are \emph{$k$-admissible} if the following hold, where $d$ is the reference function of Definition~\ref{defn:reference-function}:
\begin{enumerate}
\item[(1)] $r_0,\pi_0$ are comparable to $d$: there exist constants $0<\mfc_1\le \mfc_2$ such that
\[
    \mfc_1d \le r_0 \le \mfc_2d, \qquad \mfc_1d \le \pi_0 \le \mfc_2d \qquad\text{on }\Omega_0.
\]
\item[(2)] There exist constants $0<\mfc\le \mfC$ such that
\[
    \mfc \le \frac{\pi_0}{r_0} \le \mfC \qquad\text{on }\Omega_0.
\]
\item[(3)] The initial $4$-velocity satisfies the normalization condition $u_0^\mu (u_0)_\mu=-1$.
\item[(4)] The data have the regularity
\begin{equation*}
	r_0\in H_d^{2k,k+\frac1{2\kappa}-\frac12},\qquad \vec u_0\in H_d^{2k,k+\frac1{2\kappa}},\qquad \frac{\pi_0}{r_0}\in H^{2k,k+\frac1{2\kappa}-\frac12}_d.
\end{equation*}
\item[(5)] The initial density is small:
\begin{equation}\label{eq:uniform-smallness-r0}
	\norm{r_0}_{L^\infty}\leq \varepsilon_0,
\end{equation}
where $\varepsilon_0>0$ is the constant fixed in Remark~\ref{rmk:choice-of-epsilon0}.
\end{enumerate}
\end{defn}

\begin{rmk}\label{rmk:choice-of-epsilon0}
Condition (1) is the physical vacuum boundary condition at $t=0$. Since it makes the weights $d$ and $r_0$ comparable, we may replace $H_d$ by $H_{r_0}$ in condition (4), and we use the two notations interchangeably.

The causality condition \eqref{eq:causality} is satisfied whenever $r$ is sufficiently small. We fix $\varepsilon_0\ll1$, depending only on $\kappa$ and on $\sup_{[0,1]}\hat\zeta$, such that \eqref{eq:causality} holds at every state with $r\le2\varepsilon_0$. Condition (5) then implies \eqref{eq:causality} at $t=0$. We show in Section~\ref{sec:linearized-eqs} that, for $T$ sufficiently small, $r\le2\varepsilon_0$ on $[0,T]$, so that \eqref{eq:causality} holds on $[0,T]$ (see Remark~\ref{rmk:persistence-of-causality}).

Condition (5) is not restrictive. Near the boundary, it holds automatically by the physical vacuum boundary condition. Away from the boundary, if the causality condition \eqref{eq:causality} holds, the system \eqref{eq:Eulerian-Main-System} is a non-degenerate hyperbolic system \cite{bemfica2019causality}, and the standard local theory applies. By finite speed of propagation, the solution for general data satisfying \eqref{eq:causality} is obtained by gluing the two.
\end{rmk}

\subsection{Approximate solution and modified iteration map}\label{subsec:approximate-solution-modified-iteration-map}
To define the iteration scheme, it is convenient to recast the problem in terms of equations with vanishing initial data. In this way, the problem can be viewed as inverting a nonlinear map between Banach spaces. For this purpose, we introduce an approximate solution.

We take $k_0$-admissible data $(r_0,u_0,\pi_0)$, where $k_0$ will be chosen later on:
\begin{equation*}
		r_0\in H_{r_0}^{2k_0,k_0+\frac1{2\kappa}-\frac12},\quad u_0\in H_{r_0}^{2k_0,k_0+\frac1{2\kappa}},\quad \frac{\pi_0}{r_0}\in H^{2k_0,k_0+\frac1{2\kappa}-\frac12}_{r_0}.
\end{equation*}

We now define $\p_t\vec U|_{t=0}$. At $t=0$ one has $\eta=\mathrm{id}$, $F=I$, $(R,A_0,P)=(r_0,\pi_0/r_0,\pi_0)$,
and \eqref{eq:R-A0-defining}, \eqref{eq:R-P-defining} express $\p_tR|_{t=0}$ and $\p_tP|_{t=0}$ as affine functions of $\p_t\vec U|_{t=0}$ with coefficients depending on the data. Substituting these into \eqref{eq:iteration-map-equation} at $t=0$ yields a linear system for $\p_t\vec U|_{t=0}$. Eliminating $\p_tR|_{t=0}$, $\p_tP|_{t=0}$ in this way is equivalent to solving the system
\eqref{eq:Lagrangian-formulation}(a)--(c) at $t=0$ for the time derivatives in terms of the spatial derivatives. We define $\p_t\vec U|_{t=0}$ as the resulting vector field, for which \eqref{eq:iteration-map-equation} holds at $t=0$. Then, we can solve for $\p_t^j\vec U|_{t=0}$, for $j=1,2,\ldots$, recursively in terms of the initial data by taking the time derivative of the formulation \eqref{eq:iteration-map-equation} and restricting to $t=0$.

We introduce the following notation for the bounds of the initial data:
\begin{equation}\label{eq:initial-data-bounds}
	\begin{aligned}
		N^k_\infty &:= \norm{r_0}_{W^{k,\infty}} + \norm{\frac{\pi_0}{r_0}}_{W^{k,\infty}} + \norm{\frac{r_0}{\pi_0}}_{L^\infty},\\
		N^{l}_w &:= \norm{r_0}_{H^{2l,l+\frac{1}{2\kappa}-\frac12}_{r_0}} + \norm{\frac{\pi_0}{r_0}}_{H^{2l,l+\frac{1}{2\kappa}-\frac12}_{r_0}}.
	\end{aligned}
\end{equation}
Immediately by Lemma~\ref{lem:Morrey}, we have
\begin{equation*}
	N^1_\infty\lesssim N^{l_0}_w + \mfc^{-1},
\end{equation*} 
where $\mfc$ is the lower bound for $\frac{\pi_0}{r_0}$ as in Definition~\ref{defn:admissible-data} and $l_0$ as in Notation~\ref{nota:l0}.

\begin{rmk}
	Since the causality condition \eqref{eq:causality} is satisfied, we can solve for the time derivatives of $(R,\vec U,P,\eta)$ in terms of spatial derivatives using the Lagrangian formulation \eqref{eq:Lagrangian-formulation}, which agrees with the problem \eqref{eq:iteration-map-equation} at $t=0$. Schematically, we have the following
	\begin{equation}\label{eq:Lagrangian-time-derivative}
		\begin{aligned}
    		\p_t R &= R\calF(R,A_0,\vec U)\nabla R
      				+ R^2\calF(R,A_0,\vec U)\nabla P
      				+ R\calF(R,A_0,\vec U)\nabla \vec U
      				+ R^3\calF(R,A_0,\vec U),\\
    		\p_t \vec U &= \calF(R,A_0,\vec U)\nabla R
      				+ R\calF(R,A_0,\vec U)\nabla P
      				+ R\calF(R,A_0,\vec U)\nabla \vec U
      				+ R^2\calF(R,A_0,\vec U),\\
    		\p_t P &= R\calF(R,A_0,\vec U)\nabla R
      				+ R^2\calF(R,A_0,\vec U)\nabla P
      				+ R\calF(R,A_0,\vec U)\nabla \vec U
      				+ R\calF(R,A_0,\vec U),\\
    		\p_t\eta &= \calF(\vec U),
		\end{aligned}
	\end{equation}
	where $\calF$ is a smooth function in a neighborhood of the range of its arguments subject to change line-by-line and term-by-term.
\end{rmk}

We define the approximate solution by truncating the Taylor expansion in time at order $N_{\mathrm{app}}$:
\begin{equation}\label{eq:vec-U-approximate-solution}
	\vec U_{\mathrm{app}} = \sum_{j=0}^{N_{\mathrm{app}}}\frac{t^j}{j!}\p_t^j \vec U|_{t=0},
\end{equation}
and $R_{\mathrm{app}},A_{0,{\mathrm{app}}}$ and $\eta_{\mathrm{app}}$ are defined by
\begin{equation*}
	R_{\mathrm{app}}:=R(\vec U_{\mathrm{app}}),\quad A_{0,\mathrm{app}}:=A_0(\vec U_{\mathrm{app}}),\quad \eta_{\mathrm{app}}:=\eta(\vec U_{\mathrm{app}}),
\end{equation*}
through solving \eqref{eq:R-A0-defining} and \eqref{eq:eta-defining-eq}. We set $P_{\mathrm{app}} :=R_{\mathrm{app}}A_{0,\mathrm{app}}$.

Then, by construction,
\[
R_{\mathrm{app}}(0)=r_0,\quad
\vec U_{\mathrm{app}}(0)=\vec u_0,\quad
A_{0,\mathrm{app}}(0) = \frac{\pi_0}{r_0},\quad
P_{\mathrm{app}}(0)=\pi_0,\quad
\eta_{\mathrm{app}}(0)=\Id_{\Omega_0},
\]
and $R_{\mathrm{app}},P_{\mathrm{app}}$ vanish on $\Gamma_0$.

We now introduce a new unknown $\vec V$ by
\[
\vec U = \vec V + \vec U_{\mathrm{app}},
\]
and thus
\[
\p_t^j \vec V|_{t=0}=0,
\qquad j=0,\dots,N_{\mathrm{app}}.
\]

Define the modified nonlinear map
\[
\tilde{\mfF}(\vec V)
:= \mfF(\vec V+\vec U_{\mathrm{app}})
   - \mfF(\vec U_{\mathrm{app}}),
\qquad
\mfF^a := -\mfF(\vec U_{\mathrm{app}}).
\]
Then solving the original equation is equivalent to solving
\[
\tilde{\mfF}(\vec V)=\mfF^a,
\]
where the modified iteration map $\tilde{\mfF}$ is a nonlinear map between weighted Sobolev spaces that
vanishes to order $N_{\mathrm{app}}-1$ at the initial time $t=0$.

In our problem, we take $N_\mathrm{app}=1$, i.e.,
\begin{equation*}
	\vec U_{\mathrm{app}}:=\vec u_0 + t\p_t\vec U|_{t=0},
\end{equation*} 
since the spaces $\tilde X^{l,\frac{1}{2\kappa}}_{r_0}$ used in the Nash--Moser iteration control only one time derivative.

\begin{lem}\label{lem:approximate-solution-estimate}
Define
\begin{equation*}
	\mfN^l:=N^{l}_w + \norm{\vec u_0}_{H^{2l,l+\frac{1}{2\kappa}}_{r_0}} + \norm{\frac{r_0}{\pi_0}}_{L^\infty}.
\end{equation*}
Let $T\in(0,1]$. For any $l>l_0$ and $t\in[0,T]$
\begin{equation*}
	\norm{\vec U_{\mathrm{app}}(t)}_{\tilde X^{l,\frac{1}{2\kappa}}_{r_0}}\lesssim \const{\mfN^{l+1}}.
\end{equation*}	
\end{lem}
\begin{proof}
	Restricting the second equation of \eqref{eq:Lagrangian-time-derivative} to $t=0$, we have
	\begin{equation*}
		\begin{aligned}
			\p_t\vec U|_{t=0}&=\calF(r_0,\frac{\pi_0}{r_0},\vec u_0)\p r_0
      				+ r_0\calF(r_0,\frac{\pi_0}{r_0},\vec u_0)\p \pi_0
      				+ r_0\calF(r_0,\frac{\pi_0}{r_0},\vec u_0)\p \vec u_0 + r_0^2\calF(r_0,\frac{\pi_0}{r_0},\vec u_0)\\
      				&=\calF(r_0,\frac{\pi_0}{r_0},\vec u_0)\p r_0
      				+ r_0^2\calF(r_0,\frac{\pi_0}{r_0},\vec u_0)\p \left(\frac{\pi_0}{r_0}\right)
      				+ r_0\calF(r_0,\frac{\pi_0}{r_0},\vec u_0)\p \vec u_0 + r_0^2\calF(r_0,\frac{\pi_0}{r_0},\vec u_0),
		\end{aligned}
	\end{equation*}
	where $\calF$ is a smooth function in its arguments subject to change line-by-line and term-by-term.
	By Corollary~\ref{cor:weighted-L2-end-point-equivalence}, we have for any $l>l_0$
	\begin{equation*}
		\norm{\p_t\vec U|_{t=0}}_{H_{r_0}^{2l-2,l+\frac1{2\kappa}-1}}\lesssim\norm{r_0^{l+\frac{1}{2\kappa}-1}\p_t\vec U|_{t=0}}_{L^2} + \norm{r_0^{l+\frac{1}{2\kappa}-1}\p^{2l-2}(\p_t\vec U|_{t=0})}_{L^2}.
	\end{equation*}
	Applying Lemma~\ref{lem:WeightedEmbedding} to the first term and Corollary~\ref{Cor:Multilinear-Interpolation} to the second term, we can show
	\begin{equation*}
		\norm{r_0^{l+\frac{1}{2\kappa}-1}\p_t\vec U|_{t=0}}_{L^2},\norm{r_0^{l+\frac{1}{2\kappa}-1}\p^{2l-2}(\p_t\vec U|_{t=0})}_{L^2}\lesssim\calC(\mfN^l).
	\end{equation*}
	Then, since $T\le 1$
	\begin{equation*}
		\begin{aligned}
			\norm{\vec U_{\mathrm{app}}}_{\tilde X^{l,\frac{1}{2\kappa}}_{r_0}}&\lesssim\norm{\vec u_0}_{H^{2l,l+\frac{1}{2\kappa}}_{r_0}} + T\norm{\p_t\vec U|_{t=0}}_{H_{r_0}^{2l,l+\frac1{2\kappa}}} + \norm{\p_t\vec U|_{t=0}}_{H_{r_0}^{2l-2,l+\frac1{2\kappa}-1}}\\
			&\lesssim \calC(\mfN^{l+1}).
		\end{aligned}
	\end{equation*}
\end{proof}

\begin{rmk}
Throughout, the constants $\calC(\cdot)$ appearing in the estimates also depend on the regularity level $l$ at which the estimates are carried out. We suppress this dependence in the notation, since $l$ ranges only over values bounded by $k_0$, the highest regularity level appearing in this article.
\end{rmk}

\subsection{Iteration scheme}
We now describe the Nash--Moser iteration scheme. Initializing with $\vec V_0 = 0$, we define inductively
\[
	\vec V_{n+1} = \vec V_n + \delta \vec V_n.
\]

A Taylor expansion of $\tilde{\mfF}$ at $\vec V_n$ gives
\begin{equation}\label{eq:iteration-difference}
\begin{aligned}
\tilde{\mfF}(\vec V_{n+1}) - \tilde{\mfF}(\vec V_n)
&= (D\tilde{\mfF})(\vec V_n)\,\delta\vec V_n + e_n^{(1)} \\
&= (D\mfF)(\vec V_n+\vec U_{\mathrm{app}})\,\delta\vec V_n + e_n^{(1)} \\
&= (D\mfF)(S_{\theta_n}\vec V_n+\vec U_{\mathrm{app}})\,\delta\vec V_n
+ e_n^{(1)} + e_n^{(2)},
\end{aligned}
\end{equation}
where
\begin{equation}
e_n^{(1)}
= \mfF(\vec V_{n+1}+\vec U_{\mathrm{app}})
- \mfF(\vec V_n+\vec U_{\mathrm{app}})
- (D\mfF)(\vec V_n+\vec U_{\mathrm{app}})\,\delta\vec V_n
\end{equation}
is the quadratic error, and
\begin{equation}
e_n^{(2)}
= (D\mfF)(\vec V_n+\vec U_{\mathrm{app}})\,\delta\vec V_n
- (D\mfF)(S_{\theta_n}\vec V_n+\vec U_{\mathrm{app}})\,\delta\vec V_n
\end{equation}
is the substitution error arising from replacing the background state
$\vec V_n+\vec U_{\mathrm{app}}$ by its smoothed version
$S_{\theta_n}\vec V_n+\vec U_{\mathrm{app}}$.

We set
\[
e_n := e_n^{(1)} + e_n^{(2)},
\]
which is the total error at the $n$th iteration.

Summing \eqref{eq:iteration-difference} from $0$ to $n$ yields
\begin{equation}
\tilde{\mfF}(\vec V_{n+1})
= \sum_{i=0}^{n}
(D\mfF)(S_{\theta_i}\vec V_i+\vec U_{\mathrm{app}})\,\delta\vec V_i
+ \sum_{i=0}^{n} e_i.
\end{equation}

Formally, the iteration is designed in a way that the error terms converge to $0$ and
\[
\sum_{i=0}^n e_i \to 0,
\qquad
\sum_{i=0}^n
(D\mfF)(S_{\theta_i}\vec V_i+\vec U_{\mathrm{app}})\,\delta\vec V_i
\to \mfF^a.
\]
We therefore define the source terms 
\[
\vec G_i := (D\mfF)(S_{\theta_i}\vec V_i+\vec U_{\mathrm{app}})\,\delta\vec V_i
\]
by imposing
\[
\sum_{i=0}^{n} \vec G_i
= - S_{\theta_n}\Bigl(\sum_{i=0}^{n-1} e_i\Bigr)
+ S_{\theta_n}\mfF^a.
\]
For $n=0$, the sum on the right-hand side is empty.

This determines recursively
\begin{equation}
    \left\{\begin{aligned}
        \vec G_0 &= S_{\theta_0}\mfF^a,\\
        \vec G_1 &= (S_{\theta_1} -S_{\theta_0})\mfF^a - S_{\theta_1} e_0,\\
        \vec G_n &= (S_{\theta_{n-1}}-S_{\theta_n})(\sum_{j=0}^{n-2}e_j) - S_{\theta_n}e_{n-1} + (S_{\theta_n}-S_{\theta_{n-1}})\mfF^a \qquad n\geq 2.
    \end{aligned}\right.
\end{equation}

\subsection{Linearization of the iteration map over the Eulerian domain}

To make the iteration scheme well-defined, at each step we solve the linearized equation
\begin{equation}\label{eq:linearized-operator-iteration-map}
    (D\mfF)(\vec U)\,\delta\vec U = \vec G.
\end{equation}

To achieve that, we move the problem to the Eulerian domain and use the local well-posedness result over the Eulerian domain of Theorem~\ref{thm:Modified-linearized-equations-wellposedness}, as follows:

From the vector field $\vec U$, we have defined $(R,U^0,\vec U, P,\eta)$ and $(r,u,\pi)$ in terms of $\vec U$ through \eqref{eq:R-P-defining}, \eqref{eq:eta-defining-eq} and Remark~\ref{rmk:r-u-pi-defining}. Setting
\begin{equation}\label{eq:variables-abbreviation}
	\vec \Phi = \vec \Phi(\vec U) := 
		\begin{pmatrix}
        	R \\ \vec U \\ P \\ \eta
    	\end{pmatrix},
    	\qquad 
    	\Phi = \Phi(\vec U) := 
    	\begin{pmatrix}
        	R \\ U \\ P \\ \eta
    	\end{pmatrix},
    	\qquad
    	\vec \phi = \vec \phi(\vec U) := \begin{pmatrix}
        r \\ \vec u \\ \pi
    	\end{pmatrix},
    	\qquad
    	\phi = \phi(\vec U) := \begin{pmatrix}
        r \\ u \\ \pi 
    	\end{pmatrix}
\end{equation}
for convenience and taking variations with respect to $\vec U$, we obtain
\begin{equation}\label{eq:linearized-variables-change-of-coordinate}
    \delta\vec \Phi = \begin{pmatrix}
        \delta R \\ \delta \vec U \\ \delta P \\ \delta \eta
    \end{pmatrix}   
    = \begin{pmatrix}
        \delta\vec \phi\circ\eta + (\p_i \vec \phi)\circ\eta\delta\eta^i\\
        \delta\eta
    \end{pmatrix},
\end{equation}
where
\begin{equation*}
	\delta\vec\phi = \begin{pmatrix}
        \delta r \\ \delta\vec u \\ \delta \pi
    \end{pmatrix}.
\end{equation*}

Define
\begin{equation}\label{eq:abstract-Eulerian-operator}
	\calN=\calN(\phi)
\end{equation}
as the abstract equation operator associated to the Eulerian formulation \eqref{eq:Eulerian-Main-System} and
\begin{equation}\label{eq:abstract-spatial-Eulerian-operator}
	\vec \calN = \vec\calN(\phi)
\end{equation}
 as the operator leaving out the part associated to $u^0$ equation. We can then rewrite the operator $\mfF$ in \eqref{eq:iteration-map} as 
\[
\mfF(\vec U) = \big(\vec\calN_{\mathrm{mo}}(\phi(\vec U))\big)\circ\eta,
\]
where
\begin{equation}\label{eq:abstract-Eulerian-momentum-operator}
	\vec \calN_{\mathrm{mo}}(\phi):=u^\mu\p_\mu u^i
+\big(1+\tfrac1\kappa\big)a_2\Delta^{i\mu}\p_\mu r
+\big(2+\tfrac1\kappa\big)ra_3\Delta^{i\mu}\p_\mu\pi
\end{equation}
denotes the spatial part of the momentum part of the operator $\calN(\phi)$.

Taking variations with respect to $\vec U$ yields
\begin{equation}\label{eq:linearization-chain-rule}
    (D\mfF)(\vec U)\delta\vec U = ((D\vec\calN_{\mathrm{mo}})(\phi)\delta\phi)\circ\eta + (\p_i\vec\calN_{\mathrm{mo}}(\phi))\circ\eta\delta\eta^i.
\end{equation}
So, the equation \eqref{eq:linearized-operator-iteration-map} can be rewritten as
\[
((D\vec\calN_{\mathrm{mo}})(\phi)\delta\phi)\circ\eta + (\p_i\vec\calN_{\mathrm{mo}}(\phi))\circ\eta\delta\eta^i = \vec G.
\]

Define the variables $\xi$ and $\vec g$
\[\xi := \delta\eta\circ\eta^{-1},\qquad\vec g := \vec G \circ \eta^{-1}.\]
Composing with the inverse flow map $\eta^{-1}$ yields equations in Eulerian coordinates:
\[
D\vec\calN_{\mathrm{mo}}(\phi)\delta\phi + \p_i\vec\calN_{\mathrm{mo}}(\phi)\xi^i = \vec g,
\]
which, written out explicitly, reads
\begin{equation}\label{eq:delta-u-equation}
    \begin{aligned}
        &u^\mu\p_\mu \delta u^i + \left(1+\frac{1}{\kappa}\right)a_2\Delta^{i\mu}\p_\mu \delta r + \left(2+\frac{1}{\kappa}\right)ra_3\Delta^{i\mu}\p_\mu\delta \pi + V_2^i \delta r + (W_2)_\beta^i \delta{u}^\beta + Z_2^i\delta \pi \\
        &\quad= \vec g^i- \p_j(u^\mu \p_\mu u^i + \left(1+\frac{1}{\kappa}\right)a_2\Delta^{i\mu}\partial_\mu r + \left(2+\frac{1}{\kappa}\right)ra_3\Delta^{i\mu}\partial_\mu\pi)\xi^j
    \end{aligned}
\end{equation}
for $i=1,2,3$.

Notice that the functions $(r,\pi,\eta)$ are defined through solving \eqref{eq:eta-defining-eq} and \eqref{eq:r-pi-defining} with $\vec U$ fixed, so their variations $\delta r,\delta\pi,\delta\eta$ satisfy the linearizations of the equations \eqref{eq:eta-defining-eq} and \eqref{eq:r-pi-defining}. In Eulerian coordinates, they read
\begin{equation}
    \begin{aligned}
        u^\mu\p_\mu\delta r + \kappa r a_1\p_\mu\delta u^\mu + \p_\mu r \delta u^\mu + V_1 \delta r + r^2 Z_1\delta \pi = 0,\\
        u^\mu\p_\mu \delta \pi + r a_4\p_\mu \delta u^\mu + V_3 \delta r + W_3^\alpha\delta{u}_\alpha + Z_3 \delta \pi = 0,\\
        u^\mu\p_\mu \xi - \delta\vec u + \frac{\vec u}{u^0}\delta u^0 - u^0\p_i\left(\frac{\vec u}{u^0}\right) \xi^i=0.
    \end{aligned}
\end{equation}

Combining the above equations for $(\delta r,\delta\vec u,\delta \pi,\xi)$, we have a closed system:
\begin{subequations}\label{eq:modified-linearized-equation-with-xi}
    \begin{align}
        &u^\mu\p_\mu\delta r + \kappa r a_1\p_\mu\delta u^\mu + \p_\mu r \delta u^\mu + V_1 \delta r + r^2 Z_1\delta \pi = 0 \quad \text{in } \sptdomain, \\
        &u^\mu\p_\mu \delta u^i + \left(1+\frac{1}{\kappa}\right)a_2\Delta^{i\mu}\p_\mu \delta r + \left(2+\frac{1}{\kappa}\right)ra_3\Delta^{i\mu}\p_\mu\delta \pi + V_2^i \delta r + (W_2)_\beta^i \delta{u}^\beta + Z_2^i\delta \pi \nonumber\\
        &\quad= \vec g^i- \p_j(u^\mu \p_\mu u^i + \left(1+\frac{1}{\kappa}\right)a_2\Delta^{i\mu}\partial_\mu r + \left(2+\frac{1}{\kappa}\right)ra_3\Delta^{i\mu}\partial_\mu\pi)\xi^j \quad \text{in } \sptdomain, \\
        &u^\mu\p_\mu \delta \pi + r a_4\p_\mu \delta u^\mu + V_3 \delta r + W_3^\alpha\delta{u}_\alpha + Z_3 \delta \pi = 0 \quad \text{in } \sptdomain,\\
        &u^\mu\p_\mu \xi - \delta\vec u + \frac{\vec u}{u^0}\delta u^0 - u^0\p_i\left(\frac{\vec u}{u^0}\right)\xi^i=0 \quad \text{in } \sptdomain,\label{eq:xi-equation}\\
        &u_\mu \delta u^\mu = 0 \quad \text{in } \sptdomain, \\
        &\sptdomain = \bigcup_{0\leq t\leq T}\{t\}\times\Omega_t,\\
        &\Gamma = \bigcup_{0\leq t\leq T}\{t\}\times\Gamma_t,\\
        &\delta r(0,\cdot)=0,\qquad\delta u(0,\cdot)= 0,\qquad \delta \pi(0,\cdot) =0,\qquad \xi(0,\cdot)=0,
\end{align}
\end{subequations}
which we solve in Section~\ref{sec:linearized-eqs} by means of the modified linearized equations.

Conversely, now let $(\delta r,\delta \vec u,\delta\pi,\xi)$ be a solution of the system \eqref{eq:modified-linearized-equation-with-xi} and define
\begin{equation*}
    \delta\vec U:=\delta \vec u\circ\eta + (\p_i\vec u\xi^i)\circ\eta.
\end{equation*}
Then $\delta\vec U$ solves the linearized equation
\eqref{eq:linearized-operator-iteration-map}. Indeed, we can show that the variation of $\eta$ along $\delta\vec U$ defined as above is exactly equal to $\xi\circ\eta$. Thus, the variation of $\vec u$ along $\delta\vec U$ is $\delta\vec U\circ\eta^{-1}-\p_i\vec u\xi^i$, which is exactly $\delta\vec u$. Hence, the variations of $r$ and $\pi$ satisfy the same transport equations, with the same zero initial data, as $(\delta r,\delta\pi)$ and coincide with them by uniqueness. By the second equation of the system \eqref{eq:modified-linearized-equation-with-xi} and \eqref{eq:linearization-chain-rule}, it follows that $\delta\vec U$ solves \eqref{eq:linearized-operator-iteration-map}.

\section{Invertibility of the linearization of the iteration map over the Eulerian domain}\label{sec:linearized-eqs}
In this section, we begin by analyzing the linearized Israel--Stewart system in the Eulerian coordinates around a given background state. Throughout this section, the terms ``linearized'' or ``linearization'' refer to taking the G\^ateaux derivative of the nonlinear operator associated to the system \eqref{eq:Eulerian-Main-System} (or any equivalent formulation) at a given background state. In particular, we do not assume the existence of steady states or any other special class of background solutions.

\subsection{Linearized equations and local well-posedness}
Let the variables $\tilde\phi:=(\tilde{r},\tilde{u},\tilde \pi)$ be the linearized variables associated to the background states $\phi=(r,u,\pi)$ that solve the system \eqref{eq:Eulerian-Main-System}. Below we record the linearized equations $D\calN(\phi)\tilde\phi = \tilde F$ with $\tilde F:=(\tilde f,\tilde g,\tilde h)^T$ in explicit form:

\begin{subequations}\label{eq:Eulerian-Linearized-Equation}
    \begin{align}
        &u^\mu\p_\mu\tilde r + \kappa r a_1\p_\mu\tilde u^\mu + \p_\mu r \tilde u^\mu + V_1 \tilde r + r^2 Z_1\tilde \pi = \tilde f \quad \text{in } \sptdomain, \label{eq:tilde-r-equation}\\
        &u^\mu\p_\mu \tilde u^\alpha + \left(1+\frac{1}{\kappa}\right)a_2\Delta^{\alpha\mu}\p_\mu \tilde r + \left(2+\frac{1}{\kappa}\right)ra_3\Delta^{\alpha\mu}\p_\mu\tilde \pi + V_2^\alpha \tilde r + (W_2)_\beta^\alpha \tilde{u}^\beta + Z_2^\alpha\tilde \pi = \tilde g^\alpha \quad \text{in } \sptdomain,\label{eq:tilde-u-Equation} \\
        &u^\mu\p_\mu \tilde \pi + r a_4\p_\mu \tilde u^\mu + V_3 \tilde r + W_3^\alpha\tilde{u}_\alpha + Z_3 \tilde \pi = \tilde h \quad \text{in } \sptdomain,\label{eq:tilde-pi-equation}\\
        &u_\mu \tilde u^\mu = 0 \quad \text{in } \sptdomain, \label{eq:orthogonality-constraint}\\
        &\sptdomain= \bigcup_{0\leq t\leq T}\{t\}\times\Omega_t,\\
        &\Gamma = \bigcup_{0\leq t\le T}\{t\}\times\Gamma_t,\\
        &\tilde r(0,\cdot)=\tilde r_0(\cdot),\qquad\tilde u(0,\cdot)= \tilde u_0(\cdot),\qquad \tilde \pi(0,\cdot) =\tilde \pi_0(\cdot),
\end{align}
\end{subequations}
where $V,W,Z$ are coefficients given as the following,
\begin{equation*}
    \begin{aligned}
        &V_1 = \kappa \left[2r + 1 + \left(1 - \frac{1}{\kappa}\right) r^2 a_0^{2 + \frac{1}{\kappa}}\right] \partial_\mu u^\mu, \\
        &V_2^\alpha = 
            \left\{\left(\frac{1}{\kappa} r a_0^{2 + \frac{1}{\kappa}} - 1\right) a_2 
            \left[\left(1 + \frac{1}{\kappa}\right) \partial_\mu r + \left(2 + \frac{1}{\kappa}\right) r a_0^{1 + \frac{1}{\kappa}} \partial_\mu \pi \right] - \frac{1}{\kappa}\left(2 + \frac{1}{\kappa}\right) a_0^{1 + \frac{1}{\kappa}} \partial_\mu \pi \right\} a_2 \Delta^{\alpha \mu}, \\
        &V_3 = a_0^{3+\frac{1}{\kappa}}\left(\hat\lambda(r)+r\hat\lambda'(r)-(3+\frac{1}{\kappa})\hat\lambda(r)\right)+\frac{\hat\zeta(r)}{a_0^{1+\frac{1}{\kappa}}}\p_\mu u^\mu + r\frac{\hat{\zeta}'(r)}{a_0^{1+\frac{1}{\kappa}}}\p_\mu u^\mu+\left(1+\frac{1}{\kappa}\right)\frac{\hat\zeta(r)}{a_0^{1+\frac{1}{\kappa}}}\p_\mu u^\mu,\\
        &W_1 = 0, \\
        &(W_2)^{\alpha \beta} = \left\{ \partial_\mu u^\alpha g^{\mu \beta} + a_2 \left(u^\mu g^{\alpha \beta} + u^\alpha g^{\mu \beta}\right) 
        \left[\left(1 + \frac{1}{\kappa}\right) \partial_\mu r + \left(2 + \frac{1}{\kappa}\right) r a_0^{1 + \frac{1}{\kappa}} \partial_\mu \pi\right] \right\}, \\
        &(W_3)^\beta = g^{\mu \beta} \partial_\mu \pi, \\
        &Z_1 = \kappa \left(2 + \frac{1}{\kappa}\right) a_0^{1 + \frac{1}{\kappa}} \partial_\mu u^\mu, \\
        &Z_2^\alpha = 
            -\left(2 + \frac{1}{\kappa}\right) a_0^{\frac{1}{\kappa}} a_2 \Delta^{\alpha \mu} \Biggl\{\left(1 + \frac{1}{\kappa}\right) r a_0 a_2 \partial_\mu r + \left(2 + \frac{1}{\kappa}\right) r^2 a_0^{2 + \frac{1}{\kappa}} a_2 \partial_\mu \pi - \left(1 + \frac{1}{\kappa}\right) \partial_\mu \pi \Biggr\}, \\
        &Z_3 = 1 + (3+\frac{1}{\kappa})\hat\lambda a_0^{2+\frac{1}{\kappa}} - \left(1 + \frac{1}{\kappa}\right) \frac{\hat\zeta(r)}{a_0^{2 + \frac{1}{\kappa}}} \partial_\mu u^\mu.
    \end{aligned}
\end{equation*}
Above, the source terms $\tilde F$ are included for the purpose of the Nash--Moser iteration.

\begin{rmk}
	When linearizing the system, each member of the one-parameter family  $(r_\lambda, u_\lambda, \pi_\lambda)$ does not necessarily share the same domain, and thus there are no natural boundary conditions imposed on the linearized variables $(\tilde{r}, \tilde{u}, \tilde{\pi})$. However, the boundary behavior of the linearized variables is encoded into our choice of weighted Sobolev spaces. This is reflected in the energy estimates, where all the boundary terms vanish due to the weight vanishing on the boundary.
\end{rmk}

We now specify the function spaces for the solution of the linearized equations. The linearized $4$-velocity $\tilde u$ satisfies the orthogonality constraint \eqref{eq:orthogonality-constraint}, which arises from linearizing the normalization condition \eqref{eq:Normalization-Condition}. Thus, this constraint determines $\tilde u^0$ uniquely in terms of the spatial components by
\[
	\tilde u^0 =\frac{u^i}{u^0}\tilde u_i.
\]
Hence the map $\tilde u\mapsto(\tilde u^1,\tilde u^2,\tilde u^3)$ is a linear isomorphism onto its image, and $\tilde u$ is completely determined by its spatial part. We may therefore define the norm of $\tilde u$ using only the spatial components by setting
\[
	\abs{\tilde u}_H^2:=H^{ij}\tilde u_i\tilde u_j,\qquad H^{ij}=\delta^{ij}-\frac{u^iu^j}{(u^0)^2}.
\]
Note that $H^{ij}$ is positive definite because $u$ is time-like. We then define the weighted Sobolev spaces of the $4$-vector field $\tilde u$ orthogonal to $u$ as the space of $4$-vector fields with the following norm being finite:
\begin{equation}\label{eq:4-velocity-weighted-norm}
    \norm{\tilde u}_{H^{2l,l+\frac{1}{2\kappa}}}^2
    := \sum_{\abs{\beta}\leq 2l}
    \norm{r^{l+\frac{1}{2\kappa}}\abs{\partial^\beta \tilde u}_H}_{L^2}^2.
\end{equation}

\begin{rmk}\label{rmk:u-0-play-no-role}
As defined in \eqref{eq:4-velocity-weighted-norm}, the
weighted Sobolev norm $\norm{\cdot}_{H^{2l,l+\frac{1}{2\kappa}}}$ of the linearized velocity depends on $\tilde u$ only through its spatial components $\tilde u_i$, since $\tilde
u^0$ is determined algebraically from $\tilde u_i$ via the orthogonality
condition $u^\mu\tilde u_\mu=0$ and contributes no independent information
to the norm. Thus, the weighted norm
of such a $4$-vector field is, by definition, precisely the weighted norm
of its spatial part.
\end{rmk}

We then set
\begin{equation}\label{eq:linearized-equation-state-space}
	\calH^{2l}:= H^{2l,l+\frac{1}{2\kappa}-\frac{1}{2}}\times H^{2l,l+\frac{1}{2\kappa}}\times H^{2l,l+\frac{1}{2\kappa}}.
\end{equation}

We have the following local well-posedness result for the linearized system \eqref{eq:Eulerian-Linearized-Equation} from \cite{Zhong_2025}:

\begin{thm}\label{thm:Eulerian-Linearized-Wellposedness}
	Let $l\geq 0$ be an integer and $(r,u,\pi)$ be a sufficiently regular solution to the nonlinear system \eqref{eq:Eulerian-Main-System} on a time interval $[0,T]$ such that $a_0=\frac{\pi}{r}$ is sufficiently regular on $[0,T]$. Then, for any initial data $(\tilde r_0,\tilde u_0,\tilde \pi_0)\in \calH^{2l}$ and any source terms $(\tilde f,\tilde g,\tilde h)\in L^1_t\calH^{2l}$, the linearized system \eqref{eq:Eulerian-Linearized-Equation} has a unique solution $(\tilde r,\tilde u,\tilde \pi)\in L^\infty_t\calH^{2l}$, and it satisfies
    \begin{equation}\label{eq:energy-estimate-modified}
        \norm{(\tilde r,\tilde u,\tilde \pi)}_{\calH^{2l}}\lesssim\calC\left(\norm{(\tilde r_0,\tilde u_0,\tilde \pi_0)}_{\calH^{2l}}+\int_0^t \norm{(\tilde f,\tilde g,\tilde h)}_{\calH^{2l}}ds\right)
    \end{equation}
    for any $t\in[0,T]$, where 
    \[
        \calC=\calC\left(\norm{(\p_t r,\p_t u,\p_t\pi)}_{L^\infty_tW^{2l,\infty}_x},\norm{(r,u,\pi)}_{L^\infty_tW^{2l+1,\infty}_x},\norm{a_0}_{L^\infty_tW^{2l,\infty}_x}\right).
    \]
\end{thm}

\begin{rmk}
The present theorem is an intermediate form of \cite[Theorem~34]{Zhong_2025}. Two
substitutions, both carried out in \cite{Zhong_2025} but omitted here, recover the
exact statement:

\begin{enumerate}
\item[(i)] \textbf{Elimination of time derivatives.} Since we assume $(r,u,\pi)$ solve the nonlinear system \eqref{eq:Eulerian-Main-System}, $\partial_t(r,u,\pi)$ can be solved for algebraically in terms of
$\partial(r,u,\pi)$ via the equations of motion \eqref{eq:Eulerian-Main-System}
(by inverting the coefficient matrix of the time derivatives of the associated
first-order system). Substituting throughout the energy estimate removes
$\norm{\partial_t(r,u,\pi)}_{L^\infty_tW^{2l,\infty}_x}$-dependence from the implicit constant $\mathcal C$
entirely, replacing it with the already-present norm
$\norm{(r,u,\pi)}_{L^\infty_tW^{2l+1,\infty}_x}$.

\item[(ii)] \textbf{Elimination of $a_0$.} Unlike (i), this step is not purely
algebraic: differentiating $a_0=\pi/r$ produces a factor of $r^{-1}$ that is
not controlled by $\norm{(r,u,\pi)}_{W^{2l+1,\infty}}$ alone. In \cite[Lemma~11, Lemma~13]{Zhong_2025}, a bootstrap argument shows that, under the additional assumption
that $\partial^{2k}(\pi_0/r_0)$ is bounded for $0\le k\le l$, the norm
$\norm{a_0}_{L^\infty_tW^{2l,\infty}_x}$ is controlled by
$\norm{(r,u,\pi)}_{L^\infty_tW^{2l+1,\infty}_x}$ and the initial data, at the cost of restricting to a possibly shorter interval $[0,T^*]$.
\end{enumerate}

Combining (i) and (ii) recovers \cite[Theorem~34]{Zhong_2025}, with $\mathcal C$
depending only on $\norm{(r,u,\pi)}_{L^\infty_tW^{2l+1,\infty}_x}$ and the
time interval shrunk to $T^*=T^*(\norm{(r,u,\pi)}_{L^\infty_tW^{1,\infty}_x})$.
Since in the present article we assume directly that $a_0$ is regular on the
full interval $[0,T]$ ($T$ will be chosen later), neither
substitution is necessary, and we work with the intermediate form of the
estimate given above.
\end{rmk}

\subsection{Modified linearized equations}
Theorem~\ref{thm:Eulerian-Linearized-Wellposedness}, which we use to invert the linearized operator, presupposes that $\phi=(r,u,\pi)$ are exact solutions to the nonlinear system \eqref{eq:Eulerian-Main-System}. In the present setting of Section~\ref{sec:Nash-Moser-iteration-scheme}, however, the background states we use are not exact solutions to \eqref{eq:Eulerian-Main-System}. To still be able to invert the linearized operator, we proceed as follows.

From now on, we fix $\vec U$ defined on $[0,T]\times\Omega_0$ and define the variables $(R,U,P,\eta)$ and $(r,u,\pi,\Omega_t)$ in terms of $\vec U$ through \eqref{eq:U^0-defining-eq}, \eqref{eq:eta-defining-eq}, \eqref{eq:R-P-defining} and Remark~\ref{rmk:r-u-pi-defining} as in Section~\ref{sec:Nash-Moser-iteration-scheme}. We also employ the notations in \eqref{eq:variables-abbreviation} for convenience. As a consequence, $\phi=(r,u,\pi)$ satisfies all equations of the system \eqref{eq:Eulerian-Main-System} except possibly the momentum equations \eqref{eq:u-zero-equation} and \eqref{eq:vec-u-equation}: that is, $u$ does not necessarily solve these two equations. Thus, Theorem~\ref{thm:Eulerian-Linearized-Wellposedness} is not directly applicable. To restore solvability, we replace the $\tilde u^0$ component of the linearized system \eqref{eq:Eulerian-Linearized-Equation} by a modified one.

Applying the material derivative $D_t:= u^\mu\p_\mu$ to the orthogonality constraint \eqref{eq:orthogonality-constraint} and substituting the spatial components of \eqref{eq:tilde-u-Equation} yields the modified $\tilde u^0$ equation:
\begin{equation}\label{eq:modified-u0-equation}
    D_t \tilde u^0 + \left(1+\frac{1}{\kappa}\right)a_2\Delta^{0\mu}\p_\mu \tilde r + \left(2+\frac{1}{\kappa}\right)ra_3\Delta^{0\mu}\p_\mu \tilde \pi+ V_2^0 \tilde r + (\bar W_2)_\beta^0\tilde u^\beta  + Z_2^0\tilde \pi = \tilde g^0,
\end{equation}
where
\[
    \tilde g^0:=\frac{u^i}{u^0}\tilde g^i,\qquad (\bar W_2)_\beta^0:=\frac{u^i}{u^0}(W_2)_{i\beta} - \frac{1}{u^0}D_t u_\beta.
\]
We also set $(\bar W_2)_\beta^j := (W_2)^j_\beta$ for convenience. The system \eqref{eq:Eulerian-Linearized-Equation} with the $\tilde u^0$ 
equation replaced by \eqref{eq:modified-u0-equation} is called the modified linearized system.

\begin{rmk}
We record the derivation of \eqref{eq:modified-u0-equation} as follows:

Differentiating $u^\mu\tilde u_\mu = 0$ 
with $D_t$ gives
\[
    D_t\tilde u^0 = \frac{u^i}{u^0}D_t\tilde u_i 
    + \frac{1}{u^0}(D_t u^\mu)\tilde u_\mu.
\]
Substituting $D_t\tilde u^i$ from the spatial part of 
\eqref{eq:tilde-u-Equation} and rearranging yields \eqref{eq:modified-u0-equation}.

The gradient 
terms on the left-hand side in \eqref{eq:modified-u0-equation} carry the factor $\Delta^{0\mu}$. 
This is a consequence of the identity
\[
    \frac{u^i}{u^0}\Delta_i^\mu = \Delta^{0\mu},
\]
which follows directly from $\Delta^{\alpha\mu} 
= g^{\alpha\mu} + u^\alpha u^\mu$ and the normalization 
$u^\mu u_\mu = -1$. An identical identity holds for 
$V_2^\alpha$ and $Z_2^\alpha$ since both carry a factor of $\Delta^{\alpha\mu}$.
\end{rmk}

\begin{rmk}\label{rmk:modified-linearized-eqs-time-derivative}
	The modified linearized system can be solved algebraically for
	$\p_t(\tilde r,\vec{\tilde u},\tilde \pi)$ in terms of
	$\p(\tilde r,\vec{\tilde u},\tilde \pi)$, $(\tilde r,\vec{\tilde u},\tilde \pi)$
	and $(\tilde f,\vec{\tilde g},\tilde h)$. Explicitly,
	\begin{equation}\label{eq:modified-linearized-eqs-time-derivative}
    \left\{
    \begin{aligned}
        \p_t \tilde r &= \calF\p\tilde r + r\calF\p\vec{\tilde u} + r\calF\p\tilde\pi
            + \mfL\tilde r + \mfL\vec{\tilde u} + \mfL\tilde\pi
            + \calF\tilde f + r\calF\vec{\tilde g} + r\calF\tilde h,\\
        \p_t \vec{\tilde u} &= \calF\p\tilde r + \calF\p\vec{\tilde u} + \calF\p\tilde\pi
            + \mfL\tilde r + \mfL\vec{\tilde u} + \mfL\tilde\pi
            + \calF\tilde f + \calF\vec{\tilde g} + \calF\tilde h,\\
        \p_t \tilde \pi &= \calF\p\tilde r + \calF\p\vec{\tilde u} + \calF\p\tilde\pi
            + \mfL\tilde r + \mfL\vec{\tilde u} + \mfL\tilde\pi
            + \calF\tilde f + \calF\vec{\tilde g} + \calF\tilde h,
    \end{aligned}
    \right.
\end{equation}
where $\calF=\calF(r,a_0,\vec u)$ denotes a coefficient smooth in a
neighborhood of the range of $(r,a_0,\vec u)$, and
$\mfL=\sum_X\calF_X(r,a_0,\vec u)\p X + \calF(r,a_0,\vec u)\p_t\vec u+\calF(r,a_0,\vec u)$ for
$X\in\{r,a_0,\vec u\}$. Above, both $\calF$ and $\mfL$ are subject to
change line by line and term by term. Details can be found in
Appendix~\ref{sec:solving-time-derivative}.
\end{rmk}

Let $\tilde\phi:=(\tilde r,\tilde u,\tilde\pi)$ be the linearized variables of the modified equations. We write the modified linearized equations in the abstract form as
\begin{equation}\label{eq:abstract-modified-eqs}
	\calL(\vec\phi)\tilde\phi = \tilde F,\qquad \tilde F:=(\tilde f,\tilde g,\tilde h),
\end{equation}
where $\vec\phi=(r,\vec u,\pi)$ is the background state defined as in \eqref{eq:variables-abbreviation}. Since $\pi=ra_0$, letting $\vec\varphi:=(r,\vec u,a_0)$, the states $\vec\phi$ and $\vec\varphi$ determine one another. We write $\calL(\vec\phi)$ or $\calL(\vec\varphi)$ according to which variables are used to express the same operator, with $\calL(\vec\varphi)$ understood as $\calL(\vec\phi(\vec\varphi))$ (and conversely). When we wish to emphasize the dependence on $\vec\varphi$ -- for instance when invoking the estimates of Lemma~\ref{lem:coefficient-bounds-merged}, which are stated in terms of $\vec\varphi$ -- we write the same operator, acting on the same $\tilde\phi$, as
\begin{equation}\label{eq:abstract-modified-eqs2}
	\calL(\vec\varphi)\tilde\phi = \tilde F.
\end{equation}
Throughout, $\calL(\vec\phi)$ and $\calL(\vec\varphi)$ denote the identical linear operator; only the choice of variables used to express its coefficients differs.

Likewise, we write $\calN(\phi)$, $\vec\calN(\phi)$, $\vec\calN_{\mathrm{mo}}(\phi)$ defined in \eqref{eq:abstract-Eulerian-operator},\eqref{eq:abstract-spatial-Eulerian-operator},\eqref{eq:abstract-Eulerian-momentum-operator} or $\calN(\varphi)$, $\vec\calN(\varphi)$, $\vec\calN_{\mathrm{mo}}(\varphi)$ according to which variables are used to express the same operator, with each understood via composition with the change of variables $\phi=\phi(\varphi)$.

\subsection{Coefficient estimates}
Although $r,\pi,\vec u$ are, in our construction, all obtained simultaneously from the Lagrangian velocity field $\vec U$, the $r$- and $\pi$-equations \eqref{eq:r-equation}, \eqref{eq:pi-equation} themselves can be seen as transport equations along a given velocity field $\vec u$. Therefore, in this subsection we treat $\vec u$ as an independent variable with an associated flow map $\eta$ that satisfies $\eta(t,\Omega_0)=\Omega_t$, and establish bounds on $(r,\pi)$ and hence on the coefficient $a_0=\pi/r$ purely in terms of $\vec u$. The resulting estimates apply, in particular, to the specific $\vec u$ arising from $\vec U$ in our construction. We establish both $L^\infty$ bounds and weighted $L^2$ bounds for $(r,a_0)$ in this subsection.

\begin{lem}\label{lem:basic-Linfty-coefficient-estimate}
	For any $\Lambda>0$, there exists $T=T(\Lambda,N^0_\infty)\in(0,1]$ such that if 
	\begin{equation*}
		\norm{\vec u}_{\tilde Z^1_x}<\Lambda,
	\end{equation*}
	then
	\begin{equation}\label{eq:basic-Linfty-coefficient-estimate}
		\norm{(r,a_0,a_0^{-1})(t,\cdot)}_{L^\infty_x}\leq 2N^0_\infty
	\end{equation}	
	for $t\in[0,T]$.
\end{lem}
\begin{proof}
    We prove this statement by bootstrap. Assume that, for some $T>0$,
    \[\norm{(r,a_0,a_0^{-1})(t,\cdot)}_{L^\infty_x}\leq 4N^0_\infty, \qquad t\in[0,T].\]

    Using \eqref{eq:r-pi-defining}, we compute the evolution of $a_0,a_0^{-1}$:
    \begin{equation}\label{eq:a_0-evolution}
        \begin{aligned}
            D_t a_0 &= \kappa a_0(r+1+r^2a_0^{\frac{1}{\kappa}+2})\p_\mu u^\mu-\hat\zeta(r)(a_0^{-1})^{1+\frac{1}{\kappa}}\p_\mu u^\mu-a_0 - \hat\lambda(r)a_0^{\frac{1}{\kappa}+3},\\
            D_ta_0^{-1} &= \hat\zeta(r)(a_0^{-1})^{\frac{1}{\kappa}+3}\p_\mu u^\mu - \kappa(r+1+r^2a_0^{\frac{1}{\kappa}+2})a_0^{-1}\p_\mu u^\mu + a_0^{-1} + \hat\lambda(r) a_0^{\frac{1}{\kappa}+1}.
        \end{aligned}
    \end{equation}
    Define the vector $Y:=(r,a_0,a_0^{-1})$, and it satisfies a transport-type system of the form
    \begin{equation}\label{eq:Y-transport}
        D_t Y = \calF\left(Y,\vec u, \p_t \vec u, \p\vec u\right), \qquad t\in[0,T],
    \end{equation}
    where $\calF$ is smooth in its arguments and affine in $(\p_t\vec u,\p\vec u)$. Under the bootstrap assumption, there exists a constant $\calC>0$ non-decreasing in its argument such that
    \[\abs{\calF(Y,\vec u,\p_t\vec u,\p\vec u)}\leq \const{N^0_\infty,\Lambda},\qquad t\in[0,T].\]

    Integrating the equation \eqref{eq:Y-transport} along the flow $\eta$ and by a direct estimate, we obtain
    \begin{equation*}
    	\abs{Y}\leq N^0_\infty + T\const{N^0_\infty,\Lambda}.
    \end{equation*}

    To close the bootstrap, we choose
    \[T := \min\left\{ \frac{N^0_\infty}{\const{N^0_\infty,\Lambda}},1\right\},\]
    and thus $\norm{Y(t,\cdot)}_{L^\infty_x}\leq 2N^0_\infty$ for any $t\in[0,T]$, which improves the bootstrap bound and completes the proof.
\end{proof}

\begin{rmk}\label{rmk:persistence-of-causality}
The proof of Lemma~\ref{lem:basic-Linfty-coefficient-estimate} also shows that $r(t,\eta(t,y))\le r_0(y)+T\const{N^0_\infty,\Lambda}$ for $(t,y)\in[0,T]\times\Omega_0$. Hence, if $T\const{N^0_\infty,\Lambda}\le\varepsilon_0$, then by \eqref{eq:uniform-smallness-r0} we have $r\le2\varepsilon_0$ on $[0,T]$, and the causality condition \eqref{eq:causality} holds on $[0,T]$ by the choice of $\varepsilon_0$ in Remark~\ref{rmk:choice-of-epsilon0}.
\end{rmk}

From now on, we fix $\Lambda>0$ and choose $T=T(\Lambda,N^0_\infty)\in(0,1]$ such that \eqref{eq:basic-Linfty-coefficient-estimate} and the conclusion of Remark~\ref{rmk:persistence-of-causality} hold.

\begin{nota}
	For convenience, we introduce new notations: for any integer $k\geq1$
	\begin{equation*}
		\begin{aligned}
			\calC^k_E&:=\const{\norm{\vec u}_{\tilde Z^k_x},N^{k-1}_\infty},\\
			\calC^k_L&:=\const{\norm{\vec U}_{\tilde Z^k_y},N^{k-1}_\infty}.
		\end{aligned}
	\end{equation*}
\end{nota}

The following lemmas record estimates for a generic quantity $q$ transported along the flow of the velocity $\vec u/u^0$. They are used repeatedly below.

\begin{lem}\label{lem:generic-Linfty-transport-estimate}
Let $q,\psi$ be functions on $\sptdomain$ such that
\begin{equation*}
	\norm{q}_{L^\infty_tL^\infty_x} +\norm{\psi}_{L^\infty_tL^\infty_x}<\infty
\end{equation*}
and $q$ solves
\begin{equation}\label{eq:generic-transport}
    \left\{
    \begin{aligned}
        D_t q &= f(q,\psi),\\
        q(0,\cdot)&=q_0,
    \end{aligned}
    \right.
\end{equation}
for some $f=f(q,\psi)$ smooth in a neighborhood of the range of
$(q,\psi)$. Then, for every integer $k\ge0$,
\begin{equation}\label{eq:generic-Linfty-transport-estimate}
    \norm{q}_{W^{k,\infty}}\lesssim\const{\norm{\vec u}_{L^\infty_tW^{k,\infty}_x},\norm{\psi}_{L^\infty_tW^{k,\infty}_x},\norm{q}_{L^\infty_tL^\infty_x},\norm{q_0}_{W^{k,\infty}}}
\end{equation}
for any $t\in[0,T]$.
\end{lem}
 
\begin{proof}
We prove the estimate by induction on $k$. The case $k=0$ is automatic by the assumption. We assume the estimate \eqref{eq:generic-Linfty-transport-estimate} holds for any $k'<k$, and we now show that it holds for $k$. Applying $\p^k$ to the first equation of \eqref{eq:generic-transport} and commuting with $D_t$,
\begin{equation}\label{eq:generic-transport-commuted}
    D_t\p^kq = [u^0\p_t,\p^k]q + [u^i\p_i,\p^k]q + \p^k\big(f(q,\psi)\big).
\end{equation}
We claim the right-hand side can be written as $A(t,x)\p^kq+G(t,x)$,
where
\begin{equation*}
	\begin{aligned}
		\norm{A}_{L^\infty_tL^\infty_x}&\lesssim\const{\norm{\vec u}_{L^\infty_tW^{1,\infty}},\norm{\psi}_{L^\infty_tL^\infty_x},
    \norm{q}_{L^\infty_tL^\infty_x}},\\
    \norm{G}_{L^\infty_tL^\infty_x}&\lesssim\const{\norm{\vec u}_{L^\infty_tW^{k,\infty}},
    \norm{\psi}_{L^\infty_t W^{k,\infty}_x},
    \norm{q}_{L^\infty_tL^\infty_x},\norm{q_0}_{W^{k-1,\infty}}}.
	\end{aligned}
\end{equation*}
For the first commutator term $[u^0\p_t,\p^k]q$, we apply the Leibniz rule and it expands into a sum of terms of the form
\begin{equation*}
	\prod_{i_1=1}^{I_1} \p^{a_{i_1}} \vec u \p^c\p_tq,\quad\sum_{i_1}a_{i_1} + c = k,\quad a_{i_1}\geq 1,\quad I_1\geq 1
\end{equation*}
with coefficients smooth in $\vec u$ omitted. From \eqref{eq:generic-transport}, we can substitute 
\begin{equation*}
	\p_tq = \frac{1}{u^0}\left(f(q,\psi)-u^i\p_iq\right)
\end{equation*}
for $\p_tq$. Expanding by the Leibniz rule again refines each such term into $2$ categories of terms:

\noindent \emph{Type (i)}:
\begin{equation*}
	\prod_{i_1=1}^{I_1}\p^{a_{i_1}}\vec u\prod_{i_2=1}^{I_2}\p^{a_{i_2}}\psi\prod_{i_3=1}^{I_3}\p^{a_{i_3}}q,
	\quad \sum_{i_1} a_{i_1} + \sum_{i_2} a_{i_2} + \sum_{i_3} a_{i_3} = k,
	\quad a_{i_1}, a_{i_2}, a_{i_3} \geq 1,
	\quad I_1 \geq 1, I_2,I_3\geq 0
\end{equation*}
with coefficients smooth in $(\vec u,q,\psi)$ omitted.

\noindent\emph{Type (ii)}:
\begin{equation*}
	\prod_{i_1=1}^{I_1}\p^{a_{i_1}}\vec u\p^cq,
	\quad \sum_{i_1} a_{i_1} + c = k+1,
	\quad a_{i_1},c\geq 1,
	\quad I_1\geq 1
\end{equation*}
with coefficients smooth in $\vec u$ omitted.

Type (i) terms, by a direct estimate with the induction hypothesis, can be absorb into $G$. For Type (ii) terms, if $k$ derivatives land on $q$, then we obtain a term that belongs to $A\p^kq$. The rest of the terms belong to $G$.

The second term $[u^i\p_i,\p^k]q$ expands into a sum of terms of the same form as Type (ii) terms above and thus is treated the same way.

The third term $\p^k(f(q,\psi))$ expands into a sum of terms of the form
\begin{equation*}
	\prod_{i_1=1}^{I_1}\p^{a_{i_1}}\psi\prod_{i_2=1}^{I_2} \p^{a_{i_2}}q,
	\quad \sum_{i_1} a_{i_1}+\sum_{i_2} a_{i_2} = k,
	\quad a_{i_1},a_{i_2}\geq 1,
	\quad I_1,I_2\geq 0
\end{equation*}
with coefficients smooth in $(q,\psi)$ omitted. If $k$ derivatives hit a single $q$, it generates terms belonging to $A\p^kq$. The rest of the terms belong to $G$ with a direct estimate using the induction hypothesis.

Integrating \eqref{eq:generic-transport-commuted} along the flow $\eta$ of $\frac{\vec u}{u^0}$, we obtain
\begin{equation*}
	(\p^kq)(t,\eta(t,y)) = (\p^kq_0)(y) + \int_0^t\left(\frac{1}{u^0}(A\p^kq + G)\right)\circ\eta ds
\end{equation*}
so that
\begin{equation*}
	\abs{(\p^kq)(t,\eta(t,y))}\lesssim \norm{q_0}_{W^{k,\infty}} + \norm{G}_{L^\infty_tL^\infty_x} + \norm{A}_{L^\infty_tL^\infty_x}\int_0^t\abs{(\p^kq)(s,\eta(s,y))}ds.
\end{equation*}
An application of Gr\"onwall's inequality yields
\begin{equation*}
	\abs{\p^kq}\lesssim\const{\norm{\vec u}_{L^\infty_tW^{k,\infty}_x},\norm{\psi}_{L^\infty_tW^{k,\infty}_x},\norm{q}_{L^\infty_tL^\infty_x},\norm{q_0}_{W^{k,\infty}}}
\end{equation*}
for $t\in[0,T]$.
\end{proof}

\begin{lem}\label{lem:generic-weighted-transport-estimate}
Let $q,\psi$ be functions on $\sptdomain$ with
\begin{equation*}
	\norm{q}_{L^\infty_tL^\infty_x}+\norm{\psi}_{L^\infty_tL^\infty_x}<\infty
\end{equation*}
and $q$ solves \eqref{eq:generic-transport}. Then, for any integer $k\ge0$ and $\sigma>-\frac12$ with $k-\sigma>\frac32$, and for $t\in[0,T]$,
\begin{equation}\label{eq:generic-weighted-transport-estimate}
\norm{r^\sigma\p^kq}_{L^2}\lesssim \calC\left(\norm{r_0^\sigma\p^kq_0}_{L^2} + \int_0^t\left(\norm{\vec u(s)}_{H^{k+1,\sigma}} + \norm{r^\sigma(\p^k\psi)(s)}_{L^2}\right)ds\right),
\end{equation}
where $\calC=\const{\norm{\p_t\vec u}_{L^\infty_tL^\infty_x},\norm{\vec u}_{L^\infty_tW^{1,\infty}_x},N^0_\infty,\norm{\psi}_{L^\infty_tL^\infty_x},\norm{q}_{L^\infty_tL^\infty_x}}$.
\end{lem}

\begin{proof}
Applying $\p^k$ to \eqref{eq:generic-transport} and commuting with $D_t$, we obtain \eqref{eq:generic-transport-commuted}. By the expansion of the right-hand side of \eqref{eq:generic-transport-commuted} obtained in the proof of Lemma~\ref{lem:generic-Linfty-transport-estimate} and Corollary~\ref{Cor:Multilinear-Interpolation}, we have
\begin{equation}\label{eq:weighted-material-k-derivative-estimate}
	\norm{r^\sigma D_t\p^kq}_{L^2}\lesssim\calC\left(\norm{r^\sigma\p^kq}_{L^2} + \norm{\vec u}_{H^{k+1,\sigma}}+\norm{r^\sigma\p^k\psi}_{L^2}\right).
\end{equation}
Next, using the Reynolds transport theorem,
\begin{equation*}
	\p_t\norm{r^\sigma\p^kq}_{L^2}^2
    = \int_{\Omega_t}\frac1{u^0}D_t\big(r^{2\sigma}\abs{\p^kq}^2\big)\,dx
    + \int_{\Omega_t}\p_i\Big(\frac{u^i}{u^0}\Big)r^{2\sigma}\abs{\p^kq}^2\,dx
    =: (i)+(ii).
\end{equation*}

The second term is bounded directly:

\begin{equation*}
	(ii)\lesssim\calC\norm{r^\sigma\p^kq}_{L^2}^2.
\end{equation*}
For the first term, expanding the material derivative,
\begin{equation*}
	(i) = \int_{\Omega_t}\frac{2\sigma}{u^0}r^{2\sigma-1}(D_tr)\abs{\p^kq}^2\,dx
    + \int_{\Omega_t}\frac2{u^0}r^{2\sigma}\p^kq\cdot D_t\p^kq\,dx.
\end{equation*}
Using \eqref{eq:r-pi-defining} to trade $D_tr$ for one more power of $r$, the first term can be bounded with Lemma~\ref{lem:basic-Linfty-coefficient-estimate}:
\begin{equation*}
	\const{\norm{\p_t\vec u}_{L^\infty_tL^\infty_x},\norm{\vec u}_{L^\infty_tW^{1,\infty}_x},N^0_\infty}\norm{r^\sigma\p^kq}_{L^2}^2.
\end{equation*}
The second term is bounded by
Cauchy--Schwarz and \eqref{eq:weighted-material-k-derivative-estimate}:
\begin{equation*}
	    \lesssim \norm{r^\sigma\p^kq}_{L^2}\norm{r^\sigma D_t\p^kq}_{L^2}
    \lesssim \calC\left(\norm{r^\sigma\p^kq}_{L^2} + \norm{\vec u}_{H^{k+1,\sigma}}+\norm{r^\sigma\p^k\psi}_{L^2}\right)\norm{r^\sigma\p^kq}_{L^2}.
\end{equation*}
Combining,
\begin{equation*}
	\p_t\norm{r^\sigma\p^kq}_{L^2}\lesssim\calC\left(\norm{r^\sigma\p^kq}_{L^2} + \norm{\vec u}_{H^{k+1,\sigma}}+\norm{r^\sigma\p^k\psi}_{L^2}\right).
\end{equation*}
Integrating in time and applying Gr\"onwall's inequality yields
\begin{equation*}
	\norm{r^\sigma\p^kq}_{L^2}\lesssim\calC\left(\norm{r_0^\sigma\p^kq_0}_{L^2} + \int_0^t\left(\norm{\vec u(s)}_{H^{k+1,\sigma}} + \norm{r^\sigma(\p^k\psi)(s)}_{L^2}\right)ds\right).
\end{equation*}
\end{proof}

\begin{rmk}
	Although Lemma~\ref{lem:generic-Linfty-transport-estimate} and Lemma~\ref{lem:generic-weighted-transport-estimate} are stated and proved for scalar functions $q,\psi$, the same proofs apply verbatim when $q,\psi$ are vector- or matrix-valued, by writing these functions component-wise, carrying out the same computation on each component and taking the sum over all components.
\end{rmk}

A variant of Lemma~\ref{lem:generic-Linfty-transport-estimate} and Lemma~\ref{lem:generic-weighted-transport-estimate} is often useful. The proof follows from similar arguments as the lemmas above.

\begin{lem}\label{lem:generic-linear-transport-estimate}
Let $q,\tilde q,\tilde\psi$ be functions defined on $\sptdomain$ such that
$\norm{q}_{L^\infty_tL^\infty_x}<\infty$ and $\tilde q,\tilde\psi$
satisfy
\begin{equation}\label{eq:generic-linear-transport}
    \left\{
    \begin{aligned}
        D_t\tilde q &= f(q)\,\tilde q + \tilde\psi,\\
        \tilde q(0,\cdot) &= \tilde q_0,
    \end{aligned}
    \right.
\end{equation}
for some $f=f(q)$ smooth in a neighborhood of the range of its
argument. Then:
 
\begin{enumerate}
\item[(a)] for any integer $k\ge0$ and $\sigma>-\frac12$,
\[
    \norm{\tilde q}_{H^{k,\sigma}}
    \lesssim \const{\norm{\vec u}_{L^\infty_tW^{k,\infty}_x},
    \norm{\p_t\vec u}_{L^\infty_tL^\infty_x},N^0_\infty,
    \norm{q}_{L^\infty_tW^{k,\infty}_x}}
    \Big(\norm{\tilde q_0}_{H^{k,\sigma}_{r_0}}
    +\int_0^t\norm{\tilde\psi}_{H^{k,\sigma}}\,ds\Big);
\]
\item[(b)] if $k-\sigma>\frac32$ and $\tilde q_0=0$,
\begin{equation*}
	\begin{aligned}
		\norm{r^\sigma\p^k\tilde q}_{L^2}&\lesssim\calC\Big(\int_0^t\int_0^s\norm{\tilde\psi(\tau)}_{L^\infty_x}d\tau\left(\norm{r^\sigma \p^k\vec u}_{L^2} + \norm{r^\sigma\p^{k+1}\vec u}_{L^2} + \norm{r^\sigma\p^k q}_{L^2}\right)ds\\
		&\qquad + \int_0^t\norm{r^\sigma\p^k\tilde\psi}_{L^2} + \norm{\tilde\psi}_{L^\infty_x}\norm{r^\sigma\p^k\vec u}_{L^2}ds\Big),
	\end{aligned}
\end{equation*}
where $\calC=\const{\norm{\vec u}_{L^\infty_tW^{1,\infty}_x},
\norm{\p_t\vec u}_{L^\infty_tL^\infty_x},N^0_\infty,
\norm{q}_{L^\infty_tL^\infty_x}}$.
\end{enumerate}
\end{lem}

\begin{proof}
	Applying $\p^k$ to \eqref{eq:generic-linear-transport} and commuting with $D_t$, we have
	\begin{equation*}
	\begin{aligned}
				D_t\p^k\tilde q &= [D_t,\p^k]\tilde q + \p^k(f(q)\tilde q) + \p^k\tilde\psi\\
				&=[u^0\p_t,\p^k]\tilde q + [u^i\p_i,\p^k]\tilde q + \p^k(f(q)\tilde q) + \p^k\tilde\psi.
	\end{aligned}
	\end{equation*}

	For the first term, we replace $\p_t\tilde q$ by the equation \eqref{eq:generic-linear-transport}, i.e., substituting
	\begin{equation*}
		\p_t\tilde q = \frac{1}{u^0}\left(f(q)\tilde q + \tilde\psi - u^i\p_i\tilde q\right). 
	\end{equation*}
	By the Leibniz rule, it gives rise to $3$ types of terms:

	\noindent \emph{Type (i)}:
	\begin{equation*}
		\prod_{i_1=1}^{I_1}\p^{a_{i_1}}\vec u\prod_{i_2=1}^{I_2}\p^{a_{i_2}}q\p^c\tilde q,
		\quad \sum_{i_1} a_{i_1} + \sum_{i_2} a_{i_2} + c = k,
		\quad a_{i_1}, a_{i_2} \geq 1,
		\quad I_1 \geq 1, I_2\geq 0
	\end{equation*}
	with coefficients smooth in $(\vec u,q)$ omitted.

	\noindent\emph{Type (ii)}:
	\begin{equation*}
		\prod_{i_1=1}^{I_1}\p^{a_{i_1}}\vec u\p^c\tilde\psi,
		\quad \sum_{i_1} a_{i_1} + c = k,
		\quad a_{i_1}\geq 1,
		\quad I_1\geq 1
	\end{equation*}
	with coefficients smooth in $\vec u$ omitted.

	\noindent\emph{Type (iii)}:
	\begin{equation*}
		\prod_{i_1=1}^{I_1}\p^{a_{i_1}}\vec u\p^c\tilde q,
		\quad \sum_{i_1} a_{i_1} + c = k+1,
		\quad a_{i_1},c\geq 1,
		\quad I_1\geq 1
	\end{equation*}
	with coefficients smooth in $\vec u$ omitted.

	Expanding the second term $[u^i\p_i,\p^k]\tilde q$ yields terms that can be absorbed into \emph{Type (iii)}. The third term $\p^k(f(q)\tilde q)$ gives rise to:

	\noindent\emph{Type (iv)}:
	\begin{equation*}
		\prod_{i_1=1}^{I_1}\p^{a_{i_1}}q\p^c\tilde q,
		\quad \sum_{i_1} a_{i_1} + c = k,
		\quad a_{i_1}\geq 1,
		\quad I_1\geq 0
	\end{equation*}
	with coefficients smooth in $q$ omitted.

	Combining the expansion above and by a direct estimate, we have
	\begin{equation*}
		\norm{r^\sigma D_t\p^k\tilde q}_{L^2}\lesssim\const{\norm{\vec u}_{L^\infty_t W^{k,\infty}_x},\norm{q}_{L^\infty_tW^{k,\infty}}}\norm{\tilde q}_{H^{k,\sigma}} + \norm{r^\sigma \p^k\tilde \psi}_{L^2}.
	\end{equation*}

	By the Reynolds transport theorem,
	\begin{equation*}
		\begin{aligned}
			\p_t\norm{r^\sigma \p^k\tilde q}_{L^2}^2 &= \int_{\Omega_t}\frac{1}{u^0}D_t\left(r^{2\sigma}\abs{\p^k\tilde q}^2\right)dx + \int_{\Omega_t}\p_i\left(\frac{u^i}{u^0}\right)r^{2\sigma}\abs{\p^k\tilde q}^2dx\\
			&=\int_{\Omega_t}\frac{2}{u^0}r^{2\sigma}\p^k\tilde q\cdot D_t\p^k\tilde qdx + \int_{\Omega_t}\frac{2\sigma}{u^0}r^{2\sigma-1}D_t r \abs{\p^k\tilde q}^2dx + \int_{\Omega_t}\p_i\left(\frac{u^i}{u^0}\right)r^{2\sigma}\abs{\p^k\tilde q}^2dx\\
			&\lesssim \const{\norm{\vec u}_{L^\infty_x}}\norm{r^\sigma\p^k\tilde q}_{L^2}\norm{r^\sigma D_t\p^k\tilde q}_{L^2} + \calC^1_E\norm{r^\sigma\p^k\tilde q}_{L^2}^2.
		\end{aligned}
	\end{equation*}
	
	Combining the estimates above, we obtain
	\begin{equation*}
		\p_t\norm{r^\sigma\p^k\tilde q}_{L^2}\lesssim\const{\norm{\vec u}_{L^\infty_tW^{k,\infty}_x},\norm{\p_t\vec u}_{L^\infty_tL^\infty_x},N^0_\infty,\norm{q}_{L^\infty_tW^{k,\infty}_x}}\norm{\tilde q}_{H^{k,\sigma}}+\norm{r^\sigma\p^k\tilde\psi}_{L^2}.
	\end{equation*}
	By the same argument applied to \eqref{eq:generic-linear-transport} without differentiating the equation, we obtain
	\begin{equation*}
		\p_t\norm{r^\sigma\tilde q}_{L^2}\lesssim\const{\norm{\vec u}_{L^\infty_tW^{1,\infty}_x},\norm{\p_t\vec u}_{L^\infty_tL^\infty_x},N^0_\infty,\norm{q}_{L^\infty_tL^\infty_x}}\norm{r^\sigma\tilde q}_{L^2}+\norm{r^\sigma\tilde\psi}_{L^2}.
	\end{equation*}
	Integrating above in time and applying Corollary~\ref{cor:weighted-L2-end-point-equivalence}, we have
	\begin{equation*}
		\begin{aligned}
			\norm{\tilde q}_{H^{k,\sigma}}&\lesssim \norm{r^\sigma\tilde q}_{L^2} + \norm{r^\sigma\p^k\tilde q}_{L^2}\\
			&\lesssim \norm{\tilde q_0}_{H^{k,\sigma}_{r_0}} + \int_0^t\const{\norm{\vec u}_{L^\infty_tW^{k,\infty}_x},\norm{\p_t\vec u}_{L^\infty_tL^\infty_x},N^0_\infty,\norm{q}_{L^\infty_tW^{k,\infty}_x}}\norm{\tilde q}_{H^{k,\sigma}}+\norm{\tilde \psi}_{H^{k,\sigma}}ds.
		\end{aligned}
	\end{equation*}
	By Gr\"onwall's inequality,
	\begin{equation*}
		\norm{\tilde q}_{H^{k,\sigma}}\lesssim\const{\norm{\vec u}_{L^\infty_tW^{k,\infty}_x},\norm{\p_t\vec u}_{L^\infty_tL^\infty_x},N^0_\infty,\norm{q}_{L^\infty_tW^{k,\infty}_x}}\left(\norm{\tilde q_0}_{H^{k,\sigma}_{r_0}}+\int_0^t\norm{\tilde\psi}_{H^{k,\sigma}}ds\right).
	\end{equation*}

	Now, if $k-\sigma>3/2$, by Corollary~\ref{Cor:Multilinear-Interpolation}, we have
	\begin{equation*}
		\begin{aligned}
			\norm{r^\sigma D_t\p^k\tilde q}_{L^2}&\lesssim \const{\norm{\vec u}_{L^\infty_tW^{1,\infty}_x},\norm{q}_{L^\infty_tL^\infty_x}}\\
			&\qquad\times\left(\norm{r^\sigma\p^k\tilde q}_{L^2} + \norm{\tilde q}_{L^\infty}\left(\norm{r^\sigma\p^{k+1}\vec u}_{L^2} + \norm{r^\sigma\p^k\vec u}_{L^2} + \norm{r^\sigma\p^k q}_{L^2}\right)\right) \\
			&\qquad+ \const{\norm{\vec u}_{L^\infty_tL^\infty_x}}\left( \norm{r^\sigma\p^k\tilde\psi}_{L^2}+\norm{\tilde\psi}_{L^\infty}\norm{r^\sigma\p^k\vec u}_{L^2}\right).
		\end{aligned}
	\end{equation*}
	
	We derive an $L^\infty$ bound on $\tilde q$ as follows: consider
	\begin{equation*}
		\p_t(\tilde q\circ\eta) = (\frac{1}{u^0}D_t\tilde q)\circ\eta= \left(\frac{1}{u^0}(f(q)\tilde q+\tilde\psi)\right)\circ\eta,
	\end{equation*}
	where the second equality follows from \eqref{eq:generic-linear-transport}. Integrating in time, and by a direct estimate,
	\begin{equation*}
		\abs{\tilde q\circ\eta}\leq \underbrace{\norm{\tilde q_0}_{L^\infty}}_{=0} + \const{\norm{q}_{L^\infty_tL^\infty_x}}\int_0^t\abs{\tilde q\circ\eta} + \abs{\tilde\psi\circ\eta}ds.
	\end{equation*}
	By Gr\"onwall's inequality,
	\begin{equation*}
		\abs{\tilde q\circ\eta}\leq \const{\norm{q}_{L^\infty_tL^\infty_x}}\int_0^t\norm{\tilde\psi}_{L^\infty_x}ds.
	\end{equation*}
	Taking $L^\infty$ norm, we have
	\begin{equation*}
		\norm{\tilde q}_{L^\infty_x}=\norm{\tilde q\circ\eta}_{L^\infty_y}\leq \const{\norm{q}_{L^\infty_tL^\infty_x}}\int_0^t\norm{\tilde\psi}_{L^\infty_x}ds.
	\end{equation*}
	Combining the above estimates, we have
	\begin{equation*}
		\begin{aligned}
			\p_t\norm{r^\sigma \p^k\tilde q}_{L^2}&\lesssim \const{\norm{\vec u}_{L^\infty_x}}\norm{r^\sigma D_t\p^k\tilde q}_{L^2} + \calC^1_E\norm{r^\sigma\p^k\tilde q}_{L^2}\\
			&\lesssim\const{\norm{\vec u	}_{L^\infty_tW^{1,\infty}_x},\norm{\p_t\vec u}_{L^\infty_tL^\infty_x},N^0_\infty,\norm{q}_{L^\infty_tL^\infty_x}}\\
			&\qquad\times \left(\norm{r^\sigma\p^k\tilde q}_{L^2}+\int_0^t\norm{\tilde\psi(\tau)}_{L^\infty_x}d\tau\left(\norm{r^\sigma\p^k\vec u}_{L^2} + \norm{r^\sigma\p^{k+1}\vec u}_{L^2} + \norm{r^\sigma\p^k q}_{L^2}\right)\right)\\
			&\qquad + \const{\norm{\vec u}_{L^\infty_tL^\infty_x}}\left(\norm{r^\sigma\p^k\tilde\psi}_{L^2}+\norm{\tilde\psi}_{L^\infty_x}\norm{r^\sigma\p^k\vec u}_{L^2}\right).
		\end{aligned}
	\end{equation*}
	Integrate in time and applying Gr\"onwall's inequality, since $\tilde q_0=0$, we have
	\begin{equation*}
	\begin{aligned}
		\norm{r^\sigma\p^k\tilde q}_{L^2}&\lesssim\calC\Big(\int_0^t\int_0^s\norm{\tilde\psi(\tau)}_{L^\infty_x}d\tau\left(\norm{r^\sigma \p^k\vec u}_{L^2} + \norm{r^\sigma\p^{k+1}\vec u}_{L^2}+\norm{r^\sigma\p^k q}_{L^2}\right)ds\\
		&\qquad + \int_0^t\norm{r^\sigma\p^k\tilde\psi}_{L^2} + \norm{\tilde\psi}_{L^\infty_x}\norm{r^\sigma\p^k\vec u}_{L^2}ds\Big).
	\end{aligned}
\end{equation*}
\end{proof}

\begin{lem}\label{lem:coefficient-bounds-merged}
For any $t\in[0,T]$,
\begin{enumerate}
\item[(a)] for any integer $k\geq 0$,
\begin{equation}\label{eq:higher-Linfty-coefficient-estimate}
	\norm{(r,a_0)(t,\cdot)}_{W^{k,\infty}_x}
    \leq \calC^{k+1}_E;
\end{equation}
\item[(b)] for any integer $l>\frac{1}{2\kappa}+1$,
\begin{equation}\label{eq:weighted-sobolev-coefficient-estimate}
	\norm{r}_{H^{2l,l+\frac{1}{2\kappa}-\frac12}}+\norm{\frac{\pi}{r}}_{H^{2l,l+\frac{1}{2\kappa}-\frac12}}
    \lesssim \calC^1_E
    \Big(N^l_w+\norm{\vec u}_{\tilde X^{l+2,\frac{1}{2\kappa}}_r}\Big).
\end{equation}
\end{enumerate}
\end{lem}
 
\begin{proof}

Write $Y:=(r,a_0)$. From \eqref{eq:r-pi-defining} and \eqref{eq:a_0-evolution}, $Y$ solves a transport equation of the
form
\begin{equation}\label{eq:Y-transport2}
	D_tY = \calF(Y,\vec u,\p_t\vec u,\p\vec u).
\end{equation}
By Lemma~\ref{lem:generic-Linfty-transport-estimate}, we obtain
\begin{equation*}
	\norm{Y}_{W^{k,\infty}_x}\lesssim\calC^{k+1}_E.
\end{equation*}
By Lemma~\ref{lem:generic-weighted-transport-estimate}, Lemma~\ref{lem:basic-Linfty-coefficient-estimate} and Lemma~\ref{lem:WeightedEmbedding},

\begin{equation*}
	\norm{r^{l+\frac{1}{2\kappa}-\frac12}\p^{2l} Y}_{L^2_x}\lesssim\calC^1_E\left(N^l_w + \norm{\vec u}_{\tilde X^{l+2,\frac{1}{2\kappa}}_r}\right).
\end{equation*}
By the exact form of the equation \eqref{eq:Y-transport2} and Lemma~\ref{lem:generic-linear-transport-estimate}, we have
\begin{equation*}
	\norm{r^{l+\frac{1}{2\kappa}-\frac12}Y}_{L^2_x}\lesssim\calC^1_E N^l_w.
\end{equation*}
Combining above using Corollary~\ref{cor:weighted-L2-end-point-equivalence},
\begin{equation*}
	\begin{aligned}
		\norm{Y}_{H^{2l,l+\frac{1}{2\kappa}-\frac12}}&\lesssim \norm{r^{l+\frac{1}{2\kappa}-\frac12}Y}_{L^2_x} + \norm{r^{l+\frac{1}{2\kappa}-\frac12}\p^{2l} Y}_{L^2_x}\\
		&\lesssim\calC^1_E\left(N^l_w + \norm{\vec u}_{\tilde X^{l+2,\frac{1}{2\kappa}}_r}\right),
	\end{aligned}
\end{equation*}
which completes the proof.
\end{proof}

\begin{rmk}\label{rmk:higher-Linfty-time-derivative-bounds}
    We can use \eqref{eq:Y-transport2} to solve $\p_t Y$ in terms of the spatial derivatives of $Y$ and $\vec u$, and thus have the estimates for the time derivatives
    \[\norm{(\p_tr,\p_ta_0)(t,\cdot)}_{W^{k,\infty}}\lesssim \calC^{k+2}_E\]
    for any $t\in[0,T]$ and integer $k\geq 0$.
\end{rmk}

\subsection{Energy estimate and the invertibility of the linearization of the iteration map over the Eulerian domain}
The modified system differs from \eqref{eq:Eulerian-Linearized-Equation} only through the replacement of the zeroth-order coefficient $(W_2)^0_\beta$ by $(\bar W_2)^0_\beta$ in the $\tilde u^0$ equation, and no other coefficients, in particular no top-order coefficients, are affected. In the proof of \cite[Theorem~15]{Zhong_2025}, the term $(W_2)^\alpha_\gamma\tilde u^\gamma\tilde u_\alpha$ is estimated by a direct estimate and plays no role in the special cancellation used to control $\partial_\mu r\,\tilde u^\mu$, and therefore $(\bar W_2)^0_\beta\tilde u^\beta\tilde u^0$ can be dealt with the same way by a direct estimate after the substitution. Hence the basic energy estimate of \cite[Theorem~15]{Zhong_2025}, the higher-order estimate of \cite[Theorem~33]{Zhong_2025}, and the duality argument of \cite[Theorem~34]{Zhong_2025} carry over verbatim with $(W_2)^0_\beta$ replaced by $(\bar W_2)^0_\beta$ throughout. So, we have the exact same local well-posedness and energy estimates as in Theorem~\ref{thm:Eulerian-Linearized-Wellposedness} for the modified linearized equations. An application of Lemma~\ref{lem:coefficient-bounds-merged} yields the following theorem for the modified system.

\begin{thm}\label{thm:Modified-linearized-equations-wellposedness}
For any integer $l\geq 0$, any initial data $(\tilde r_0,\tilde u_0,\tilde \pi_0)\in\calH^{2l}$ and any source terms $(\tilde f,\tilde g,\tilde h)\in L^1_t\calH^{2l}$, the modified
linearized equations have a unique solution $(\tilde r,\tilde u,\tilde \pi)\in L^\infty_t\calH^{2l}$, and it satisfies
\begin{equation}
    \norm{(\tilde r,\tilde u,\tilde\pi)}_{\calH^{2l}}
    \lesssim \calC^{2l+2}_E\Big(\norm{(\tilde r_0,\tilde u_0,\tilde\pi_0)}_{\calH^{2l}}
    +\int_0^t\norm{(\tilde f,\tilde g,\tilde h)}_{\calH^{2l}}\,ds\Big),
    \qquad t\in[0,T].
\end{equation}
\end{thm}

\begin{rmk}\label{rmk:estimates-only-involve-spatial}
	Since the norm $\norm{\cdot}_{\calH^{2l}}$ depends on $\tilde u$ only
	through the spatial part (see Remark~\ref{rmk:u-0-play-no-role}), the
	estimate of Theorem~\ref{thm:Modified-linearized-equations-wellposedness},
	although stated in terms of $(\tilde r,\tilde u,\tilde\pi)$ and
	$(\tilde f,\tilde g,\tilde h)$, in fact involves only the spatial
	components $(\tilde r,\vec{\tilde u},\tilde\pi)$ and $(\tilde f,\vec {\tilde g},\tilde h)$.
\end{rmk}

Now we return to the system \eqref{eq:modified-linearized-equation-with-xi}. Although the system does not explicitly prescribe an equation for $\delta u^0$, one may derive such an equation exactly as $\tilde u^0$'s equation \eqref{eq:modified-u0-equation} was derived for the modified linearized equations: define $\delta u^0:=\frac{u^i}{u^0}\delta u_i$ via the orthogonality constraint, and take the material derivative $D_t$ of this identity, substituting the equation already satisfied by $\delta\vec u$. With $\delta u^0$'s equation so obtained, the equations satisfied by $(\delta r,\delta u^0,\delta\vec u,\delta\pi)$ coincide term-by-term with the modified linearized equations for $(\tilde r,\tilde u,\tilde\pi)$, with the source terms $(\tilde f,\tilde g,\tilde h)$ replaced by
\begin{equation}\label{eq:modified-linearized-equations-with-xi-source}
	\Big(0,\frac{u^i}{u^0}\big(\vec g_i-\p_j(\vec\calN_{\mathrm{mo}})_i(\vec \varphi)\xi^j\big),
    \ \vec g-\p_j\vec\calN_{\mathrm{mo}}(\vec \varphi)\xi^j,\ 0\Big).
\end{equation}
The system \eqref{eq:modified-linearized-equation-with-xi} is thus the modified linearized equations themselves, supplemented with an additional transport equation for $\xi$.

\begin{thm}\label{thm:modified-linearized-with-xi-local-wellposedness}
    For any integer $l\geq 0$ and $\vec g\in L^1_tH^{2l,l+\frac{1}{2\kappa}}$, the linearized system \eqref{eq:modified-linearized-equation-with-xi} has a unique solution $(\delta r,\delta u,\delta \pi, \xi)\in L^\infty_t\left(\calH^{2l}\times H^{2l,l+\frac{1}{2\kappa}}\right)$, and it satisfies 
    \begin{equation}\label{eq:Linfty-energy-bound-velocity-variation}
        \norm{(\delta r,\delta \vec u,\delta \pi,\xi)}_{\calH^{2l}\times H^{2l,l+\frac{1}{2\kappa}}}\lesssim\calC^{2l+3}_E\int_0^t \norm{\vec g}_{H^{2l,l+\frac{1}{2\kappa}}}ds
    \end{equation}
    for any $t\in[0,T]$.
\end{thm}
\begin{proof}
	It is sufficient to close the estimate. By Lemma~\ref{lem:coefficient-bounds-merged}, Remark~\ref{rmk:higher-Linfty-time-derivative-bounds} and Theorem~\ref{thm:Modified-linearized-equations-wellposedness}, we have
	\begin{equation}\label{eq:basic-Linfty-xi-estimate-pre1}
		\begin{aligned}
			\norm{(\delta r,\delta \vec u,\delta \pi)}_{\calH^{2l}}&\lesssim\calC^{2l+2}_E\int_0^t \norm{\vec g-\p_j\vec\calN_{\mathrm{mo}}(\varphi)\xi^j}_{H^{2l,l+\frac{1}{2\kappa}}}ds\\
			&\lesssim \calC^{2l+2}_E\int_0^t \norm{\vec g}_{H^{2l,l+\frac{1}{2\kappa}}}ds + \calC^{2l+3}_E\int_0^t\norm{\xi}_{H^{2l,l+\frac{1}{2\kappa}}}ds.
		\end{aligned}
	\end{equation}

	By Lemma~\ref{lem:generic-weighted-transport-estimate}, we have
	\begin{equation}\label{eq:basic-Linfty-xi-estimate-pre2}
		\norm{\xi}_{H^{2l,l+\frac{1}{2\kappa}}}\lesssim\const{\norm{\vec u}_{L^\infty_tW^{2l+1,\infty}_x},\norm{\p_t\vec u}_{L^\infty_tL^\infty_x},N^0_\infty}\int_0^t\norm{\delta \vec u}_{H^{2l,l+\frac{1}{2\kappa}}}ds.
	\end{equation}
	Summing \eqref{eq:basic-Linfty-xi-estimate-pre1} and \eqref{eq:basic-Linfty-xi-estimate-pre2}, we have
	\begin{equation*}
		\begin{aligned}
			\norm{(\delta r,\delta \vec u,\delta \pi)}_{\calH^{2l}}+\norm{\xi}_{H^{2l,l+\frac{1}{2\kappa}}}&\lesssim \calC^{2l+2}_E\int_0^t \norm{\vec g}_{H^{2l,l+\frac{1}{2\kappa}}}ds + \calC^{2l+3}_E\int_0^t\left(\norm{\delta \vec u}_{H^{2l,l+\frac{1}{2\kappa}}}+\norm{\xi}_{H^{2l,l+\frac{1}{2\kappa}}}\right)ds\\
			&\lesssim \calC^{2l+2}_E\int_0^t \norm{\vec g}_{H^{2l,l+\frac{1}{2\kappa}}}ds + \calC^{2l+3}_E\int_0^t\left(\norm{(\delta r,\delta \vec u,\delta \pi)}_{\calH^{2l}}+\norm{\xi}_{H^{2l,l+\frac{1}{2\kappa}}}\right)ds.
		\end{aligned}
	\end{equation*}
	By Gr\"onwall's inequality,
	\begin{equation*}
		\norm{(\delta r,\delta \vec u,\delta \pi)}_{\calH^{2l}}+\norm{\xi}_{H^{2l,l+\frac{1}{2\kappa}}}\lesssim \calC^{2l+3}_E\int_0^t \norm{\vec g}_{H^{2l,l+\frac{1}{2\kappa}}}ds.
	\end{equation*}
\end{proof}

\section{First-order tame estimate of the linearization of the iteration map over the Eulerian domain}\label{sec:tame-estimate-linearization-Eulerian}

The estimate in Theorem~\ref{thm:modified-linearized-with-xi-local-wellposedness} involves a loss of derivatives through the dependence of the constant $\calC^{2l+3}_E$ on high norms of the background state. This dependence is not tame. We therefore establish a refined estimate for the modified linearized system, in which the dependence on the background coefficients is of tame type. This is the estimate required by the Nash--Moser iteration defined in Section~\ref{sec:Nash-Moser-iteration-scheme}.

We first introduce the weighted elliptic operators used throughout this section. Applying the material derivative $D_t$ to \eqref{eq:tilde-r-equation} and substituting \eqref{eq:tilde-u-Equation}, and vice versa, we obtain at leading order
\[
    D_t^2\tilde r\sim\tilde L_1\tilde r,\qquad D_t^2\tilde u\sim\tilde L_2\tilde u,
\]
where the second-order operators $\tilde L_1$ and $\tilde L_2$ are defined by
\begin{subequations}
    \begin{align}
        \tilde{L}_1 \tilde{r}&=(\kappa+1) H^{i j}\left(r \partial_i \partial_j \tilde{r}+\frac{1}{\kappa} \partial_i r \partial_j \tilde{r}\right),\label{op:weighted-laplacian}\\
        \left(\tilde{L}_2 \tilde{u}\right)_\alpha&=(\kappa+1) A_\alpha^i H^{j k}\left(\partial_i\left(r \partial_j \tilde{u}_k\right)+\frac{1}{\kappa} \partial_j r \partial_i \tilde{u}_k\right),\label{op:weighted-div}
    \end{align}
\end{subequations}
with
\begin{equation*}
    A_\alpha^i=g_\alpha^i-\frac{u^i}{u^0} g_\alpha^0.
\end{equation*}
The operator $\tilde L_1$ is a weighted Laplacian acting on the scalar $\tilde r$, and $\tilde L_2$ is a weighted divergence operator acting on the spatial components of a $4$-vector field $\tilde u$ orthogonal to $u$. Since $\tilde L_2$ controls only the divergence of $\tilde u$, we pair it with the operator
\begin{equation}\label{op:weighted-curl}
    \left(\tilde{L}_3 \tilde{u}\right)_\alpha=(\kappa+1) A_\alpha^i H^{jk} r^{-\frac{1}{\kappa}} \partial_j\left(r^{\frac{1}{\kappa}+1}\left(\partial_k \tilde{u}_i-\partial_i \tilde{u}_k\right)\right),
\end{equation}
associated with the curl, so that $\tilde L_2+\tilde L_3$ is elliptic on the spatial components of $\tilde u$. Since $u^\alpha A^i_\alpha=0$, the outputs $\tilde L_2\tilde u$ and $\tilde L_3\tilde u$ are again $4$-vector fields orthogonal to $u$.

The variable $\tilde\pi$ is measured in the same weighted space as $\tilde u$, with half a power of $r$ more than $\tilde r$, and in this space $\tilde L_1$ is not coercive. We therefore use for $\tilde\pi$ the scalar operator
\begin{equation}\label{op:weighted-laplacian-pi}
    \tilde L_4\tilde\pi=(\kappa+1)H^{ij}\Big(r\p_i\p_j\tilde\pi+\Big(1+\frac1\kappa\Big)\p_ir\,\p_j\tilde\pi\Big).
\end{equation}

We commute the modified linearized equation \eqref{eq:abstract-modified-eqs2} with the elliptic operator defined by
\begin{equation}\label{Op:VecEllipticOperator}
    L :=
    \begin{pmatrix}
        \tilde L_1 & 0 & 0\\
        0 & (\tilde L_2+\tilde L_3)& 0\\
        0 & 0 &\tilde L_4
        \end{pmatrix},
\end{equation}
i.e., $L(\tilde r, \tilde u, \tilde\pi) := 
(\tilde L_1 \tilde r,\, (\tilde L_2+\tilde L_3)\tilde u,\, \tilde L_4\tilde\pi)$. This gives
\[
    \calL(\vec \phi)L\tilde\phi = L\tilde F + [\calL(\vec\phi),L]\tilde\phi.
\]

Written out componentwise, this equation reads:
\begin{subequations}\label{eq:commuted-system}
\begin{align}
D_t \tilde L_1\tilde r
&+ \kappa r a_1 \p_\mu [(\tilde L_2+\tilde L_3)\tilde u]^\mu + \p_\mu r [(\tilde L_2+\tilde L_3)\tilde u]^\mu
+ V_1\tilde L_1\tilde r
+ r^2Z_1\tilde L_4\tilde\pi
=
\tilde L_1 \tilde f + \calR_{\tilde r},
\\
D_t [(\tilde L_2+\tilde L_3)\tilde u]^\alpha
&+ (1+\tfrac1\kappa)a_2\Delta^{\alpha\mu}\p_\mu \tilde L_1\tilde r
+ (2+\tfrac1\kappa)ra_3\Delta^{\alpha\mu}\p_\mu \tilde L_4\tilde \pi \notag
\\
&\quad+ (V_2)^\alpha\tilde L_1\tilde r
+ (\bar W_2)^\alpha_\beta[(\tilde L_2+\tilde L_3)\tilde u]^\beta
+ (Z_2)^\alpha \tilde L_4\tilde \pi
=
[(\tilde L_2+\tilde L_3)\tilde g]^\alpha + \calR_{\tilde u}^\alpha,
\\
D_t\tilde L_4 \tilde \pi
&+ r a_4\p_\mu[(\tilde L_2+\tilde L_3)\tilde u]^\mu
+ V_3 \tilde L_1\tilde r
+ (W_3)_\alpha[(\tilde L_2+\tilde L_3)\tilde u]^\alpha
+ Z_3\tilde L_4\tilde \pi
=
\tilde L_4\tilde h + \calR_{\tilde\pi},
\end{align}
\end{subequations}
where the remainder terms $\calR_{\tilde r},\mathcal R_{\tilde u},\mathcal R_{\tilde\pi}$
collect all commutator terms and are given by
\begin{equation}\label{eq:commutator-modified-linearized-equation-elliptic-operator}
    \begin{aligned}
        \calR_{\tilde r} &:= [D_t,\tilde L_1]\tilde r + \kappa ra_1\p_\mu [(\tilde L_2+\tilde L_3)\tilde u]^\mu - \tilde L_1(\kappa r a_1\p_\mu \tilde u^\mu) + \p_\mu r [(\tilde L_2+\tilde L_3)\tilde u]^\mu-\tilde L_1(\p_\mu r\tilde u^\mu) \\
        &\quad + V_1\tilde L_1\tilde r - \tilde L_1(V_1\tilde r)+ r^2Z_1\tilde L_4\tilde \pi - \tilde L_1(r^2Z_1\tilde \pi),\\
        \calR_{\tilde u} &:= [D_t,(\tilde L_2+\tilde L_3)]\tilde u + \left(1+\frac1\kappa\right)a_2\Delta^{\cdot\mu}\p_\mu(\tilde L_1\tilde r) - (\tilde L_2 + \tilde L_3)\left[\left(1+\frac1\kappa\right)a_2\Delta^{\cdot\mu}\p_\mu \tilde r\right]\\
        &\quad + \left( 2 + \frac1\kappa \right)ra_3\Delta^{\cdot\mu}\p_\mu\tilde L_4\tilde\pi - (\tilde L_2+\tilde L_3)\left[\left(2+\frac1\kappa\right)ra_3\Delta^{\cdot\mu}\p_\mu\tilde\pi\right] \\
        &\quad + V_2^\cdot \tilde L_1\tilde r - (\tilde L_2 + \tilde L_3)(V_2^\cdot\tilde r) + (\bar W_2)^\cdot_\beta [(\tilde L_2+\tilde L_3)\tilde u]^\beta - (\tilde L_2 + \tilde L_3)\left((\bar W_2)^\cdot_\beta \tilde u^\beta\right) \\
        &\quad + Z^\cdot_2 \tilde L_4\tilde \pi - (\tilde L_2 + \tilde L_3)(Z_2^\cdot\tilde \pi),\\
        \calR_{\tilde\pi} &:= [D_t,\tilde L_4]\tilde \pi + ra_4\p_\mu[(\tilde L_2+\tilde L_3)\tilde u]^\mu - \tilde L_4(ra_4\p_\mu\tilde u^\mu) + V_3 \tilde L_1\tilde r - \tilde L_4 (V_3\tilde r)\\
        &\quad + (W_3)_\beta [(\tilde L_2+\tilde L_3)\tilde u]^\beta - \tilde L_4\left((W_3)_\beta \tilde u^\beta\right) + Z_3\tilde L_4 \tilde \pi - \tilde L_4 (Z_3\tilde \pi).
    \end{aligned}
\end{equation}

Thus, the commutator operator $[\calL(\vec\varphi),L]$ acting on the linearized variables $(\tilde r,\tilde u,\tilde \pi)$ can be written as
\begin{equation}\label{eq:commutator-operator-one-elliptic}
    [\calL(\vec\varphi),L]\tilde\phi=\begin{pmatrix}
        \calR_{\tilde r} \\ \calR_{\tilde u} \\ \calR_{\tilde\pi}
    \end{pmatrix}.
\end{equation}

The proof of the tame estimates breaks into four steps. First, we introduce a bookkeeping scheme, which assigns to each term produced by commuting the weighted elliptic operators and the linearized operator a class and an order. We show that every term that is lower-order in the sense of the bookkeeping scheme, can be controlled. Second, we show by the tame elliptic estimates that $\norm{L^l\tilde\phi}_{\calH^0}$ controls $\norm{\tilde\phi}_{\calH^{2l}}$. Third, we verify that all the terms of the commutator $[\calL(\vec\varphi),L^l]$ are lower-order in the sense of the bookkeeping scheme, and thus the commuator can be controlled. Finally, the basic energy estimate applied to $L^l\tilde\phi$ closes the argument.

\subsection{Bookkeeping scheme}\label{subsec:bookkeeping}
Commuting the modified linearized equations with $L^l$ produces a large number of terms. Each of them is a product of a power of the weight $r$, derivatives of the background states $\vec\varphi$, and derivatives of one linearized variable or one source term. To estimate such a term we only need to know the power $\sigma$ of the weight, the number $N_r$ of factors of derivatives of the background states of the form $\p^a r$, the total count $K$ of derivatives, and the number $c$ of derivatives falling on the linearized variables. Throughout this section, we write
\begin{equation*}
    \alpha:=\frac{1}{2\kappa}-\frac{1}{2}.
\end{equation*}
Then, we can write $\calH^0 = H^{0,\alpha}\times H^{0,\alpha+\frac{1}{2}}\times H^{0,\alpha+\frac{1}{2}}$. We refer to the three components of $\calH^0$ as the $r$-, $\vec u$- and $\pi$-slots.

We set, for a slot $s\in\{r,u,\pi\}$ and for $X\in\{\tilde r,\tilde u,\tilde\pi,\tilde f,\tilde g,\tilde h\}$,
\begin{equation*}
    \varsigma_s:=
    \begin{cases}
        0, & s=r,\\
        \tfrac12, & s\in\{u,\pi\},
    \end{cases}
    \qquad
    \varsigma_X:=
    \begin{cases}
        0, & X\in\{\tilde r,\tilde f\},\\
        \tfrac12, & X\in\{\tilde u,\tilde\pi,\tilde g,\tilde h\},
    \end{cases}
\end{equation*}
and
\begin{equation*}
    \alpha_s:=\alpha+\varsigma_s,\qquad \alpha_X:=\alpha+\varsigma_X.
\end{equation*}
Thus the slot $s$ of $\calH^0$ is measured in $H^{0,\alpha_s}$, and the
norm of $X$ in $\calH^{2l}$ is $\norm{X}_{H^{2l,l+\alpha_X}}$.

\begin{defn}\label{defn:elementary-differential-operators}
Let $\vec\varphi:=(r,\vec u,a_0)$ be the background state. An elementary operator is a differential operator of the form
\begin{equation}\label{eq:elementary-operator}
    D=D(\vec\varphi)=\calF(\vec\varphi)\,r^{\sigma}\prod_{j=1}^{J}\p^{a_j}\varphi_j\;\p^{c},
    \qquad \varphi_j\in\{r,\vec u,a_0\},\quad a_j\ge1,\quad \sigma,c,J\ge0,
\end{equation}
where $\calF$ is a smooth function in a neighborhood of the range of $\vec\varphi$, and at most one of the derivatives appearing in \eqref{eq:elementary-operator} is $\p_t$. We set
\begin{equation*}
    \sigma(D) = \sigma,\qquad
    N_r(D):=\#\{j:\varphi_j=r\},\qquad
    K(D):=\sum_{j=1}^{J}a_j+c,\qquad c(D):=c,\qquad
    \operatorname{ord}(D):=c-\sigma.
\end{equation*}
For integers $l,k\ge0$, we say that $D$ is \emph{of class
$\mathcal{M}_{l,k}$} if
\begin{equation*}
    \sigma+N_r(D)\ge l\qquad\text{and}\qquad K(D) \le k.
\end{equation*}
A differential operator is of class $\mathcal{M}_{l,k}$ if it is a finite sum of elementary operators of class $\mathcal{M}_{l,k}$. The order of a differential operator is the largest order among the summands, and $c(D)$ is the largest value of $c$ among its elementary operators.
\end{defn}

\begin{defn}\label{defn:monomial}
Let $X\in\{\tilde r,\tilde u,\tilde\pi,\tilde f,\tilde g,\tilde h\}$. A monomial in $X$ is an expression $DX$ with $D$ an elementary operator, and we set
\begin{equation*}
    \sigma(DX):=\sigma(D),\qquad N_r(DX):=N_r(D),\qquad
    K(DX):=K(D),\qquad c(DX):=c(D).
\end{equation*}
Given a slot $s\in\{r,u,\pi\}$, the effective weight and the order of
$DX$ in the slot $s$ are
\begin{equation*}
    \bar\sigma_s(DX):=\sigma(D)+\varsigma_s-\varsigma_X,
    \qquad
    \ord_s(DX):=c(D)-\bar\sigma_s(DX)
    =\ord(D)+\varsigma_X-\varsigma_s.
\end{equation*}
For integers $l,k,j\ge0$, we say that $DX$ is of class
$\calM^s_{l,k}$ if
\begin{equation*}
    \bar\sigma_s(DX)+N_r(DX)\ge l\qquad\text{and}\qquad K(DX)\le k ,
\end{equation*}
and that $DX$ is lower order at level $2j$ in the slot $s$ if
$\ord_s(DX)\le j$ and $c(DX)\le 2j$. For a finite sum $\Sigma$ of monomials, possibly in different variables, $\ord_s(\Sigma)$ and $c(\Sigma)$ are the largest values of
$\ord_s$ and $c$ among its summands, and $\Sigma$ is of class $\mathcal{M}^s_{l,k}$, and lower order at level $2j$ in the slot $s$, if each of its summands is.
\end{defn}

\begin{lem}\label{lem:composition-elementary-operators}
Let $D_1$ and $D_2$ be elementary operators of class $\mathcal{M}_{l_1,k_1}$ and $\mathcal{M}_{l_2,k_2}$, respectively, and assume that at most one $\p_t$ appears in $D_1$ and $D_2$ combined. Then $D_1D_2$ is a finite sum of elementary operators, each of class $\mathcal{M}_{l_1+l_2,\,k_1+k_2}$, with
\begin{equation*}
    \ord{(D_1D_2)}\le \ord{(D_1)} + \ord{(D_2)},
    \qquad
    c(D_1D_2)\le c(D_1)+c(D_2).
\end{equation*}
By linearity, the same holds when $D_1$ and $D_2$ are differential operators of class $\mathcal{M}_{l_1,k_1}$ and $\mathcal{M}_{l_2,k_2}$.
\end{lem}

\begin{proof}
Write $D_i = \calF_i(\vec\varphi)r^{\sigma_i}P_i \p^{c_i}$, where $P_i$ is the product of the factors involving the background variables $\vec\varphi$. Since $D_1,D_2$ has at most one time derivative, by the Leibniz rule, $D_1D_2$ is a finite sum of elementary operators, each obtained by derivatives of $\p^{c_1}$ in $D_1$ either fall on $\p^{c_2}$, the weight $r^{\sigma_2}$, or the coefficient $\calF_2(\vec\varphi)$ or on one of the background factors $P_2$ of $D_2$. In the first case, the number $c$ of derivatives of each summand increases by one. In the second, the power $\sigma$ of weights $r^\sigma$ of each summand decreases by one and the number $N_r$ of differentiated background $r$ increases by one. In the third, a background factor is created or lengthened. Among all cases, the total number $K$ of derivatives and the number of time derivatives are unchanged, $\sigma+N_r$ is non-decreasing, and $c-\sigma$ is non-increasing. Therefore, each resulting summand has $\sigma+N_r\ge l_1+l_2$, $K=K(D_1)+K(D_2)\le
k_1+k_2$, $c\le c_1+c_2$ and order $\le \ord(D_1)+\ord(D_2)$. 
\end{proof}

\begin{rmk}\label{rmk:class-of-L}
    The weighted elliptic operators $\tilde L_1$, $\tilde L_2+\tilde L_3$ and $\tilde L_4$ are all sums of purely spatial elementary operators of class $\calM_{1,2}$ of order $1$. By Lemma~\ref{lem:composition-elementary-operators}, $\tilde L_1^l$, $(\tilde L_2+\tilde L_3)^l$ and $(\tilde L_4)^l$ are of class $\calM_{l,2l}$ of order $l$. Thus, $L^l\tilde\phi$ is lower-order in each slot at level $2l$.
\end{rmk}

The monomials in the commutator $[\calL(\vec\varphi),L^l]\tilde\phi$ contain at most one time derivative. Notice that the estimate in Lemma~\ref{lem:multilinear-tame-estimate} below is stated only for monomials in which the time derivative, if present, falls on $\vec u$. The following lemma reduces to this case. It uses the substitution formulas
\begin{align}
    \p_t r &= \calF(\vec\varphi)\,r\,\p\vec u+\calF(\vec\varphi)\,r\,\p_t\vec u
        +\calF(\vec\varphi)\,\p r,\label{eq:dt-r-schematic}\\
    \p_t a_0 &= \calF(\vec\varphi)\,\p\vec u+\calF(\vec\varphi)\,\p_t\vec u
        +\calF(\vec\varphi)\,\p a_0+\calF(\vec\varphi),
        \label{eq:dt-a0-schematic}
\end{align}
which follow from writing the first equations of \eqref{eq:r-pi-defining} and \eqref{eq:a_0-evolution} schematically and writing $\p_t=\frac{1}{u^0}(D_t-u^i\p_i)$ and $\p_\mu u^\mu=\p_iu^i+\frac{u_i\p_tu^i}{u^0}$. Here and below $\calF$ denotes a generic smooth function in a neighborhood of the range of $\vec\varphi$, subject to change from term to term. It also uses the following variant of the formulas in Remark~\ref{rmk:modified-linearized-eqs-time-derivative}, which is
\begin{equation}\label{eq:dt-phi-schematic}
    \left\{
    \begin{aligned}
        \p_t \tilde r &= \calF\p\tilde r + r\calF\p\vec{\tilde u} + r\calF \p\tilde \pi + \mfL\tilde r + r\mfL\vec{\tilde u} + r\mfL\tilde\pi + \calF\p r \vec{\tilde u} + \calF \tilde f+ r\calF\tilde g + r\calF\tilde h,\\
        \p_t \vec{\tilde u} &= \calF\p\tilde r + \calF\p\vec{\tilde u} + \calF\p\tilde\pi
            + \mfL\tilde r + \mfL\vec{\tilde u} + \mfL\tilde\pi
            + \calF\tilde f + \calF\vec{\tilde g} + \calF\tilde h,\\
        \p_t \tilde \pi &= \calF\p\tilde r + \calF\p\vec{\tilde u} + \calF\p\tilde\pi
            + \mfL\tilde r + \mfL\vec{\tilde u} + \mfL\tilde\pi
            + \calF\tilde f + \calF\vec{\tilde g} + \calF\tilde h,
    \end{aligned}
    \right.
\end{equation}
where $\mfL=\sum_Y\calF\,\p Y+\calF\,\p_t\vec u+\calF$ with $Y\in\{r,\vec u,a_0\}$. The only difference between \eqref{eq:dt-phi-schematic} and
Remark~\ref{rmk:modified-linearized-eqs-time-derivative} is that, in the first line, the zeroth-order terms in $\tilde u$ and $\tilde\pi$ are written as $r\mfL\vec{\tilde u}+\calF\,\p r\,\vec{\tilde u}$ and
$r\mfL\tilde\pi$ rather than $\mfL\vec{\tilde u}$ and $\mfL\tilde\pi$.
This can be directly read off from \eqref{eq:final}: the coefficients $\mfa_{12}$, $\frac{\kappa ra_1u_i}{(u^0)^2}$ and $\frac{ra_4u_i}{(u^0)^2}$ of $\mfR_{\tilde\pi}$ and $\mfR^i_{\vec{\tilde u}}$ carry a factor $r$, and in $\mfR_{\tilde r}$ the coefficients of $\tilde u_i$ and $\tilde\pi$ are
$\frac{\kappa ra_1}{u^0}H^{ij}\p_tu_j+\frac{u_i}{u^0}\p_tr+\p_ir$ and
$r^2Z_1$, the former carrying $r$ or $\p r$ in each term once $\p_tr$ is replaced by \eqref{eq:dt-r-schematic}.

\begin{lem}[Elimination of time derivatives]\label{lem:time-derivative-elimination}
Let $\tilde\phi=(\tilde r,\tilde u,\tilde\pi)$ solve the modified linearized equations with source terms $\tilde F=(\tilde f,\tilde g,\tilde h)$. Let $DX$ be a monomial in $X\in\{\tilde r,\tilde u,\tilde\pi\}$ of class $\mathcal{M}^{s}_{l,k}$ with $\ord_s{(DX)}=m$, containing one time derivative which does not fall on $\vec u$. Then $DX$ is a finite sum of monomials of the following two types:
\begin{enumerate}
    \item[(a)] monomials $D'X'$ in some $X'\in\{\tilde r,\tilde u,\tilde\pi\}$, of class $\mathcal{M}^{s}_{l,k}$ with order $\ord_s{(D'X')}\leq m$ and $K(DX)\ge K(D'X')\ge K(DX)-1$ and $c(D'X')\le c(DX)$, in which the time derivative, if present, falls on $\vec u$,
    \item[(b)] monomials $D''X''$ in some $X''\in\{\tilde f,\tilde g,\tilde h\}$, of class $\mathcal{M}^{s}_{l,k-1}$ with order $\ord_s{(D''X'')} \le m$ and $K(D''X'')=K(DX)-1$ and $c(D''X'')\le c(DX)-1$, containing no time derivative.
\end{enumerate}
Monomials of type (b) occur only if the time derivative of $DX$ falls on $X$.
\end{lem}
\begin{proof}
    We write $D$ as in \eqref{eq:elementary-operator}. Let $\bar\sigma:=\bar\sigma_s(DX)$ be the effective weight of $DX$, $N_r:=N_r(DX)$ and $c=c(DX)$. Since $DX$ is a monomial in $X$ of class $\calM^s_{l,k}$ with order $\ord_s(DX)=m$, we have
    \begin{equation*}
        \bar\sigma + N_r \geq l, \qquad 
        c-\bar\sigma = m.
    \end{equation*}
    We divide the proof into cases according to which factor carries the time derivative $\p_t$. In each case, we will substitute \eqref{eq:dt-r-schematic}, \eqref{eq:dt-a0-schematic} and \eqref{eq:dt-phi-schematic} and expand the remaining spatial derivatives of the same factor by the Leibniz rule. All monomials thus produced contain at most one time derivative, and it can only fall on $\vec u$. It then remains to check the class, $K$ and $c$ of each generated term.

    If the time derivative falls on a factor $\p^{a-1}\p_tr$ with $a\ge1$, then by \eqref{eq:dt-r-schematic}, this factor is replaced by $\p^{a-1}(\calF r\p\vec u)$, $\p^{a-1}(\calF r\p_t\vec u)$ and $\p^{a-1}(\calF\p r)$, which generates monomials in the same $X$. In the terms generated from the first two factors, if the factor $r$ is not differentiated, it joins the power $r^\sigma$, so that $\sigma$ increases by one while $N_r$ could decrease by one. If it is differentiated, $N_r$ is unchanged. In the terms generated from the third, $N_r$ is unchanged. In all cases, the derivatives falling on $\calF$ can only increase $N_r$. Hence, $\sigma(D'X')\geq \sigma(DX)$, $\bar\sigma_s(D'X')+N_r(D'X')\ge\bar\sigma_s(DX)+N_r(DX)$, and $K(D'X')=K(DX)$, $c(D'X')=c(DX)$, so that $\ord_s(D'X')\le\ord_s(DX)$.

    If the time derivative falls on a factor $\p^{a-1}\p_ta_0$ with $a\ge1$, then by \eqref{eq:dt-a0-schematic} this factor is replaced by $\p^{a-1}(\calF\p\vec u)$, $\p^{a-1}(\calF\p_t\vec u)$, $\p^{a-1}(\calF\p a_0)$ and $\p^{a-1}\calF$. The first three factors are treated as in the previous case: $\sigma$, $c$ and $K$ are unchanged and $N_r$ does not decrease. In the fourth factor, $\calF$ carries no derivative: if $a\ge2$, $\p^{a-1}\calF$ is a sum of products of derivatives of $r$, $\vec u$ and $a_0$ of total count $a-1$ with coefficients smooth in $\vec\varphi$, and if $a=1$ it is absorbed into the coefficient of $D$. In both cases $K$ decreases by one, $N_r$ does not decrease, and $\sigma$, $c$ are unchanged. Hence, $\bar\sigma_s(D'X')+N_r(D'X')\ge\bar\sigma_s(DX)+N_r(DX)$, $\ord_s(D'X')=\ord_s(DX)$, $c(D'X')=c(DX)$, and these monomials are of type (a) with $K(D'X')\ge K(DX)-1$.

    If the time derivative falls on $X$, i.e., $DX$ contains $\p^{c-1}\p_t X$ with $c\ge 1$, we substitute \eqref{eq:dt-phi-schematic} and distribute $\p^{c-1}$ by the Leibniz rule. Since the generated monomials may be in a variable $X'\in\{\tilde r,\tilde u,\tilde\pi\}$ or $X''\in\{\tilde f,\tilde g,\tilde h\}$ different from $X$, we now write $\sigma',N_r',K',c'$ and $\bar\sigma'_s:=\bar\sigma_s(D'X')=\sigma'+\varsigma_s-\varsigma_{X'}$ for the quantities of a generated monomial $D'X'$, and similarly with double primes for a generated monomial $D''X''$ in a source term. First, we assume $X=\tilde r$, and substitute the first line of \eqref{eq:dt-phi-schematic}. The term $\p^{c-1}(\calF\p\tilde r)$ produces monomials in $\tilde r$ with $c'\leq c$, $\sigma'=\sigma$, $N'_r\geq N_r$ and $K'=K$. The term $\p^{c-1}(\mfL\tilde r)$ produces monomials in $\tilde r$ with $c'\leq c-1$, $\sigma'=\sigma$, $N'_r\geq N_r$ and $K'\geq K-1$. The term $\p^{c-1}(\calF\tilde f)$ produces spatial monomials in $\tilde f$ with $c'' \le c-1$, $\sigma''=\sigma$, $N_r''\geq N_r$ and $K'' = K-1$. Since $\varsigma_{\tilde r} = \varsigma_{\tilde f}$, $\bar\sigma'' = \bar\sigma$. These terms above are terms of type (a) and (b), respectively.

    The remaining terms of the first line of \eqref{eq:dt-phi-schematic} contain a variable $X'\in\{\tilde u,\tilde\pi\}$ or $X''\in\{\tilde g,\tilde h\}$, for which $\varsigma_{X'}=\varsigma_{X''}=\frac12$, so that $\bar\sigma'=\sigma'+\varsigma_s-\frac12=\bar\sigma+(\sigma'-\sigma)-\frac12$ and likewise for $\bar\sigma''$. They are of two kinds.

    For terms with a factor $r$, i.e., $r\calF\p X'$, $r\mfL X'$ and $r\calF X''$, let $b\ge0$ be the number of derivatives of $\p^{c-1}$ falling on this factor $r$. If $b=0$, the factor $r$ joins the power $r^\sigma$, so $\sigma' = \sigma + 1$ and $N_r'\geq N_r$. Hence, $\bar\sigma'=\bar\sigma+\frac12$ and $c'\leq c$. If $b\geq 1$, the factor becomes $\p^b r$, so $\sigma'=\sigma$, $N_r'\geq N_r+1$. Hence, $\bar\sigma'=\bar\sigma -\frac{1}{2}$, and $c'\leq c-b\leq c-1$. In both cases,
    \begin{equation*}
        \bar\sigma'+N'_r\geq \bar\sigma + N_r +\frac{1}{2}, \qquad c'-\bar\sigma'\leq c-\bar\sigma-\frac12.
    \end{equation*}
    As for $K$ and $c$, the terms $r\calF\p X'$ have $K'=K$ and $c'\le c$, and the terms $r\mfL X'$ have $K'\ge K-1$ and $c'\le c-1$. For the terms $r\calF X''$, the same computation gives $\bar\sigma''+N_r''\ge\bar\sigma+N_r+\frac12$, $c''-\bar\sigma''\le c-\bar\sigma-\frac12$, $K''=K-1$ and $c''\le c-1$. So, these terms above are of type (a) and (b), respectively.

    For the term $\calF\p r\tilde u$, we have $\sigma'=\sigma$, $N'_r\geq N_r+1$, $K'=K$ and $c'\leq c-1$. Hence, $\bar\sigma'=\bar\sigma-\frac12$ and 
    \begin{equation*}
        \bar\sigma'+N_r'\geq \bar\sigma + N_r +\frac{1}{2},\qquad
        c'-\bar\sigma'\leq c-\bar\sigma-\frac{1}{2}.
    \end{equation*}
    These terms are of type (a).

    In all cases, the generated monomials $D'X'$ are of class $\calM^s_{l,k}$ with $\ord_s(D'X')\leq m$, $K - 1\leq K'\leq K$ and $c'\leq c$, which are of type (a). The generated monomials $D''X''$ are of the class $\calM^s_{l,k-1}$ with $\ord_s(D''X'')\leq m$ and $K'' = K-1$ and $c''\leq c-1$, which are of type (b).

    If $X\in\{\tilde u,\tilde\pi\}$, we substitute the second or third line of \eqref{eq:dt-phi-schematic}. Since $\varsigma_X=\frac12$, the effective weight does not decrease when $X$ is replaced by any $X'$ or $X''$, and the same counting as above shows that the generated monomials are of type (a) and (b).
\end{proof}

In Section~\ref{subsec:commutator}, we show that, for each $s\in\{r,u,\pi\}$, the $s$-component of the commutator $[\calL(\vec\varphi),L^l]\tilde\phi$ is a finite sum of monomials of class $\mathcal{M}^s_{l,2l+1}$ which are lower order at level $2l$ in the slot $s$, and that the $s$-component of $L^l\tilde\phi$ is a finite sum of monomials of the same kind. The following lemma bounds all such monomials after elimination of the time derivative using Lemma~\ref{lem:time-derivative-elimination}.

\begin{lem}[Tame estimate for lower-order monomials]\label{lem:multilinear-tame-estimate}
Let $l>\frac{1}{2\kappa}+2$, $t\in[0,T]$ and $s\in\{r,u,\pi\}$. Let $DX$ be a monomial in $X\in\{\tilde r,\tilde u,\tilde\pi,\tilde f,\tilde g,\tilde h\}$ of class $\mathcal{M}^{s}_{l,2l+1}$ which is lower order at level $2l$ in the slot $s$, with $K(DX)\ge 2l$, and in which the time derivative, if present, falls on $\vec u$. Then
\begin{equation}\label{eq:multilinear-tame-1}
\begin{aligned}
    \norm{DX}_{H^{0,\alpha_s}}
    &\lesssim \calC_1\Big(\norm{X}_{H^{2l,l+\alpha_X}}
    +\norm{X}_{L^\infty}\Big(
    \norm{r}_{H^{2l+4,\,l+\frac{1}{2\kappa}+\frac{3}{2}}}
    +\norm{\vec u}_{H^{2l+4,\,l+\frac{1}{2\kappa}+2}}\\
    &\qquad\qquad\qquad
    +\norm{\p_t\vec u}_{H^{2l+2,\,l+\frac{1}{2\kappa}+1}}
    +\norm{a_0}_{H^{2l+2,\,l+\frac{1}{2\kappa}+\frac{1}{2}}}
    \Big)\Big),
\end{aligned}
\end{equation}
with $\calC_1=\calC_1\big(\norm{r}_{W^{2,\infty}},\norm{\vec u}_{W^{1,\infty}},\norm{\p_t\vec u}_{L^\infty},\norm{a_0}_{W^{1,\infty}}\big)$. If moreover $K(DX)\le 2l$ and $DX$ contains no time derivative, then
\begin{equation}\label{eq:multilinear-tame-2}
\begin{aligned}
    \norm{DX}_{H^{0,\alpha_s}}
    &\lesssim \calC_2\Big(\norm{X}_{H^{2l,l+\alpha_X}}
    +\norm{X}_{L^\infty}\Big(
    \norm{r}_{H^{2l+2,\,l+\frac{1}{2\kappa}+\frac{1}{2}}}
    +\norm{\vec u}_{H^{2l+2,\,l+\frac{1}{2\kappa}+1}}
    +\norm{a_0}_{H^{2l,\,l+\frac{1}{2\kappa}-\frac{1}{2}}}
    \Big)\Big),
\end{aligned}
\end{equation}
with $\calC_2=\calC_2\big(\norm{r}_{W^{1,\infty}},\norm{\vec u}_{L^\infty},\norm{a_0}_{L^\infty}\big)$. In both estimates, the weighted norm of $a_0$ (respectively of $\vec u$, $\p_t\vec u$) on the right-hand side is absent if no factor $\p^{a_j}\varphi_j$ of $D$ involves $a_0$ (respectively $\vec u$, $\p_t\vec u$). The same holds for the arguments of $\calC_1$ and $\calC_2$ provided, in addition, that the coefficient $\calF$ does not depend on $a_0$ (respectively $\vec u$).
\end{lem}

\begin{proof}
Write $D$ as in \eqref{eq:elementary-operator} and set $\bar\sigma:=\bar\sigma_s(DX)=\sigma+\varsigma_s-\varsigma_X$. The coefficient $\calF(\vec\varphi)$ is bounded by a constant depending on $\norm{r}_{L^\infty},\norm{\vec u}_{L^\infty},\norm{a_0}_{L^\infty}$, so we may assume $\calF\equiv1$ for simplicity, and then
\begin{equation}\label{eq:product-to-estimate}
    \norm{DX}_{H^{0,\alpha_s}}=\Big\|r^{\sigma+\alpha_s}\prod_{j=1}^{J}\p^{a_j}\varphi_j\;\p^{c}X\Big\|_{L^2}.
\end{equation}
The hypotheses of the lemma read
\begin{equation}\label{eq:hypotheses-in-numbers}
    \bar\sigma+N_r\ge l,\qquad 2l\le K\le 2l+1,\qquad c\le l+\bar\sigma,\qquad c\le 2l.
\end{equation}
Let $k_*:=\lfloor l+\bar\sigma\rfloor$. Since $N_r$, $c$ and $l$ are integers, the first and third conditions give
\begin{equation}\label{eq:Nr-and-c}
    N_r\ge 2l-k_*,\qquad c\le k_*.
\end{equation}
Note that $\bar\sigma\ge-\frac12$, hence $k_*\ge l-1$, and $k_*\ge l$ if $\varsigma_s\ge\varsigma_X$.

\noindent\textbf{The plan.}
For functions $g_1,\dots,g_n$ and integers $k_1,\dots,k_n\ge1$ with $k_1+\dots+k_n=k$, Corollary~\ref{Cor:Multilinear-Interpolation} gives
\begin{equation}\label{eq:cor-recalled}
    \Big\|r^{\tau}\prod_{i=1}^n\p^{k_i}g_i\Big\|_{L^2}\lesssim\sum_{i=1}^n\Big(\prod_{i'\ne i}\norm{g_{i'}}_{L^\infty}\Big)\norm{r^{\tau}\p^{k}g_i}_{L^2}
\end{equation}
if $\tau>-\frac12$ and $k-\tau>\frac32$. We shall write \eqref{eq:product-to-estimate} in the form of the left-hand side of \eqref{eq:cor-recalled}, with $X$ among the $g_i$, and then bound each term on the right by Lemma~\ref{lem:WeightedEmbedding}. For the term with $g_i=X$, we need $\norm{r^{\sigma+\alpha_s}\p^kX}_{L^2}\lesssim\norm{X}_{H^{2l,l+\alpha_X}}$, and by Lemma~\ref{lem:WeightedEmbedding} this holds if $k\le 2l$ and $l+\alpha_X\le\sigma+\alpha_s+2l-k$, i.e., if
\begin{equation}\label{eq:budget}
    k\le \min\{k_*,2l\}.
\end{equation}
If regarded as a product of derivatives of $r,\vec u,a_0$ and $X$, the product \eqref{eq:product-to-estimate} has $K$ derivatives, and $K$ may exceed $\min\{k_*,2l\}$. We therefore begin by regarding some of the background factors as derivatives of $\p r$, $\p^2r$, $\p\vec u$ or $\p a_0$ instead of $r$, $\vec u$ or $a_0$, so as to lower the total number of derivatives in the product. This is the reason why $\calC_1$ depends on $\norm{r}_{W^{2,\infty}}$, $\norm{\vec u}_{W^{1,\infty}}$ and $\norm{a_0}_{W^{1,\infty}}$ rather than on $L^\infty$ norms alone.

\noindent\textbf{Step 1: lowering the number of derivatives.}
Observe that the arguments of $\calC_1$ bound the $L^\infty$ norms of the eight functions
\begin{equation}\label{eq:bounded-functions}
    r,\ \p r,\ \p^2r,\ \vec u,\ \p\vec u,\ a_0,\ \p a_0,\ \p_t\vec u.
\end{equation}
A factor $\p^{a}r$ of the product is the derivative of order $a$ of $r$, but it is also the derivative of order $a-1$ of $\p r$ and, if $a\ge2$, the derivative of order $a-2$ of $\p^2r$. Likewise, $\p^{a}\vec u$ and $\p^a a_0$ are the derivatives of order $a-1$ of $\p\vec u$ and $\p a_0$, respectively, and $\p^{a-1}\p_t\vec u$ is the derivative of order $a-1$ of $\p_t\vec u$. Thus, the product in \eqref{eq:product-to-estimate} can be written as
\begin{equation}\label{eq:reduced-product}
    \prod_{j=1}^{J}\p^{b_j}f_j\;\p^{c}X,\qquad f_j\text{ among }\eqref{eq:bounded-functions},\quad 0\le b_j\le a_j ,
\end{equation}
and the total number of derivatives in this representation is $\sum_jb_j+c$.

Let $J_t\in\{0,1\}$ be the number of factors containing $\p_t$. If every factor is read as a derivative of $r$, $\vec u$, $a_0$ or $\p_t\vec u$, then the number of derivatives is
\begin{equation*}
    S:=K-J_t ,
\end{equation*}
and this is the largest possible value of $\sum_jb_j+c$ among all such representations. Reading a factor $\p^{a_j}r$ as a derivative of $\p r$ or $\p^2r$ lowers the number by one or two, the latter only if $a_j\ge2$, and reading $\p^{a_j}\vec u$ or $\p^{a_j}a_0$ as a derivative of $\p\vec u$ or $\p a_0$ lowers it by one. Therefore, the number of derivatives can be lowered by any amount between $0$ and $R$ given by
\begin{equation*}
    R:=N_r+\#\{j:\varphi_j=r,\ a_j\ge2\}+\#\{j:\varphi_j\in\{\vec u,a_0\}\text{ without }\p_t\}.
\end{equation*}
We choose the number of derivatives to be
\begin{equation}\label{eq:choice-of-k}
    k:=\min\{S,\,k_*,\,2l\},
\end{equation}
which satisfies \eqref{eq:budget}. This requires lowering the number of derivatives by $S-k$, which is possible only if $S-k\le R$. We now show that this inequality holds.

If $k=S$, there is nothing to prove. If $k=2l<S$, then $S=2l+1$ because $S\le K\le 2l+1$, so $S-k=1$. To show the inequality, we need $R\ge1$, i.e. a background factor without $\p_t$. Such a factor exists, since if there were none, then either $J=0$ and thus $S=c\le 2l$, or the only background factor contains $\p_t$ and $S=K-J_t = K-1\le 2l$. Both contradict $S=2l+1$. Finally if $k=k_*<S$, by \eqref{eq:Nr-and-c},
\begin{equation*}
   S-k = S-k_*=K-J_t-k_*\le(2l+1)-J_t-k_*=(2l-k_*)+1-J_t\le N_r+1-J_t.
\end{equation*}
By definition, $R\ge N_r$. If some background factor without $\p_t$ is not of the form $\p r$, then $R\ge N_r+1$ and $S-k\le N_r+1\le R$. Otherwise every background factor without $\p_t$ is $\p r$ and $R=N_r$. If $J_t=1$, the display above gives $S-k\le N_r=R$. If $J_t=0$, then $K = N_r+c$ and $S-k=K-k_*=N_r+(c-k_*)\le N_r=R$ by \eqref{eq:Nr-and-c}. This proves $S-k\le R$.

\noindent\textbf{Step 2: adjusting the weight.}
To apply \eqref{eq:cor-recalled} we need $k-\tau>\frac32$, which may fail for $\tau=\sigma+\alpha_s$ when $\sigma$ is large. Since $r$ is bounded, we may replace $r^{\sigma}$ by $r^{\sigma'}$ with integer $\sigma'$ such that $0\le\sigma'\le\sigma$. We choose $\sigma'$ as small as possible subject to $k\le l+\bar\sigma'$, where $\bar\sigma':=\sigma'+\varsigma_s-\varsigma_X$, so that $\norm{r^{\sigma'+\alpha_s}\p^kX}_{L^2}\lesssim\norm{X}_{H^{2l,l+\alpha_X}}$ still holds by Lemma~\ref{lem:WeightedEmbedding}. This means
\begin{equation}\label{eq:choice-of-sigma-prime}
    \lfloor l+\bar\sigma'\rfloor=k.
\end{equation} 
Such a $\sigma'$ exists. Indeed, as $\sigma'$ increases from $0$ to $\sigma$, the integer $\lfloor l+\bar\sigma'\rfloor$ increases by one at each step, from
\begin{equation*}
    k_0:=\lfloor l+\varsigma_s-\varsigma_X\rfloor=
    \begin{cases} l, & \varsigma_s\ge\varsigma_X,\\ l-1, & \varsigma_s<\varsigma_X,\end{cases}
\end{equation*}
to $k_*$, so it suffices to show $k_0\le k\le k_*$. The upper bound is part of \eqref{eq:choice-of-k}. For the lower bound, recall $k=\min\{S,k_*,2l\}$, where $S\ge K-1\ge 2l-1\ge l$ and $2l\ge l$. If $\varsigma_s\ge\varsigma_X$, then $\bar\sigma\ge0$ and $k_*\ge l=k_0$; if $\varsigma_s<\varsigma_X$, then $k_*\ge l-1=k_0$. In both cases $k\ge k_0$.

Set $\tau':=\sigma'+\alpha_s$, then $\tau'>-\frac12$. Since $l+\bar\sigma'\in\frac12\mathbb{Z}$, we have $l+\bar\sigma'-\frac12\le k\le l+\bar\sigma'$, and therefore
\begin{align}
    k&\le l+\sigma'+\varsigma_s-\varsigma_X ,\label{eq:R1}\\
    k-\tau'&\ge l+\bar\sigma'-\tfrac12-\sigma'-\alpha_s=l-\tfrac12-\varsigma_X-\alpha
    \ge l-\tfrac{1}{2\kappa}-\tfrac12>\tfrac32 ,\label{eq:R2}
\end{align}
where the last inequality is the assumption $l>\frac{1}{2\kappa}+2$. Using $r^{\sigma+\alpha_s}\le\norm{r}_{L^\infty}^{\sigma-\sigma'}r^{\tau'}$ and \eqref{eq:reduced-product} with $\sum_jb_j+c=k$, we have
\begin{equation}\label{eq:after-weight-adjustment}
    \norm{DX}_{H^{0,\alpha_s}}\lesssim\Big\|r^{\tau'}\prod_{j=1}^{J}\p^{b_j}f_j\;\p^{c}X\Big\|_{L^2},\qquad\sum_jb_j+c=k.
\end{equation}

\noindent\textbf{Step 3: interpolation.}
We apply Corollary~\ref{Cor:Multilinear-Interpolation} to the right-hand side of \eqref{eq:after-weight-adjustment}. First, every factor $\p^{b_j}f_j$ with $b_j=0$, and the factor $X$ if $c=0$, is bounded in $L^\infty$ and taken out of the norm. This produces $L^\infty$ norms of the functions \eqref{eq:bounded-functions}, which are arguments of $\calC_1$, and possibly one factor $\norm{X}_{L^\infty}$. Let $g_1,\dots,g_n$ be the distinct functions among the remaining $f_j$ and, if $c\ge1$, $X$. In the notation of Corollary~\ref{Cor:Multilinear-Interpolation}, $g_i$ is the function $f_i$, $J_i\ge1$ is the number of factors of the form $\p^{\alpha}g_i$ with $\alpha\ge1$, and the total number of derivatives is $k$. By \eqref{eq:R2} and $\tau'>-\frac12$, the corollary applies with weight $\tau'$ and gives
\begin{equation}\label{eq:after-interpolation}
    \Big\|r^{\tau'}\prod_{j=1}^{J}\p^{b_j}f_j\;\p^{c}X\Big\|_{L^2}
    \lesssim\calC\sum_{i=1}^{n}\prod_{i'\ne i}\norm{g_{i'}}_{L^\infty}\,\norm{r^{\tau'}\p^{k}g_i}_{L^2},
    \qquad\calC=\prod_{i=1}^{n}\norm{g_i}_{L^\infty}^{J_i-1}.
\end{equation}
Since $X$ occurs in at most one factor, $J_i\le1$ for $g_i=X$, so $\calC$ involves only the functions \eqref{eq:bounded-functions}. If $c\ge1$, then $X$ is one of the $g_i$, and $\norm{X}_{L^\infty}$ appears exactly once in every term of \eqref{eq:after-interpolation} except the term with $g_i=X$, where it does not appear; if $c=0$, then $\norm{X}_{L^\infty}$ was taken out of the norm and no $g_i$ equals $X$. In both cases,
\begin{equation}\label{eq:tame-structure}
    \norm{DX}_{H^{0,\alpha_s}}\lesssim\calC_1\Big(\norm{r^{\tau'}\p^{k}X}_{L^2}+\norm{X}_{L^\infty}\sum_{g}\norm{r^{\tau'}\p^{k}g}_{L^2}\Big),
\end{equation}
with the sum running over those $g$ in \eqref{eq:bounded-functions} which occur among the $g_i$. It remains to bound $\norm{r^{\tau'}\p^{k}g}_{L^2}$ for $g=X$ and for each $g$ in \eqref{eq:bounded-functions}.

\noindent\textbf{Step 4: Hardy inequality.}
Let $g$ be a function, $N\ge k$ an integer and $\beta>-\frac12$ with $\beta\le\tau'+N-k$. Applying Lemma~\ref{lem:WeightedEmbedding} $N-k$ times and using that $r$ is bounded,
\begin{equation}\label{eq:hardy-step}
    \norm{r^{\tau'}\p^{k}g}_{L^2}\le\norm{g}_{H^{k,\tau'}}
    \lesssim\norm{g}_{H^{N,\tau'+N-k}}\lesssim\norm{g}_{H^{N,\beta}}.
\end{equation}
We apply \eqref{eq:hardy-step} with the following choices of $N$ and $\beta$. In each case $N\ge k$ holds because $k\le 2l$, and the condition $\beta\le\tau'+N-k$, after inserting $\tau'=\sigma'+\alpha+\varsigma_s$, reduces to the inequality stated after the choice of $N$ and $\beta$, which follows from \eqref{eq:R1}.
\begin{itemize}
    \item $g=X$, $N=2l$, $\beta=l+\alpha_X$: the condition is $k\le l+\sigma'+\varsigma_s-\varsigma_X$, which is \eqref{eq:R1}; the bound is $\norm{X}_{H^{2l,l+\alpha_X}}$.
    \item $g=\p^2r$, $N=2l+2$, $\beta=l+2+\alpha$: the condition is $k\le l+\sigma'+\varsigma_s$; the bound is $\norm{\p^2r}_{H^{2l+2,l+2+\alpha}}\le\norm{r}_{H^{2l+4,l+\frac{1}{2\kappa}+\frac32}}$.
    \item $g=\p r$, $N=2l+1$, $\beta=l+1+\alpha$: the condition is $k\le l+\sigma'+\varsigma_s$; the bound is $\norm{\p r}_{H^{2l+1,l+1+\alpha}}\le\norm{r}_{H^{2l+2,l+\frac{1}{2\kappa}+\frac12}}$.
    \item $g=r$, $N=2l+2$, $\beta=l+1+\alpha$: the condition is $k\le l+1+\sigma'+\varsigma_s$; the bound is $\norm{r}_{H^{2l+2,l+\frac{1}{2\kappa}+\frac12}}$.
    \item $g=\p\vec u$, $N=2l+3$, $\beta=l+\alpha+\frac52$: the condition is $k\le l+\frac12+\sigma'+\varsigma_s$; the bound is $\norm{\p\vec u}_{H^{2l+3,l+\alpha+\frac52}}\le\norm{\vec u}_{H^{2l+4,l+\frac{1}{2\kappa}+2}}$.
    \item $g=\vec u$, $N=2l+2$, $\beta=l+\alpha+\frac32$: the condition is $k\le l+\frac12+\sigma'+\varsigma_s$; the bound is $\norm{\vec u}_{H^{2l+2,l+\frac{1}{2\kappa}+1}}$.
    \item $g=\p a_0$, $N=2l+1$, $\beta=l+1+\alpha$: the condition is $k\le l+\sigma'+\varsigma_s$; the bound is $\norm{\p a_0}_{H^{2l+1,l+1+\alpha}}\le\norm{a_0}_{H^{2l+2,l+\frac{1}{2\kappa}+\frac12}}$.
    \item $g=a_0$, $N=2l$, $\beta=l+\alpha$: the condition is $k\le l+\sigma'+\varsigma_s$; the bound is $\norm{a_0}_{H^{2l,l+\frac{1}{2\kappa}-\frac12}}$.
    \item $g=\p_t\vec u$, $N=2l+2$, $\beta=l+\alpha+\frac32$: the condition is $k\le l+\frac12+\sigma'+\varsigma_s$; the bound is $\norm{\p_t\vec u}_{H^{2l+2,l+\frac{1}{2\kappa}+1}}$.
\end{itemize}
By Lemma~\ref{lem:WeightedEmbedding}, we have $\norm{r}_{H^{2l+2,l+\frac{1}{2\kappa}+\frac12}}\lesssim\norm{r}_{H^{2l+4,l+\frac{1}{2\kappa}+\frac32}}$, $\norm{\vec u}_{H^{2l+2,l+\frac{1}{2\kappa}+1}}\lesssim\norm{\vec u}_{H^{2l+4,l+\frac{1}{2\kappa}+2}}$ and $\norm{a_0}_{H^{2l,l+\frac{1}{2\kappa}-\frac12}}\lesssim\norm{a_0}_{H^{2l+2,l+\frac{1}{2\kappa}+\frac12}}$. Inserting these bounds into \eqref{eq:tame-structure} gives \eqref{eq:multilinear-tame-1}.

\noindent\textbf{Step 5: proof of \eqref{eq:multilinear-tame-2}.}
Assume $K\le 2l$ and that $D$ contains no time derivative. Then $S=K\le 2l$, so $k=\min\{K,k_*\}$ and $S-k=\max\{0,K-k_*\}\le \max\{0,2l-k_*\}\le N_r$ by \eqref{eq:Nr-and-c}. Since $S-k\le N_r$, we can lower the number of derivatives from $S$ to $k$ by choosing $S-k$ of the factors $\p^{a_j}r$ and reading each of them as the derivative of order $a_j-1$ of $\p r$, while every other factor is read as a derivative of $r$, $\vec u$ or $a_0$. Then the functions $f_j$ in \eqref{eq:reduced-product} are among $r,\p r,\vec u,a_0$, whose $L^\infty$ norms are the arguments of $\calC_2$. Steps 2--4 go through unchanged, with $\calC_1$ replaced by $\calC_2$ and only the bounds of Step 4 for $g\in\{r,\p r,\vec u,a_0\}$ being used; these are the norms in \eqref{eq:multilinear-tame-2}.

We note that the functions $g$ in \eqref{eq:after-interpolation} are among $X$ and the functions $f_j$ of \eqref{eq:reduced-product}, and that $f_j\in\{r,\p r,\p^2r\}$ if $\varphi_j=r$, $f_j\in\{\vec u,\p\vec u\}$ if $\varphi_j=\vec u$, $f_j\in\{a_0,\p a_0\}$ if $\varphi_j=a_0$, and $f_j=\p_t\vec u$ if the factor contains $\p_t$. Hence a norm of $a_0$ (respectively $\vec u$, $\p_t\vec u$) occurs in \eqref{eq:after-interpolation} only if some factor of $D$ involves $a_0$ (respectively $\vec u$, $\p_t\vec u$). The $L^\infty$ norms collected in $\calC_1$ arise from the coefficient $\calF$, from $\norm{r}_{L^\infty}$ in Steps 2 and 4, from the factors with $b_j=0$ taken out in Step 3, and from the constant of Corollary~\ref{Cor:Multilinear-Interpolation} and the products $\prod_{i'\ne i}\norm{g_{i'}}_{L^\infty}$ in \eqref{eq:after-interpolation}; the last three involve only functions occurring in $D$, and the first only the arguments of $\calF$. This proves the final assertion of the lemma.
\end{proof}

\subsection{Weighted elliptic estimates}\label{subsec:elliptic-estimates}

In this subsection we record the elliptic estimates satisfied by the operators $\tilde L_1$, $\tilde L_2+\tilde L_3$ and $\tilde L_4$ introduced at the beginning of Section~\ref{sec:tame-estimate-linearization-Eulerian}. The basic estimates for $\tilde L_1$ and $\tilde L_2+\tilde L_3$ are essentially those of Lemma~25 in \cite{Zhong_2025} (see also \cite{ifrim2023compressible,disconzi2022relativistic}). The estimate for $\tilde L_4$ follows from the same integration by parts argument as that for $\tilde L_1$, with the weight $r^{\frac{1}{2\kappa}-\frac12}$ replaced by $r^{\frac{1}{2\kappa}}$. We refer to \cite{Zhong_2025} for the proofs.

\begin{lem}\label{Lem:BasicEllipticEstimates}
Let $t\in[0,T]$. For any sufficiently smooth scalar $\tilde r$,
\begin{equation}\label{Eq:BasicEllipticEstimate}
    \norm{\tilde r}_{H^{2,\frac{1}{2\kappa}+\frac12}}\lesssim\norm{\tilde{L}_1\tilde{r}}_{H^{0,\frac{1}{2\kappa}-\frac{1}{2}}}+\norm{\tilde r}_{H^{0,\frac{1}{2\kappa}-\frac{1}{2}}},
\end{equation}
for any sufficiently smooth $4$-vector field $\tilde u$ orthogonal to $u$,
\begin{equation}\label{Eq:BasicDivCurlEstimate}
    \norm{\tilde u}_{H^{2,\frac{1}{2\kappa}+1}}\lesssim\norm{(\tilde{L}_2+\tilde L_3)\tilde{u}}_{H^{0,\frac{1}{2\kappa}}}+\norm{\tilde u}_{H^{0,\frac{1}{2\kappa}}},
\end{equation}
and for any sufficiently smooth scalar $\tilde\pi$,
\begin{equation}\label{Eq:BasicEllipticEstimatePi}
    \norm{\tilde \pi}_{H^{2,\frac{1}{2\kappa}+1}}\lesssim\norm{\tilde L_4\tilde{\pi}}_{H^{0,\frac{1}{2\kappa}}}+\norm{\tilde \pi}_{H^{0,\frac{1}{2\kappa}}}.
\end{equation}
The implicit constants depend on $\norm{r}_{W^{1,\infty}}$ and $\norm{\vec u}_{W^{1,\infty}}$.
\end{lem}

\begin{rmk}\label{Rmk:ModifiedBasicEllipticEstimates}
Given an integer $b\geq0$, define
\begin{equation}\label{op:b-families}
\begin{aligned}
    \tilde L_{1,b}\tilde r&:=\tilde L_1\tilde r+b(\kappa+1)H^{ij}\p_ir\,\p_j\tilde r,\\
    \big(\tilde L_{2,b}\tilde u\big)_\alpha&:=\big(\tilde L_2\tilde u\big)_\alpha+b(\kappa+1)A^i_\alpha H^{jk}\p_jr\,\p_i\tilde u_k,\\
    \big(\tilde L_{3,b}\tilde u\big)_\alpha&:=(\kappa+1)A^i_\alpha H^{jk}r^{-\frac1\kappa-b}\p_j\Big(r^{\frac1\kappa+1+b}\big(\p_k\tilde u_i-\p_i\tilde u_k\big)\Big),\\
    \tilde L_{4,b}\tilde\pi&:=\tilde L_4\tilde\pi+b(\kappa+1)H^{ij}\p_ir\,\p_j\tilde\pi.
\end{aligned}
\end{equation}
Then, the operators $\tilde L_{1,b}$, $\tilde L_{2,b}+\tilde L_{3,b}$ and $\tilde L_{4,b}$ satisfy \eqref{Eq:BasicEllipticEstimate}, \eqref{Eq:BasicDivCurlEstimate} and \eqref{Eq:BasicEllipticEstimatePi}, respectively, by the same proof. 
\end{rmk}

The proof of the first-order tame estimate of Theorem~\ref{thm:modified-linearized-with-xi-tame-estimate} requires a tame version of the elliptic estimates above.
\begin{lem}\label{lem:tame-elliptic-estimate}
    Let $l>\frac{1}{2\kappa}+2$ be an integer and $t\in[0,T]$. For any sufficiently smooth scalar $\tilde r$,
    \begin{equation}\label{eq:tame-elliptic-esimtate}
    \norm{\tilde{r}}_{H^{2l,l+\frac{1}{2\kappa}-\frac{1}{2}}}\lesssim\norm{\tilde{L}_1^l\tilde{r}}_{H^{0,\frac{1}{2\kappa}-\frac{1}{2}}}+\norm{\tilde r}_{H^{0,\frac{1}{2\kappa}-\frac{1}{2}}}+\norm{\tilde r}_{L^\infty}\norm{(r,\vec u)}_{H^{2l+2,l+\frac{1}{2\kappa}+\frac{1}{2}}\times H^{2l+2,l+1+\frac{1}{2\kappa}}},
    \end{equation}
    for any $4$-vector field $\tilde u$ orthogonal to $u$,
    \begin{equation}\label{eq:tame-div-curl-estimate}
        \norm{\tilde{u}}_{H^{2l,l+\frac{1}{2\kappa}}}\lesssim\norm{(\tilde{L}_2+\tilde L_3)^l\tilde{u}}_{H^{0,\frac{1}{2\kappa}}}+\norm{\tilde u}_{H^{0,\frac{1}{2\kappa}}}+\norm{\tilde u}_{L^\infty}\norm{(r,\vec u)}_{H^{2l+2,l+\frac{1}{2\kappa}+\frac{1}{2}}\times H^{2l+2,l+1+\frac{1}{2\kappa}}},
    \end{equation}
    and for any sufficiently smooth scalar $\tilde\pi$,
    \begin{equation}\label{eq:tame-elliptic-estimate-pi}
        \norm{\tilde{\pi}}_{H^{2l,l+\frac{1}{2\kappa}}}\lesssim\norm{\tilde{L}_4^l\tilde{\pi}}_{H^{0,\frac{1}{2\kappa}}}+\norm{\tilde \pi}_{H^{0,\frac{1}{2\kappa}}}+\norm{\tilde \pi}_{L^\infty}\norm{(r,\vec u)}_{H^{2l+2,l+\frac{1}{2\kappa}+\frac{1}{2}}\times H^{2l+2,l+1+\frac{1}{2\kappa}}}.
    \end{equation}
Here, the implicit constants depend on $\norm{r}_{W^{2,\infty}}$ and $\norm{\vec u}_{W^{1,\infty}}$.
\end{lem}
\begin{proof}
    The estimates \eqref{eq:tame-div-curl-estimate} and \eqref{eq:tame-elliptic-estimate-pi} are proved with a similar argument as \eqref{eq:tame-elliptic-esimtate}. Therefore, we only prove \eqref{eq:tame-elliptic-esimtate} in detail. We set $\alpha:=\frac{1}{2\kappa}-\frac{1}{2}$ throughout the proof and let $C$ be the implicit constants depending on $\norm{r}_{W^{2,\infty}_x}$ and $\norm{\vec u}_{W^{1,\infty}_x}$, but not $\varepsilon$ chosen below. By a partition of unity, we can localize the estimate to a small ball near the boundary $\Gamma_t$ and assume $\tilde r$ is supported in the small ball. The estimate is standard otherwise.

    Let $x_0\in \Gamma_t$. Near $x_0$, we choose the coordinate $e_3 = \frac{(\nabla r)(t,x_0)}{|(\nabla r)(t,x_0)|}$ and $e_1,e_2$ an orthonormal basis of the tangent space of the free boundary $\Gamma_t$ at $x_0$. For sufficiently small $\delta >0$, on the ball $B_\delta(x_0)$, we have
    \begin{equation}\label{eq:tame-localization}
    r\le\varepsilon,\qquad\abs{(\nabla r)(t,\cdot)-(\nabla r)(t,x_0)}<\varepsilon,\qquad\text{hence}\quad\abs{\p'r}\le\varepsilon,\quad \p_3r\ge\tfrac12c_0
    \end{equation}
    for some $c_0>0$, where we let $\p'\in\{\p_1,\p_2\}$ and $\varepsilon>0$ be sufficiently small. For any integer $\lambda \geq 1$ and $0\leq a \leq \lambda$, we set 
    \begin{equation}\label{eq:E-lambda-a}
        E^{\lambda, a} := \norm{r^a\p^{\lambda + a}\tilde r}_{H^{0,\alpha}}.
    \end{equation}

    We will make use of the following $3$ preliminary results in the proof:

    \noindent\textbf{Preliminary 1: The remainder class.} Let $\lambda\ge1$ be an integer. Consider a differential operator of the form
    \begin{equation}\label{Eq:TameRemainderStructure}
        D=\calF\,r^\sigma\prod_{i_1=1}^{I_1}\p^{a_{i_1}}r\prod_{i_2=1}^{I_2}\p^{a_{i_2}}\vec u\;\p^c,
    \end{equation}
    where $\calF$ is a smooth function in a neighborhood of the range of $(r,\vec u)$, and $\sigma,c,I_1,I_2\ge0$, $a_{i_1},a_{i_2}\ge1$ are integers with $\sigma\le\lambda$. Set
    \begin{equation*}
        K:=\sum_{i_1}a_{i_1}+\sum_{i_2}a_{i_2}+c,\qquad R:=I_1+\#\{i_1:a_{i_1}\ge2\}+I_2.
    \end{equation*}
    We say that a differential operator belongs to the class $\calR^\lambda$ if it is a finite sum of operators $D$ of the form \eqref{Eq:TameRemainderStructure} satisfying
    \begin{itemize}
        \item[(i)] $c\le\lambda+\sigma-1$,
        \item[(ii)] $\lambda + \sigma -1\leq K\le\lambda+\sigma-1+R$.
    \end{itemize}
    For a function $f$, we use $\calR^\lambda(f)$ to denote a finite sum of norms $\norm{Df}_{H^{0,\alpha}}$ with $D\in\calR^\lambda$. The number of terms and the operators may change from line to line. Condition (i) ensures that $D\tilde r$ is of lower order than the quantities $E^{\lambda,a}$ of level $\lambda$: compared with $E^{\lambda, a}$, it carries at least one derivative less, or equivalently one power of $r$ to spare, which is traded for a factor $\varepsilon$ through \eqref{eq:tame-localization}. Condition (ii) says that the excess $K-(\lambda+\sigma-1)$ of derivatives can be hidden by regarding $\p^{a_{i_1}}r$ as a derivative of $\p r$ or of $\p^2r$, and $\p^{a_{i_2}}\vec u$ as a derivative of $\p\vec u$, so that at the end of the proof Corollary~\ref{Cor:Multilinear-Interpolation} can be applied with exactly $\lambda+\sigma-1$ derivatives.

    We record one property of the class, used in the inductive step: if $D\in\calR^{\lambda}$, then $D\tilde L_1\in\calR^{\lambda+1}$. Indeed, in each term of the Leibniz expansion of $D\tilde L_1$, the total count $K$ increases by $2$, and either $(\sigma,c)$ increases by $(1,2)$ with $R$ unchanged, when $\p^c$ falls entirely on $\tilde r$ in the term $rH^{ij}\p_i\p_j$, or $\sigma$ and $c$ increase by at most $1$ each and $R$ by at least $1$, since a factor $\p^ar$ or $\p^a\vec u$ is produced. In both cases (i)--(ii) for $\lambda+1$ follow from $D\in\calR^\lambda$.

    \noindent\textbf{Preliminary 2: Homogeneous weighted embedding.} Since $r\le\varepsilon$ on the support of $\tilde r$, by Lemma~\ref{lem:WeightedEmbedding}, we obtain
    \begin{equation*}
        \norm{r^\sigma\p^b\tilde r}_{H^{0,\alpha}}\lesssim\norm{r^{\sigma+1}\p^b\tilde r}_{H^{0,\alpha}}+\norm{r^{\sigma+1}\p^{b+1}\tilde r}_{H^{0,\alpha}}.
    \end{equation*}
    The first term on the right is at most $\varepsilon\norm{r^\sigma\p^b\tilde r}_{H^{0,\alpha}}$ and can be absorbed into the left for $\varepsilon$ sufficiently small. Hence, for any $\sigma\geq 0$ and integer $b\ge0$,
    \begin{equation}\label{Eq:HomogeneousWeightedEmbedding}
        \norm{r^\sigma\p^b\tilde r}_{H^{0,\alpha}}\lesssim\norm{r^{\sigma+1}\p^{b+1}\tilde r}_{H^{0,\alpha}},
    \end{equation}
    and in particular
    \begin{equation}\label{Eq:HardyChain}
        E^{\lambda,a}\lesssim E^{\lambda,a+1}\qquad(0\le a<\lambda).
    \end{equation}
    Now, let $0\le\sigma\le\lambda$ and $0\le c\le\lambda+\sigma-1$ be integers. Applying \eqref{Eq:HomogeneousWeightedEmbedding} and $r\le\varepsilon$ gives
    \begin{equation}\label{Eq:LowerOrderBound}
        \norm{r^\sigma\p^c\tilde r}_{H^{0,\alpha}}\lesssim \varepsilon\,E^{\lambda,\lambda}.
    \end{equation}

    \noindent\textbf{Preliminary 3: Commutator formula.} We set $\p_N:=H^{ij}\p_ir\,\p_j$ for the derivative along the $H$-gradient of $r$, so that $\tilde L_{1,b}=\tilde L_1+b(\kappa+1)\p_N$. For integers $k\ge1$ and $0\le b\le k$, the Leibniz rule gives
    \begin{equation}\label{Eq:TameCommutatorFormula}
        \p'^{k-b}\p_3^b(\tilde L_1\tilde r)
        =\tilde L_{1,b}(\p'^{k-b}\p_3^b\tilde r)+\p'r\,\p'^{k-b-1}\p_3^b\p^2\tilde r+\p'r\,\p'^{k-b}\p_3^b\p\tilde r+\p_3r\,\p'^{k-b+1}\p_3^{b-1}\p\tilde r+S^k\tilde r,
    \end{equation}
    up to coefficients depending on $\vec u$. Above, the second, third and fourth terms come from the terms in which the factor $r$ of $rH^{ij}\p_i\p_j$ receives one derivative and the other coefficients of $\tilde L_1$ receive none (the third and fourth terms are absent when $b=0$), and $S^k\tilde r$ collects the remaining terms in which some derivative falls on the coefficients of $\tilde L_1$. In more detail, the terms with exactly one derivative on the coefficients, falling on the factor $r$ of $rH^{ij}\p_i\p_j$, are
    \begin{equation*}
        (\kappa+1)\big[(k-b)\,\p'r\,H^{ij}\p_i\p_j(\p'^{k-b-1}\p_3^b\tilde r)+b\,\p_3r\,H^{ij}\p_i\p_j(\p'^{k-b}\p_3^{b-1}\tilde r)\big].
    \end{equation*}
    The first term above is the second term in \eqref{Eq:TameCommutatorFormula}. In the second term above, the terms with $i\in\{1,2\}$ belong to the fourth term in \eqref{Eq:TameCommutatorFormula}, and by $\p_3r\,H^{3j}\p_j=\p_N-\p'_\gamma r\,H^{\gamma j}\p_j$, the terms with $i=3$ equal the sum of $b(\kappa+1)\p_N(\p'^{k-b}\p_3^b\tilde r)$, which combines with $\tilde L_1(\p'^{k-b}\p_3^b\tilde r)$ into $\tilde L_{1,b}(\p'^{k-b}\p_3^b\tilde r)$, and of the third term in \eqref{Eq:TameCommutatorFormula}.

    Every term of $S^k\tilde r$ is of the form $D\tilde r$ with $D$ as in \eqref{Eq:TameRemainderStructure}, where $\sigma\in\{0,1\}$, $K=k+2\geq(k+1)+\sigma-1$, $c\le k+\sigma$ and $R\ge2-\sigma$: if $\sigma=1$, the factor $r$ of $rH^{ij}$ is not differentiated and at least one derivative falls on $H^{ij}$, producing a factor $\p^{a_{i_2}}\vec u$; if $\sigma=0$, either a derivative falls on $H^{ij}\p_ir$, producing a factor $\p^{a_{i_1}}r$ with $a_{i_1}\ge2$ or a factor $\p^{a_{i_2}}\vec u$ next to $\p r$, or the factor $r$ of $rH^{ij}$ is differentiated and at least one further derivative falls on $rH^{ij}$. In all cases (i)--(ii) hold with $\lambda=k+1$, so that $S^k\in\calR^{k+1}$.

    \noindent\textbf{The induction.} We show by induction on $\lambda\ge1$ that
    \begin{equation}\label{Eq:TameEllipticInductionHypothesis}
        \mathbf H(\lambda):\qquad E^{\lambda,\lambda}\lesssim\norm{\tilde L_1^\lambda\tilde r}_{H^{0,\alpha}}+\calR^\lambda(\tilde r).
    \end{equation}
    The base case $\lambda=1$ follows from \eqref{Eq:BasicEllipticEstimate}. Let $\lambda\ge2$ and assume $\mathbf H(\lambda-1)$. Since $\tilde L_1\tilde r$ has the same support as $\tilde r$, for $0\le a\le\lambda-1$, the inequality \eqref{Eq:HardyChain} and the inductive hypothesis $\mathbf H(\lambda-1)$ apply to it and give:
    \begin{equation*}
        \norm{r^a\p^{\lambda-1+a}\tilde L_1\tilde r}_{H^{0,\alpha}}\lesssim \norm{r^{\lambda-1}\p^{2\lambda-2}\tilde L_1\tilde r}_{H^{0,\alpha}}\lesssim \norm{\tilde L_1^\lambda\tilde r}_{H^{0,\alpha}}+\calR^{\lambda-1}(\tilde L_1\tilde r) \lesssim  \norm{\tilde L_1^\lambda\tilde r}_{H^{0,\alpha}}+\calR^\lambda(\tilde r),
    \end{equation*}
    where the last inequality is by Preliminary~1. Hence, for $0\le a\le\lambda-1$,
    \begin{equation}\label{Eq:InductionTails}
        \norm{r^a\p^{\lambda-1+a}\tilde L_1\tilde r}_{H^{0,\alpha}}\lesssim \norm{\tilde L_1^\lambda\tilde r}_{H^{0,\alpha}}+\calR^\lambda(\tilde r).
    \end{equation}

    We now prove successively for $a=0,1,\ldots,\lambda$,
    \begin{equation}\label{Eq:TameEllipticInductionOnWeight}
        E^{\lambda,a}\lesssim \varepsilon\,E^{\lambda,\lambda}+ \norm{\tilde L_1^\lambda\tilde r}_{H^{0,\alpha}}+\calR^\lambda(\tilde r).
    \end{equation}
    Then, choosing $a=\lambda$ and $\varepsilon$ so small that the first term can be absorbed into the left-hand side, this gives $\mathbf H(\lambda)$.

    \smallskip
    \noindent\emph{Step 1: the case $a=0$, proved by a sub-induction on the number $b$ of normal derivatives.} For $0\le b\le\lambda$, set $F^b:=\norm{\p_3^b\p'^{\lambda-b}\tilde r}_{H^{0,\alpha}}$, then $E^{\lambda,0}\le\sum_bF^b$. We show by induction on $b$ that
    \begin{equation}\label{Eq:TameEllipticInductionHypothesisOnB}
        F^b\lesssim \varepsilon\,E^{\lambda,0} + \varepsilon\,E^{\lambda,\lambda} + \norm{\tilde L_1^\lambda\tilde r}_{H^{0,\alpha}}+\calR^\lambda(\tilde r).
    \end{equation}
    Summing over $b$ and absorbing the term $\varepsilon E^{\lambda,0}$ for $\varepsilon$ sufficiently small gives \eqref{Eq:TameEllipticInductionOnWeight} for $a=0$.

    \emph{Cases $b=0,1$.} Let $w:=\p'^{\lambda-1}\tilde r$. By \eqref{Eq:BasicEllipticEstimate} and Lemma~\ref{lem:WeightedEmbedding}, $F^0+F^1\lesssim \norm{\p w}_{H^{0,\alpha}}\lesssim \norm{w}_{H^{2,\alpha+1}}\lesssim \norm{\tilde L_1w}_{H^{0,\alpha}}+\norm{w}_{H^{0,\alpha}}$, and $\norm{w}_{H^{0,\alpha}}\lesssim \varepsilon E^{\lambda,\lambda}$ by \eqref{Eq:LowerOrderBound}. Applying \eqref{Eq:TameCommutatorFormula} with $k=\lambda-1$, $b=0$,
    \begin{equation*}
        \tilde L_1w=\p'^{\lambda-1}\tilde L_1\tilde r-\p'r\,\p'^{\lambda-2}\p^2\tilde r-S^{\lambda-1}\tilde r.
    \end{equation*}
    The first term is bounded by \eqref{Eq:InductionTails} with $a=0$. The second carries $\lambda$ derivatives of $\tilde r$ and the factor $\p'r$, so its norm is at most $\varepsilon E^{\lambda,0}$. The third has norm bounded by $\calR^\lambda(\tilde r)$, since $S^{\lambda-1}\in\calR^\lambda$ by Preliminary~3. This gives \eqref{Eq:TameEllipticInductionHypothesisOnB} for $b=0,1$.

    \emph{Case $b+1$, $1\le b\le\lambda-1$.} Assume \eqref{Eq:TameEllipticInductionHypothesisOnB} for all normal derivative counts up to $b$, and set $w:=\p_3^{b}\p'^{\lambda-1-b}\tilde r$. The basic elliptic estimate for $\tilde L_{1,b}$ (Remark~\ref{Rmk:ModifiedBasicEllipticEstimates}) gives $F^{b+1}\lesssim \norm{\p w}_{H^{0,\alpha}}\lesssim \norm{\tilde L_{1,b}w}_{H^{0,\alpha}}+\norm{w}_{H^{0,\alpha}}$, and $\norm{w}_{H^{0,\alpha}}\lesssim\varepsilon E^{\lambda,\lambda}$ is bounded as before. Applying \eqref{Eq:TameCommutatorFormula} with $k=\lambda-1$,
    \begin{equation*}
        \tilde L_{1,b}w=\p'^{\lambda-1-b}\p_3^b\tilde L_1\tilde r-\p'r\,\p'^{\lambda-2-b}\p_3^b\p^2\tilde r-\p'r\,\p'^{\lambda-1-b}\p_3^b\p\tilde r-\p_3r\,\p'^{\lambda-b}\p_3^{b-1}\p\tilde r-S^{\lambda-1}\tilde r.
    \end{equation*}
    Taking the norm $\norm{\cdot}_{H^{0,\alpha}}$, on the right hand side,
    \begin{itemize}
        \item \textit{First term.} Bounded by \eqref{Eq:InductionTails} with $a=0$.
        \item \textit{Second and third terms.} Both carry $\lambda$ derivatives of $\tilde r$ and the factor $\p'r$; since $\abs{\p'r}\le\varepsilon$, they are bounded by $\varepsilon E^{\lambda,0}$.
        \item \textit{Fourth term.} It contains at most $b$ normal derivatives, which is bounded by the induction hypothesis on $b$.
        \item \textit{Fifth term.} Bounded by $\calR^\lambda(\tilde r)$, by Preliminary~3.
    \end{itemize}
    This gives \eqref{Eq:TameEllipticInductionHypothesisOnB} for $b+1$, which closes the sub-induction.

    \smallskip
    \noindent\textit{Step 2: from the weight $a$ to the weight $a+1$.} Let $0\le a\le\lambda-1$ and assume \eqref{Eq:TameEllipticInductionOnWeight} for $E^{\lambda,a}$. Fix a multi-index of length $m:=\lambda-1+a$ and let $g:=\p^m\tilde r$, $v:=r^ag$, then $E^{\lambda,a+1}$ is the sum over these multi-indices of $\norm{r^{a+1}\p^2g}_{H^{0,\alpha}}$. The basic estimate \eqref{Eq:BasicEllipticEstimate} applied to $v$ gives
    \begin{equation}\label{Eq:WeightStepBasic}
        \norm{r\p^2v}_{H^{0,\alpha}}\lesssim \norm{\tilde L_1v}_{H^{0,\alpha}}+\norm{v}_{H^{0,\alpha}},
    \end{equation}
    and $\norm{v}_{H^{0,\alpha}}=\norm{r^a\p^{\lambda-1+a}\tilde r}_{H^{0,\alpha}}\lesssim \varepsilon E^{\lambda,\lambda}$ by \eqref{Eq:LowerOrderBound}. Expanding by the Leibniz rule on the left, 
    \begin{equation*}
        r\p^2 v = r\p^2(r^ag)=r^{a+1}\p^2g+2a\,r^a\p r\,\p g+a(a-1)\,r^{a-1}\p r\,\p r\,g+a\,r^a\p^2r\,g.
    \end{equation*}
    The last three terms have norms bounded by $CE^{\lambda,a}$, $CE^{\lambda,a-1}\le CE^{\lambda,a}$ (by \eqref{Eq:HardyChain}) and $\calR^\lambda(\tilde r)$, respectively. On the right, a direct computation gives
    \begin{equation}\label{Eq:WeightCommutator}
        \tilde L_1v=r^a\tilde L_1g+(\kappa+1)\Big(2a\,r^a\p_Ng+a\Big(a-1+\frac1\kappa\Big)\p_Nr\;r^{a-1}g+a\,r^aH^{ij}\p_i\p_jr\,g\Big),
    \end{equation}
    with $\p_N$ defined as in Preliminary~3, and the three terms in the bracket are bounded by $CE^{\lambda,a}$, $CE^{\lambda,a-1}\le CE^{\lambda,a}$ (by \eqref{Eq:HardyChain}) and $\calR^\lambda(\tilde r)$ as the last three terms on the left. Finally, $r^a\tilde L_1g=r^a\p^m\tilde L_1\tilde r-r^a[\p^m,\tilde L_1]\tilde r$. The first term is bounded by \eqref{Eq:InductionTails}. Expanding the commutator $[\p^m,\tilde L_1]\tilde r$ by the Leibniz rule, the terms with exactly one derivative on the coefficients of $\tilde L_1$ have bounded coefficients and carry at most $\lambda+a$ derivatives of $\tilde r$, so after multiplication by $r^a$, their norms are bounded by $CE^{\lambda,a}$. The terms with at least two derivatives on the coefficients are among those described in Preliminary~3 for $S^k\tilde r$, whose count does not depend on the direction of the derivatives. With $k=m$, they are of the form $D\tilde r$, $D$ as in \eqref{Eq:TameRemainderStructure}, with $\sigma\in\{0,1\}$, $K=\lambda+a+1$, $c\le\lambda+a-1$ and $R\ge2-\sigma$. Multiplication by $r^a$ replaces $\sigma$ by $\sigma+a\le\lambda$ and leaves $c$, $K$, $R$ unchanged, so that (i)--(ii) hold for $r^aD$ with this $\lambda$, and thus the norms of these terms are bounded by $\calR^\lambda(\tilde r)$.


    \noindent\textbf{Interpolation of $\calR^l(\tilde r)$.} Since $\calR^l(\tilde r)$ is a finite sum of norms $\norm{D\tilde r}$ with $D\in\calR^l$, and the number of terms depends only on $l$, it suffices to bound a single such term, with $D$ as in \eqref{Eq:TameRemainderStructure}. The coefficient $\calF$ is bounded and we drop it for simplicity. We claim that
    \begin{equation}\label{Eq:RemainderBound}
        \norm{D\tilde r}_{H^{0,\alpha}}\lesssim \varepsilon\norm{r^l\p^{2l}\tilde r}_{H^{0,\alpha}}+\norm{\tilde r}_{L^\infty}\Big(\norm{r}_{H^{2l+2,l+1+\alpha}}+\norm{\vec u}_{H^{2l+2,l+\frac32+\alpha}}\Big).
    \end{equation}
    Setting $k:=l+\sigma-1$, by (ii), we have $0\le K-k\le R$, so we may read $K-k$ of the factors $\p^ar$, $\p^a\vec u$ as derivatives of $\p r$, $\p^2r$ or $\p\vec u$, which lowers the number of derivatives by one or two each, until the count is exactly $k$. The remaining factors are read as they stand. Then $D\tilde r$ takes the form
    \begin{equation}\label{Eq:RemainderReduced}
        r^{\sigma}\prod_{j}\p^{b_j}g_j\;\p^c\tilde r,\qquad g_j\in\{r,\p r,\p^2r,\vec u,\p\vec u\},\quad b_j\ge0,\quad\sum_jb_j+c=k,
    \end{equation}
    and $c\le k$ by (i). Set $\tau:=\sigma+\alpha>-\frac12$, then $k-\tau=l-1-\alpha>\frac32$ by the assumption on $l$. The factors with $b_j=0$ are bounded by $\norm{r}_{W^{2,\infty}}$ or $\norm{\vec u}_{W^{1,\infty}}$, and $\tilde r$ by $\norm{\tilde r}_{L^\infty}$ if $c=0$. We take them out of the norm. Letting $g$ range over the distinct functions among the remaining $g_j$ and $\tilde r$, Corollary~\ref{Cor:Multilinear-Interpolation} gives
    \begin{equation}\label{Eq:RemainderInterpolated}
        \norm{D\tilde r}_{H^{0,\alpha}}\lesssim\sum_g\Big(\prod_{g'\ne g}\norm{g'}_{L^\infty}\Big)\norm{r^\tau\p^kg}_{L^2},
    \end{equation}
    where the implicit constant depends on $\norm{r}_{W^{2,\infty}}$ and $\norm{\vec u}_{W^{1,\infty}}$ but not on $\tilde r$, since $\tilde r$ occurs in at most one factor. The factor $\norm{\tilde r}_{L^\infty}$ appears in every term except the one with $g=\tilde r$. For $g=\tilde r$, the term is $\norm{r^\sigma\p^{l+\sigma-1}\tilde r}_{H^{0,\alpha}}\lesssim \varepsilon E^{l,l}=\varepsilon \norm{r^l\p^{2l}\tilde r}_{H^{0,\alpha}}$ by \eqref{Eq:LowerOrderBound}. For the background states, we use $\norm{r^\tau\p^kg}_{L^2}\lesssim\norm{g}_{H^{N,\tau+N-k}}$ for $N\ge k$ by Lemma~\ref{lem:WeightedEmbedding}, where $\tau+N-k=\alpha+N-l+1$: for $g=\p^2r$ ($N=2l$), $g=\p r$ ($N=2l+1$) and $g=r$ ($N=2l+2$), it gives $\norm{r}_{H^{2l+2,l+1+\alpha}}$, using $r\lesssim1$ in the last two cases; for $g=\p\vec u$ ($N=2l+1$) and $g=\vec u$ ($N=2l+2$) it gives $\norm{\vec u}_{H^{2l+2,l+2+\alpha}}$ and $\norm{\vec u}_{H^{2l+2,l+3+\alpha}}$, both bounded by $\norm{\vec u}_{H^{2l+2,l+\frac32+\alpha}}$. This proves \eqref{Eq:RemainderBound}.

    Inserting \eqref{Eq:RemainderBound} into $\mathbf H(l)$ and choosing $\varepsilon$ small to absorb $C\varepsilon E^{l,l}$ into the left-hand side, we obtain
    \begin{equation}\label{Eq:BoundaryBallEstimate}
        \norm{r^l\p^{2l}\tilde r}_{H^{0,\alpha}}\lesssim \norm{\tilde L_1^l\tilde r}_{H^{0,\alpha}} + \norm{\tilde r}_{L^\infty_x}\Big(\norm{r}_{H^{2l+2,l+1+\alpha}}+\norm{\vec u}_{H^{2l+2,l+\frac32+\alpha}}\Big),
    \end{equation}
    which finishes the proof.
\end{proof}

We also have the boundedness of the elliptic operators.
\begin{lem}\label{lem:boundedness-elliptic-operator}
For any integer $l>\frac{1}{2\kappa}+2$, $t\in[0,T]$, and any sufficiently smooth scalar $\tilde r$,
\begin{equation}\label{eq:boundedness-elliptic-operator1}
    \norm{\tilde{L}_1^l\tilde{r}}_{H^{0,\frac{1}{2\kappa}-\frac{1}{2}}}\lesssim \norm{\tilde{r}}_{H^{2l,l+\frac{1}{2\kappa}-\frac{1}{2}}}+\norm{\tilde r}_{L^\infty}\norm{(r,\vec u)}_{H^{2l+2,l+\frac{1}{2\kappa}+\frac{1}{2}}\times H^{2l+2,l+1+\frac{1}{2\kappa}}},
\end{equation}
for any $4$-vector field $\tilde u$ orthogonal to $u$,
\begin{equation}\label{eq:boundedness-elliptic-operator2}
    \norm{(\tilde{L}_2+\tilde L_3)^l\tilde{u}}_{H^{0,\frac{1}{2\kappa}}} \lesssim \norm{\tilde{u}}_{H^{2l,l+\frac{1}{2\kappa}}} + \norm{\tilde u}_{L^\infty}\norm{(r,\vec u)}_{H^{2l+2,l+\frac{1}{2\kappa}+\frac{1}{2}}\times H^{2l+2,l+1+\frac{1}{2\kappa}}},
\end{equation}
and for any sufficiently smooth scalar $\tilde\pi$,
\begin{equation}\label{eq:boundedness-elliptic-operator3}
    \norm{\tilde{L}_4^l\tilde{\pi}}_{H^{0,\frac{1}{2\kappa}}}\lesssim\norm{\tilde{\pi}}_{H^{2l,l+\frac{1}{2\kappa}}}+\norm{\tilde \pi}_{L^\infty}\norm{(r,\vec u)}_{H^{2l+2,l+\frac{1}{2\kappa}+\frac{1}{2}}\times H^{2l+2,l+1+\frac{1}{2\kappa}}}.
\end{equation}
The implicit constants depend on $\norm{r}_{W^{1,\infty}}$ and $\norm{\vec u}_{L^\infty}$.
\end{lem}

\begin{proof}
By Remark~\ref{rmk:class-of-L}, $\tilde L_1^l$, $(\tilde L_2+\tilde L_3)^l$ and $\tilde L_4^l$ are of class $\calM_{l,2l}$ of order $l$, with purely spatial elementary operators whose coefficients are smooth functions of $\vec u$ and whose background factors involve only $r$ and $\vec u$. Hence, $\tilde L_1^l\tilde r$, $(\tilde L_2+\tilde L_3)^l\tilde u$ and $\tilde L_4^l\tilde\pi$ are finite sums of monomials in $\tilde r$, $\tilde u$, $\tilde\pi$ of class $\calM^s_{l,2l}$ in the slots $s=r,u,\pi$ respectively, lower order at level $2l$, with $K=2l$ and no time derivative. The three estimates are \eqref{eq:multilinear-tame-2}, with the norms of $a_0$ absent by the last assertion of Lemma~\ref{lem:multilinear-tame-estimate}.
\end{proof}

\subsection{The commutator estimate}\label{subsec:commutator}
To obtain the first-order tame estimate, we commute the modified linearized equations with the elliptic operator $L^l$, and apply the energy estimate \eqref{eq:energy-estimate-modified} at the basic level together with the tame elliptic estimates (Lemma~\ref{lem:tame-elliptic-estimate}). In this subsection, we show that each component of the commutator $[\calL(\vec\varphi),L^l]\tilde\phi$ is a finite sum of monomials that are lower-order at level $2l$ in the corresponding slot, in the sense of Section~\ref{subsec:bookkeeping}. We then estimate the commutator by Lemma~\ref{lem:multilinear-tame-estimate}.

\begin{lem}\label{lem:commutator-operator-scheme-basic}
Let $\tilde\phi=(\tilde r,\tilde u,\tilde \pi)$ be a solution to the modified linearized equations with source terms $\tilde F=(\tilde f,\tilde g,\tilde h)$. Then
\begin{equation}\label{op:commutator-operator-basic}
    [\calL(\vec\varphi),L]\tilde\phi
    =\begin{pmatrix}
        D_{11} & D_{12} & D_{13}\\
        D_{21} & D_{22} & D_{23}\\
        D_{31} & D_{32} & D_{33}
    \end{pmatrix}\tilde\phi
    + \begin{pmatrix}
        D_{11}' & D_{12}' & D_{13}'\\
        D_{21}' & D_{22}' & D_{23}'\\
        D_{31}' & D_{32}' & D_{33}'
    \end{pmatrix}\tilde F,
\end{equation}
where, writing $s_1=r$, $s_2=u$, $s_3=\pi$ for the slots and $\tilde\phi_j$, $\tilde F_j$ for the components of $\tilde\phi$, $\tilde F$,
\begin{enumerate}
    \item[(a)] for all $i,j$, $D_{ij}\tilde\phi_j$ is a finite sum of monomials in $\tilde\phi_j$ of class $\calM^{s_i}_{1,3}$ which are lower order at level $2$ in the slot $s_i$; every such monomial has $K\ge2$, and every such monomial containing a time derivative has $K=3$;
    \item[(b)] for all $i,j$, $D'_{ij}\tilde F_j$ is a finite sum of monomials in $\tilde F_j$ of class $\calM^{s_i}_{1,2}$ which are lower order at level $2$ in the slot $s_i$, contain no time derivative, and have $K=2$;
    \item[(c)] $D'_{11}=D'_{13}=D'_{22}=D'_{31}=D'_{33}=0$.
\end{enumerate}
\end{lem}

\begin{nota}\label{nota:simeq}
In the proof of Lemma~\ref{lem:commutator-operator-scheme-basic}, in order to avoid keeping track of a large number of terms, we fix the slot $s$ under consideration and write $A\simeq B$ if $A-B$ is a finite sum of monomials as in (a) and (b) of the lemma for that slot, and such terms are absorbed into the right-hand side of \eqref{op:commutator-operator-basic}.
\end{nota}

\begin{proof}
    Recall that $[\calL(\vec\varphi),L]\tilde\phi$ is given by \eqref{eq:commutator-operator-one-elliptic}, with $\calR_{\tilde r}$, $\calR_{\tilde u}$, $\calR_{\tilde\pi}$ defined in \eqref{eq:commutator-modified-linearized-equation-elliptic-operator}. We first explain the counts in (a) and (b), which hold for every term produced below. Each of $\calR_{\tilde r}$, $\calR_{\tilde u}$, $\calR_{\tilde\pi}$ is a component of the commutator $[\calL(\vec\varphi),L]\tilde\phi=\calL(\vec\varphi)L\tilde\phi-L\calL(\vec\varphi)\tilde\phi$. In the two expressions on the right, the terms in which no derivative of $\calL(\vec\varphi)$ falls on the coefficients of $L$ and no derivative of $L$ falls on the coefficients of $\calL(\vec\varphi)$ are the same and cancel. Hence every term of the commutator carries at least one derivative of the background state or of the weight, and $c\le2$. Replacing $D_t\tilde r$, $D_t\tilde u$ and $D_t\tilde \pi$ by the modified linearized equations does not increase $c$. Since $\calL(\vec\varphi)$ has one derivative and $L$ has two, every term has $K\le3$ before any substitution. When $D_t\tilde r$, $D_t\tilde u$ or $D_t\tilde\pi$ is replaced by the modified linearized equations, a derivative is lost only in the terms generated by a source term or by the derivative-free parts of the coefficients $V_3$ and $Z_3$ of the $\tilde\pi$-equation: these terms contain no time derivative, and $K=2$, while the other terms introduced have $K=3$. Hence, every term in $\tilde\phi$ has $K\ge2$, and those containing a time derivative have $K=3$. Every term in $\tilde F$ arises from a source term, has $K=2$ and no time derivative. It remains to check the classes and the orders, which we do slot by slot.

    \noindent\textbf{The $r$-slot.} Consider the first term $[D_t,\tilde L_1]\tilde r=D_t(\tilde L_1\tilde r)-\tilde L_1(D_t\tilde r)$ of $\calR_{\tilde r}$ in \eqref{eq:commutator-modified-linearized-equation-elliptic-operator}. Substituting \eqref{op:weighted-laplacian}, the terms in which $D_t$ falls on $\tilde r$ cancel, and
    \[
        \begin{aligned}
            [D_t,\tilde L_1]\tilde r
            &=(\kappa + 1)D_t H^{ij}\Big(r\p_i\p_j\tilde r + \frac{1}{\kappa}\p_i r\p_j\tilde r\Big)
            + (\kappa+1)H^{ij}\Big(D_t r\,\p_i\p_j\tilde r + \frac{1}{\kappa}D_t\p_ir\,\p_j\tilde r\Big)\\
            &\quad + (\kappa + 1)H^{ij}\Big(r[D_t,\p_i\p_j]\tilde r + \frac{1}{\kappa}\p_i r[D_t,\p_j]\tilde r\Big).
        \end{aligned}
    \]
    Using $[D_t,\p_i]=-\p_iu^\mu\p_\mu$ and \eqref{eq:r-pi-defining} to trade $D_tr$ for one more power of $r$, every term is of one of the forms $r\,\p\vec u\,\p^2\tilde r$, $r\,\p^2\vec u\,\p\tilde r$, $\p r\,\p\vec u\,\p\tilde r$, $r\,\p r\,\p\vec u\,\p\tilde r$ up to smooth coefficients, where one of the derivatives may be $\p_t$. 
    Each has $\bar\sigma_r+N_r\ge1$, $c\le2$, $c-\bar\sigma_r\le1$ and $K=3$, so it is lower order in the $r$-slot and $[D_t,\tilde L_1]\tilde r\simeq0$. 

    We consider the second term of $\calR_{\tilde r}$ given as the following:
    \begin{equation}\label{Term:SecondinRtilder}
        \kappa r a_1\p_\mu[(\tilde L_2+\tilde L_3)\tilde u]^\mu - \tilde L_1(\kappa r a_1\p_\mu\tilde u^\mu).
    \end{equation}
    Here, by a simple counting, we see that a cancellation is needed: both expressions in \eqref{Term:SecondinRtilder} contain terms with three derivatives on $\tilde u$, which have $c=3$, and the terms of the forms $r\,\p r\,\p^2\tilde u$ and $\p r\,\p r\,\p\tilde u$, which have $\bar\sigma_r+N_r\ge1$ but $c-\bar\sigma_r=\frac32$. Neither kind is lower-order in the $r$-slot, so we must show that these terms cancel.

    To get the cancellation, we write the two four-divergences in terms of the spatial components of $\tilde u$: for $\p_\mu\tilde u^\mu$, we use the orthogonality condition $u^\mu\tilde u_\mu=0$, i.e. $\tilde u^0=\frac{u^i}{u^0}\tilde u_i$, and for $\p_\mu[(\tilde L_2+\tilde L_3)\tilde u]^\mu$ the identity $[(\tilde L_2+\tilde L_3)\tilde u]^0=\frac{u^i}{u^0}[(\tilde L_2+\tilde L_3)\tilde u]_i$, which follows from $A^i_0=-\frac{u^i}{u^0}$. Setting
    \begin{equation}\label{eq:hat-partial}
        \hat\p_i:=\p_i+\frac{u_i}{u^0}\p_t=H^{im}\p_m+\frac{u_i}{(u^0)^2}D_t ,
    \end{equation}
    which is the operator that arises when the four-divergence of a vector field orthogonal to $u$ is written in terms of its spatial components, we have
    \begin{equation}\label{eq:four-divergence}
        \p_\mu\tilde u^\mu=\hat\p_i\tilde u_i+\p_t\Big(\frac{u^i}{u^0}\Big)\tilde u_i,\qquad
        \p_\mu[(\tilde L_2+\tilde L_3)\tilde u]^\mu=\hat\p_i[(\tilde L_2+\tilde L_3)\tilde u]_i+\p_t\Big(\frac{u^i}{u^0}\Big)[(\tilde L_2+\tilde L_3)\tilde u]_i ,
    \end{equation}
    and the terms with $\p_t(\frac{u^i}{u^0})$ contribute admissible terms to \eqref{Term:SecondinRtilder}, of the forms $r\,\p_t\vec u\,r\p^2\tilde u$, $r\,\p_t\vec u\,\p r\,\p\tilde u$ and $r\,\p^2(r\,\p_t\vec u\,\tilde u)$.

    From \eqref{op:weighted-div} and \eqref{op:weighted-curl}, for a spatial index $i$,
    \begin{equation}\label{eq:L2L3-explicit}
        [(\tilde L_2+\tilde L_3)\tilde u]_i = (\kappa+1)H^{jk} \Big(r\p_j\partial_k\tilde u_i + \p_i r\p_j\tilde u_k + \Big(1+\frac1\kappa\Big)\p_j r\p_k\tilde u_i - \p_j r\partial_i\tilde u_k\Big).
    \end{equation}

    Applying $\hat\p_i$ to \eqref{eq:L2L3-explicit}, we retain only the terms in which the derivative falls on the weight $r$ or on $\tilde u$. The others, in which it falls on $H^{jk}$ or on $\p r$, are of the forms $r\,\p\vec u\,r\p^2\tilde u$, $r\,\p\vec u\,\p r\,\p\tilde u$, $r\,\p^2r\,\p\tilde u$, where one derivative may be $\p_t$, and are lower order in the $r$-slot. Among the retained terms, the one in which $D_t$ falls on the weight, $r\,D_tr\,H^{jk}\p_j\p_k\tilde u_i$, has order $\frac32$ in the $r$-slot. Replacing $D_tr$ by \eqref{eq:r-pi-defining}, it becomes $r^2\p\vec u\,\p^2\tilde u$, which is lower order. Using \eqref{eq:hat-partial}, we obtain
    \begin{equation}\label{eq:first-term-second}
        \kappa r a_1\p_\mu[(\tilde L_2+\tilde L_3)\tilde u]^\mu\simeq (\kappa+1)\kappa r a_1\Big(H^{im}H^{jk}T_1 + \frac{u^i}{(u^0)^2}H^{jk}T_2\Big),
    \end{equation}
    where
    \[
        \begin{aligned}
            T_1 &= r\p_j\p_k\p_m\tilde u_i + \Big(2+\frac1\kappa\Big)\p_j r\p_k\p_m\tilde u_i,\\
            T_2 & =r\p_j\p_kD_t\tilde u_i + \p_ir \p_jD_t\tilde u_k +\Big(1+\frac1\kappa\Big)\p_j r\p_kD_t\tilde u_i - \p_j r\p_i D_t\tilde u_k.
        \end{aligned}
    \]
    Here $T_1$ collects the terms in which the spatial part $H^{im}\p_m$ of $\hat\p_i$ falls on $r$ or $\tilde u$: the terms $H^{im}\p_mr\,H^{jk}\p_j\p_k\tilde u_i$ and $-H^{im}H^{jk}\p_jr\,\p_i\p_m\tilde u_k$, coming from the first and fourth terms of \eqref{eq:L2L3-explicit}, cancel after relabeling $(i,m)\leftrightarrow(j,k)$, and the second and third terms of \eqref{eq:L2L3-explicit} combine into the coefficient $2+\frac1\kappa$. All terms of $T_1$ are of top order. $T_2$ collects the terms in which $D_t$ falls on $\tilde u$, with $D_t$ commuted past the spatial derivatives. The commutators $[D_t,\p_j]$ and $[D_t,\p_j\p_k]$ produce terms of the forms $r^2\p\vec u\,\p^2\tilde u$, $r^2\p^2\vec u\,\p\tilde u$ and $r\,\p r\,\p\vec u\,\p\tilde u$, where one derivative of $\tilde u$ may be $\p_t$, and these are lower order in the $r$-slot. In $T_2$ itself, only the first term is of top order. The remaining three are of the form $r\,\p r\,\p D_t\tilde u$, which has order $\frac32$ in the $r$-slot. Replacing $D_t\tilde u$ by the modified linearized equation, they become terms of the forms $r\,\p r\,\p^2\tilde r$, $r\,\p r\,\p(r\p\tilde\pi)$ and $r\,\p r\,\p\tilde g$, together with terms in $\tilde r$, $\tilde u$, $\tilde\pi$ with one derivative less. All of these are lower order in the $r$-slot. The terms in $\tilde g$ have class $\calM^r_{1,2}$, order $\frac12$ and $K=2$, and are absorbed into $D_{12}'\tilde g$. Thus
    \[
        T_2\simeq r\p_j\p_kD_t\tilde u_i.
    \]

    Next, we turn to the second term of \eqref{Term:SecondinRtilder}. By \eqref{eq:four-divergence} and \eqref{eq:hat-partial},
    \[
        \tilde L_1(\kappa ra_1\p_\mu\tilde u^\mu) \simeq \tilde L_1\big(\kappa ra_1H^{im}\p_m\tilde u_i\big) + \tilde L_1\Big(\kappa ra_1\frac{u^i}{(u^0)^2}D_t\tilde u_i\Big)=: (i) +(ii).
    \]
    Retaining only the top-order contributions, we obtain
    \[
        (i)\simeq (\kappa + 1)\kappa a_1H^{im}H^{jk}T_3,\qquad
        T_3 = r^2\p_j\p_k\p_m\tilde u_i + \Big(2+\frac1\kappa\Big)r\p_j r\p_k\p_m\tilde u_i + \frac{1}{\kappa}\p_j r\p_k r\p_m\tilde u_i ,
    \]
    and similarly
    \[
        (ii)\simeq (\kappa + 1)\kappa r^2a_1\frac{u^i}{(u^0)^2}H^{jk}\p_j\p_kD_t\tilde u_i.
    \]
    Above, in $(i)$, $T_3$ collects the terms of order $\frac32$, and the others are lower order in the $r$-slot. In $(ii)$, the terms of the forms $r\,\p r\,\p D_t\tilde u$ and $\p r\,\p r\,D_t\tilde u$ have order $\frac32$ and become lower order after $D_t\tilde u$ is replaced by the modified linearized equation, exactly as for $T_2$. Other than that, the terms in which a derivative falls on a coefficient are lower order. Hence the only term of $(ii)$ which is not lower order is the one in which both derivatives fall on $D_t\tilde u_i$.

    Subtracting $(i)+(ii)$ from \eqref{eq:first-term-second}, the terms containing $D_t\tilde u$ cancel and $rT_1-T_3=-\frac1\kappa\p_jr\p_kr\p_m\tilde u_i$. Writing $\abs{\nabla r}_H^2:=H^{jk}\p_jr\p_kr$ and $\operatorname{div}_H\tilde u:=H^{im}\p_m\tilde u_i$, we obtain
    \begin{equation}\label{eq:residual-term-2}
        \kappa r a_1\p_\mu[(\tilde L_2+\tilde L_3)\tilde u]^\mu - \tilde L_1(\kappa r a_1\p_\mu\tilde u^\mu)
        \simeq -(\kappa+1)\,a_1\,\abs{\nabla r}_H^2\,\operatorname{div}_H\tilde u.
    \end{equation}
    The right-hand side has order $\frac32$ in the $r$-slot and is not admissible by itself.

    Now we turn to the third term of $\calR_{\tilde r}$,
    \[
        \p_\mu r\big[(\tilde L_2+\tilde L_3)\tilde u\big]^\mu - \tilde L_1(\p_\mu r\,\tilde u^\mu)
        = \hat\p_ir\,\big[(\tilde L_2+\tilde L_3)\tilde u\big]^i - \tilde L_1(\hat\p_ir\,\tilde u^i),
    \]
    where we used $\p_\mu r\,\tilde u^\mu=\hat\p_ir\,\tilde u^i$ and $\p_\mu r[(\tilde L_2+\tilde L_3)\tilde u]^\mu=\hat\p_ir\,[(\tilde L_2+\tilde L_3)\tilde u]^i$ follows from a similar computation as in \eqref{eq:four-divergence}. Substituting \eqref{eq:L2L3-explicit} into the first term and expanding the second by \eqref{op:weighted-laplacian}, the terms $r\,\hat\p_ir\,H^{jk}\p_j\p_k\tilde u^i$ cancel. Replacing $D_tr$ in $\hat\p_ir$ by \eqref{eq:r-pi-defining}, the terms of order $\frac32$ that remain are
    \[
        (\kappa+1)H^{jk}H^{im}\p_mr\Big(\p_ir\,\p_j\tilde u_k+\p_jr\,\p_k\tilde u_i-\p_jr\,\p_i\tilde u_k\Big),
    \]
    where the last two terms cancel after relabeling $j\leftrightarrow m$. Hence,
    \begin{equation}\label{eq:residual-term-3}
        \p_\mu r\big[(\tilde L_2+\tilde L_3)\tilde u\big]^\mu - \tilde L_1(\p_\mu r\,\tilde u^\mu)
        \simeq (\kappa+1)\,\abs{\nabla r}_H^2\,\operatorname{div}_H\tilde u.
    \end{equation}
    Adding \eqref{eq:residual-term-2} and \eqref{eq:residual-term-3}, the second and third terms of $\calR_{\tilde r}$ together satisfy
    \[
        \text{(second term)}+\text{(third term)}\simeq -(\kappa+1)(a_1-1)\abs{\nabla r}_H^2\operatorname{div}_H\tilde u\simeq 0,
    \]
    since $a_1-1=r\big(1+ra_0^{2+\frac1\kappa}\big)$.

    The remaining two terms of $\calR_{\tilde r}$, $V_1\tilde L_1\tilde r-\tilde L_1(V_1\tilde r)$ and $r^2Z_1\tilde L_4\tilde \pi-\tilde L_1(r^2Z_1\tilde \pi)$, are lower order in the $r$-slot by a direct counting. Finally, the only terms of $\calR_{\tilde r}$ involving $\tilde F$ arise from $\tilde g$ when $D_t\tilde u$ is replaced, so $D'_{11}=D'_{13}=0$.

    \noindent\textbf{The $u$-slot.} We consider the first term of $\calR_{\tilde u}$ in \eqref{eq:commutator-modified-linearized-equation-elliptic-operator}, i.e., the commutator $[D_t,\tilde L_2+\tilde L_3]\tilde u$, which is treated exactly as $[D_t,\tilde L_1]\tilde r$: expanding the material derivative in \eqref{eq:L2L3-explicit}, using \eqref{eq:r-pi-defining} for $D_tr$, all terms are of the forms $r\,\p\vec u\,\p^2\tilde u$, $r\,\p^2\vec u\,\p\tilde u$, $\p r\,\p\vec u\,\p\tilde u$, $r\,\p_t\vec u\,\p^2\tilde u$, $r\,\p\vec u\,\p r\,\p\tilde u$, where one of the derivatives of $\tilde u$ possibly being $\p_t$, which are lower-order in the $u$-slot. Hence, $[D_t,\tilde L_2+\tilde L_3]\tilde u\simeq0$.

    We set $c_1:=(1+\frac1\kappa)a_2$ for convenience and consider the second term of $\calR_{\tilde u}$,
    \begin{equation}\label{eq:commutator-u-second}  
        c_1\Delta^{\alpha\mu}\p_\mu(\tilde L_1\tilde r)
        - \big[(\tilde L_2+\tilde L_3)\big(c_1\Delta^{\cdot\mu}\p_\mu\tilde r\big)\big]^\alpha.
    \end{equation}
    Both expressions in \eqref{eq:commutator-u-second} contain terms of the forms $r\,\p^3\tilde r$ and $\p r\,\p^2\tilde r$, which are not lower order in the $u$-slot. So, cancellations are needed here. We use the two identities followed from a direct computation,
    \begin{equation}\label{eq:projection-identities}
        A^i_\alpha\Delta_i^{\ \mu}=\Delta_\alpha^{\ \mu},\qquad
        H^k_j\Delta_k^{\ \mu}\p_\mu=\hat\p_j.
    \end{equation}
    Substituting \eqref{op:weighted-laplacian} into the first expression in \eqref{eq:commutator-u-second} and \eqref{eq:L2L3-explicit} into the second, and retaining only the terms of order $\frac32$ in the $u$-slot, we obtain
    \[
        \begin{aligned}
        c_1\Delta_\alpha{}^{\mu}\p_\mu(\tilde L_1\tilde r)
        &\simeq(\kappa+1)c_1\Delta_\alpha{}^{\mu}H^{jk}\Big(r\p_j\p_k\p_\mu\tilde r+\p_\mu r\,\p_j\p_k\tilde r+\tfrac1\kappa\p_jr\,\p_k\p_\mu\tilde r\Big),\\
        \big[(\tilde L_2+\tilde L_3)\big(c_1\Delta^{\cdot\mu}\p_\mu\tilde r\big)\big]_\alpha
        &\simeq(\kappa+1)c_1\Big(\Delta_\alpha{}^{\mu}H^{jk}r\p_j\p_k\p_\mu\tilde r
        +\big(1+\tfrac1\kappa\big)\Delta_\alpha{}^{\mu}H^{jk}\p_jr\,\p_k\p_\mu\tilde r\\
        &\qquad\qquad\qquad
        +A^i_\alpha\p_ir\,\hat\p^{\,j}\p_j\tilde r
        -A^i_\alpha\p^{\,j}r\,\hat\p_j\p_i\tilde r\Big),
        \end{aligned}
    \]
    where the first identity in \eqref{eq:projection-identities} was used in the first two terms of the second expression, and the second identity in the last two. The terms with three derivatives on $\tilde r$ cancel, and the difference of the two expressions is, up to the factor $(\kappa+1)c_1$,
    \[
        \Delta_\alpha{}^{\mu}\p_\mu r\,H^{jk}\p_j\p_k\tilde r-\Delta_\alpha{}^{\mu}H^{jk}\p_jr\,\p_k\p_\mu\tilde r
        -A^i_\alpha\p_ir\,\hat\p^{\,j}\p_j\tilde r+A^i_\alpha\p^{\,j}r\,\hat\p_j\p_i\tilde r.
    \]
    In the first two terms we write $\Delta_\alpha{}^{\mu}\p_\mu=\p_\alpha+u_\alpha D_t$. Up to lower order terms, the contributions of $u_\alpha D_t$ are $u_\alpha D_tr\,H^{jk}\p_j\p_k\tilde r$, which becomes $r\,\p\vec u\,\p^2\tilde r$ after \eqref{eq:r-pi-defining}, and $u_\alpha H^{jk}\p_jr\,\p_kD_t\tilde r$, which has order $\frac32$ and becomes lower order after $D_t\tilde r$ is replaced by the modified linearized equation. This is where $\tilde f$ enters, through terms of the form $\p r\,\p\tilde f$, of class $\calM^u_{1,2}$, order $\frac12$ and $K=2$, which are absorbed into $D_{21}'\tilde f$. Since both expressions in \eqref{eq:commutator-u-second} are orthogonal to $u$, it suffices to consider a spatial index $\alpha=a$, for which $A^i_a=\delta^i_a$. What remains is
    \[
        \begin{aligned}
        &\p_ar\,H^{jk}\p_j\p_k\tilde r-H^{jk}\p_jr\,\p_k\p_a\tilde r
        -\p_ar\,\hat\p^{\,j}\p_j\tilde r+\p^{\,j}r\,\hat\p_j\p_a\tilde r\\
        &\qquad\simeq
        -\p_ar\,\tfrac{u^j}{u^0}\p_j\Big(\tfrac1{u^0}D_t\tilde r\Big)
        +\tfrac{u^j}{u^0}\p_jr\,\p_a\Big(\tfrac1{u^0}D_t\tilde r\Big).
        \end{aligned}
    \]
    Replacing $D_t\tilde r$ by the modified linearized equation as above, the right-hand side is lower order in the $u$-slot. Hence, the second term of $\calR_{\tilde u}$ is $\simeq0$.

    We set $c_2:=(2+\frac1\kappa)ra_3$ for convenience and consider the third term of $\calR_{\tilde u}$,
    \[
        c_2\Delta^{\alpha\mu}\p_\mu\tilde L_4\tilde\pi - \big[(\tilde L_2+\tilde L_3)\big(c_2\Delta^{\cdot\mu}\p_\mu\tilde\pi\big)\big]^\alpha.
    \]
    Its terms with three derivatives on $\tilde\pi$ are $(\kappa+1)c_2\Delta_\alpha^{\ \mu}H^{jk}r\p_j\p_k\p_\mu\tilde\pi$ on both sides, by the first identity in \eqref{eq:projection-identities}, and cancel. Because of the factor $r$ in $c_2$, all remaining terms are lower order in the $u$-slot. The terms containing $D_t\tilde\pi$ are replaced by the modified linearized equation. This produces terms of the forms $r\,\p r\,\p\tilde h$ and $\p r\,\p r\,\tilde h$, which have class $\calM^u_{1,2}$ and $K=2$ and go into $D_{23}'\tilde h$.

    The remaining terms of $\calR_{\tilde u}$, involving $V_2$, $\bar W_2$ and $Z_2$, are lower order in the $u$-slot by a direct counting, as for $\calR_{\tilde r}$. Finally, the only terms of $\calR_{\tilde u}$ involving $\tilde F$ arise from $\tilde f$ and $\tilde h$ when $D_t\tilde r$ and $D_t\tilde\pi$ are replaced, so $D'_{22}=0$.

    \noindent\textbf{The $\pi$-slot.} The terms of $\calR_{\tilde\pi}$ are treated as the corresponding terms of $\calR_{\tilde r}$, with $\tilde L_1$ replaced by $\tilde L_4$ and $\kappa a_1$ by $a_4$. In the second term $ra_4\p_\mu[(\tilde L_2+\tilde L_3)\tilde u]^\mu-\tilde L_4(ra_4\p_\mu\tilde u^\mu)$, the terms with three derivatives on $\tilde u$ cancel as in the $r$-slot. All remaining terms are lower-order due to the extra $r^{1/2}$ in the $\pi$-slot. Finally, the only terms of $\calR_{\tilde\pi}$ involving $\tilde F$ arise from $\tilde g$ when $D_t\tilde u$ is replaced, so $D'_{31}=D'_{33}=0$.

    This verifies the classes and orders in each slot and completes the proof of (a), (b) and (c).
\end{proof}

\begin{rmk}
Since $\tilde\phi$ solves the modified linearized equations \eqref{eq:abstract-modified-eqs2} and thus $\calL(\vec\varphi)\tilde \phi = \tilde F$, it is sometimes convenient to re-write \eqref{op:commutator-operator-basic} in the form
\begin{equation}\label{op:modified-commutator-operator}
    [\calL(\vec\varphi),L]\tilde\phi
    =\begin{pmatrix}
        D_{11} & D_{12} & D_{13}\\
        D_{21} & D_{22} & D_{23}\\
        D_{31} & D_{32} & D_{33}
    \end{pmatrix}\tilde \phi
    + \begin{pmatrix}
        D_{11}' & D_{12}' & D_{13}'\\
        D_{21}' & D_{22}' & D_{23}'\\
        D_{31}' & D_{32}' & D_{33}'
    \end{pmatrix}\calL(\vec\varphi)\tilde \phi.
\end{equation}
\end{rmk}

We now consider the higher-order version of the previous lemma.

\begin{lem}\label{lem:higher-order-commutator-operator-scheme}
Let $\tilde\phi=(\tilde r,\tilde u,\tilde \pi)$ be a solution to the modified linearized equations with source terms $\tilde F=(\tilde f,\tilde g,\tilde h)$. Then, for every $l\ge1$,
\begin{equation}\label{op:higher-order-commutator-operator}
    [\calL(\vec\varphi),L^l]\tilde\phi
    =D^l\tilde\phi+D^{l'}\tilde F,
    \qquad
    D^l=\big(D^l_{ij}\big),\quad D^{l'}=\big(D^{l'}_{ij}\big),
\end{equation}
where, with $s_1=r$, $s_2=u$, $s_3=\pi$,
\begin{enumerate}
    \item[(a)] for all $i,j$, $D^l_{ij}\tilde\phi_j$ is a finite sum of monomials in $\tilde\phi_j$ of class $\calM^{s_i}_{l,2l+1}$ which are lower order at level $2l$ in the slot $s_i$; every such monomial has $K\ge2l$, and every such monomial containing a time derivative has $K=2l+1$;
    \item[(b)] for all $i,j$, $D^{l'}_{ij}\tilde F_j$ is a finite sum of monomials in $\tilde F_j$ of class $\calM^{s_i}_{l,2l}$ which are lower order at level $2l$ in the slot $s_i$, contain no time derivative, and have $K=2l$.
\end{enumerate}
\end{lem}

\begin{proof}
We argue by induction on $l$. The case $l=1$ is Lemma~\ref{lem:commutator-operator-scheme-basic}. Assume that \eqref{op:higher-order-commutator-operator} holds for some $l$, and write $D,D'$ for the matrices of Lemma~\ref{lem:commutator-operator-scheme-basic}. Using the induction hypothesis, we compute
\[
    \begin{aligned}
        [\calL(\vec\varphi),L^{l+1}]\tilde \phi &= [\calL(\vec\varphi),L^l]L\tilde\phi + L^l[\calL(\vec\varphi),L]\tilde\phi\\
        &= D^lL\tilde\phi + D^{l'}\big(L\tilde F + [\calL(\vec\varphi),L]\tilde\phi\big) +L^l[\calL(\vec \varphi),L]\tilde\phi\\
        &= \big(D^lL + D^{l'}D + L^lD\big)\tilde \phi + \big(D^{l'}L + D^{l'}D' + L^lD'\big)\tilde F ,
    \end{aligned}
\]
where we used $\calL(\vec\varphi)L\tilde\phi=L\tilde F+[\calL(\vec\varphi),L]\tilde\phi$ and Lemma~\ref{lem:commutator-operator-scheme-basic}. By Remark~\ref{rmk:class-of-L}, $L$ is diagonal, of class $\calM_{1,2}$ and order $1$, with $K=2$ and $c\le2$ in every entry, and $L^l$ is of class $\calM_{l,2l}$ and order $l$, with $K=2l$ and $c\le2l$. The six products are products of matrices of operators, and we apply Lemma~\ref{lem:composition-elementary-operators} entrywise. We obtain: the entries of $D^lL$, $D^{l'}D$, $L^lD$ applied to $\tilde\phi_j$ are finite sums of monomials of class $\calM^{s_i}_{l+1,2l+3}$ which are lower order at level $2l+2$ in the slot $s_i$; the entries of $D^{l'}L$, $D^{l'}D'$, $L^lD'$ applied to $\tilde F_j$ are finite sums of monomials of class $\calM^{s_i}_{l+1,2l+2}$ which are lower order at level $2l+2$ in the slot $s_i$ and contain only spatial derivatives. Since the total number $K$ of derivatives adds under composition, every monomial of $D^lL$, $D^{l'}D$, $L^lD$ has $K\ge 2l+2$, those containing a time derivative have $K=2l+3$, and every monomial of $D^{l'}L$, $D^{l'}D'$, $L^lD'$ has $K=2l+2$. This closes the induction.
\end{proof}

We can now estimate the commutator $[\calL(\vec\varphi),L^l]\tilde\phi$ in $\calH^0$: by Lemma~\ref{lem:higher-order-commutator-operator-scheme}, each of its components is a finite sum of monomials of class $\calM^s_{l,2l+1}$ which are lower order at level $2l$ in their slot, and after the elimination of the time derivatives by Lemma~\ref{lem:time-derivative-elimination}, each of them is controlled by Lemma~\ref{lem:multilinear-tame-estimate}.

\begin{lem}\label{Lem:CommutatorEstimate}
Let $\tilde\phi := (\tilde r,\tilde u,\tilde \pi)$ be a solution to the modified linearized system with source terms $\tilde F:=(\tilde f,\tilde g,\tilde h)$. Then, for any $l>\frac{1}{2\kappa}+2$ and $t\in[0,T]$,
\begin{equation}\label{eq:commutator-estimate-1}
\begin{aligned}
\norm{[\calL(\vec\varphi),L^l]\tilde\phi}_{\calH^{0}}
&\lesssim \calC_1\Bigl(
\norm{\tilde\phi}_{\calH^{2l}}
+\norm{\tilde\phi}_{L^\infty}
\Bigl(
\norm{r}_{H^{2l+4,\,l+\frac{1}{2\kappa}+\frac{3}{2}}}
+\norm{\vec u}_{H^{2l+4,\,l+\frac{1}{2\kappa}+2}}\\
&\qquad\qquad\qquad
+\norm{\p_t\vec u}_{H^{2l+2,\,l+\frac{1}{2\kappa}+1}}
+\norm{a_0}_{H^{2l+2,\,l+\frac{1}{2\kappa}+\frac{1}{2}}}
\Bigr)\Bigr)\\
&\quad+\calC_2\Bigl(
\norm{\tilde F}_{\calH^{2l}}
+\norm{\tilde F}_{L^\infty}
\Bigl(
\norm{r}_{H^{2l+2,\,l+\frac{1}{2\kappa}+\frac{1}{2}}}
+\norm{\vec u}_{H^{2l+2,\,l+\frac{1}{2\kappa}+1}}
+\norm{a_0}_{H^{2l,\,l+\frac{1}{2\kappa}-\frac{1}{2}}}
\Bigr)\Bigr),
\end{aligned}
\end{equation}
where
\[
\begin{aligned}
\calC_1 &= \calC_1\!\left(
\norm{r}_{W^{2,\infty}},
\norm{\vec u}_{W^{1,\infty}},
\norm{\p_t\vec u}_{L^\infty},
\norm{a_0}_{W^{1,\infty}}
\right),\\
\calC_2 &= \calC_2\!\left(
\norm{r}_{W^{1,\infty}},
\norm{\vec u}_{L^\infty},
\norm{a_0}_{L^\infty}
\right).
\end{aligned}
\]
\end{lem}
\begin{proof}
Fix a slot $s\in\{r,u,\pi\}$, say $s=s_i$. The $i$-th component of $[\calL(\vec\varphi),L^l]\tilde\phi$ is $\sum_j D^l_{ij}\tilde\phi_j+\sum_jD^{l'}_{ij}\tilde F_j$, and we estimate it in $H^{0,\alpha_s}$.

By Lemma~\ref{lem:higher-order-commutator-operator-scheme}(a), each $D^l_{ij}\tilde\phi_j$ is a finite sum of monomials $DX$ in $X=\tilde\phi_j$ of class $\calM^{s}_{l,2l+1}$, lower order at level $2l$ in the slot $s$, with $K(DX)\ge 2l$ and at most one time derivative. If the time derivative is absent or falls on $\vec u$, Lemma~\ref{lem:multilinear-tame-estimate} applies directly, and \eqref{eq:multilinear-tame-1} bounds $\norm{DX}_{H^{0,\alpha_s}}$ by the $\calC_1$ term of \eqref{eq:commutator-estimate-1}. Otherwise $K(DX)=2l+1$, and by Lemma~\ref{lem:time-derivative-elimination}, $DX$ is a finite sum of monomials $D'X'$ in $X'\in\{\tilde r,\tilde u,\tilde\pi\}$ of class $\calM^{s}_{l,2l+1}$, lower order at level $2l$ in the slot $s$, with $K(D'X')\ge 2l$ and time derivative, if any, on $\vec u$, and of monomials $D''X''$ in $X''\in\{\tilde f,\tilde g,\tilde h\}$ of class $\calM^{s}_{l,2l}$, lower order at level $2l$ in the slot $s$, with $K(D''X'')=2l$ and no time derivative. Lemma~\ref{lem:multilinear-tame-estimate} applies to both: \eqref{eq:multilinear-tame-1} bounds the $D'X'$ by the $\calC_1$ term of \eqref{eq:commutator-estimate-1}, and \eqref{eq:multilinear-tame-2} bounds the $D''X''$ by the $\calC_2$ term.

By Lemma~\ref{lem:higher-order-commutator-operator-scheme}(b), each $D^{l'}_{ij}\tilde F_j$ is a finite sum of monomials in $\tilde F_j$ of class $\calM^{s}_{l,2l}$, lower order at level $2l$ in the slot $s$, with $K=2l$ and no time derivative, and \eqref{eq:multilinear-tame-2} bounds each of them by the $\calC_2$ term of \eqref{eq:commutator-estimate-1}.

Summing over the finitely many monomials, whose number depends only on $l$, and over the three slots, we obtain \eqref{eq:commutator-estimate-1}.
\end{proof}

\subsection{First-order tame estimate}

Now, we can prove the first-order tame estimate of the linearization of iteration map over the Eulerian domain.

\begin{thm}\label{thm:modified-linearized-with-xi-tame-estimate}
   Let $l_0$ be as in Notation~\ref{nota:l0}. For any $l>l_0$, the solution $(\delta r,\delta u,\delta \pi,\xi)$
    to the linearized system \eqref{eq:modified-linearized-equation-with-xi}
    satisfies:
    \begin{equation}\label{eq:tame-bound-velocity-variation}
            \norm{(\delta r,\delta \vec u,\delta \pi,\xi)}_{\calH^{2l}\times H^{2l,l+\frac{1}{2\kappa}}}
            \lesssim\calC^{2l_0+3}_E\Big(\norm{\vec g}_{L^1_tH^{2l,l+\frac{1}{2\kappa}}}+\norm{\vec g}_{L^1_tH^{2l_0,l_0+\frac{1}{2\kappa}}}
            \Big(N^{l+2}_w + \norm{\vec u}_{\tilde X^{l+4,\frac{1}{2\kappa}}_r}\Big)\Big)        
    \end{equation}
    for any $t\in[0,T]$.
    \end{thm}

\begin{proof}
Let $\delta\phi:=(\delta r,\delta u,\delta\pi)$. We can write the system \eqref{eq:modified-linearized-equation-with-xi} as
\begin{subequations}
	\begin{align}
		\calL(\vec\varphi)\delta\phi &= \tilde F\label{eq:modified-linearized-equations-with-xi-abstract1},\\
		u^\mu \p_\mu\xi^i &= u^0\p_j\left(\frac{u^i}{u^0}\right)\xi^j + H^{ij}\delta u_j\label{eq:modified-linearized-equations-with-xi-abstract2}
	\end{align}
\end{subequations}
with the source terms $\tilde F$ given by \eqref{eq:modified-linearized-equations-with-xi-source}.

Applying the weighted elliptic operator \eqref{Op:VecEllipticOperator} to \eqref{eq:modified-linearized-equations-with-xi-abstract1} and commuting with $\calL(\vec\varphi)$ gives
\begin{equation*}
	\calL(\vec\varphi)L^l\delta\phi = L^l\tilde F +[\calL(\vec\varphi),L^l]\delta\phi.
\end{equation*}

Then, by Theorem~\ref{thm:Modified-linearized-equations-wellposedness} at the basic level, we obtain
\begin{equation}\label{eq:tame-estimate-with-xi-pre1}
	\norm{L^l\delta\phi}_{\calH^0}\lesssim \calC^2_E\int_0^t\left(\norm{L^l\tilde F}_{\calH^0} + \norm{[\calL(\vec\varphi),L^l]\delta\phi}_{\calH^0}\right)ds.
\end{equation}

By Lemma~\ref{lem:boundedness-elliptic-operator} and Lemma~\ref{lem:coefficient-bounds-merged},
\begin{equation*}
	\begin{aligned}
		\norm{L^l\tilde F}_{\calH^0}&\lesssim \calC^2_E\Big(\norm{\tilde F}_{\calH^{2l}}+ \norm{\tilde F}_{L^\infty_x}\big(N_w^{l+1}+ \norm{\vec u}_{\tilde X^{l+3,\frac{1}{2\kappa}}_r}\big)\Big)\\
		&\lesssim \calC^2_E\left(\norm{\vec g - \p_j\vec\calN_{\mathrm{mo}}(\vec\varphi)\xi^j}_{H^{2l,l+\frac{1}{2\kappa}}}+\norm{\vec g-\p_j\vec\calN_{\mathrm{mo}}(\vec \varphi)\xi^j}_{L^\infty_x}\left(N^{l+1}_w+\norm{\vec u}_{\tilde X^{l+3,\frac{1}{2\kappa}}_r}\right)\right)\\
		&\lesssim\calC^3_E\left(\norm{\vec g}_{H^{2l,l+\frac{1}{2\kappa}}}+\norm{\p_j\vec\calN_{\mathrm{mo}}(\vec\varphi)\xi^j}_{H^{2l,l+\frac{1}{2\kappa}}}+\left(\norm{\vec g}_{L^\infty_x}+\norm{\xi}_{L^\infty_x}\right)\left(N^{l+1}_w+\norm{\vec u}_{\tilde X^{l+3,\frac{1}{2\kappa}}_r}\right)\right).
	\end{aligned}
\end{equation*}

For the weighted $L^2$ bound above, by Corollary~\ref{cor:weighted-L2-end-point-equivalence}, we obtain
\begin{equation*}
	\norm{\p_j\vec\calN_\mathrm{mo}(\vec\varphi)\xi^j}_{H^{2l,l+\frac1{2\kappa}}}
    \lesssim \norm{r^{l+\frac1{2\kappa}}\p_j\vec\calN_\mathrm{mo}(\vec\varphi)\xi^j}_{L^2}
    + \norm{r^{l+\frac1{2\kappa}}\p^{2l}\big(\p_j\vec\calN_\mathrm{mo}(\vec\varphi)\xi^j\big)}_{L^2}.
\end{equation*}
The first term is estimated directly with Lemma~\ref{lem:coefficient-bounds-merged}:
\begin{equation*}
	\norm{r^{l+\frac1{2\kappa}}\p_j\vec\calN_\mathrm{mo}(\vec\varphi)\xi^j}_{L^2}
    \lesssim \const{\norm{\p_t\vec u}_{L^\infty_tW^{2,\infty}_x},
    \norm{\vec u}_{L^\infty_tW^{3,\infty}_x},N_\infty^2}\norm{\xi}_{H^{2l,l+\frac1{2\kappa}}}.
\end{equation*}

For the second term, we first use \eqref{eq:r-pi-defining} and \eqref{eq:a_0-evolution} to express $\p_tr,\p_ta_0$ in terms of spatial derivatives of $(r,a_0,\vec u)$ and $\p_t\vec u$, so that $\p^{2l}\big(\p_j\vec\calN_{\mathrm{mo}}(\varphi)\xi^j\big)$ expands, via the Leibniz rule, into a finite sum of multilinear expressions of total order $2l+2$ in $\varphi$ and $\xi$, with the number of derivatives falling on $\xi$ not exceeding $2l$. Each such term is estimated by promoting a single factor to top order via Corollary~\ref{Cor:Multilinear-Interpolation}, with the resulting coefficients bounded using Lemma~\ref{lem:coefficient-bounds-merged},
\begin{equation}\label{eq:xi-coupling-highorder}
    \norm{r^{l+\frac1{2\kappa}}\p^{2l}\big(\p_j\vec\calN_\mathrm{mo}(\vec\varphi)\xi^j\big)}_{L^2}
    \lesssim \calC^3_E\Big(\norm{\xi}_{H^{2l,l+\frac1{2\kappa}}}
        + \norm{\xi}_{L^\infty_x}\big(N_w^{l+2}
        + \norm{\vec u}_{\tilde X^{l+4,\frac{1}{2\kappa}}_r}\big)\Big).
\end{equation} 

By Lemma~\ref{Lem:CommutatorEstimate} and Lemma~\ref{lem:coefficient-bounds-merged},
\begin{equation}\label{eq:commutator-bound}
\begin{aligned}
    \norm{[\calL(\vec\varphi),L^l]\delta\phi}_{\calH^0}
    &\lesssim \calC^3_E\Big(\norm{\delta\phi}_{\calH^{2l}}
        + \norm{\delta\phi}_{L^\infty_x}\big(N_w^{l+2}
        +\norm{\vec u}_{\tilde X^{l+4,\frac{1}{2\kappa}}_r}\big)\Big)\\
    &\quad+\calC^2_E
    \Big(\norm{\tilde F}_{\calH^{2l}}
        + \norm{\tilde F}_{L^\infty_x}\big(N_w^{l+1}
        +\norm{\vec u}_{\tilde X^{l+3,\frac{1}{2\kappa}}_r}\big)\Big).
\end{aligned}
\end{equation}

By Lemma~\ref{lem:tame-elliptic-estimate} and Lemma~\ref{lem:coefficient-bounds-merged},
\begin{equation}\label{eq:elliptic-bound}
\begin{aligned}
    \norm{\delta\phi}_{\calH^{2l}}\lesssim \calC^3_E\Big(\norm{L^l\delta\phi}_{\calH^0}
    	+ \norm{\delta\phi}_{\calH^0}
        + \norm{\delta\phi}_{L^\infty_x}\big(N_w^{l+1}
        + \norm{\vec u}_{\tilde X^{l+3,\frac{1}{2\kappa}}_r}\big)\Big).
\end{aligned}
\end{equation}
Combining all the bounds above,
\begin{equation}\label{eq:combined-1}
\begin{aligned}
    \norm{\delta\phi}_{\calH^{2l}}
    &\lesssim \calC^3_E\int_0^t \norm{\delta\phi}_{\calH^{2l}}
        + \norm{\xi}_{H^{2l,l+\frac1{2\kappa}}}
        + \big(\norm{\delta\phi}_{L^\infty_x}+\norm{\xi}_{L^\infty_x}\big)\big(N_w^{l+2}
        + \norm{\vec u}_{\tilde X^{l+4,\frac{1}{2\kappa}}_r}\big)\,ds\\
    &\quad+ \calC^3_E\int_0^t \norm{\vec g}_{H^{2l,l+\frac1{2\kappa}}}
        + \norm{\vec g}_{L^\infty_x}\big(N_w^{l+1}
        +\norm{\vec u}_{\tilde X^{l+3,\frac{1}{2\kappa}}_r}\big)\,ds\\
    &\quad+ \calC^3_E\norm{\delta\phi}_{\calH^0}\\
    &\quad+ \calC^3_E\norm{\delta\phi}_{L^\infty_x}\big(N_w^{l+1}
        +\norm{\vec u}_{\tilde X^{l+3,\frac{1}{2\kappa}}_r}\big).
\end{aligned}
\end{equation}

By Theorem~\ref{thm:modified-linearized-with-xi-local-wellposedness} and Lemma~\ref{lem:WeightedEmbedding},
\[
    \norm{\delta\phi}_{\calH^0}
    \lesssim \calC^3_E\int_0^t \norm{\vec g}_{H^{0,\frac1{2\kappa}}}\,ds,
\]
and 
\[
    \norm{\delta\phi}_{L^\infty_x} + \norm{\xi}_{L^\infty_x}
    \lesssim \norm{\delta\phi}_{\calH^{2l_0}} + \norm{\xi}_{H^{2l_0,l_0+\frac{1}{2\kappa}}}
    \lesssim \calC^{2l_0+3}_E
    \int_0^t \norm{\vec g}_{H^{2l_0,l_0+\frac1{2\kappa}}}\,ds.
\]
Substituting these bounds back,
\begin{equation}\label{eq:combined-2}
\begin{aligned}
    \norm{\delta\phi}_{\calH^{2l}}
    &\lesssim \calC^{2l_0+3}_E\int_0^t \norm{\delta\phi}_{\calH^{2l}}
        + \norm{\xi}_{H^{2l,l+\frac1{2\kappa}}}ds\\
    &\quad+ \calC^{2l_0+3}_E\Big(\norm{\vec g}_{L^1_tH^{2l,l+\frac1{2\kappa}}}
        +\norm{\vec g}_{L^1_tH^{2l_0,l_0+\frac1{2\kappa}}}\big(N_w^{l+2}
        +\norm{\vec u}_{\tilde X^{l+4,\frac{1}{2\kappa}}_r}\big)\Big),
\end{aligned}
\end{equation}

By Lemma~\ref{lem:generic-linear-transport-estimate}, we have
\begin{equation}\label{eq:xi-final}
	\norm{\xi}_{H^{2l,l+\frac1{2\kappa}}}\lesssim \calC^{2l_0+3}_E\left(\int_0^t \norm{\delta\vec u}_{H^{2l,l+\frac1{2\kappa}}}ds
        + \norm{\vec g}_{L^1_tH^{2l_0,l_0+\frac{1}{2\kappa}}}\norm{\vec u}_{L^\infty_tH^{2l+2,l+\frac{1}{2\kappa}+1}}\right).
\end{equation}
Summing \eqref{eq:combined-2} and \eqref{eq:xi-final}, 
\begin{equation}\label{eq:sum}
\begin{aligned}
    \norm{\delta\phi}_{\calH^{2l}}+\norm{\xi}_{H^{2l,l+\frac1{2\kappa}}}
    &\lesssim \calC^{2l_0+3}_E\int_0^t \norm{\delta\phi}_{\calH^{2l}}
        + \norm{\xi}_{H^{2l,l+\frac1{2\kappa}}}ds\\
    &\quad+ \calC^{2l_0+3}_E\Big(\norm{\vec g}_{L^1_tH^{2l,l+\frac1{2\kappa}}}
        +\norm{\vec g}_{L^1_tH^{2l_0,l_0+\frac1{2\kappa}}}\big(N_w^{l+2}
        +\norm{\vec u}_{\tilde X^{l+4,\frac{1}{2\kappa}}_r}\big)\Big).
\end{aligned}
\end{equation}

By Gr\"onwall's inequality,
\[
    \norm{\delta\phi}_{\calH^{2l}}+\norm{\xi}_{H^{2l,l+\frac1{2\kappa}}}
    \lesssim\calC^{2l_0+3}_E\Big(\norm{\vec g}_{L^1_tH^{2l,l+\frac1{2\kappa}}}
        +\norm{\vec g}_{L^1_tH^{2l_0,l_0+\frac1{2\kappa}}}\big(N_w^{l+2}
        +\norm{\vec u}_{\tilde X^{l+4,\frac{1}{2\kappa}}_r}\big)\Big),
\]
which finishes the proof.
  
\end{proof}

\section{Tame estimates in Lagrangian coordinates}\label{sec:Lagrangian-tame}
In this section, we transfer the energy estimate and the first-order tame estimate for the linearized iteration map, established over the Eulerian domain in Section~\ref{sec:linearized-eqs} and Section~\ref{sec:tame-estimate-linearization-Eulerian}, to the Lagrangian domain, where the Nash--Moser iteration is defined. As before, we fix $\vec U$ defined on $[0,T]\times\Omega_0$ with $T\in(0,1]$, and define the variables $(R,U,P,\eta)$ and $(r,u,\pi,\Omega_t)$ in terms of $\vec U$ as in Section~\ref{sec:Nash-Moser-iteration-scheme}, via \eqref{eq:U^0-defining-eq}, \eqref{eq:eta-defining-eq}, \eqref{eq:R-P-defining}, and Remark~\ref{rmk:r-u-pi-defining}. We also employ the notation of \eqref{eq:variables-abbreviation} for convenience.

\subsection{Change-of-variables estimates between the Eulerian and the Lagrangian domain}\label{subsec:change-of-coordinate}
In this subsection, we establish the estimates needed for the change of coordinates between the Eulerian and Lagrangian domains. Specifically, we need to show how to bound the weighted Sobolev norm of $\mfG:=\mfg\circ\eta$ in Lagrangian coordinates by the norm $\norm{\cdot}_{H^{2l,l+\frac{1}{2\kappa}}_r}$ over the Eulerian domain in a tame way, and vice versa. By the chain rule, such bounds reduce to controlling the deformation matrix $F\circ\eta^{-1}$ pulled back to the Eulerian domain and its inverse $F^{-1}$ in the appropriate norms. We will estabish the $L^\infty$ and weighted $L^2$ control of $F\circ\eta^{-1}$ and $F^{-1}$ and with them eventually establish the tame bounds for the change of coordinates between the Eulerian and Lagrangian domain.

First, we show the $L^\infty$ bound of $F^{-1}$:

\begin{lem}\label{lem:deformation-inverse-Linfty-estimate}
	For any $M>0$, there exists $T=T(M)\in(0,1]$ such that if
	\begin{equation*}
		\norm{\vec U}_{L^\infty_tW_y^{1,\infty}}<M,
	\end{equation*}
	then for any $t\in[0,T]$,
	\begin{equation}\label{eq:deformation-det-estimate}
		\frac12\leq\det\nabla_y\eta\leq\frac32,
	\end{equation}
	and
	\begin{equation}\label{eq:deformation-inverse-Linfty-estimate}
		\norm{F^{-1}}_{L^\infty_y}\leq \const{\norm{\vec U}_{L^\infty_tW^{1,\infty}_y}}.
	\end{equation}
\end{lem}

\begin{proof}
	Differentiating \eqref{eq:eta-defining-eq} gives
	\begin{equation*}
		F=\nabla_y\eta = I + \int_0^t\left(\nabla_y\left(\frac{\vec U}{U^0}\right)\right)(s,y)\,ds,
	\end{equation*}
	which implies
	\begin{equation*}
		\norm{F}_{L^\infty_y}\leq 1 + T\norm{\nabla_y\left(\frac{\vec U}{U^0}\right)}_{L^\infty_tL^\infty_y}\leq\const{\norm{\vec U}_{L^\infty_tW^{1,\infty}_y}}.
	\end{equation*}

	Now we show that $T=T(M)$ can be chosen so that $\det F$ stays bounded
	away from $0$. The determinant $\det F$ satisfies
	\begin{equation}\label{eq:transport-det-deformation-matrix}
		\p_t(\det F) = \tr\big(\adj F\,\p_tF\big)=\tr\left(\adj F\,\nabla_y\left(\frac{\vec U}{U^0}\right)\right).
	\end{equation}
	Each entry of the cofactor matrix $\adj F$ is a polynomial in the
	entries of $F$, and thus
	\begin{equation*}
		\norm{\adj F}_{L^\infty_y}\leq\const{\norm{F}_{L^\infty_y}}\leq\const{\norm{\vec U}_{L^\infty_tW^{1,\infty}_y}}.
	\end{equation*}
	Integrating \eqref{eq:transport-det-deformation-matrix} and estimating
	directly, we obtain
	\begin{equation*}
		\begin{aligned}
			\abs{\det F-1}&\leq \abs{\int_0^t\tr\left(\adj F\,\nabla_y\left(\frac{\vec U}{U^0}\right)\right)ds}\\
			&\leq \int_0^t \norm{\adj F}_{L^\infty_y}\norm{\nabla_y\left(\frac{\vec U}{U^0}\right)}_{L^\infty_y}\,ds\leq T\const{\norm{\vec U}_{L^\infty_tW^{1,\infty}_y}}.
		\end{aligned}
	\end{equation*}
	Thus,
	\begin{equation*}
		\det F\geq 1-T\const{\norm{\vec U}_{L^\infty_tW^{1,\infty}_y}} \geq 1-T\const{M}\geq \frac12
	\end{equation*}
	if $T$ is chosen sufficiently small depending on $M$, which gives
	\eqref{eq:deformation-det-estimate}.

	By Cramer's rule,
	\begin{equation*}
		\norm{F^{-1}}_{L^\infty_y} = \norm{\frac{\adj F}{\det F}}_{L^\infty_y}\leq\frac{\const{\norm{\vec U}_{L^\infty_tW^{1,\infty}_y}}}{1/2}=\const{\norm{\vec U}_{L^\infty_tW^{1,\infty}_y}}
	\end{equation*}
	on $[0,T]$, which gives \eqref{eq:deformation-inverse-Linfty-estimate}.
\end{proof}

\begin{rmk}\label{rmk:eta-diffeo}
By Lemma~\ref{lem:deformation-inverse-Linfty-estimate}, we can show that, for $T=T(M)$ small, $\eta(t,\cdot)$ is a bi-Lipschitz map of $\Omega_0$ onto $\Omega_t$ with $\eta(t,\Gamma_0)=\Gamma_t$. In
particular $\eta^{-1}$ is well defined.
\end{rmk}

The following lemma is a Lagrangian version of Lemma~\ref{lem:generic-Linfty-transport-estimate}.

\begin{lem}\label{lem:Lagrangian-generic-transport-estimate}
	Let $Q,\Psi$ be functions such that
	\begin{equation*}
		\norm{Q}_{L^\infty_tL^\infty_y} + \norm{\Psi}_{L^\infty_tL^\infty_y}<\infty
	\end{equation*}
	and $Q$ solves
	\begin{equation}\label{eq:Lagrangian-generic-transport}
    \left\{
    \begin{aligned}
        \p_t Q &= f(Q,\Psi),\\
        Q(0,\cdot)&=Q_0,
    \end{aligned}
    \right.
\end{equation}
for some function $f=f(Q,\Psi)$ smooth in a neighborhood of the range of $(Q,\Psi)$. Then, for every integer $k\geq0$
\begin{equation}
	\norm{Q}_{W^{k,\infty}_y}\lesssim\const{\norm{\Psi}_{L^\infty_tW^{k,\infty}_y},\norm{Q}_{L^\infty_tL^\infty_y},\norm{Q_0}_{W^{k,\infty}}}
\end{equation}
for any $t\in[0,T]$.
\end{lem}

\begin{proof}
    The proof is the same as that of Lemma~\ref{lem:generic-Linfty-transport-estimate}, and simpler: since $\p_t$ and $\p_y^k$ commute, there are no commutator terms. We apply $\p_y^k$ to $\p_tQ=f(Q,\Psi)$, expand $\p_y^k\big(f(Q,\Psi)\big)$ by the Leibniz rule as for the term $\p^k(f(q,\psi))$ in that proof, integrate in $t$ and apply Gr\"onwall's inequality. We omit the details.
\end{proof}

\begin{lem}\label{lem:deformation-inverse-Wkinfty-estimate}
For any $M>0$, there exists $T=T(M)\in(0,1]$ such that if $\norm{\vec U}_{L^\infty_tW_y^{1,\infty}}<M$, then for any integer $k\geq 0$ and $t\in[0,T]$,
	\begin{equation}\label{eq:deformation-inverse-estimate}
		\norm{F^{-1}}_{W^{k,\infty}_y}\leq\const{\norm{\vec U}_{L^\infty_tW^{k+1,\infty}_y}}.
	\end{equation}
\end{lem}
\begin{proof}
	Fix $M>0$, and pick $T=T(M)$ as in Lemma~\ref{lem:deformation-inverse-Linfty-estimate}. Notice that $F^{-1}$ solves
	\begin{equation}\label{eq:transport-deformation-inverse}
		\p_t(F^{-1}) = -F^{-1}\nabla_y\left(\frac{\vec U}{U^0}\right) F^{-1}.
	\end{equation}
	Applying Lemma~\ref{lem:Lagrangian-generic-transport-estimate} with $Q:=F^{-1}$, $\Psi:=\nabla_y\left(\frac{\vec U}{U^0}\right)$, and $Q_0=F^{-1}(0,\cdot)=I$ (since $\eta(0,\cdot)=\Id_{\Omega_0}$) yields
\eqref{eq:deformation-inverse-estimate}.
\end{proof}

\begin{thm}\label{thm:Eulerian-Lagrangian-weighted-Sobolev-comparison}
For any $M>0$, there exists $T=T(M)\in(0,1]$ such that if
$\norm{\vec U}_{L^\infty_tW^{1,\infty}_y}<M$, then for any $t\in[0,T]$ and
integer $l\geq0$,
\begin{equation}\label{eq:Eulerian-Lagrangian-weighted-Sobolev-comparison}
	\norm{\mfg}_{H^{2l,l+\frac{1}{2\kappa}}_r}\lesssim\const{\norm{\vec U}_{L^\infty_tW^{2l,\infty}_y}}\norm{\mfG}_{H^{2l,l+\frac{1}{2\kappa}}_R}.
\end{equation}
\end{thm}
\begin{proof}
By Corollary~\ref{cor:weighted-L2-end-point-equivalence},
\begin{equation*}
	\norm{\mfg}_{H^{2l,l+\frac{1}{2\kappa}}_r}\lesssim \norm{r^{l+\frac{1}{2\kappa}}\mfg}_{L^2_x} + \norm{r^{l+\frac{1}{2\kappa}}\p^{2l}\mfg}_{L^2_x}.
\end{equation*}
For the first term, changing variables via $\eta$ yields
\begin{equation*}
	\begin{aligned}
		\norm{r^{l+\frac{1}{2\kappa}} \mfg}_{L^2_x}^2 &= \int_{\Omega_t}r^{2l+\frac{1}{\kappa}}\abs{\mfG\circ\eta^{-1}}^2\,dx\\
		&=\int_{\Omega_0}R^{2l+\frac{1}{\kappa}}\abs{\mfG}^2\det F\,dy\\
		&\lesssim\norm{R^{l+\frac{1}{2\kappa}}\mfG}_{L^2_y}^2,
	\end{aligned}
\end{equation*}
where the last inequality uses the bound $\det F\leq\frac32$ from
\eqref{eq:deformation-det-estimate}.
For the second term, similarly,
\begin{equation*}
	\begin{aligned}
		\norm{r^{l+\frac{1}{2\kappa}}\p^{2l}\mfg}_{L^2_x}^2 &= \int_{\Omega_t}r^{2l+\frac{1}{\kappa}}\abs{\p^{2l}\mfg}^2\,dx\\
		&=\int_{\Omega_0}R^{2l+\frac{1}{\kappa}}\abs{(\p^{2l}\mfg)\circ\eta}^2\det F\,dy\\
		&\lesssim\norm{R^{l+\frac{1}{2\kappa}}(\p^{2l}\mfg)\circ\eta}_{L^2_y}^2,
	\end{aligned}
\end{equation*}
again using $\det F\leq\frac32$ from \eqref{eq:deformation-det-estimate}.
Over the Lagrangian domain, $(\p^{2l}\mfg)\circ\eta=(F^{-1}\nabla_y)^{2l}\mfG$,
which expands, by the Leibniz rule, into a finite sum of multilinear
expressions
\[
	\Big(\prod_{i_1=1}^{I_1}\p_y^{a_{i_1}}F^{-1}\Big)\p_y^c\mfG,\qquad \sum_{i_1} a_{i_1}+c=2l,\quad a_{i_1},c\ge1,
\]
in derivatives of $F^{-1}$ and $\mfG$ with coefficients smoothly depending
on $F^{-1}$. Note that the innermost application of $\nabla_y$ in
$(F^{-1}\nabla_y)^{2l}$ acts directly on $\mfG$, and thus every such term
has $c\ge1$. In particular, no single factor of $F^{-1}$ carries more than
$2l-1$ derivatives. Placing the coefficient smoothly depending on $F^{-1}$
in $L^\infty$, bounded via Lemma~\ref{lem:deformation-inverse-Wkinfty-estimate},
and leaving the single factor $\p_y^c\mfG$ ($c\le2l$), multiplied by
$R^{l+\frac{1}{2\kappa}}$, as one of the terms summed in
$\norm{\mfG}_{H^{2l,l+\frac{1}{2\kappa}}_R}$, we obtain
\begin{equation*}
	\norm{R^{l+\frac{1}{2\kappa}}(\p^{2l}\mfg)\circ\eta}_{L^2_y}\lesssim\const{\norm{\vec U}_{L^\infty_tW^{2l,\infty}_y}}\norm{\mfG}_{H^{2l,l+\frac{1}{2\kappa}}_R}.
\end{equation*}
Combining the estimates for the first and second terms yields
\eqref{eq:Eulerian-Lagrangian-weighted-Sobolev-comparison}.
\end{proof}

The following lemma controls the $L^\infty$ bounds of velocity under the change of coordinates between Eulerian and Lagrangian domains.

\begin{lem}\label{lem:velocity-Linfty-from-Eulerian-to-Lagrangian}
For any $M>0$, there exists $T=T(M)\in(0,1]$ such that, if
$\norm{\vec U}_{L^\infty_tW^{1,\infty}_y}<M$, then for any integer $k\ge0$
and $t\in[0,T]$,
\begin{equation}\label{eq:velocity-Linfty-from-Eulerian-to-Lagrangian}
    \begin{aligned}
        \norm{\vec u}_{W^{k,\infty}(\Omega_t)} &\lesssim \const{\norm{\vec U}_{L^\infty_tW^{k,\infty}_y}},\\
        \norm{\p_t\vec u}_{W^{k,\infty}(\Omega_t)} &\lesssim \const{\norm{\p_t\vec U}_{L^\infty_tW^{k,\infty}_y},
        \norm{\vec U}_{L^\infty_tW^{k+1,\infty}_y}}.
    \end{aligned}
\end{equation}
\end{lem}

\begin{proof}
Throughout, fix $\beta$ with $\abs\beta\le k$, and recall from the chain
rule that for any Eulerian quantity $\mfg$,
\begin{equation}\label{eq:pullback-operator-recall}
    (\p_x^\beta\mfg)\circ\eta = (F^{-1}\nabla_y)^{\abs\beta}(\mfg\circ\eta),
\end{equation}
exactly as in the proof of
Theorem~\ref{thm:Eulerian-Lagrangian-weighted-Sobolev-comparison}.

Since $\vec u\circ\eta=\vec U$,
applying \eqref{eq:pullback-operator-recall} with $\mfg=\vec u$ gives
$(\p_x^\beta\vec u)\circ\eta=(F^{-1}\nabla_y)^{\abs\beta}\vec U$, which
expands, by the Leibniz rule, into a finite sum of terms
\[
    \Big(\prod_{i_1=1}^{I_1}\p_y^{a_{i_1}}F^{-1}\Big)\p_y^c\vec U,\qquad \sum_{i_1} a_{i_1}+c=\abs\beta,\quad a_{i_1}\ge1,
\]
with coefficients smoothly depending on $F^{-1}$. As in the proof of Theorem~\ref{thm:Eulerian-Lagrangian-weighted-Sobolev-comparison}, the innermost application of $\nabla_y$ acts directly on $\vec U$, so $c\ge1$ and no single factor of $F^{-1}$ carries more than $\abs\beta-1\le k-1$ derivatives. By Lemma~\ref{lem:deformation-inverse-Wkinfty-estimate} (at order $\leq k-1$),
\[
    \norm{(\p_x^\beta\vec u)\circ\eta}_{L^\infty_y}
    \lesssim\const{\norm{\vec U}_{L^\infty_tW^{k,\infty}_y}}\norm{\vec U}_{W^{\abs\beta,\infty}_y}
    \lesssim\const{\norm{\vec U}_{L^\infty_tW^{k,\infty}_y}}.
\]
Since $\norm{\p_x^\beta\vec u}_{L^\infty_x}=\norm{(\p_x^\beta\vec u)\circ\eta}_{L^\infty_y}$ and $\beta$ was arbitrary with $\abs\beta\le k$, this proves the first estimate in \eqref{eq:velocity-Linfty-from-Eulerian-to-Lagrangian}.

For the bound on $\p_t\vec u$, we differentiate the identity
$\vec u\circ\eta=\vec U$ in $t$, at fixed $y$. By the chain rule,
\[
    \p_t\vec U = \p_t(\vec u\circ\eta)=(\p_t\vec u)\circ\eta+\big((\nabla_x\vec u)\circ\eta\big)\,\p_t\eta.
\]
Substituting $\p_t\eta=(\vec U/U^0)$ and $(\nabla_x\vec u)\circ\eta=(\nabla_y\vec U)F^{-1}$, we
obtain
\begin{equation}\label{eq:dt-velocity-identity}
    (\p_t\vec u)\circ\eta = \p_t\vec U - (\nabla_y\vec U)F^{-1}\frac{\vec U}{U^0}.
\end{equation}

By \eqref{eq:pullback-operator-recall} with $\mfg=\p_t\vec u$, together with \eqref{eq:dt-velocity-identity},
\[
    (\p_x^\beta\p_t\vec u)\circ\eta = (F^{-1}\nabla_y)^{\abs\beta}\p_t\vec U
    \ -\ (F^{-1}\nabla_y)^{\abs\beta}\Big[(\nabla_y\vec U)F^{-1}\frac{\vec U}{U^0}\Big].
\]
Taking the $L^\infty$ norm yields
\[
    \norm{\p_x^\beta\p_t\vec u}_{L^\infty_x}=\norm{(\p_x^\beta\p_t\vec u)\circ\eta}_{L^\infty_y}
    \lesssim\const{\norm{\vec U}_{L^\infty_tW^{k,\infty}_y}}\norm{\p_t\vec U}_{W^{k,\infty}_y}
    +\const{\norm{\vec U}_{L^\infty_tW^{k+1,\infty}_y}},
\]
which is controlled by
$\const{\norm{\p_t\vec U}_{L^\infty_tW^{k,\infty}_y},\norm{\vec U}_{L^\infty_tW^{k+1,\infty}_y}}$,
proving the second estimate in
\eqref{eq:velocity-Linfty-from-Eulerian-to-Lagrangian}.
\end{proof}

\begin{conv}[Choice of $M$ and $T$]\label{conv:choice-of-T}
Fix $M>0$. Let $T_1=T_1(M)\in(0,1]$ be the time of Lemma~\ref{lem:deformation-inverse-Linfty-estimate};
since Lemmas~\ref{lem:deformation-inverse-Wkinfty-estimate} and \ref{lem:velocity-Linfty-from-Eulerian-to-Lagrangian} rely only on that lemma, they hold on $[0,T_1]$ as
well, and in particular $\norm{\vec u}_{\tilde Z^1_x}<\Lambda:=\calC(M)$ whenever $\norm{\vec U}_{\tilde Z^1_y}<M$.
By Lemma~\ref{lem:basic-Linfty-coefficient-estimate} there is $T_2=T_2(\Lambda,N^0_\infty)\in(0,1]$ such that
\begin{equation}\label{eq:basic-Linfty-R-A0-bound}
    \norm{(R,A_0)}_{L^\infty_y}=\norm{(r,a_0)}_{L^\infty_x}\le2N^0_\infty\qquad\text{on }[0,T_2].
\end{equation}
We fix $T\in(0,T_*]$, $T_*:=\min\{T_1,T_2\}=T_*(M,N^0_\infty)\le1$. All results below hold on $[0,T]$ for any such $T$, so $T$ may be decreased whenever needed.
\end{conv}

For the rest of the paper we assume $\norm{\vec U}_{\tilde Z^1_y}<M$, with $M$ and $T$ as in
Convention~\ref{conv:choice-of-T}.

The following lemma is a Lagrangian version of
Lemma~\ref{lem:generic-weighted-transport-estimate}.

\begin{lem}\label{lem:Lagrangian-generic-weighted-transport-estimate}
In the setting of Lemma~\ref{lem:Lagrangian-generic-transport-estimate}, for any integer $k\ge0$ and $\sigma>-\frac12$ with $k-\sigma>\frac32$,
and any $t\in[0,T]$,
\begin{equation}
    \norm{R^\sigma\p^kQ}_{L^2_y}
    \lesssim \calC\Big(\norm{r_0^\sigma\p^kQ_0}_{L^2_y}
    +\int_0^t\norm{R^\sigma\p^k\Psi(s)}_{L^2_y}\,ds\Big),
\end{equation}
where $\calC=\const{\norm{\vec U}_{L^\infty_tW^{1,\infty}_y},
\norm{\p_t\vec U}_{L^\infty_tL^\infty_y},N^0_\infty,
\norm{Q}_{L^\infty_tL^\infty_y},\norm{\Psi}_{L^\infty_tL^\infty_y}}$.
\end{lem}

\begin{proof}
Taking time derivative of the norm,
\[
    \p_t\norm{R^\sigma\p^kQ}_{L^2_y}^2
    = \p_t\int_{\Omega_0}R^{2\sigma}\abs{\p^kQ}^2\,dy
    = \int_{\Omega_0}2\sigma R^{2\sigma-1}(\p_tR)\abs{\p^kQ}^2\,dy
    + \int_{\Omega_0}2R^{2\sigma}\p^kQ\,\p_t\p^kQ\,dy.
\]

To bound the first term, we use the first equation of
\eqref{eq:R-P-defining} to trade $\p_tR$ for one more power of $R$, which
introduces coefficients in $R,A_0,\p_t\vec U,\p\vec U$, and $F^{-1}$. By a direct estimate with \eqref{eq:basic-Linfty-R-A0-bound}, the
first term is controlled by
\[
    \int_{\Omega_0}2\sigma R^{2\sigma-1}(\p_tR)\abs{\p^kQ}^2\,dy
    \lesssim \const{\norm{\vec U}_{L^\infty_tW^{1,\infty}_y},
    \norm{\p_t\vec U}_{L^\infty_tL^\infty_y},N^0_\infty}
    \norm{R^\sigma\p^kQ}_{L^2_y}^2.
\]

To bound the second term, by Cauchy--Schwarz, it suffices to bound $\norm{R^\sigma\p_t\p^kQ}_{L^2_y}$. Applying $\p^k$ to \eqref{eq:Lagrangian-generic-transport} gives directly $\p_t\p^kQ=\p^k\big(f(Q,\Psi)\big)$. Expanding the right-hand side by the Leibniz rule, exactly as for the term $\p^k(f(q,\psi))$ in the proof of Lemma~\ref{lem:generic-weighted-transport-estimate}, and applying
Corollary~\ref{Cor:Multilinear-Interpolation} to bound each resulting
term as in that proof, we obtain
\[
    \norm{R^\sigma\p_t\p^kQ}_{L^2_y}
    \lesssim \const{\norm{Q}_{L^\infty_tL^\infty_y},\norm{\Psi}_{L^\infty_tL^\infty_y}}
    \Big(\norm{R^\sigma\p^kQ}_{L^2_y}
    +\norm{R^\sigma\p^k\Psi}_{L^2_y}\Big).
\]

Putting the two bounds together,
\[
    \p_t\norm{R^\sigma\p^kQ}_{L^2_y}
    \lesssim \calC\Big(\norm{R^\sigma\p^kQ}_{L^2_y}
    +\norm{R^\sigma\p^k\Psi}_{L^2_y}\Big),
\]
where $\calC=\const{\norm{\vec U}_{L^\infty_tW^{1,\infty}_y},
\norm{\p_t\vec U}_{L^\infty_tL^\infty_y},N^0_\infty,
\norm{Q}_{L^\infty_tL^\infty_y},\norm{\Psi}_{L^\infty_tL^\infty_y}}$.
By Gr\"onwall's inequality,
\[
    \norm{R^\sigma\p^kQ(t)}_{L^2_y}
    \lesssim \calC\Big(\norm{r_0^\sigma\p^kQ_0}_{L^2_y}
    +\int_0^t\norm{R^\sigma\p^k\Psi(s)}_{L^2_y}\,ds\Big),
\]
which is the claimed estimate.
\end{proof}

A slight variation of
Lemma~\ref{lem:Lagrangian-generic-weighted-transport-estimate}, which is
the Lagrangian version of
Lemma~\ref{lem:generic-linear-transport-estimate}, is recorded as
follows.

\begin{lem}\label{lem:Lagrangian-generic-linear-transport-estimate}
	Let $Q,\tilde Q,\tilde\Psi$ be
	functions such that $\norm{Q}_{L^\infty_tL^\infty_y}<\infty$ and
	$\tilde Q,\tilde\Psi$ satisfy
	\begin{equation}\label{eq:Lagrangian-generic-linear-transport}
    \left\{
    \begin{aligned}
        \p_t\tilde Q &= f(Q)\,\tilde Q + \tilde\Psi,\\
        \tilde Q(0,\cdot) &= \tilde Q_0,
    \end{aligned}
    \right.
\end{equation}
	for some $f=f(Q)$ smooth in a neighborhood of the range of its
	argument. Then,
	\begin{enumerate}
	\item[(a)] for any integer $k\ge0$ and $\sigma>-\frac12$,
	\[
	    \norm{\tilde Q}_{H^{k,\sigma}_R}
	    \lesssim \const{\norm{\vec U}_{L^\infty_tW^{1,\infty}_y},
	    \norm{\p_t\vec U}_{L^\infty_tL^\infty_y},N^0_\infty,
	    \norm{Q}_{L^\infty_tW^{k,\infty}_y}}
	    \Big(\norm{\tilde Q_0}_{H^{k,\sigma}_{r_0}}
	    +\int_0^t\norm{\tilde\Psi}_{H^{k,\sigma}_R}\,ds\Big);
	\]
	\item[(b)] if $\tilde Q_0=0$, then we have an $L^\infty$ estimate:
	\begin{equation}\label{eq:tildeQ-Linfty-bound}
    	\norm{\tilde Q}_{L^\infty_y}\lesssim\const{\norm{Q}_{L^\infty_tL^\infty_y}}\int_0^t\norm{\tilde\Psi(s)}_{L^\infty_y}\,ds;
	\end{equation}
	
	\item[(c)] if $k-\sigma>\frac32$ and $\tilde Q_0=0$,
	\[
	    \norm{R^\sigma\p^k\tilde Q}_{L^2_y}
	    \lesssim \calC\Big(\int_0^t\int_0^s\norm{\tilde\Psi(\tau)}_{L^\infty_y}\,d\tau\,
	    \norm{R^\sigma\p^kQ(s)}_{L^2_y}\,ds
	    +\int_0^t\norm{R^\sigma\p^k\tilde\Psi(s)}_{L^2_y}\,ds\Big),
	\]
	where $\calC=\const{\norm{\vec U}_{L^\infty_tW^{1,\infty}_y},
	\norm{\p_t\vec U}_{L^\infty_tL^\infty_y},N^0_\infty,
	\norm{Q}_{L^\infty_tL^\infty_y}}$.
	\end{enumerate}
\end{lem}

\begin{proof}
Applying $\p_y^k$ to
\eqref{eq:Lagrangian-generic-linear-transport} gives
\begin{equation}\label{eq:k-derivative-Lagrangian-generic-transport}
	\p_t\p^k\tilde Q=\p^k\big(f(Q)\tilde Q\big)+\p^k\tilde\Psi.
\end{equation}
Expanding $\p^k(f(Q)\tilde Q)$ by the Leibniz rule and bounding every
factor of $Q$ in $L^\infty$, we obtain
\begin{equation}\label{eq:dt-dk-tildeQ-bound}
    \norm{R^\sigma\p_t\p^k\tilde Q}_{L^2_y}
    \lesssim\const{\norm{Q}_{L^\infty_tW^{k,\infty}_y}}\norm{\tilde Q}_{H^{k,\sigma}_R}
    +\norm{R^\sigma\p^k\tilde\Psi}_{L^2_y}.
\end{equation}
Differentiating the norm in time,
\[
    \p_t\norm{R^\sigma\p^k\tilde Q}_{L^2_y}^2
    =\int_{\Omega_0}2\sigma R^{2\sigma-1}(\p_tR)\abs{\p^k\tilde Q}^2\,dy
    +\int_{\Omega_0}2R^{2\sigma}\p^k\tilde Q\cdot\p_t\p^k\tilde Q\,dy.
\]

The first term is bounded by
\[
\const{\norm{\vec U}_{L^\infty_tW^{1,\infty}_y},\norm{\p_t\vec U}_{L^\infty_tL^\infty_y},N^0_\infty}\norm{R^\sigma\p^k\tilde Q}_{L^2_y}^2,
\]
using the first equation of
\eqref{eq:R-P-defining} to trade $\p_tR$ for one more power of $R$, which
introduces coefficients in $R,A_0,\p_t\vec U,\p\vec U$, and $F^{-1}$.

The second, by Cauchy--Schwarz and \eqref{eq:dt-dk-tildeQ-bound}, is
bounded by
$2\norm{R^\sigma\p^k\tilde Q}_{L^2_y}\norm{R^\sigma\p_t\p^k\tilde Q}_{L^2_y}$.
Combining,
\[
    \p_t\norm{R^\sigma\p^k\tilde Q}_{L^2_y}
    \lesssim\const{\norm{\vec U}_{L^\infty_tW^{1,\infty}_y},\norm{\p_t\vec U}_{L^\infty_tL^\infty_y},N^0_\infty,\norm{Q}_{L^\infty_tW^{k,\infty}_y}}\norm{\tilde Q}_{H^{k,\sigma}_R}
    +\norm{R^\sigma\p^k\tilde\Psi}_{L^2_y}.
\]
The identical argument applied to \eqref{eq:Lagrangian-generic-linear-transport}
without differentiating (i.e.\ $k=0$) gives the same bound for
$\p_t\norm{R^\sigma\tilde Q}_{L^2_y}$. Integrating both in time and
applying Corollary~\ref{cor:weighted-L2-end-point-equivalence} and
Gr\"onwall's inequality, we obtain part (a).

For part (b) and (c), assume $k-\sigma>\frac32$ and $\tilde Q_0=0$. We first derive an $L^\infty$ bound for $\tilde Q$. Integrating
\eqref{eq:Lagrangian-generic-linear-transport} directly in $t$,
\[
    \tilde Q(t,y)=\underbrace{\tilde Q_0(y)}_{=0}+\int_0^t\big(f(Q)\tilde Q+\tilde\Psi\big)(s,y)\,ds,
\]
so that
\[
    \abs{\tilde Q(t,y)}\le\const{\norm{Q}_{L^\infty_tL^\infty_y}}\int_0^t\abs{\tilde Q(s,y)}\,ds+\int_0^t\abs{\tilde\Psi(s,y)}\,ds.
\]
Gr\"onwall's inequality (in $t$, uniformly in $y$) gives \eqref{eq:tildeQ-Linfty-bound}.

Now, we derive a tame bound on $\p_t\p^k\tilde Q$ using Corollary~\ref{Cor:Multilinear-Interpolation}. Expanding the first term on the right hand side of \eqref{eq:k-derivative-Lagrangian-generic-transport} by the Leibniz rule, each resulting term is a multilinear expression of total derivative count $k$ in factors of $Q$ and $\tilde Q$. Applying Corollary~\ref{Cor:Multilinear-Interpolation}, we obtain
\[
    \norm{R^\sigma\p^k\big(F(Q)\tilde Q\big)}_{L^2_y}
    \lesssim\const{\norm{Q}_{L^\infty_tL^\infty_y}}
    \Big(\norm{R^\sigma\p^k\tilde Q}_{L^2_y}+\norm{\tilde Q}_{L^\infty_y}\norm{R^\sigma\p^kQ}_{L^2_y}\Big),
\]
so that
\[
    \norm{R^\sigma\p_t\p^k\tilde Q}_{L^2_y}
    \lesssim\const{\norm{Q}_{L^\infty_tL^\infty_y}}
    \Big(\norm{R^\sigma\p^k\tilde Q}_{L^2_y}+\norm{\tilde Q}_{L^\infty_y}\norm{R^\sigma\p^kQ}_{L^2_y}\Big)
    +\norm{R^\sigma\p^k\tilde\Psi}_{L^2_y}.
\]

Combining the above estimate,
\begin{equation*}
	\begin{aligned}
		\p_t\norm{R^\sigma\p^k\tilde Q}_{L^2_y}
    &\lesssim\const{\norm{\vec U}_{L^\infty_tW^{1,\infty}_y},\norm{\p_t\vec U}_{L^\infty_tL^\infty_y},N^0_\infty,\norm{Q}_{L^\infty_tL^\infty_y}}\\
    &\quad \times\Big(\norm{R^\sigma\p^k\tilde Q}_{L^2_y}+\norm{\tilde Q}_{L^\infty_y}\norm{R^\sigma\p^kQ}_{L^2_y}\Big)
    +\norm{R^\sigma\p^k\tilde\Psi}_{L^2_y}.
	\end{aligned}
\end{equation*}
Substituting \eqref{eq:tildeQ-Linfty-bound} for $\norm{\tilde Q}_{L^\infty_y}$,
integrating in time, and applying Gr\"onwall's inequality (using
$\tilde Q_0=0$),
\[
    \norm{R^\sigma\p^k\tilde Q(t)}_{L^2_y}
    \lesssim\calC\Big(\int_0^t\int_0^s\norm{\tilde\Psi(\tau)}_{L^\infty_y}\,d\tau\,\norm{R^\sigma\p^kQ(s)}_{L^2_y}\,ds
    +\int_0^t\norm{R^\sigma\p^k\tilde\Psi(s)}_{L^2_y}\,ds\Big),
\]
which is the claimed estimate.
\end{proof}

\begin{lem}\label{lem:weighted-L2-deformation-inverse-estimate}
For any integer $k\ge1$ and $\sigma>-\frac12$ with $k-\sigma>\frac32$, and
any $t\in[0,T]$,
\begin{equation}\label{eq:weighted-L2-deformation-inverse-estimate}
    \norm{R^\sigma\p^kF^{-1}}_{L^2_y}
    \lesssim \calC^1_L \int_0^t\norm{\vec U}_{H^{k+1,\sigma}_R}\,ds.
\end{equation}
\end{lem}
\begin{proof}
Equation \eqref{eq:transport-deformation-inverse} is of the form \eqref{eq:Lagrangian-generic-transport} with
$Q:=F^{-1}$ and $\Psi:=\nabla_y(\vec U/U^0)$, and $\norm{F^{-1}}_{L^\infty_y}\le\calC^1_L$ by
Lemma~\ref{lem:deformation-inverse-Linfty-estimate}. Since $F^{-1}(0,\cdot)=I$ and $k\ge1$, the initial term in
Lemma~\ref{lem:Lagrangian-generic-weighted-transport-estimate} vanishes, and
\begin{equation*}
		\norm{R^\sigma\p^kF^{-1}}_{L^2_y}\lesssim\calC^1_L\int_0^t\norm{R^\sigma\p^k\nabla_y\left(\frac{\vec U}{U^0}\right)}_{L^2_y}\,ds.
\end{equation*}
By Corollary~\ref{Cor:Multilinear-Interpolation}, the integrand is bounded by
$\calC^1_L\big(\norm{R^\sigma\p^{k+1}\vec U}_{L^2_y}+\norm{R^\sigma\p^k\vec U}_{L^2_y}\big)\le\calC^1_L\norm{\vec U}_{H^{k+1,\sigma}_R}$.
\end{proof}

The following lemma controls the weighted $L^2$ bounds of velocity
under the change of coordinates between Eulerian and Lagrangian
domains.

\begin{lem}\label{lem:velocity-weighted-L2-from-Eulerian-to-Lagrangian}
For any $t\in[0,T]$ and any integer $l>l_0$,
\begin{equation}\label{eq:velocity-weighted-L2-from-Eulerian-to-Lagrangian}
    \begin{aligned}
        \norm{\vec u}_{H^{2l,l+\frac{1}{2\kappa}}_r} &\lesssim \calC^1_L\norm{\vec U}_{L^\infty_tH^{2l,l+\frac{1}{2\kappa}}_R},\\
        \norm{\p_t\vec u}_{H^{2l,l+\frac{1}{2\kappa}}_r} &\lesssim \calC^1_L\norm{\vec U}_{\tilde X^{l+1,\frac{1}{2\kappa}}_R}.
    \end{aligned}
\end{equation}
\end{lem}

\begin{proof}
By Corollary~\ref{cor:weighted-L2-end-point-equivalence},
\begin{equation*}
	\norm{\vec u}_{H^{2l,l+\frac{1}{2\kappa}}_r}\lesssim \norm{r^{l+\frac{1}{2\kappa}}\vec u}_{L^2} + \norm{r^{l+\frac{1}{2\kappa}}\p^{2l}\vec u}_{L^2}.
\end{equation*}

For the first term, changing variables $x=\eta(t,y)$ and
using \eqref{eq:deformation-det-estimate} yields
\begin{equation*}
	\norm{r^{l+\frac{1}{2\kappa}}\vec u}_{L^2}^2=\int_{\Omega_t}r^{2l+\frac1\kappa}\abs{\vec u}^2dx
	=\int_{\Omega_0}R^{2l+\frac1\kappa}\abs{\vec U}^2\det Fdy
	\lesssim\norm{R^{l+\frac{1}{2\kappa}}\vec U}_{L^2}^2.
\end{equation*}

For the second term, similarly,
\begin{equation*}
	\norm{r^{l+\frac{1}{2\kappa}}\p^{2l}\vec u}_{L^2}^2=\int_{\Omega_t}r^{2l+\frac1\kappa}\abs{\p^{2l}\vec u}^2dx
	=\int_{\Omega_0}R^{2l+\frac1\kappa}\abs{(\p^{2l}\vec u)\circ\eta}^2\det Fdy
	\lesssim\norm{R^{l+\frac{1}{2\kappa}}(F^{-1}\nabla_y)^{2l}\vec U}_{L^2}^2.
\end{equation*}
Expanding $R^{l+\frac{1}{2\kappa}}(F^{-1}\nabla_y)^{2l}\vec U$ by the
Leibniz rule gives a sum of terms of the form
\[
	R^{l+\frac{1}{2\kappa}}\prod_{i=1}^{I}\p^{a_i}F^{-1}\p^c\vec U,\qquad\sum_ia_i+c=2l,\qquad a_i,c\ge1,
\]
with coefficients smoothly depending on $F^{-1}$. By
Corollary~\ref{Cor:Multilinear-Interpolation},
\begin{equation*}
	\begin{aligned}
		\norm{R^{l+\frac{1}{2\kappa}}\prod_{i=1}^{I}\p^{a_i}F^{-1}\p^c\vec U}_{L^2}
		&\lesssim\const{\norm{F^{-1}}_{L^\infty_y},\norm{\p\vec U}_{L^\infty_y}}\Big(\norm{R^{l+\frac{1}{2\kappa}}\p^{2l}\vec U}_{L^2}+\norm{R^{l+\frac{1}{2\kappa}}\p^{2l-1}F^{-1}}_{L^2}\Big)\\
		&\lesssim\calC^1_L\Big(\norm{\vec U}_{H^{2l,l+\frac{1}{2\kappa}}_R}+\int_0^t\norm{\vec U(s)}_{H^{2l,l+\frac{1}{2\kappa}}_R}\,ds\Big)\lesssim \calC^1_L\norm{\vec U}_{L^\infty_tH^{2l,l+\frac{1}{2\kappa}}_R},
	\end{aligned}
\end{equation*}
where the second inequality uses Lemma~\ref{lem:deformation-inverse-Linfty-estimate} and Lemma~\ref{lem:weighted-L2-deformation-inverse-estimate} with $k=2l-1,\sigma=l+\frac{1}{2\kappa}$, and the last inequality uses $T\leq 1$. Combining the above yields
\begin{equation*}
	\norm{\vec u}_{H^{2l,l+\frac{1}{2\kappa}}_r}\lesssim\calC^1_L\norm{\vec U}_{L^\infty_tH^{2l,l+\frac{1}{2\kappa}}_R}.
\end{equation*}

For the second inequality, by Corollary~\ref{cor:weighted-L2-end-point-equivalence},
\begin{equation*}
	\norm{\p_t\vec u}_{H^{2l,l+\frac{1}{2\kappa}}_r}\lesssim\norm{r^{l+\frac{1}{2\kappa}}\p_t\vec u}_{L^2}+\norm{r^{l+\frac{1}{2\kappa}}\p^{2l}(\p_t\vec u)}_{L^2}.
\end{equation*}
By the identity \eqref{eq:dt-velocity-identity} in the proof of Lemma~\ref{lem:velocity-Linfty-from-Eulerian-to-Lagrangian} and changing coordinates as above, for the first term, we have
\begin{equation*}
	\norm{r^{l+\frac{1}{2\kappa}}\p_t\vec u}_{L^2}^2=\int_{\Omega_0}R^{2l+\frac1\kappa}\Big|\p_t\vec U-(\nabla_y\vec U)F^{-1}\tfrac{\vec U}{U^0}\Big|^2\det F\,dy
	\lesssim\calC^1_L\Big(\norm{R^{l+\frac{1}{2\kappa}}\p_t\vec U}_{L^2}^2+\norm{R^{l+\frac{1}{2\kappa}}\p\vec U}_{L^2}^2\Big).
\end{equation*}

For the second term, similarly,
\begin{equation*}
	\begin{aligned}
		\norm{r^{l+\frac{1}{2\kappa}}\p^{2l}(\p_t\vec u)}_{L^2}^2&=\int_{\Omega_0}R^{2l+\frac1\kappa}\abs{(\p^{2l}(\p_t\vec u))\circ\eta}^2\det Fdy\\
		&\lesssim\norm{R^{l+\frac{1}{2\kappa}}(F^{-1}\nabla_y)^{2l}\p_t\vec U}_{L^2}^2+\norm{R^{l+\frac{1}{2\kappa}}(F^{-1}\nabla_y)^{2l}\Big((\nabla_y\vec U)F^{-1}\tfrac{\vec U}{U^0}\Big)}_{L^2}^2.
	\end{aligned}
\end{equation*}
Expanding $R^{l+\frac{1}{2\kappa}}(F^{-1}\nabla_y)^{2l}\p_t\vec U$ by the Leibniz rule and a direct estimate yield:
\begin{equation*}
	\begin{aligned}
		\norm{R^{l+\frac{1}{2\kappa}}\Big(\prod_i\p^{a_i}F^{-1}\Big)\p^c\p_t\vec U}_{L^2}
		&\lesssim\const{\norm{F^{-1}}_{L^\infty_y},\norm{\p_t\vec U}_{L^\infty_y}}\Big(\norm{R^{l+\frac{1}{2\kappa}}\p^{2l}\p_t\vec U}_{L^2}+\norm{R^{l+\frac{1}{2\kappa}}\p^{2l}F^{-1}}_{L^2}\Big)\\
		&\lesssim\calC^1_L\Big(\norm{\p_t\vec U}_{L^\infty_tH^{2l,l+\frac{1}{2\kappa}}_R}+\norm{\vec U}_{L^\infty_tH^{2l+2,l+\frac{1}{2\kappa}+1}_R}\Big),
	\end{aligned}
\end{equation*}
where the last inequality again is by Lemma~\ref{lem:deformation-inverse-Linfty-estimate}, Lemma~\ref{lem:weighted-L2-deformation-inverse-estimate}, Lemma~\ref{lem:WeightedEmbedding} and $T\le1$ as above.

Expanding $R^{l+\frac{1}{2\kappa}}(F^{-1}\nabla_y)^{2l}\big((\nabla_y\vec U)F^{-1}(\vec U/U^0)\big)$
by the Leibniz rule yields a sum of terms of the form
\begin{equation*}
	R^{l+\frac{1}{2\kappa}}\prod_{i_1=1}^{I_1}\p^{a_{i_1}}F^{-1}\prod_{i_2=1}^{I_2}\p^{a_{i_2}}\vec U\p^c(\p\vec U),\qquad\sum_{i_1}a_{i_1}+\sum_{i_2}a_{i_2}+c=2l,\qquad a_{i_1},a_{i_2}\ge1,
\end{equation*}
with coefficients smoothly depending on $F^{-1}$ and $\vec U$. Using Corollary~\ref{Cor:Multilinear-Interpolation},
\begin{equation*}
	\begin{aligned}
		\norm{R^{l+\frac{1}{2\kappa}}\prod_{i_1=1}^{I_1}\p^{a_{i_1}}F^{-1}\prod_{i_2=1}^{I_2}\p^{a_{i_2}}\vec U\p^c(\p\vec U)}_{L^2}
		&\lesssim\const{\norm{F^{-1}}_{L^\infty_y},\norm{\vec U}_{W^{1,\infty}_y}}\Big(\norm{R^{l+\frac{1}{2\kappa}}\p^{2l+1}\vec U}_{L^2}\\
		&\qquad+\norm{R^{l+\frac{1}{2\kappa}}\p^{2l}\vec U}_{L^2}+\norm{R^{l+\frac{1}{2\kappa}}\p^{2l}F^{-1}}_{L^2}\Big)\\
		&\lesssim\calC^1_L\norm{\vec U}_{L^\infty_tH^{2l+2,l+\frac{1}{2\kappa}+1}_R}.
	\end{aligned}
\end{equation*}
Combining the above estimates,
\begin{equation}\label{eq:vel-weighted-second-bound}
	\norm{\p_t\vec u}_{H^{2l,l+\frac{1}{2\kappa}}_r}\lesssim\calC^1_L\left(\norm{\p_t\vec U}_{L^\infty_t H^{2l,l+\frac{1}{2\kappa}}_R}+\norm{\vec U}_{L^\infty_t H^{2l+2,l+\frac{1}{2\kappa}+1}_R}\right).
\end{equation}
\end{proof}

\begin{thm}\label{thm:Eulerian-Lagrangian-tame-bound}
For any $t\in[0,T]$ and $l>\frac{1}{2\kappa}+\frac32$,
\begin{equation}
	\norm{\mfg}_{H^{2l,l+\frac{1}{2\kappa}}_r}\lesssim \calC^1_L\left(\norm{\mfG}_{H^{2l,l+\frac{1}{2\kappa}}_R}+\norm{\mfG}_{L^\infty_y}\int_0^t\norm{\vec U}_{H^{2l+2,l+\frac{1}{2\kappa}+1}_R}\,ds\right).
\end{equation}
\end{thm}
\begin{proof}
	By Corollary~\ref{cor:weighted-L2-end-point-equivalence}, we have
	\begin{equation*}
		\norm{\mfg}_{H^{2l,l+\frac{1}{2\kappa}}_r}\lesssim \norm{r^{l+\frac{1}{2\kappa}}\mfg}_{L^2_x} + \norm{r^{l+\frac{1}{2\kappa}}\p^{2l}\mfg}_{L^2_x}.
	\end{equation*}
	For the first term,
	\begin{equation*}
		\begin{aligned}
			\norm{r^{l+\frac{1}{2\kappa}} \mfg}_{L^2}^2 &= \int_{\Omega_t}r^{2l+\frac{1}{\kappa}}\abs{\mfG\circ\eta^{-1}}^2\,dx\\
			&=\int_{\Omega_0}R^{2l+\frac{1}{\kappa}}\abs{\mfG}^2\det F\,dy\\
			&\lesssim\norm{R^{l+\frac{1}{2\kappa}}\mfG}_{L^2_y}^2,
		\end{aligned}
	\end{equation*}
	where the last inequality uses the $L^\infty$ bound on $\det F$ from
	\eqref{eq:deformation-det-estimate}.

	For the second term,
	\begin{equation*}
		\begin{aligned}
			\norm{r^{l+\frac{1}{2\kappa}}\p^{2l}\mfg}_{L^2_x}^2 &= \int_{\Omega_t}r^{2l+\frac{1}{\kappa}}\abs{\p^{2l}\mfg}^2\,dx\\
			&=\int_{\Omega_0}R^{2l+\frac{1}{\kappa}}\abs{(\p^{2l}\mfg)\circ\eta}^2\det F\,dy\\
			&\lesssim \norm{R^{l+\frac{1}{2\kappa}}(\p^{2l}\mfg)\circ\eta}_{L^2_y}^2.
		\end{aligned}
	\end{equation*}
	Over the Lagrangian domain, $(\p^{2l}\mfg)\circ\eta$ can be written as
	$(F^{-1}\nabla_y)^{2l}\mfG$, which expands into multilinear
	expressions in derivatives of $F^{-1}$ and $\mfG$ with a total
	derivative count $2l$. By Corollary~\ref{Cor:Multilinear-Interpolation},
	promoting either $\mfG$ or a single factor of $F^{-1}$ to top order,
	\begin{equation*}
		\begin{aligned}
			\norm{R^{l+\frac{1}{2\kappa}}(\p^{2l}\mfg)\circ\eta}_{L^2_y}
			&\lesssim \const{\norm{F^{-1}}_{L^\infty_y}}\left(\norm{R^{l+\frac{1}{2\kappa}}\p^{2l}\mfG}_{L^2_y}
			+ \norm{\mfG}_{L^\infty_y}\norm{R^{l+\frac{1}{2\kappa}}\p^{2l}F^{-1}}_{L^2_y}\right).
		\end{aligned}
	\end{equation*}
	Applying Lemma~\ref{lem:deformation-inverse-Wkinfty-estimate}, Lemma~\ref{lem:weighted-L2-deformation-inverse-estimate} (with $k=2l$ and $\sigma=l+\frac{1}{2\kappa}$) and then Lemma~\ref{lem:WeightedEmbedding} yields
	\begin{equation*}
			\norm{R^{l+\frac{1}{2\kappa}}(\p^{2l}\mfg)\circ\eta}_{L^2_y}\lesssim \calC^1_L \left(\norm{R^{l+\frac{1}{2\kappa}}\p^{2l}\mfG}_{L^2_y}+\norm{\mfG}_{L^\infty_y}\int_0^t\norm{\vec U}_{H^{2l+2,l+\frac{1}{2\kappa}+1}_R}\,ds\right).
	\end{equation*}
	Combining the estimates for the first and second terms, and absorbing
	the constant into $\calC^1_L$, yields the desired inequality.
\end{proof}

For the reverse direction, in which norms of Lagrangian variables are bounded by those of Eulerian variables, we have the following analogous results.

\begin{lem}\label{lem:basic-Linfty-estimate-deformation-Eulerian}
The deformation matrix $F=\nabla_y\eta$ satisfies
\begin{equation}\label{eq:Eulerian-deformation-det-estimate}
	\const{\norm{\vec u}_{L^\infty_t W^{1,\infty}_x}}^{-1}\leq \det F(t,y)\leq \const{\norm{\vec u}_{L^\infty_tW^{1,\infty}_x}}, \qquad  t\in[0,T],\ y\in\Omega_0,
\end{equation}
and its pushforward to the Eulerian domain, $F\circ\eta^{-1}$, satisfies
\begin{equation}\label{eq:Eulerian-deformation-Linfty-estimate}
	\norm{F\circ\eta^{-1}}_{L^\infty_x}\lesssim\const{\norm{\vec u}_{L^\infty_tW^{1,\infty}_x}},\qquad t\in[0,T].
\end{equation}
\end{lem}

\begin{proof}
	From the equation \eqref{eq:eta-defining-eq}, we can deduce the
	equation for $\det F$:
	\begin{equation*}
		\p_t\det F = \left(\mathrm{div}_x\left(\frac{\vec u}{u^0}\right)\right)\circ\eta\,\det F,\qquad F(0,y)= I.
	\end{equation*}
	We can directly solve the equation for $\det F$:
	\begin{equation*}
		\det F = \exp\left(\int_0^t\left(\mathrm{div}_x\left(\frac{\vec u}{u^0}\right)\right)\circ\eta\,ds\right).
	\end{equation*}
	By a direct estimate, we have
	\begin{equation*}
		\calC(\norm{\vec u}_{L^\infty_tW^{1,\infty}_x})^{-1}\leq\exp\left(-T\norm{\frac{\vec u}{u^0}}_{L^\infty_tW^{1,\infty}_x}\right)\leq \det F\leq\exp\left(T\norm{\frac{\vec u}{u^0}}_{L^\infty_tW^{1,\infty}_x}\right)\leq\calC(\norm{\vec u}_{L^\infty_tW^{1,\infty}_x}),
	\end{equation*}
	which is \eqref{eq:Eulerian-deformation-det-estimate}.

	Taking $\nabla_y$ of \eqref{eq:eta-defining-eq} we have
	\begin{equation*}
		F = I + \int_0^t\nabla_x\left(\frac{\vec u}{u^0}\right)\circ\eta\,F\,ds.
	\end{equation*}
	By a direct estimate and Gr\"onwall's inequality,
	\begin{equation*}
		\abs{F}\lesssim\const{\norm{\vec u}_{L^\infty_tW^{1,\infty}_x}}
	\end{equation*}
	for any $t\in[0,T]$. Since composition with the diffeomorphism $\eta^{-1}$
	does not change the $L^\infty_x$ norm, this is
	\eqref{eq:Eulerian-deformation-Linfty-estimate}.
\end{proof}

\begin{lem}\label{lem:Eulerian-deformation-Wkinfty-estimate}
The pushforward $F\circ\eta^{-1}$ of the deformation matrix to the
Eulerian domain satisfies
\begin{equation*}
	\norm{F\circ\eta^{-1}}_{W^{k,\infty}_x}\lesssim\const{\norm{\vec u}_{L^\infty_tW^{k+1,\infty}_x}}
\end{equation*}
for any integer $k\geq 0$ and $t\in[0,T]$.
\end{lem}
\begin{proof}
	From \eqref{eq:eta-defining-eq}, we obtain the following equation
	for $F\circ\eta^{-1}$:
	\begin{equation}\label{eq:deformation-matrix-Eulerian-eq}
		D_t\big(F\circ\eta^{-1}\big)=u^0\,\nabla_x\left(\frac{\vec u}{u^0}\right)F\circ\eta^{-1},
		\qquad \big(F\circ\eta^{-1}\big)(0,\cdot)=I.
	\end{equation}
	Applying Lemma~\ref{lem:generic-Linfty-transport-estimate} with $q:=F\circ\eta^{-1}$, $\psi:=u^0\nabla_x(\vec u/u^0)$ yields the result.
\end{proof}

\begin{thm}\label{thm:Lagrangian-Eulerian-weighted-Sobolev-comparison}
For any $t\in[0,T]$ and integer $l\geq0$,
\begin{equation}\label{eq:Lagrangian-Eulerian-weighted-Sobolev-comparison}
	\norm{\mfG}_{H^{2l,l+\frac{1}{2\kappa}}_R(\Omega_0)}\lesssim\calC\left(\norm{\vec u}_{L^\infty_tW^{2l,\infty}_x}\right)\norm{\mfg}_{H^{2l,l+\frac{1}{2\kappa}}_r}.
\end{equation}
\end{thm}

\begin{proof}
The proof is analogous to that of Theorem~\ref{thm:Eulerian-Lagrangian-weighted-Sobolev-comparison}, with the roles of $(R,\mfG,F^{-1})$ and $(r,\mfg,F\circ\eta^{-1})$ interchanged, using \eqref{eq:Eulerian-deformation-det-estimate} in place of \eqref{eq:deformation-det-estimate} and Lemma~\ref{lem:Eulerian-deformation-Wkinfty-estimate} in place of Lemma~\ref{lem:deformation-inverse-Wkinfty-estimate}; since the constant is allowed to depend on $\norm{\vec u}_{L^\infty_t W^{2l,\infty}_x}$, all factors of $F\circ\eta^{-1}$ are bounded in $L^\infty$. We omit the details.
\end{proof}

\begin{lem}\label{lem:weighted-L2-deformation-Eulerian-estimate}
For any integer $k\ge2$ and $\sigma>-\frac12$ with $k-\sigma>\frac32$, and
any $t\in[0,T]$,
\begin{equation}\label{eq:weighted-L2-deformation-Eulerian-estimate}
    \norm{r^\sigma\p^k(F\circ\eta^{-1})}_{L^2_x}
    \lesssim \calC^1_E \int_0^t\norm{\vec u}_{H^{k+1,\sigma}_r}\,ds.
\end{equation}
\end{lem}
\begin{proof}
Recall that $F\circ\eta^{-1}$ satisfies \eqref{eq:deformation-matrix-Eulerian-eq}. By Lemma~\ref{lem:Eulerian-deformation-Wkinfty-estimate}, $\norm{F\circ\eta^{-1}}_{L^\infty_tL^\infty_x} \lesssim\calC(\norm{\vec u}_{L^\infty_tW^{1,\infty}_x})<\infty$, so Lemma~\ref{lem:generic-weighted-transport-estimate} applies with $q:=F\circ\eta^{-1}$ and $\psi:=u^0\nabla_x(\vec u/u^0)$. Since $(F\circ\eta^{-1})(0,\cdot)=I$ is constant and $k\ge2$, the initial-value term vanishes, $\norm{r_0^\sigma\p^k(F\circ\eta^{-1})(0,\cdot)}_{L^2_x}=\norm{r_0^\sigma\p^kI}_{L^2_x}=0$, giving
\begin{equation*}
	\norm{r^\sigma\p^k(F\circ\eta^{-1})}_{L^2_x}
	\lesssim\calC^1_E\left(\int_0^t\norm{\vec u}_{H^{k+1,\sigma}_r}
	+\norm{r^\sigma\p^k\nabla_x\left(\frac{\vec u}{u^0}\right)}_{L^2_x}\,ds\right).
\end{equation*}
By Corollary~\ref{Cor:Multilinear-Interpolation}, the second term inside
the integral is controlled by $\const{\norm{\vec u}_{L^\infty_tW^{1,\infty}_x}}\norm{\vec u}_{H^{k+1,\sigma}_r}$,
so the right-hand side is bounded by
$\calC^1_E\int_0^t\norm{\vec u}_{H^{k+1,\sigma}_r}\,ds$, which is
\eqref{eq:weighted-L2-deformation-Eulerian-estimate}.
\end{proof}

\begin{thm}\label{thm:Lagrangian-Eulerian-tame-comparison}
For any $t\in[0,T]$ and $l>\frac{1}{2\kappa}+\frac32$,
\begin{equation}\label{eq:Lagrangian-Eulerian-tame-comparison}
	\norm{\mfG}_{H^{2l,l+\frac{1}{2\kappa}}_R(\Omega_0)}\lesssim\calC^1_E\left(\norm{\mfg}_{H^{2l,l+\frac{1}{2\kappa}}_r}+\norm{\mfg}_{L^\infty_x}\int_0^t\norm{\vec u(s)}_{H^{2l+2,l+\frac{1}{2\kappa}+1}}\,ds\right).
\end{equation}
\end{thm}
\begin{proof}
The proof is identical to that of Theorem~\ref{thm:Eulerian-Lagrangian-tame-bound}, with the roles of $(\mfG,\mfg,F^{-1})$ and $(\mfg,\mfG,F\circ\eta^{-1})$ interchanged. The only difference is the treatment of term $\norm{r^{l+\frac{1}{2\kappa}}\p^{2l}(F\circ\eta^{-1})}_{L^2_x}$: applying Lemma~\ref{lem:weighted-L2-deformation-Eulerian-estimate} (with $k=2l$ and $\sigma=l+\frac{1}{2\kappa}$) gives the bound $\calC^1_E\int_0^t\norm{\vec u}_{H^{2l+1,l+\frac{1}{2\kappa}}_r}\,ds$, which is then embedded into $\int_0^t\norm{\vec u(s)}_{H^{2l+2,l+\frac{1}{2\kappa}+1}}\,ds$ exactly as in that proof.
\end{proof}

\subsection{First-order tame estimate of the linearization of the iteration map}\label{sec:first-order-tame-estimate-iteration-map}
Now, we move the energy estimate in Theorem~\ref{thm:modified-linearized-with-xi-local-wellposedness} over to the Lagrangian domain:
\begin{thm}\label{thm:Lagrangian-energy-bound}
	For any integer $l\geq 1$ and any $\vec G\in L^\infty_tH^{2l,l+\frac{1}{2\kappa}}_R$, there exists a solution $\delta \vec U$ to the equation
	\begin{equation*}
		(D\mfF)(\vec U)\delta\vec U = \vec G
	\end{equation*}
	such that
	\begin{equation}\label{eq:Lagrangian-energy-bound}
		\norm{\delta\vec U}_{\tilde X^{l,\frac{1}{2\kappa}}_R}\lesssim\calC^{2l+3}_L\norm{\vec G}_{\tilde Y^{l,\frac{1}{2\kappa}}_R}.
	\end{equation}
\end{thm}

\begin{proof}
	Let $\vec g:=\vec G\circ\eta^{-1}$. Taking the variation $\delta$ of $\vec U=\vec u\circ\eta$, we obtain
	\begin{equation}\label{eq:relation-between-delta-vec-U-and-delta-vec-u}
	\delta\vec U = \delta\vec u\circ\eta + (\p_i\vec u)\circ\eta\,\delta\eta^i.
	\end{equation}
	Taking the norm $\norm{\cdot}_{H^{2l,l+\frac{1}{2\kappa}}_R}$ in \eqref{eq:relation-between-delta-vec-U-and-delta-vec-u} yields
	\begin{equation*}
		\norm{\delta\vec U}_{H^{2l,l+\frac{1}{2\kappa}}_R}\lesssim \norm{\delta\vec u\circ\eta }_{H^{2l,l+\frac{1}{2\kappa}}_R} + \norm{(\p_i\vec u)\circ\eta\,\delta\eta^i}_{H^{2l,l+\frac{1}{2\kappa}}_R}.
	\end{equation*}
	For the first term, we apply Theorem~\ref{thm:Lagrangian-Eulerian-weighted-Sobolev-comparison}, Theorem~\ref{thm:modified-linearized-with-xi-local-wellposedness} and then Theorem~\ref{thm:Eulerian-Lagrangian-weighted-Sobolev-comparison} and obtain
	\begin{equation*}
		\norm{\delta\vec u\circ\eta }_{H^{2l,l+\frac{1}{2\kappa}}_R}\lesssim\calC^{2l}_E\norm{\delta\vec u}_{H^{2l,l+\frac{1}{2\kappa}}_r}\lesssim \calC^{2l+3}_E\int_0^t\norm{\vec g(s)}_{H^{2l,l+\frac{1}{2\kappa}}_r}ds\lesssim \calC^{2l+3}_L\int_0^t\norm{\vec G(s)}_{H^{2l,l+\frac{1}{2\kappa}}_R}\,ds.
	\end{equation*}
	For the second term, using the same theorems as for the first term,
	\begin{equation*}
		\begin{aligned}
			\norm{(\p_i\vec u)\circ\eta\delta\eta^i}_{H^{2l,l+\frac{1}{2\kappa}}_R}&\lesssim\calC^{2l}_E\norm{(\p_i\vec u)\xi^i}_{H^{2l,l+\frac{1}{2\kappa}}_r}\\
			&\lesssim \calC^{2l+1}_E\norm{\xi}_{H^{2l,l+\frac{1}{2\kappa}}_r}\\
			&\lesssim\calC^{2l+3}_E\int_0^t\norm{\vec g(s)}_{H^{2l,l+\frac{1}{2\kappa}}_r}ds\\
			&\lesssim \calC^{2l+3}_L\int_0^t\norm{\vec G(s)}_{H^{2l,l+\frac{1}{2\kappa}}_R}ds.
		\end{aligned}
	\end{equation*}
	Combining the two terms and taking $L^\infty$ in time, since $T\leq 1$, we have
	\begin{equation}\label{eq:Lagrangian-energy-bound1}
	\norm{\delta\vec U}_{L^\infty_tH^{2l,l+\frac{1}{2\kappa}}_R}\lesssim\calC^{2l+3}_L\norm{\vec G}_{\tilde Y^{l,\frac1{2\kappa}}_R}.
	\end{equation}
	
	Now, we deal with the time derivative part $\norm{\p_t\delta\vec U}_{L^\infty_tH^{2l-2,l+\frac{1}{2\kappa}-1}_R}$. Differentiating \eqref{eq:relation-between-delta-vec-U-and-delta-vec-u} in time yields
	\begin{equation}\label{eq:time-differentiated-relation}
		\p_t\delta\vec U = \left(\frac{1}{u^0}D_t\delta\vec u\right)\circ\eta + \left(\frac{1}{u^0}D_t\p_i\vec u\xi^i\right)\circ\eta + \left(\p_i\vec u\frac{1}{u^0}D_t\xi^i\right)\circ\eta.
	\end{equation}
	Taking the norm $\norm{\cdot}_{H^{2l-2,l+\frac{1}{2\kappa}-1}_R}$, we then bound each of the three terms separately. For the term $\left(\frac1{u^0}D_t\delta\vec u\right)\circ\eta$, by Theorem~\ref{thm:Lagrangian-Eulerian-weighted-Sobolev-comparison}, we have
	\begin{equation*}
		\norm{\left(\frac{1}{u^0}D_t\delta\vec u\right)\circ\eta}_{H^{2l-2,l+\frac{1}{2\kappa}-1}_R}\lesssim \calC^{2l-2}_E\norm{\frac{1}{u^0}D_t\delta\vec u}_{H^{2l-2,l+\frac{1}{2\kappa}-1}_r}.
	\end{equation*}
	By Remark~\ref{rmk:modified-linearized-eqs-time-derivative}, we have
	\begin{equation*}
		\frac{1}{u^0}D_t\delta\vec u = \calF\p\delta r + \calF\p\delta\vec u + \calF\p\delta\pi
            + \mfL\delta r + \mfL\delta\vec u + \mfL\delta\pi + \calF\left(\vec g-\p_j\vec\calN_{\mathrm{mo}}(\varphi)\xi^j\right),
	\end{equation*}
	where $\calF=\calF(r,a_0,\vec u)$ denotes a coefficient smooth in a neighborhood of the range of $(r,a_0,\vec u)$, and $\mfL=\sum_X\calF_X(r,a_0,\vec u)\p X+\calF(r,a_0,\vec u)\p_t\vec u+\calF(r,a_0,\vec u)$ for $X\in\{r,a_0,\vec u\}$. Above, both $\calF$ and $\mfL$ are subject to change term by term. 
	
	We bound $\calF\p\delta r$ in detail. The other terms can be bounded in a similar way. By Corollary~\ref{cor:weighted-L2-end-point-equivalence},
	\begin{equation*}
		\norm{\calF\p\delta r}_{H^{2l-2,l+\frac{1}{2\kappa}-1}_r}\lesssim \norm{r^{l+\frac{1}{2\kappa}-1}(\calF\p\delta r)}_{L^2} + \norm{r^{l+\frac{1}{2\kappa}-1}\p^{2l-2}(\calF\p\delta r)}_{L^2}.
	\end{equation*}
	The first term can be bounded directly with Lemma~\ref{lem:basic-Linfty-coefficient-estimate} and Lemma~\ref{lem:WeightedEmbedding},
	\begin{equation*}
		\norm{r^{l+\frac{1}{2\kappa}-1}(\calF\p\delta r)}_{L^2}\lesssim \calC^1_E\norm{r^{l+\frac{1}{2\kappa}-1}\p\delta r}_{L^2}\lesssim \calC^1_E\norm{\delta r}_{H^{2l,l+\frac{1}{2\kappa}}_r}\lesssim \calC^1_E\norm{\delta r}_{H^{2l,l+\frac{1}{2\kappa}-\frac{1}{2}}_r}.
	\end{equation*}
	For the second term, apply Lemma~\ref{lem:coefficient-product-estimate}(a) with $k=2l-2,\sigma=l+\frac{1}{2\kappa}-1$ and then Lemma~\ref{lem:basic-Linfty-coefficient-estimate} and Lemma~\ref{lem:WeightedEmbedding}
	\begin{equation*}
		\norm{r^{l+\frac{1}{2\kappa}-1}\p^{2l-2}\left(\calF\p\delta r\right)}_{L^2}\lesssim\calC^{2l-1}_E\norm{r^{l+\frac{1}{2\kappa}-1}\p^{2l-1}\delta r}_{L^2}\lesssim \calC^{2l-1}_E\norm{\delta r}_{H^{2l,l+\frac{1}{2\kappa}}_r}\lesssim \calC^{2l-1}_E\norm{\delta r}_{H^{2l,l+\frac{1}{2\kappa}-\frac{1}{2}}_r}.
	\end{equation*}
	Thus, we obtain
	\begin{equation*}
		\norm{\calF\p\delta r}_{H^{2l-2,l+\frac{1}{2\kappa}-1}_r}\lesssim\calC^{2l-1}_E\norm{\delta r}_{H^{2l,l+\frac{1}{2\kappa}-\frac12}_r},
	\end{equation*}
	and similarly
	\begin{equation*}
	\begin{aligned}
		\norm{\calF\p\delta \vec u}_{H^{2l-2,l+\frac{1}{2\kappa}-1}_r}&\lesssim\calC^{2l-1}_E\norm{\delta \vec  u}_{H^{2l,l+\frac{1}{2\kappa}}_r},\\
		\norm{\calF\p\delta \pi}_{H^{2l-2,l+\frac{1}{2\kappa}-1}_r}&\lesssim\calC^{2l-1}_E\norm{\delta \pi}_{H^{2l,l+\frac{1}{2\kappa}}_r},\\
		\norm{\mfL\delta r}_{H^{2l-2,l+\frac{1}{2\kappa}-1}_r}&\lesssim \calC^{2l}_E\norm{\delta r}_{H^{2l-2,l+\frac{1}{2\kappa}-\frac{3}{2}}_r}\lesssim \calC^{2l}_E\norm{\delta r}_{H^{2l,l+\frac{1}{2\kappa}-\frac{1}{2}}_r},\\
		\norm{\mfL\delta \vec u}_{H^{2l-2,l+\frac{1}{2\kappa}-1}_r}&\lesssim \calC^{2l}_E\norm{\delta \vec u}_{H^{2l-2,l+\frac{1}{2\kappa}-1}_r}\lesssim \calC^{2l}_E\norm{\delta \vec u}_{H^{2l,l+\frac{1}{2\kappa}}_r},\\
		\norm{\mfL\delta \pi}_{H^{2l-2,l+\frac{1}{2\kappa}-1}_r}&\lesssim \calC^{2l}_E\norm{\delta \pi}_{H^{2l-2,l+\frac{1}{2\kappa}-1}_r}\lesssim \calC^{2l}_E\norm{\delta \pi}_{H^{2l,l+\frac{1}{2\kappa}}_r},\\
		\norm{\calF\big(\vec g-\p_j\vec\calN_{\mathrm{mo}}(\varphi)\xi^j\big)}_{H^{2l-2,l+\frac{1}{2\kappa}-1}_r}&\lesssim \calC^{2l-1}_E\norm{\vec g}_{H^{2l-2,l+\frac{1}{2\kappa}-1}_r} + \calC^{2l+1}_E\norm{\xi}_{H^{2l-2,l+\frac{1}{2\kappa}-1}_r}\\
		&\qquad \lesssim \calC^{2l-1}_E\norm{\vec g}_{H^{2l,l+\frac{1}{2\kappa}}_r} + \calC^{2l+1}_E\norm{\xi}_{H^{2l,l+\frac{1}{2\kappa}}_r}.
	\end{aligned}
	\end{equation*}
	Combining the above estimates and using Theorem~\ref{thm:modified-linearized-with-xi-local-wellposedness} and Theorem~\ref{thm:Eulerian-Lagrangian-weighted-Sobolev-comparison},
 	\begin{equation*}
		\begin{aligned}
			\norm{\left(\frac{1}{u^0}D_t\delta\vec u\right)\circ\eta}_{H^{2l-2,l+\frac{1}{2\kappa}-1}_R}&\lesssim\calC^{2l+1}_E\left(\norm{\delta\vec\phi}_{\calH^{2l}}+\norm{\xi}_{H^{2l,l+\frac{1}{2\kappa}}_r}+\norm{\vec g}_{H^{2l,l+\frac{1}{2\kappa}}_r}\right)\\
			&\lesssim\calC^{2l+3}_E\left(\int_0^t\norm{\vec g}_{H^{2l,l+\frac{1}{2\kappa}}_r}ds + \norm{\vec g}_{H^{2l,l+\frac{1}{2\kappa}}_r}\right)\\
			&\lesssim\calC^{2l+3}_L\left(\int_0^t\norm{\vec G(s)}_{H^{2l,l+\frac{1}{2\kappa}}_R}ds+\norm{\vec G}_{H^{2l,l+\frac{1}{2\kappa}}_R}\right).
		\end{aligned}
	\end{equation*}
	Taking the $L^\infty$ in time,
	\begin{equation*}
		\norm{\left(\frac{1}{u^0}D_t\delta\vec u\right)\circ\eta}_{L^\infty_t H^{2l-2,l+\frac{1}{2\kappa}-1}_R}\lesssim\calC^{2l+3}_L\norm{\vec G}_{\tilde Y^{l,\frac1{2\kappa}}_R}.
	\end{equation*}
	
	By a similar argument as above, we can also show that
	\begin{equation*}
	 	\norm{\left(\frac{1}{u^0}D_t\p_i\vec u\xi^i\right)\circ\eta}_{L^\infty_tH^{2l-2,l+\frac{1}{2\kappa}-1}_R},\norm{\left(\p_i\vec u\frac{1}{u^0}D_t\xi^i\right)\circ\eta}_{L^\infty_tH^{2l-2,l+\frac{1}{2\kappa}-1}_R}\lesssim\calC^{2l+3}_L\norm{\vec G}_{\tilde Y^{l,\frac1{2\kappa}}_R}.
	\end{equation*}
	Thus, 
	\begin{equation}\label{eq:Lagrangian-energy-bound2}
	\norm{\p_t\delta\vec U}_{L^\infty_tH^{2l-2,l+\frac{1}{2\kappa}-1}_R}\lesssim\calC^{2l+3}_L\norm{\vec G}_{\tilde Y^{l,\frac1{2\kappa}}_R}.
	\end{equation}
	Summing \eqref{eq:Lagrangian-energy-bound1} and \eqref{eq:Lagrangian-energy-bound2} yields the desired bound \eqref{eq:Lagrangian-energy-bound}.
\end{proof}

Next, we move the first-order tame estimate in Theorem~\ref{thm:modified-linearized-with-xi-tame-estimate} to the Lagrangian domain:
\begin{thm}\label{thm:Lagrangian-tame-bound}
	For any $l> l_0$, the solution $\delta \vec U$ to the equation
	\begin{equation*}
		(D\mfF)(\vec U)\delta\vec U = \vec G
	\end{equation*}
	satisfies the estimate
	\begin{equation}\label{eq:Lagrangian-tame-bound}
		\norm{\delta\vec U}_{\tilde X^{l,\frac{1}{2\kappa}}_R}\lesssim\calC^{2l_0+3}_L\left(\norm{\vec G}_{\tilde Y^{l,\frac{1}{2\kappa}}_R} + \norm{\vec G}_{\tilde Y^{l_0,\frac{1}{2\kappa}}_R}\left(N^{l+2}_w+\norm{\vec U}_{\tilde X^{l+4,\frac{1}{2\kappa}}_R}\right)\right).
	\end{equation}
\end{thm}

\begin{proof}
Taking the norm $\norm{\cdot}_{H^{2l,l+\frac1{2\kappa}}_R}$ of \eqref{eq:relation-between-delta-vec-U-and-delta-vec-u} yields 

\begin{equation*}
	\norm{\delta\vec U}_{H^{2l,l+\frac1{2\kappa}}_R}\lesssim \norm{\delta\vec u\circ\eta}_{H^{2l,l+\frac1{2\kappa}}_R} + \norm{(\p_i\vec u)\circ\eta\,\delta\eta^i}_{H^{2l,l+\frac1{2\kappa}}_R}.
\end{equation*}
We bound both terms on the right hand side of the above inequality by the same argument. For a quantity $\Xi\in\{\delta\vec u,(\p_i\vec u)\xi^i\}$ over the Eulerian domain, applying Theorem~\ref{thm:Lagrangian-Eulerian-tame-comparison} and since $T\leq 1$,
\begin{equation*}
	\norm{\Xi\circ\eta}_{H^{2l,l+\frac1{2\kappa}}_R}\lesssim\calC^1_E\left(\norm{\Xi}_{H^{2l,l+\frac1{2\kappa}}_r}+\norm{\Xi}_{L^\infty_x}\norm{\vec u}_{L^\infty_t H^{2l+2,l+\frac1{2\kappa}+1}_r}\right).
\end{equation*}

The quantity $\norm{\Xi}_{H^{2l,l+\frac{1}{2\kappa}}_r}$ can be controlled by an application of Theorem~\ref{thm:modified-linearized-with-xi-tame-estimate}. Note that if $\Xi=\p_i\vec u\xi^i$, we apply Corollary~\ref{cor:weighted-L2-end-point-equivalence} and Corollary~\ref{Cor:Multilinear-Interpolation} to reduce to the case $\Xi=\xi$. We obtain the following:
\begin{equation*}
	\norm{\Xi}_{H^{2l,l+\frac{1}{2\kappa}}_r}\lesssim\calC^{2l_0+3}_E\left(\norm{\vec g}_{L^1_tH^{2l,l+\frac{1}{2\kappa}}_r}+\norm{\vec g}_{L^1_tH^{2l_0,l_0+\frac{1}{2\kappa}}_r}\left(N^{l+2}_w+\norm{\vec u}_{\tilde X^{l+4,\frac{1}{2\kappa}}_r}\right)\right).
\end{equation*}

Then, we apply Lemma~\ref{lem:WeightedEmbedding} and Theorem~\ref{thm:modified-linearized-with-xi-local-wellposedness} at level $l_0$ and obtain
\begin{equation*}
	\norm{\Xi}_{L^\infty_x}\lesssim\norm{\Xi}_{H^{2l_0,l_0+\frac{1}{2\kappa}}_r}\lesssim\calC^{2l_0+3}_E\int_0^t\norm{\vec g}_{H^{2l_0,l_0+\frac{1}{2\kappa}}_r}ds.
\end{equation*}

Next, we apply the bounds in Theorem~\ref{thm:Eulerian-Lagrangian-tame-bound}, Theorem~\ref{thm:Eulerian-Lagrangian-weighted-Sobolev-comparison} and Lemma~\ref{lem:velocity-weighted-L2-from-Eulerian-to-Lagrangian} to bound the norms of $\vec g$ and $\vec u$ by the norms of $\vec G$ and $\vec U$. Taking the $L^\infty$ norm in time and due to $T\leq1$, we finally obtain
\begin{equation*}
	\norm{\Xi\circ\eta}_{H^{2l,l+\frac1{2\kappa}}_R}\lesssim \calC^{2l_0+3}_L\left(\norm{\vec G}_{\tilde Y^{l,\frac1{2\kappa}}_R}+\norm{\vec G}_{\tilde Y^{l_0,\frac1{2\kappa}}_R}\left(N^{l+2}_w+\norm{\vec U}_{\tilde X^{l+4,\frac1{2\kappa}}_R}\right)\right)
\end{equation*}
and thus
\begin{equation*}
	\norm{\delta\vec U}_{H^{2l,l+\frac1{2\kappa}}_R}\lesssim\calC^{2l_0+3}_L\left(\norm{\vec G}_{\tilde Y^{l,\frac1{2\kappa}}_R}+\norm{\vec G}_{\tilde Y^{l_0,\frac1{2\kappa}}_R}\left(N^{l+2}_w+\norm{\vec U}_{\tilde X^{l+4,\frac1{2\kappa}}_R}\right)\right).
\end{equation*}

Now we differentiate \eqref{eq:relation-between-delta-vec-U-and-delta-vec-u} in time and obtain \eqref{eq:time-differentiated-relation}. Taking the norm $\norm{\cdot}_{H^{2l-2,l+\frac{1}{2\kappa}-1}_R}$, we obtain
\begin{equation*}
	\begin{aligned}
		\norm{\p_t\delta\vec U}_{H^{2l-2,l+\frac{1}{2\kappa}-1}_R}&\lesssim\norm{\left(\frac{1}{u^0}D_t\delta\vec u\right)\circ\eta}_{H^{2l-2,l+\frac{1}{2\kappa}-1}_R} + \norm{\left(\frac{1}{u^0}(D_t\p_i\vec u)\xi^i\right)\circ\eta}_{H^{2l-2,l+\frac{1}{2\kappa}-1}_R}\\
		&\qquad + \norm{\left(\p_i\vec u\frac{1}{u^0}D_t\xi^i\right)\circ\eta}_{H^{2l-2,l+\frac{1}{2\kappa}-1}_R}\\
		&= (i) + (ii) + (iii).
	\end{aligned}
\end{equation*}

Applying Corollary~\ref{cor:weighted-L2-end-point-equivalence} to $(i)$, we obtain
\begin{equation*}
	\begin{aligned}
		\norm{\Big(\frac{1}{u^0}D_t\delta\vec u\Big)\circ\eta}_{H^{2l-2,l+\frac{1}{2\kappa}-1}_R}&\lesssim\norm{R^{l+\frac{1}{2\kappa}-1}\Big(\frac{1}{u^0}D_t\delta\vec u\Big)\circ\eta}_{L^2} +\norm{R^{l+\frac{1}{2\kappa}-1}\p_y^{2l-2}\Big(\Big(\frac{1}{u^0}D_t\delta\vec u\Big)\circ\eta\Big)}_{L^2}\\
		&=(i,i) + (i,ii).
	\end{aligned}
\end{equation*}

To bound $(i,i)$, we change coordinates and apply Lemma~\ref{lem:basic-Linfty-estimate-deformation-Eulerian} to obtain
\begin{equation*}
	\norm{R^{l+\frac{1}{2\kappa}-1}\Big(\frac{1}{u^0}D_t\delta\vec u\Big)\circ\eta}_{L^2}
	\lesssim\calC^1_E\norm{r^{l+\frac{1}{2\kappa}-1}\frac{1}{u^0}D_t\delta\vec u}_{L^2}.
\end{equation*}

From Remark~\ref{rmk:modified-linearized-eqs-time-derivative}, we have
\begin{equation}\label{eq:Dt-delta-vec-u-bookkeeping}
		\frac{1}{u^0}D_t\delta\vec u = \calF\p\delta r + \calF\p\delta\vec u + \calF\p\delta\pi
            + \mfL\delta r + \mfL\delta\vec u + \mfL\delta\pi + \calF\left(\vec g-\p_j\vec\calN_{\mathrm{mo}}(\varphi)\xi^j\right),
\end{equation}
where $\calF=\calF(r,a_0,\vec u)$ denotes a coefficient smooth in a neighborhood of the range of $(r,a_0,\vec u)$, and $\mfL=\sum_X\calF_X(r,a_0,\vec u)\p X+\calF(r,a_0,\vec u)\p_t\vec u+\calF(r,a_0,\vec u)$ for $X\in\{r,a_0,\vec u\}$. 

By a direct estimate of bounding all coefficients in $L^\infty$ and applying Lemma~\ref{lem:WeightedEmbedding},
\begin{equation}\label{eq:(i,i)-bound}
	\norm{r^{l+\frac{1}{2\kappa}-1}\frac{1}{u^0}D_t\delta\vec u}_{L^2}
	\lesssim\calC^3_E\Big(\norm{\delta r}_{H^{2l,l+\frac{1}{2\kappa}-\frac12}_r}+\norm{\delta\vec u}_{H^{2l,l+\frac{1}{2\kappa}}_r}+\norm{\delta\pi}_{H^{2l,l+\frac{1}{2\kappa}}_r}+\norm{\vec g}_{H^{2l,l+\frac{1}{2\kappa}}_r}+\norm\xi_{H^{2l,l+\frac{1}{2\kappa}}_r}\Big).
\end{equation}
For $(i,ii)$, similarly as above,
\begin{equation*}
	\norm{R^{l+\frac{1}{2\kappa}-1}\p_y^{2l-2}\Big(\big(\frac{1}{u^0}D_t\delta\vec u\big)\circ\eta\Big)}_{L^2}
	\lesssim\calC^1_E\norm{r^{l+\frac{1}{2\kappa}-1}(F\circ\eta^{-1}\nabla_x)^{2l-2}\left(\frac{1}{u^0}D_t\delta\vec u\right)}_{L^2},
\end{equation*}
and we bound the right-hand side term by term with \eqref{eq:Dt-delta-vec-u-bookkeeping}.

We bound the term from $\calF\p\delta r$ in \eqref{eq:Dt-delta-vec-u-bookkeeping} in detail. The terms from $\calF\p\delta\vec u$ and $\calF\p\delta\pi$ are bounded in a similar way. Expanding the term by the Leibniz rule and applying Corollary~\ref{Cor:Multilinear-Interpolation} with $\sigma = l+\frac{1}{2\kappa}-1$ and $k=2l-1$, together with Lemma~\ref{lem:basic-Linfty-estimate-deformation-Eulerian}, Lemma~\ref{lem:weighted-L2-deformation-Eulerian-estimate}, Lemma~\ref{lem:coefficient-bounds-merged} and Lemma~\ref{lem:Weighted-Interpolation}, we obtain
\begin{equation*}
	\norm{r^{l+\frac{1}{2\kappa}-1}(F\circ\eta^{-1}\nabla_x)^{2l-2}(\calF\p\delta r)}_{L^2}
	\lesssim\calC^1_E\Big(\norm{\delta r}_{H^{2l,l+\frac{1}{2\kappa}-\frac12}_r}+\norm{\delta r}_{L^\infty_x}\big(N^l_w+\norm{\vec u}_{\tilde X^{l+2,\frac{1}{2\kappa}}_r}\big)\Big)
\end{equation*}
and the bounds from $\calF\p\delta\vec u$ and $\calF\p\delta\pi$ are 
\begin{equation*}
	\begin{aligned}
		\norm{r^{l+\frac{1}{2\kappa}-1}(F\circ\eta^{-1}\nabla_x)^{2l-2}(\calF\p\delta \vec u)}_{L^2}&\lesssim\calC^1_E\left(\norm{\delta \vec u}_{H^{2l,l+\frac{1}{2\kappa}}_r} + \norm{\delta \vec u}_{L^\infty_x}\left(N^l_w + \norm{\vec u}_{\tilde X^{l+2,\frac{1}{2\kappa}}_r}\right)\right),\\
		\norm{r^{l+\frac{1}{2\kappa}-1}(F\circ\eta^{-1}\nabla_x)^{2l-2}(\calF\p\delta \pi)}_{L^2}&\lesssim\calC^1_E\left(\norm{\delta \pi}_{H^{2l,l+\frac{1}{2\kappa}}_r} + \norm{\delta \pi}_{L^\infty_x}\left(N^l_w + \norm{\vec u}_{\tilde X^{l+2,\frac{1}{2\kappa}}_r}\right)\right).
	\end{aligned}
\end{equation*}

We bound the term from $\mfL\delta r$ in \eqref{eq:Dt-delta-vec-u-bookkeeping} in detail. The terms from $\mfL\delta\vec u$ and $\mfL\delta\pi$ are bounded in a similar way. Writing $\mfL\delta r=\sum_X(\calF_X\p X)\delta r+(\calF\,\p_t\vec u)\delta r+\calF\delta r$, the summands in the first term give terms that expand into multilinear expressions with the same $\sigma=l+\frac{1}{2\kappa}-1$ and $k=2l-1$ as the term from $\calF\p\delta r$ and thus yield the same bounds. The second and the third term give terms with weight $r^{l+\frac{1}{2\kappa}-1}$ and total derivative count $2l-2$. By Corollary~\ref{Cor:Multilinear-Interpolation} and Lemma~\ref{lem:WeightedEmbedding}, they also yield the same bound. In total, we have
\begin{equation*}
	\begin{aligned}
		\norm{r^{l+\frac{1}{2\kappa}-1}(F\circ\eta^{-1}\nabla_x)^{2l-2}(\mfL\delta r)}_{L^2}&\lesssim\calC^1_E\left(\norm{\delta r}_{H^{2l,l+\frac{1}{2\kappa}}_r} + \norm{\delta r}_{L^\infty_x}\left(N^l_w + \norm{\vec u}_{\tilde X^{l+2,\frac{1}{2\kappa}}_r}\right)\right),\\
		\norm{r^{l+\frac{1}{2\kappa}-1}(F\circ\eta^{-1}\nabla_x)^{2l-2}(\mfL\delta \vec u)}_{L^2}&\lesssim\calC^1_E\left(\norm{\delta \vec u}_{H^{2l,l+\frac{1}{2\kappa}}_r} + \norm{\delta \vec u}_{L^\infty_x}\left(N^l_w + \norm{\vec u}_{\tilde X^{l+2,\frac{1}{2\kappa}}_r}\right)\right),\\
		\norm{r^{l+\frac{1}{2\kappa}-1}(F\circ\eta^{-1}\nabla_x)^{2l-2}(\mfL\delta \pi)}_{L^2}&\lesssim\calC^1_E\left(\norm{\delta \pi}_{H^{2l,l+\frac{1}{2\kappa}}_r} + \norm{\delta \pi}_{L^\infty_x}\left(N^l_w + \norm{\vec u}_{\tilde X^{l+2,\frac{1}{2\kappa}}_r}\right)\right).
	\end{aligned}
\end{equation*}
We bound the term from $\calF(\vec g-\p_j\vec\calN_{\mathrm{mo}}(\varphi)\xi^j)$ in \eqref{eq:Dt-delta-vec-u-bookkeeping}. For the term from $\calF\vec g$, we can show the bound:
\begin{equation*}
	\norm{r^{l+\frac{1}{2\kappa}-1}(F\circ\eta^{-1}\nabla_x)^{2l-2}(\calF\vec g)}_{L^2}
	\lesssim\calC^1_E\Big(\norm{\vec g}_{H^{2l-2,l+\frac1{2\kappa}-1}_r}+\norm{\vec g}_{L^\infty_x}\big(N^{l-1}_w+\norm{\vec u}_{\tilde X^{l+1,\frac1{2\kappa}}_r}\big)\Big).
\end{equation*}

For the term from $\calF\p_j\vec\calN_{\mathrm{mo}}(\varphi)\xi^j$, called $(*)$, by Corollary~\ref{Cor:Multilinear-Interpolation}, we have
\begin{equation*}
	\begin{aligned}
		(*)&\lesssim\calC^1_E\Big(\norm{r^{l+\frac{1}{2\kappa}-1}\p^{2l-1}\xi}_{L^2}+\norm\xi_{L^\infty_x}\Big(\norm{r^{l+\frac{1}{2\kappa}-1}\p^{2l-1}\vec\calN_{\mathrm{mo}}(\varphi)}_{L^2}+\norm{r^{l+\frac{1}{2\kappa}-1}\p^{2l-1}(r,\vec u,a_0)}_{L^2}\Big)\Big)\\
		&\lesssim\calC^1_E\left(\norm{\xi}_{H^{2l,l+\frac{1}{2\kappa}}_r}+\norm{\xi}_{L^\infty_x}\left(N^{l+1}_w + \norm{\vec u}_{\tilde X^{l+3,\frac{1}{2\kappa}}_r}\right)\right).
	\end{aligned}
\end{equation*}

Assembling all the estimates above and by Lemma~\ref{lem:WeightedEmbedding}, Theorem~\ref{thm:modified-linearized-with-xi-local-wellposedness}, Theorem~\ref{thm:modified-linearized-with-xi-tame-estimate} and the change of coordinate bounds in Section~\ref{subsec:change-of-coordinate}, we have
\begin{equation*}
	\begin{aligned}
		(i)&\lesssim\calC^3_E\left(\norm{\delta\vec \phi}_{\calH^{2l}}+\norm{\xi}_{H^{2l,l+\frac{1}{2\kappa}}_r}+\left(\norm{\delta\vec\phi}_{L^\infty_x}+\norm{\xi}_{L^\infty}\right)\left(N^{l+1}_w + \norm{\vec u}_{\tilde X^{l+3,\frac{1}{2\kappa}}_r}\right)\right)\\
		&+\calC^3_E\left(\norm{\vec g}_{H^{2l,l+\frac{1}{2\kappa}}}+\norm{\vec g}_{L^\infty_x}\left(N^{l+1}_w + \norm{\vec u}_{\tilde X^{l+3,\frac{1}{2\kappa}}_r}\right)\right)\\
		&\lesssim \calC^{2l_0+3}_E\left(\norm{\vec g}_{L^1_tH^{2l,l+\frac{1}{2\kappa}}_r}+\norm{\vec g}_{L^1_tH^{2l_0,l_0+\frac{1}{2\kappa}}_r}\left(N^{l+2}_w+\norm{\vec u}_{\tilde X^{l+4,\frac{1}{2\kappa}}_r}\right)\right)\\
		&+\calC^3_E\left(\norm{\vec g}_{H^{2l,l+\frac{1}{2\kappa}}}+\norm{\vec g}_{L^\infty_x}\left(N^{l+1}_w + \norm{\vec u}_{\tilde X^{l+3,\frac{1}{2\kappa}}_r}\right)\right)\\
		&\lesssim\calC^{2l_0+3}_L\left(\norm{\vec G}_{\tilde Y^{l,\frac1{2\kappa}}_R}+\norm{\vec G}_{\tilde Y^{l_0,\frac1{2\kappa}}_R}\left(N^{l+2}_w+\norm{\vec U}_{\tilde X^{l+4,\frac1{2\kappa}}_R}\right)\right).
	\end{aligned}
\end{equation*}
With a similar argument, we also have the same bound for $(ii)$ and $(iii)$
\begin{equation*}
	(ii),(iii)\lesssim\calC^{2l_0+3}_L\left(\norm{\vec G}_{\tilde Y^{l,\frac1{2\kappa}}_R}+\norm{\vec G}_{\tilde Y^{l_0,\frac1{2\kappa}}_R}\left(N^{l+2}_w+\norm{\vec U}_{\tilde X^{l+4,\frac1{2\kappa}}_R}\right)\right),
\end{equation*}
which finishes the proof.
\end{proof}

\subsection{Second-order tame estimate of the linearization of the iteration map}\label{sec:second-order-tame-estimate-iteration-map}
In this subsection, we establish the tame estimate for the second-order derivative of the iteration map $\mfF$, with which we control the quadratic error of the Nash--Moser iteration. Let $\delta \vec U$ and $\epsilon\vec U$ be sufficiently regular $3$-dimensional vector fields defined on $[0,T]\times\Omega_0$. For any $s=(s_1,s_2)$ in a neighborhood of $0\in\R^2$, define
\begin{equation*}
	\vec U^{(s)}:=\vec U + s_1\delta\vec U + s_2\epsilon\vec U.
\end{equation*}
Let $\mfG$ denote a quantity defined through $\vec U$ (e.g. $\vec U$ itself, $\eta,R,A_0$ and $\mfF(\vec U)$). We define
\begin{equation*}
	\delta \mfG:=\frac{\p}{\p s_1}\Big|_{s=0}\mfG(\vec U^{(s)}), \qquad \epsilon\mfG:=\frac{\p}{\p s_2}\Big|_{s=0}\mfG(\vec U^{(s)}),\qquad \epsilon\delta\mfG:=\frac{\p^2}{\p s_1\p s_2}\Big|_{s=0}\mfG(\vec U^{(s)}),
\end{equation*}
and write
\begin{equation*}
	(D\mfF)(\vec U)\delta\vec U:=\delta(\mfF(\vec U)), \qquad (D^2\mfF)(\vec U)[\delta\vec U,\epsilon\vec U]:=\epsilon\delta(\mfF(\vec U)).
\end{equation*}

Substituting $P=RA_0$ into the definition of $\mfF$ in \eqref{eq:iteration-map} and eliminating $\p_tR$ and $\p_tA_0$ by equation \eqref{eq:R-A0-defining}, we can write $\mfF$ schematically as
\begin{equation}\label{eq:schematic-iteration-map}
	\mfF(\vec U)=\calF\p_t\vec U + \calF F^{-1}\nabla_y\vec U + \calF F^{-1}\nabla_y R + \calF F^{-1}\nabla_y A_0 + \calF,
\end{equation}
where $\calF =\calF(R,\vec U,A_0)$ denotes a coefficient smooth in a neighborhood of the range of $(R,\vec U,A_0)$ subject to change term by term. 

\begin{nota}
	A brace $\{X_1,\cdots,X_k\}$ appearing as a factor in a schematic term denotes a sum of terms obtained by replacing the brace with each of the $X_j$'s. For example, $\calF F^{-1}\nabla_y\{\delta R,\delta\vec U,\delta A_0\}$ in \eqref{eq:iteration-map-first-variation} denotes the sum
	\begin{equation*}
		\calF F^{-1}\nabla_y\delta R + \calF F^{-1}\nabla_y\delta \vec U + \calF F^{-1}\nabla_y\delta A_0.
	\end{equation*}
\end{nota}

Taking the $\delta$-variation, we obtain
\begin{equation}\label{eq:iteration-map-first-variation}
	\begin{aligned}
		(D\mfF)(\vec U)\delta\vec U &=\calF\p_t \delta\vec U + \calF F^{-1}\nabla_y\{\delta R,\delta\vec U,\delta A_0\} + \calF\delta F^{-1}\nabla_y\{R,\vec U,A_0\}\\
		&+\left(\calF\p_t\vec U+\calF F^{-1}\nabla_y\{R,\vec U,A_0\}+\calF\right)\{\delta R,\delta\vec U,\delta A_0\}.
	\end{aligned}
\end{equation}
Taking $\epsilon$-variation of \eqref{eq:iteration-map-first-variation}, we obtain
\begin{equation}\label{eq:iteration-map-second-variation}
	\begin{aligned}
		(D^2\mfF)(\vec U)[\delta\vec U,\epsilon\vec U]&=\calF\p_t\delta\vec U\{\epsilon R,\epsilon\vec U,\epsilon A_0\}+\calF\p_t\epsilon\vec U\{\delta R,\delta\vec U,\delta A_0\}\\
	&+\calF F^{-1}\nabla_y\{\epsilon\delta R,\epsilon\delta A_0\} + \left(\calF\p_t\vec U+\calF F^{-1}\nabla_y\{R,\vec U,A_0\}+\calF\right)\{\epsilon\delta R,\epsilon\delta A_0\}\\
	&+\calF\epsilon\delta F^{-1}\nabla_y\{R,\vec U,A_0\}+\calF\delta F^{-1}\nabla_y\{\epsilon R,\epsilon\vec U,\epsilon A_0\}+\calF\epsilon F^{-1}\nabla_y\{\delta R,\delta\vec U,\delta A_0\}\\
	&+\calF\delta F^{-1}\nabla_y\{R,\vec U, A_0\}\{\epsilon R,\epsilon\vec U,\epsilon A_0\} + \calF\epsilon F^{-1}\nabla_y\{R,\vec U,A_0\}\{\delta R,\delta\vec U,\delta A_0\}\\
	&+\calF F^{-1}\nabla_y\{\delta R,\delta\vec U,\delta A_0\}\{\epsilon R,\epsilon\vec U,\epsilon A_0\}+\calF F^{-1}\nabla_y\{\epsilon R,\epsilon\vec U,\epsilon A_0\}\{\delta R,\delta\vec U,\delta A_0\}\\
	&+\left(\calF\p_t\vec U+\calF F^{-1}\nabla_y\{R,\vec U,A_0\}+\calF\right)\{\delta R,\delta \vec U,\delta A_0\}\{\epsilon R,\epsilon\vec U,\epsilon A_0\}.
	\end{aligned}
\end{equation}

Notice that the expression for $(D^2\mfF)(\vec U)[\delta\vec U,\epsilon\vec U]$ contains variations of $R,A_0,F^{-1}$, which are not independent quantities but determined by $\vec U$, $\delta\vec U$ and $\epsilon\vec U$. Since the right hand side of the desired estimate of $(D^2\mfF)(\vec U)[\delta\vec U,\epsilon\vec U]$ only involves the norms of $\vec U$, $\delta\vec U$ and $\epsilon\vec U$, we must bound the variations of $R,A_0,F^{-1}$ by the norms of $\vec U$, $\delta\vec U$ and $\epsilon\vec U$.

\begin{lem}\label{lem:second-order-estimates-variations-of-dependent-variables}
	For any $l>l_0$ and $t\in[0,T]$, the following estimates hold,
	\begin{equation*}
		\begin{aligned}
			\norm{\delta F^{-1}}_{H^{2l,l+\frac{1}{2\kappa}}_R}&\lesssim\calC^1_L\left(\norm{\delta\vec U}_{L^\infty_t H^{2l+2,l+\frac{1}{2\kappa}+1}_R}+\norm{\delta\vec U}_{L^\infty_t W^{1,\infty}_y}\norm{\vec U}_{L^\infty_tH^{2l+2,l+\frac{1}{2\kappa}+1}_R}\right),\\
			\norm{\epsilon\delta F^{-1}}_{H^{2l,l+\frac{1}{2\kappa}}_R}&\lesssim\calC^1_L\Big(\norm{\epsilon\vec U}_{L^\infty_t H^{2l+2,l+\frac{1}{2\kappa}+1}_R}\norm{\delta\vec U}_{L^\infty_tW^{1,\infty}_y} + \norm{\epsilon\vec U}_{L^\infty_t W^{1,\infty}_y}\norm{\delta \vec U}_{L^\infty_tH^{2l+2,l+\frac{1}{2\kappa}+1}_R}\\
			&+\norm{\epsilon\vec U}_{L^\infty_t W^{1,\infty}_y}\norm{\delta\vec U}_{L^\infty_t W^{1,\infty}_y}\norm{\vec U}_{L^\infty_t H^{2l+2,l+\frac{1}{2\kappa}+1}_R}\Big),\\
			\norm{\delta R}_{H^{2l,l+\frac{1}{2\kappa}}_R} + \norm{\delta A_0}_{H^{2l,l+\frac{1}{2\kappa}}_R}&\lesssim\calC^1_L\Bigg(\norm{\delta\vec U}_{\tilde X^{l+1,\frac{1}{2\kappa}}_R} + \norm{\delta \vec U}_{\tilde Z^1_y}\left(N^l_w + \norm{\vec U}_{\tilde X^{l+2,\frac{1}{2\kappa}}_R}\right)\Bigg),
		\end{aligned}
	\end{equation*}
	and
	\begin{equation*}
		\begin{aligned}
			\norm{\epsilon\delta R}_{H^{2l,l+\frac{1}{2\kappa}}_R} + \norm{\epsilon\delta A_0}_{H^{2l,l+\frac{1}{2\kappa}}_R}&\lesssim\calC^1_L\Bigg(\norm{\delta\vec U}_{\tilde X^{l+1,\frac{1}{2\kappa}}_R}\norm{\epsilon\vec U}_{\tilde Z^1_y} + \norm{\delta\vec U}_{\tilde Z^1_y}\norm{\epsilon\vec U}_{\tilde X^{l+1,\frac{1}{2\kappa}}_R} \\
			&\qquad + \norm{\delta \vec U}_{\tilde Z^1_y} \norm{\epsilon\vec U}_{\tilde Z^1_y}\left(N^l_w + \norm{\vec U}_{\tilde X^{l+2,\frac{1}{2\kappa}}_R}\right)\Bigg).
		\end{aligned}
	\end{equation*}
\end{lem}
\begin{proof}
The variations of $F^{-1}$, $R$ and $A_0$ satisfy linear equations of the form
\eqref{eq:Lagrangian-generic-linear-transport}, obtained by a long but straightforward differentiation of
\eqref{eq:eta-defining-eq} and \eqref{eq:R-A0-defining}. The estimates follow by applying
Lemma~\ref{lem:Lagrangian-generic-linear-transport-estimate} to each equation, together with
Corollary~\ref{cor:weighted-L2-end-point-equivalence} and Corollary~\ref{Cor:Multilinear-Interpolation},
as the estimates in the proof of Theorem~\ref{thm:Lagrangian-tame-bound}. We omit the details.
\end{proof}

Using Lemma~\ref{lem:second-order-estimates-variations-of-dependent-variables} to estimate every term of \eqref{eq:iteration-map-second-variation}, we obtain the second-order tame estimate:

\begin{thm}\label{thm:second-order-tame-estimate-iteration-map}
For any $l>l_0$,
\begin{equation*}
\begin{aligned}
\norm{(D^2\mfF)(\vec U)[\delta\vec U,\epsilon\vec U]}_{\tilde Y^{l,\frac1{2\kappa}}_R}
&\lesssim\calC^2_L\Bigg(\norm{\delta\vec U}_{\tilde X^{l+2,\frac1{2\kappa}}_R}\norm{\epsilon\vec U}_{\tilde X^{l_0+1,\frac1{2\kappa}}_R}
+\norm{\delta\vec U}_{\tilde X^{l_0+1,\frac1{2\kappa}}_R}\norm{\epsilon\vec U}_{\tilde X^{l+2,\frac1{2\kappa}}_R}\\
&\qquad+\norm{\delta\vec U}_{\tilde X^{l_0+1,\frac1{2\kappa}}_R}\norm{\epsilon\vec U}_{\tilde X^{l_0+1,\frac1{2\kappa}}_R}
\Big(N^{l+1}_w+\norm{\vec U}_{\tilde X^{l+3,\frac1{2\kappa}}_R}\Big)\Bigg).
\end{aligned}
\end{equation*}
\end{thm}

\begin{proof}
The estimate follows from \eqref{eq:iteration-map-second-variation} by the same computation as in the proof of
Theorem~\ref{thm:Lagrangian-tame-bound}, using Lemma~\ref{lem:second-order-estimates-variations-of-dependent-variables} for the weighted norms of the
variations of $(R,A_0,F^{-1})$ and Lemma~\ref{lem:Lagrangian-generic-linear-transport-estimate}(b) for their
$L^\infty_y$ norms. The terms containing $\nabla_y$ of a variation require that lemma at
level $l+1$, which is why $\tilde X^{l+2,\frac1{2\kappa}}_R$ and $N^{l+1}_w+\norm{\vec U}_{\tilde X^{l+3,\frac1{2\kappa}}_R}$
appear; the terms containing a time derivative are bounded at level $l$.
\end{proof}

\section{Proof of the main theorem}\label{sec:proof-main-theorem}

The estimates of Sections~\ref{sec:linearized-eqs}--\ref{sec:tame-estimate-linearization-Eulerian} are formulated in spaces weighted by the density $r$ on the moving domain $\Omega_t$, and those of Sections~\ref{subsec:change-of-coordinate}--\ref{sec:second-order-tame-estimate-iteration-map} in spaces weighted by $R=r\circ\eta$ on the Lagrangian domain $\Omega_0$. In both cases, the weight is determined by the background $\vec U$ through \eqref{eq:R-A0-defining}, and it changes from one iterate to the next. On the other hand, the Nash--Moser scheme requires a fixed scale of Banach spaces. We therefore use the scale weighted by the initial density $r_0$. The passage from $R=R(\vec U)$ to $r_0$ rests on the following lemma:

\begin{lem}[Comparability of $R$ and $r_0$]\label{lem:comparability}
For all $(t,y)\in[0,T]\times\Omega_0$,
\begin{equation}\label{eq:comparability}
    e^{-\calC^1_L\,t}\,r_0(y)\le R(t,y)\le e^{\calC^1_L\,t}\,r_0(y).
\end{equation}
In particular $R(t,\cdot)>0$ in $\Omega_0$, $R(t,\cdot)=0$ on $\Gamma_0$ and $R(t,\cdot)$ is comparable to the reference function $d$ uniformly in $t$, i.e., the physical vacuum boundary condition is preserved on $[0,T]$.
\end{lem}

\begin{proof}
Let $(R,A_0)$ be the solution to \eqref{eq:R-A0-defining} on $[0,T]$. The first equation of \eqref{eq:R-A0-defining} in the variable $R$ can be schematically written as
\[
    \p_tR=\mfM\,R,\qquad \mfM:=\calF\,\p_t\vec U+\calF\,F^{-1}\nabla_y\vec U,
\]
where $\calF=\calF(R,\vec U,A_0)$ is smooth in a neighborhood of the range of its arguments. Since
$(R,A_0)$ is already determined, $\mfM$ is a known function on $[0,T]\times\Omega_0$, bounded by
$\abs{\mfM}\le\calC^1_L$: the range of $(R,A_0)$ is controlled by \eqref{eq:basic-Linfty-R-A0-bound} and
$F^{-1}$ by Lemma~\ref{lem:deformation-inverse-Linfty-estimate}.
For each fixed $y\in\Omega_0$, $\p_tR=\mfM R$ with $R(0,y)=r_0(y)$ is a scalar linear ODE in $t$, hence
\[
    R(t,y)=r_0(y)\exp\Big(\int_0^t\mfM(s,y)\,ds\Big),
\]
and $\abs{\int_0^t\mfM\,ds}\le\calC^1_L\,t$ gives \eqref{eq:comparability}. Since $r_0$ is a non-degenerate defining function of $\Omega_0$ and $e^{\pm\calC^1_Lt}$ is bounded
above and below on $[0,T]$ with $T\leq 1$, so is $R(t,\cdot)$, which finishes the proof.
\end{proof}

Consequently, for every $N$, $\sigma\ge0$ and $t\in[0,T]$,
\begin{equation}\label{eq:change-of-weight}
    e^{-\sigma\calC^1_L}\norm{f}_{H^{N,\sigma}_{r_0}}\le\norm{f}_{H^{N,\sigma}_{R(t)}}\le e^{\sigma\calC^1_L}\norm{f}_{H^{N,\sigma}_{r_0}},
\end{equation}
and the norms of $\tilde X^{l,\frac1{2\kappa}}$, $\tilde Y^{l,\frac1{2\kappa}}$ with weights $R$ and $r_0$ are equivalent, with constants depending on $\calC^1_L$.

Let $m_0\ge 3l_0+3$ and $k\ge m_0+7$ be integers, and set
\begin{equation}\label{eq:NM-parameters}
    \alpha:=k+1,\qquad \tilde\alpha:=2\alpha+2-m_0,\qquad m_*:=\tilde\alpha+4,\qquad k_0:=\tilde\alpha+5=2k+9-m_0.
\end{equation}
Among the above parameters in the Nash--Moser scheme, $k$ is the regularity of the solution to be constructed, $m_0$ plays the role of the index $m_0$ in \cite[Assumption 2.1--2.2]{secchi2016nash}, which is the index of the space $X_{m_0}$ containing the neighborhood $U$ of the background states in their notation. The lower bound $3l_0+3$ is what the embedding $\tilde X^{m_0,\frac1{2\kappa}}_{r_0}\hookrightarrow\tilde Z^{2l_0+3}_y$ requires (see Lemma~\ref{lem:time-dependent-norms-properties}). The parameters $\alpha$ and $\tilde\alpha$ are the auxiliary indices in \cite[Lemma~3.6, 3.8]{secchi2016nash} for the derivative losses $r=2$, $r'=3$, $s=0$, $s'=4$ (in their notation) of the tame estimates. The index $m_*$ is the highest index at which the smoothing operators and the tame estimates are used (see also Remark~\ref{rmk:finite-regularity}), and $k_0$ is the regularity of the data needed such that $\vec U_{\mathrm{app}}\in\tilde X^{m_*,\frac{1}{2\kappa}}_{r_0}$ by Lemma~\ref{lem:approximate-solution-estimate}. Indeed, for $k_0$-admissible data $(r_0,\vec u_0,\pi_0)$ as defined in Definition~\ref{defn:admissible-data} such that
\begin{equation*}
    \mfN^{k_0}:=N^{k_0}_w + \norm{\vec u_0}_{H^{2k_0,k_0+\frac{1}{2\kappa}}_{r_0}} + \norm{\frac{r_0}{\pi_0}}_{L^\infty}<\infty,
\end{equation*}
Lemma~\ref{lem:approximate-solution-estimate} gives
\begin{equation}\label{eq:approximate-solution-high-regularity-bound}
    \norm{\vec U_{\mathrm{app}}}_{\tilde X^{m_*,\frac1{2\kappa}}_{r_0}}\lesssim\const{\mfN^{k_0}}.
\end{equation}

Recall the modified iteration map $\tilde\mfF(\vec V)=\mfF(\vec V+\vec U_{\mathrm{app}})-\mfF(\vec U_{\mathrm{app}})$ defined in Section~\ref{subsec:approximate-solution-modified-iteration-map}. This is the map we use as $F$ in \cite[Theorem~2.4]{secchi2016nash}, and
$\mfF^a=-\mfF(\vec U_{\mathrm{app}})$ is the source $f$. Since $\mfF(\vec U_{\mathrm{app}})$ does not depend on
$\vec V$,
\begin{equation}\label{eq:derivatives-of-modified-iteration-map}
    d\tilde\mfF(\vec V)\delta\vec V=(D\mfF)(\vec V+\vec U_{\mathrm{app}})\delta\vec V,\qquad
    d^2\tilde\mfF(\vec V)[\delta\vec V,\epsilon\vec V]=(D^2\mfF)(\vec V+\vec U_{\mathrm{app}})[\delta\vec V,\epsilon\vec V],
\end{equation}
so the estimates of $\tilde\mfF$ can be obtained by applying the estimates of $\mfF$ in Section~\ref{sec:Lagrangian-tame} at the background state $\vec V+ \vec U_{\mathrm{app}}$.

The scales $X_m$, $Y_m$ of \cite[Definition~2.3]{secchi2016nash} are, in our case, for integers $m_0\le m\le m_*$,
\begin{equation}\label{eq:NM-scales}
    X_m:=\{\vec V\in\tilde X^{m,\frac1{2\kappa}}_{r_0}:\ \vec V(0,\cdot)=0\},\qquad Y_m:=\tilde Y^{m,\frac1{2\kappa}}_{r_0}.
\end{equation}
Above, since $\tilde X^{m,\frac1{2\kappa}}_{r_0} \subset W^{1,\infty}([0,T],H^{2m-2,m-1+\frac1{2\kappa}}_{r_0})$, the trace at $t=0$ is well-defined and thus $X_m$ is a closed subspace. Both scales of spaces are decreasing in $m$ (see Lemma~\ref{lem:time-dependent-norms-properties}(i)).

The smoothing operators $S_\theta$ of
Lemma~\ref{lem:spatial-smoothing-operator} (with $w=r_0$) act on $Y_m$, and on $X_m$ because they act at each fixed time and commute with $\p_t$, so that $(S_\theta\vec V)(0)=S_\theta(\vec V(0))=0$. Hence, the smoothing hypothesis \cite[Definition~2.3]{secchi2016nash} is satisfied on both scales for indices in $[m_0,m_*]$, which contains the range used in the proof of \cite[Theorem~2.4(i)]{secchi2016nash} (see Remark~\ref{rmk:smoothing-range}).

Finally, let $\calU:=\{\vec V\in X_{m_0}:\norm{\vec V}_{X_{m_0}}<1\}$; this plays the role of the neighbourhood $U$ of \cite[Assumptions~2.1--2.2]{secchi2016nash}. For $\vec V\in\calU$, Lemma~\ref{lem:Morrey} (with $m_0\ge3l_0+3$) and \eqref{eq:approximate-solution-high-regularity-bound} give 
\begin{equation}\label{eq:choice-of-M}
    \norm{\vec V+\vec U_{\mathrm{app}}}_{\tilde Z^{2l_0+3}_y}\le\const{\mfN^{k_0}}=:M.
\end{equation}

We fix this $M$ in Convention~\ref{conv:choice-of-T}, so that $T_*=T_*(M,N^0_\infty)=T_*(\mfN^{k_0})$, and choose
$T\in(0,T_*]$. Later in this section, $T$ will be chosen smaller to make the source $\mfF^a$ sufficiently small, as the
Nash--Moser theorem requires.

\cite[Assumptions~2.1--2.2]{secchi2016nash} require tame estimates for the second-order derivative and for a
right inverse of the first-order derivative of $\tilde\mfF$, respectively. By
\eqref{eq:derivatives-of-modified-iteration-map}, these follow from
Theorem~\ref{thm:second-order-tame-estimate-iteration-map} and Theorem~\ref{thm:Lagrangian-tame-bound},
respectively, after replacing the weight $R$ by $r_0$ via Lemma~\ref{lem:comparability} and absorbing the norms of
$\vec U_{\mathrm{app}}$ into the constants.

The following is \cite[Assumption 2.1]{secchi2016nash} with the derivative losses $r=2$, $r'=3$.
\begin{thm}\label{thm:second-order-tame-estimate-modified-iteration-map}
Let $T\in(0,T_*]$ be as in Convention~\ref{conv:choice-of-T} with $M$ as chosen in \eqref{eq:choice-of-M}. For any
integer $m\in[m_0,\tilde\alpha]$, any $\vec V\in\calU\cap X_{m_*}$ and any
$\delta\vec V,\epsilon\vec V\in X_{m_*}$,
\begin{equation}\label{eq:second-order-tame-estimate-modified-iteration-map}
\begin{aligned}
    \norm{d^2\tilde\mfF(\vec V)[\delta\vec V,\epsilon\vec V]}_{Y_m}
    \le\const{\mfN^{k_0}}\Big(&\norm{\delta\vec V}_{X_{m+2}}\norm{\epsilon\vec V}_{X_{m_0}}
    +\norm{\delta\vec V}_{X_{m_0}}\norm{\epsilon\vec V}_{X_{m+2}}\\
    &+\norm{\delta\vec V}_{X_{m_0}}\norm{\epsilon\vec V}_{X_{m_0}}\big(1+\norm{\vec V}_{X_{m+3}}\big)\Big).
\end{aligned}
\end{equation}
\end{thm}
\begin{proof}
Let $\vec U:=\vec V+\vec U_{\mathrm{app}}$. From \eqref{eq:choice-of-M}, we have $\norm{\vec U}_{\tilde Z^{2}_y}\le M$ and $\calC^2_L=\const{\norm{\vec U}_{\tilde Z^2_y},N^1_\infty}\le\const{M,N^1_\infty}\le\const{\mfN^{k_0}}$.
Applying Theorem~\ref{thm:second-order-tame-estimate-iteration-map} with $l=m$ and background $\vec U$, and
replacing the weight $R$ by $r_0$ on both sides via Lemma~\ref{lem:comparability} (with constants again
$\le\const{\mfN^{k_0}}$), we obtain \eqref{eq:second-order-tame-estimate-modified-iteration-map} with the bounds $\norm{\cdot}_{\tilde X^{l_0+1,\frac1{2\kappa}}_{r_0}}\le\norm{\cdot}_{X_{m_0}}$, since $l_0+1\le m_0$, and
\[
    N^{m+1}_w+\norm{\vec U}_{\tilde X^{m+3,\frac1{2\kappa}}_{r_0}}
    \lesssim N^{k_0}_w+\norm{\vec U_{\mathrm{app}}}_{\tilde X^{m_*,\frac1{2\kappa}}_{r_0}}+\norm{\vec V}_{X_{m+3}}
    \le\const{\mfN^{k_0}}\big(1+\norm{\vec V}_{X_{m+3}}\big),
\]
since $m+3\le m_*$, $m+1\le k_0$ and \eqref{eq:approximate-solution-high-regularity-bound}.
\end{proof}

\begin{rmk}\label{rmk:mapping-property}
In \cite[Assumption~2.1]{secchi2016nash}, the map $F$ sends $X_m$ to $Y_m$, whereas $\tilde\mfF$ sends $\calU\cap X_m$ only to $Y_{m-3}$, since $R,A_0$ are two levels less regular than $\vec U$ (see Lemma~\ref{lem:coefficient-bounds-merged}(b)) and $\mfF$ contains their first derivatives, which costs a further level in the scale. This does not affect the proof of \cite[Theorem~2.4(i)]{secchi2016nash}, which uses only the estimates \cite[(2.1), (2.2)]{secchi2016nash} and Taylor expansion for $F$ and $F'$, and the latter holds regardless of the space in which $F$ takes its values.
\end{rmk}

The following is \cite[Assumption 2.2]{secchi2016nash} with the derivative losses $s=0$, $s'=4$.
\begin{thm}\label{thm:first-order-tame-estimate-modified-iteration-map}
Let $T\in(0,T_*]$ be as in Convention~\ref{conv:choice-of-T} with $M$ as chosen in \eqref{eq:choice-of-M}. For any integer $m\in[m_0,\tilde\alpha]$, $\vec V\in\calU\cap X_{m_*}$ and $\vec G\in Y_m$, the solution $\delta \vec V$ to the equation
\begin{equation*}
    (D\tilde \mfF)(\vec V)\delta \vec V = \vec G
\end{equation*}
satisfies the estimate 
\begin{equation}\label{eq:first-order-tame-estimate-modified-iteration-map}
    \norm{\delta\vec V}_{X_m}\le\const{\mfN^{k_0}}\Big(\norm{\vec G}_{Y_m}
    +\norm{\vec G}_{Y_{m_0}}\norm{\vec V}_{X_{m+4}}\Big).
\end{equation}
\end{thm}
\begin{proof}
As in the preceding proof, apply Theorem~\ref{thm:Lagrangian-tame-bound} with $l=m$ to $\vec U=\vec V+\vec U_{\mathrm{app}}$, use Lemma~\ref{lem:comparability}, and bound $\norm{\cdot}_{\tilde Y^{l_0,\frac{1}{2\kappa}}_R}\le\norm{\cdot}_{Y_{m_0}}$ and
$N^{m+2}_w+\norm{\vec U}_{\tilde X^{m+4,\frac{1}{2\kappa}}_{r_0}}\le\const{\mfN^{k_0}}\big(1+\norm{\vec V}_{X_{m+4}}\big)$.
\end{proof}

\begin{rmk}\label{rmk:finite-regularity}
In \cite{secchi2016nash}, the background states $u$ lie in $U\cap X_\infty$ and the iterates constructed are in $X_\infty$. Here, the background states $\vec V$ lie in  $\calU\cap X_{m_*}$ and the iterates constructed are  in $X_{\tilde\alpha}$. This suffices: inspection of the proof of \cite[Theorem~2.4(i)]{secchi2016nash} shows that the estimate (2.2) in \cite[Assumption~2.2]{secchi2016nash} is applied only at the background states with regularities up to $X_{\tilde\alpha+s'}$ in \cite[Lemma~3.8]{secchi2016nash}, and the estimate (2.1) in \cite[Assumption~2.1]{secchi2016nash} at the background states with regularities up to $X_{\tilde\alpha-\max(r,r')+r'}$ in \cite[Lemma~3.3,Lemma~3.4]{secchi2016nash}. In both cases, the background states are at most $X_{m_*}$. 
\end{rmk}

It remains to verify the smallness hypothesis of \cite[Theorem~2.4(i)]{secchi2016nash} on the source $f$, which in our case is
$\mfF^a=-\mfF(\vec U_{\mathrm{app}})$. This is where we potentially shrink the choice of $T$: since, as we will show, $\mfF^a$ vanishes at $t=0$, its size is controlled by $T$.

\begin{lem}\label{lem:source-control}
Let $T\in(0,T_*]$ be as in Convention~\ref{conv:choice-of-T} with $M$ as chosen in \eqref{eq:choice-of-M}. Then
$\mfF^a\in W^{1,\infty}\big([0,T];H^{2(\tilde\alpha+1),\tilde\alpha+1+\frac1{2\kappa}}_{r_0}\big)$, $\mfF^a(0,\cdot)=0$, and
\begin{equation}\label{eq:source-control-bound}
    \norm{\mfF^a}_{Y_{\tilde\alpha+1}}\le T\,\const{\mfN^{k_0}}.
\end{equation}
\end{lem}

\begin{proof}
Since $\vec U_{\mathrm{app}}(0)=\vec u_0$ and $\p_t\vec U_{\mathrm{app}}(0)=\p_t\vec U|_{t=0}$, while
$\p_tR_{\mathrm{app}}(0)$, $\p_tP_{\mathrm{app}}(0)$ are given by \eqref{eq:R-A0-defining}, \eqref{eq:R-P-defining}
at $t=0$, the expression $\mfF(\vec U_{\mathrm{app}})(0,\cdot)$ is \eqref{eq:iteration-map-equation} at $t=0$
with $\p_t\vec U=\p_t\vec U|_{t=0}$, which vanishes by the definition of $\p_t\vec U|_{t=0}$ in
Section~\ref{subsec:approximate-solution-modified-iteration-map}. Hence $\mfF^a(0)=0$.

Estimating term by term of \eqref{eq:schematic-iteration-map}, we can show $\mfF^a\in L^\infty([0,T],H^{2(\tilde\alpha+1),\tilde\alpha+1+\frac{1}{2\kappa}}_{r_0})$. For the time derivative, write $\vec U$ for $\vec U_{\mathrm{app}}$, $(R,A_0,\eta)$ for
$(R_{\mathrm{app}},A_{0,\mathrm{app}},\eta_{\mathrm{app}})$, $F:=\nabla_y\eta$, and $\calF=\calF(R,A_0,\vec U)$.

From \eqref{eq:R-A0-defining} and \eqref{eq:transport-deformation-inverse}, we have
\begin{equation*}
    \begin{aligned}
        \p_tR&=R\calF\,\p_t\vec U+R\calF\,F^{-1}\nabla_y\vec U,\\
        \p_tA_0&=\calF+\calF\,\p_t\vec U+\calF\,F^{-1}\nabla_y\vec U,\\
        \p_tF^{-1}&=F^{-1}\calF\,\nabla_y\vec U\,F^{-1}.
    \end{aligned}
\end{equation*}

Differentiating \eqref{eq:schematic-iteration-map} pointwise in $t$, using $\p_t^2\vec U=0$ and substituting these
identities, we find that $\p_t\mfF^a$ is a finite sum of terms of the four types:
\[
    \calF,\qquad \calF\,q,\qquad \calF\,q\,q',\qquad \calF\,q_2,
\]
where 
\begin{equation*}
    \begin{aligned}
        q,q'&\in\{\p_t\vec U,\ F^{-1}\nabla_y\vec U,\ F^{-1}\nabla_yR,\ F^{-1}\nabla_yA_0\},\\
        q_2&\in\{F^{-1}\nabla_y\p_t\vec U,\ F^{-1}(\nabla_yF^{-1})\nabla_y\vec U,\ F^{-1}F^{-1}\nabla_y^2\vec U\}.
    \end{aligned}
\end{equation*}

Let $m:=\tilde\alpha+1$. Each term is estimated in $H^{2m,m+\frac1{2\kappa}}_{r_0}$ by
Corollary~\ref{cor:weighted-L2-end-point-equivalence}, Lemma~\ref{lem:coefficient-product-estimate}(b) and Corollary~\ref{Cor:Multilinear-Interpolation}, with the functions $f_j$ taken among
\[
    R,\ A_0,\ \vec U,\ F^{-1},\ \p_t\vec U,\ \nabla_yR,\ \nabla_yA_0,\ \nabla_y\vec U,\ \nabla_yF^{-1},\ \nabla_y\p_t\vec U,\ \nabla_y^2\vec U
\]
and with total derivative count $2m$ and weight $r_0^{m+\frac{1}{2\kappa}}$. Their $L^\infty_y$ norms are bounded by $\calC(\mfN^{k_0})$ using Lemma~\ref{lem:WeightedEmbedding}, Lemma~\ref{lem:deformation-inverse-Wkinfty-estimate}, Lemma~\ref{lem:coefficient-bounds-merged}(a), and Lemma~\ref{lem:velocity-Linfty-from-Eulerian-to-Lagrangian}. Their weighted norms $\norm{r_0^{m+\frac1{2\kappa}}\p^{2m}f_j}_{L^2}$ are bounded by $\calC(\mfN^{k_0})$ using Lemma~\ref{lem:WeightedEmbedding}, Lemma~\ref{lem:weighted-L2-deformation-inverse-estimate}, Theorem~\ref{thm:Lagrangian-Eulerian-tame-comparison}, Lemma~\ref{lem:coefficient-bounds-merged}(b), and Lemma~\ref{lem:velocity-weighted-L2-from-Eulerian-to-Lagrangian}. Hence, we have $\norm{\p_t\mfF^a}_{L^\infty_tH^{2m,m+\frac1{2\kappa}}_{r_0}}\le\const{\mfN^{k_0}}$. Integrating in time of $\p_t\mfF^a$ and a direct estimate gives \eqref{eq:source-control-bound}.
\end{proof}

\begin{thm}[Main theorem]\label{thm:main}
Let $m_0$, $k$, $k_0$ be as in \eqref{eq:NM-parameters}. For every $k_0$-admissible data $(r_0,\vec u_0,\pi_0)$
there exist $T=T(\mfN^{k_0})\in(0,1]$ and $\vec U\in\tilde X^{k,\frac1{2\kappa}}_{r_0}$ solving
\eqref{eq:iteration-map-equation} on $[0,T]\times\Omega_0$. Moreover, with $U^0$, $\eta$, $(R,A_0)$ defined from
$\vec U$ by \eqref{eq:U^0-defining-eq}, \eqref{eq:eta-defining-eq}, \eqref{eq:R-A0-defining} and $P:=RA_0$, $(R,\vec U,P,\eta)$ solves the Lagrangian system \eqref{eq:Lagrangian-formulation} on
$[0,T]\times\Omega_0$, with $\eta-\mathrm{Id}\in W^{1,\infty}_tH^{2k,k+\frac1{2\kappa}}_{r_0}$ and $R,A_0\in W^{1,\infty}_tH^{2(k-2),k-2+\frac1{2\kappa}-\frac12}_{r_0}$.
\end{thm}

\begin{proof}
Let $M$ be as in \eqref{eq:choice-of-M} and $T\in(0,T_*]$ as in Convention~\ref{conv:choice-of-T}. By Theorems~\ref{thm:first-order-tame-estimate-modified-iteration-map},~\ref{thm:second-order-tame-estimate-modified-iteration-map} and Remarks~\ref{rmk:finite-regularity},~\ref{rmk:mapping-property}, \cite[Theorem~2.4(i)]{secchi2016nash} applies to the scales \eqref{eq:NM-scales} and the map $\tilde\mfF$ with the parameters $r=2$, $r'=3$, $s=0$, $s'=4$, $\alpha=k+1$ and $\tilde\alpha$, $m_*=\tilde\alpha+\max(r',s')$ as in \eqref{eq:NM-parameters}. The constant $\theta_0\ge1$ is chosen and fixed in \cite[Lemma~3.1]{secchi2016nash} and is independent of all parameters and of $T$. The constant $\delta>0$ is fixed in \cite[Lemmas~3.3--3.8]{secchi2016nash}. Inspection of those proofs shows that $\delta$ depends only on the parameters $m_0,\alpha,\tilde\alpha,r,r',s,s'$ and on the constants in \cite[Assumptions~2.1--2.2]{secchi2016nash} and in the smoothing estimates \cite[(2.3)]{secchi2016nash}. Since the proof of \cite[Theorem~2.4(i)]{secchi2016nash} only involves indices $m\in[m_0,m_*]$, these constants are uniformly bounded in $m$. In our case, they are the constants of
\eqref{eq:first-order-tame-estimate-modified-iteration-map},
\eqref{eq:second-order-tame-estimate-modified-iteration-map} and Lemma~\ref{lem:spatial-smoothing-operator},
all bounded in terms of the parameters \eqref{eq:NM-parameters} and $\mfN^{k_0}$. Hence, $\delta=\delta(\mfN^{k_0})$, independent of $T$. By Lemma~\ref{lem:source-control} we may decrease $T$ to
$T=T(\mfN^{k_0})$ so that $\norm{\mfF^a}_{Y_{\alpha}}/\delta$ is as small as \cite[Section~3.6]{secchi2016nash}
requires. The Nash--Moser scheme \cite[(3.1)--(3.6)]{secchi2016nash} with $f=\mfF^a$ then yields
$\vec V_n\in\calU\cap X_{\tilde\alpha}$ with
\[
    \vec V_n\to\vec V\ \text{in }X_k,\qquad \tilde\mfF(\vec V_n)\to\mfF^a\ \text{in }Y_k\subset Y_{k-3} ,
\]
by \cite[Section~3.6]{secchi2016nash} with $m'=\alpha-1=k$. Since $\tilde\mfF:X_m\to Y_{m-3}$ is twice-differentiable, we have $\tilde\mfF(\vec V_n)\to\tilde\mfF(\vec V)$ in $Y_{k-3}$. Thus, we get $\tilde\mfF(\vec V)=\mfF^a$ in $Y_{k-3}$, i.e., if we set $\vec U:=\vec V+\vec U_{\mathrm{app}}$, $\mfF(\vec U)=0$ and $\vec U(0,\cdot)=\vec u_0$, which is \eqref{eq:iteration-map-equation}. Moreover, $\norm{\vec V}_{X_{m_0}}\leq\delta<1$ since $\norm{\vec V_n}_{X_{m_0}}\leq\delta$ by \cite[Lemma~3.1]{secchi2016nash} and $\vec V_n\to \vec V$ in $X_k\subset X_{m_0}$. Hence, $\vec V\in\calU$.

In \eqref{eq:Lagrangian-formulation}, (a) and (c) are \eqref{eq:R-P-defining}, (b) is $\mfF(\vec U)=0$, (d) and (e) are \eqref{eq:eta-defining-eq} and \eqref{eq:U^0-defining-eq}; the initial conditions hold by construction. Since $\norm{\vec U}_{\tilde Z^1_y}<M$, the regularity
of $R,A_0$ follows from Theorem~\ref{thm:Lagrangian-Eulerian-tame-comparison}, Lemma~\ref{lem:coefficient-bounds-merged}(b), Lemma~\ref{lem:velocity-weighted-L2-from-Eulerian-to-Lagrangian}, and Lemma~\ref{lem:comparability}, and  that of
$\eta$ from \eqref{eq:eta-defining-eq} and a direct estimate.
\end{proof}

\section*{Acknowledgements}
The author is grateful to his Ph.D. advisor Marcelo Disconzi for suggesting the problem and for many helpful discussions and comments on earlier versions of this manuscript. The author also thanks his postdoctoral mentor Samuel Walsh for his mentorship and encouragement.

\appendix

\section{Solving the modified linearized system for the time derivatives}\label{sec:solving-time-derivative}

First, we eliminate $\tilde u^0$ via \eqref{eq:orthogonality-constraint}.
\begin{itemize}
	\item Differentiating \eqref{eq:Normalization-Condition} yields
	\begin{equation}\label{eq:dt-u0}
		\p_t u^0 = \frac{u^i}{u^0}\p_t u_i.
	\end{equation}
	\item Differentiating \eqref{eq:orthogonality-constraint} and applying \eqref{eq:dt-u0} yields
	\begin{equation}\label{eq:dt-tilde-u0}
    \p_t\tilde u^0=\frac{u_i}{u^0}\p_t\tilde u^i+\frac1{u^0}H^{ij}(\p_tu_j)\tilde u_i,
    \qquad H^{ij}:=\delta^{ij}-\frac{u^iu^j}{(u^0)^2}.
\end{equation}
\end{itemize}
 
Substituting \eqref{eq:dt-tilde-u0} and $\tilde u^0 = \frac{u^i}{u^0}\tilde u_i$ from \eqref{eq:orthogonality-constraint} into \eqref{eq:tilde-r-equation}, the spatial part of \eqref{eq:tilde-u-Equation} and \eqref{eq:tilde-pi-equation}, we obtain
\begin{align}
	u^0\p_t\tilde r+\frac{\kappa ra_1}{u^0}u_i\p_t\tilde u^i&=\mfR_{\tilde r}\label{eq:time-derivative-tilde-r},\\
	\left(1+\frac1\kappa\right)a_2\Delta^{i0}\p_t\tilde r+u^0\p_t\tilde u^i+\left(2+\frac1\kappa\right)ra_3\Delta^{i0}\p_t\tilde\pi&=\mfR_{\vec{\tilde u}}^i,\label{eq:time-derivative-tilde-u}\\
	\frac{ra_4}{u^0}u_i\p_t\tilde u^i+u^0\p_t\tilde\pi&=\mfR_{\tilde \pi}\label{eq:time-derivative-tilde-pi},
\end{align}
where
\begin{align}
	\mfR_{\tilde r}&:=\tilde f-u^i\p_i\tilde r-\kappa ra_1\p_i\tilde u^i
    -\left(\frac{\kappa ra_1}{u^0}H^{ij}\p_tu_j+\frac{u_i}{u^0}\p_tr+ \delta^{ij} \p_j r\right)\tilde u_i
    -V_1\tilde r-r^2Z_1\tilde\pi,\label{eq:time-derivative-tilde-r-remainder}\\
    \mfR_{\vec{\tilde u}}^i &:= \tilde g^i-\left(1+\frac1\kappa\right)a_2\Delta^{ij}\p_j\tilde r-u^j\p_j\tilde u^i-\Big(2+\frac1\kappa\Big)ra_3\Delta^{ij}\p_j\tilde\pi\notag\\
    &\qquad-V_2^i\tilde r-\Big((W_2)^i_j+(W_2)^i_0\frac{u_j}{u^0}\Big)\tilde u^j-Z_2^i\tilde\pi, \label{eq:time-derivative-tilde-u-remainder}\\
    \mfR_{\tilde\pi}&:=\tilde h-ra_4\p_i\tilde u^i-u^i\p_i\tilde\pi
    -\Big(\frac{ra_4}{u^0}H^{ij}\p_tu_j+W_3^i-W_3^0\frac{u_i}{u^0}\Big)\tilde u_i
    -V_3\tilde r-Z_3\tilde\pi.
\end{align}
 
From \eqref{eq:time-derivative-tilde-u}, we can obtain
\begin{equation}\label{eq:dt-u-i-subs}
    \p_t\tilde u^i=\frac1{u^0}\mfR_{\vec{\tilde u}}^i-\frac1{u^0}\Big(1+\frac1\kappa\Big)a_2\Delta^{i0}\p_t\tilde r-\frac1{u^0}\Big(2+\frac1\kappa\Big)ra_3\Delta^{i0}\p_t\tilde\pi.
\end{equation}
Substituting \eqref{eq:dt-u-i-subs} into \eqref{eq:time-derivative-tilde-r} and \eqref{eq:time-derivative-tilde-pi}, we obtain the following $2\times 2$ system
\begin{equation}\label{eq:2-by-2-system}
	\begin{pmatrix}
		\mfa_{11} & \mfa_{12}\\
		\mfa_{21} &\mfa_{22}		
	\end{pmatrix}
	\begin{pmatrix}
		\p_t\tilde r\\
		\p_t\tilde\pi
	\end{pmatrix}=
	\begin{pmatrix}
		\mfR_{\tilde r}- \frac{\kappa r a_1 u_i}{(u^0)^2}\mfR^{i}_{\vec{\tilde u}}\\
		\mfR_{\tilde \pi} - \frac{ra_4 u_i}{(u^0)^2}\mfR_{\vec{\tilde u}}^i
	\end{pmatrix},
\end{equation}
where
\begin{align*}
	\mfa_{11}&=u^0-\frac{\kappa ra_1}{(u^0)^2}\left(1+\frac1\kappa\right)a_2u_i\Delta^{i0},\\
	\mfa_{12}&=-\frac{\kappa ra_1}{(u^0)^2}\left(2+\frac1\kappa\right)ra_3u_i\Delta^{i0},\\
	\mfa_{21}&=-\frac{ra_4}{(u^0)^2}\left(1+\frac1\kappa\right)a_2u_i\Delta^{i0},\\
	\mfa_{22}&=u^0-\frac{ra_4}{(u^0)^2}\left(2+\frac1\kappa\right)ra_3u_i\Delta^{i0}.
\end{align*}

Let 
\begin{equation}\label{eq:coefficient-matrix-of-time-derivative}
    \mfA:=\begin{pmatrix}\mfa_{11} & \mfa_{12}\\ \mfa_{21} &\mfa_{22}\end{pmatrix}.
\end{equation}
Its determinant 
\begin{equation*}
	D=\det\mfA=\mfa_{11}\mfa_{22}-\mfa_{12}\mfa_{21}\geq c>0
\end{equation*}
is bounded away from $0$, since $r\le 2\varepsilon_0$ on $[0,T]$ by Remark~\ref{rmk:persistence-of-causality}.

By Cramer's rule, we obtain
\begin{equation}\label{eq:final}
    \begin{aligned}
    	\p_t\tilde r&=\frac{\mfa_{22}(\mfR_{\tilde r}- \frac{\kappa r a_1 u_i}{(u^0)^2}\mfR^{i}_{\vec{\tilde u}})-\mfa_{12}(\mfR_{\tilde \pi} - \frac{ra_4 u_i}{(u^0)^2}\mfR_{\vec{\tilde u}}^i)}{D},\\
    \p_t\tilde\pi&=\frac{\mfa_{11}(\mfR_{\tilde \pi} - \frac{ra_4 u_i}{(u^0)^2}\mfR_{\vec{\tilde u}}^i)\,-\mfa_{21}(\mfR_{\tilde r}- \frac{\kappa r a_1 u_i}{(u^0)^2}\mfR^{i}_{\vec{\tilde u}})}{D}.
    \end{aligned}
\end{equation}
Substituting \eqref{eq:final} back into \eqref{eq:dt-u-i-subs} gives $\p_t\tilde u^i$; together with
\eqref{eq:dt-tilde-u0}, this gives $\p_t\tilde u^0$.

Substituting \eqref{eq:r-pi-defining} to write $\p_t(r,\pi)$ in terms of spatial derivatives of $(r,\pi)$ and $(\p\vec u,\p_t\vec u)$, and replacing $\pi$ by $\pi=ra_0$, we have schematically
\begin{equation}
    \left\{
    \begin{aligned}
        \p_t \tilde r &= \calF\p\tilde r + r\calF\p\vec{\tilde u} + r\calF\p\tilde\pi
            + \mfL\tilde r + \mfL\tilde u + \mfL\tilde\pi
            + \calF\tilde f + r\calF\vec{\tilde g} + r\calF\tilde h,\\
        \p_t \vec{\tilde u} &= \calF\p\tilde r + \calF\p\vec{\tilde u} + \calF\p\tilde\pi
            + \mfL\tilde r + \mfL\tilde u + \mfL\tilde\pi
            + \calF\tilde f + \calF\vec{\tilde g} + \calF\tilde h,\\
        \p_t \tilde \pi &= \calF\p\tilde r + \calF\p\vec{\tilde u} + \calF\p\tilde\pi
            + \mfL\tilde r + \mfL\tilde u + \mfL\tilde\pi
            + \calF\tilde f + \calF\vec{\tilde g} + \calF\tilde h,
    \end{aligned}
    \right.
\end{equation}
where $\calF=\calF(r,a_0,\vec u)$ denotes a coefficient smooth in a neighborhood of the range of $(r,a_0,\vec u)$,
and $\mfL=\sum_X\calF_X(r,a_0,\vec u)\p X+\calF(r,a_0,\vec u)\p_t\vec u+\calF(r,a_0,\vec u)$ for $X\in\{r,a_0,\vec u\}$. Above, both
$\calF$ and $\mfL$ are subject to change line by line and term by term.

\begin{rmk}
The term $\calF(r,a_0,\vec u)$ in $\mfL$ comes from the zeroth-order terms of \eqref{eq:Eulerian-Main-System} (for instance, the term $\pi$ in \eqref{eq:pi-equation}), and the term $\calF(r,a_0,\vec u)\p_t\vec u$ from replacing $\p_tr$ and $\p_ta_0$ using \eqref{eq:r-pi-defining} and \eqref{eq:a_0-evolution}. 
\end{rmk}

\section{Smoothing operators in weighted Sobolev spaces}\label{app:smoothing}

Here, we construct a family of smoothing operators on the scale $H^{2l,l+\frac1{2\kappa}}_w(\Omega)$, where $w$ is a non-degenerate defining function for $\Omega$, satisfying estimates similar to \cite[(2.3)]{secchi2016nash}. Standard mollification does not preserve the weighted spaces near $\partial\Omega$, where the weight degenerates. Thus, we employ the good regularization kernels in \cite[Section~2.4]{ifrim2023compressible}, whose scale adapts to the distance to the boundary.

\begin{lem}[Smoothing operators in weighted Sobolev spaces]\label{lem:spatial-smoothing-operator-appendix}
Let $\Omega\subset\R^3$ be a bounded domain with smooth boundary, $w$ a non-degenerate defining function for $\Omega$, and $m_*\ge1$ a fixed integer. There exists a family of linear operators
\begin{equation*}
    S_\theta:H^{0,\frac1{2\kappa}}_w\to\bigcap_{l\ge0}H^{2l,l+\frac1{2\kappa}}_w,\qquad \theta\ge1
\end{equation*}
depending smoothly on $\theta$, such that for all integers $0\le\alpha,\beta\le m_*$ and $u\in H^{2\alpha,\alpha+\frac1{2\kappa}}_w$,
\begin{align}
    \norm{S_\theta u}_{H^{2\beta,\beta+\frac1{2\kappa}}_w}&\le C\theta^{(\beta-\alpha)_+}\norm{u}_{H^{2\alpha,\alpha+\frac1{2\kappa}}_w},\label{eq:S1}\\
    \norm{S_\theta u-u}_{H^{2\beta,\beta+\frac1{2\kappa}}_w}&\le C\theta^{\beta-\alpha}\norm{u}_{H^{2\alpha,\alpha+\frac1{2\kappa}}_w},\qquad\beta\le\alpha,\label{eq:S2}\\
    \norm{\tfrac{d}{d\theta}S_\theta u}_{H^{2\beta,\beta+\frac1{2\kappa}}_w}&\le C\theta^{\beta-\alpha-1}\norm{u}_{H^{2\alpha,\alpha+\frac1{2\kappa}}_w},\label{eq:S3}
\end{align}
with $C=C(m_*,\kappa,\Omega,w)$ independent of $\theta,\alpha,\beta,u$.
\end{lem}
\begin{proof}
    First, we notice that the spaces in the statement of the theorem form a decreasing scale. For integers $\alpha\geq\beta$, losing derivatives and applying Lemma~\ref{lem:WeightedEmbedding} gives the following embedding of spaces:
    \begin{equation}\label{eq:Decreasing-weighted-scale}
        H^{2\alpha,\alpha+\frac{1}{2\kappa}}_w \subset H^{\alpha+\beta,\alpha+\frac{1}{2\kappa}}_w \subset H^{2\beta,\beta+\frac{1}{2\kappa}}_w.
    \end{equation}

    Choose a sufficiently large integer $N_*$ such that $N_*\geq m_*$. By \cite[Definition~2.8, Lemma~2.9]{ifrim2023compressible} there exist good regularization kernels $K^h$, $h\in\N\cup\{0\}$, adapted to $w$, with the integer $N=N_*$ in (ii), (iii) of \cite[Definition~2.8]{ifrim2023compressible}. They are defined on an enlargement $\tilde\Omega^{[h]}$ of $\Omega$, where
    \begin{equation*}
        \tilde\Omega^{[h]}:=\{x\in\R^3: \dist(x,\Omega)\leq c2^{-2h}\}
    \end{equation*}
    for some small universal constant $c$. In our case, we restrict them to $x\in\Omega$. We use the following properties from \cite[Definition~2.8]{ifrim2023compressible}: 
    \begin{itemize}
        \item [(i)]The support condition:
        \begin{equation*}
            K^h:\tilde\Omega^{[h]}\times\Omega\to \R
        \end{equation*}
        with support properties
        \begin{equation*}
            \supp K^h\subset \{(x,y)\in \tilde\Omega^{[h]}\times\Omega^{[<h]}:|x-y|\lesssim \delta^h(y):= 2^{-2h}+2^{-h}w(y)^{\tfrac12}\},
        \end{equation*}
        where
        \begin{equation*}
            \Omega^{[<h]}:=\{x\in\Omega:w(x)\gtrsim 2^{-2h}\}.
        \end{equation*}
        \item [(ii)] The derivative bounds:
        \begin{equation*}
            \abs{\p^\mu_x\p^\nu_yK^h(x,y)}\le C\big(2^{-2h}+2^{-h}w(y)^{1/2}\big)^{-N_*-\abs\mu-\abs\nu},\qquad \abs\mu+\abs\nu\le4N_*.
        \end{equation*}
        \item [(iii)]The approximate identity and moment conditions:
        \begin{equation*}
            \begin{aligned}
                \int_\Omega K^h(x,y)\,dy&=1,\\
                \int_\Omega K^h(x,y)(x-y)^\gamma\,dy&=0, \quad 1\le\abs\gamma\le2N_*.
            \end{aligned}
        \end{equation*}
    \end{itemize}
    Define
    \begin{equation*}
        (\Psi_hu)(x):=\int_\Omega K^h(x,y)\,u(y)\,dy,\qquad x\in\Omega.
    \end{equation*}
    For $u\in H^{0,\frac1{2\kappa}}_w$, $\Psi_hu$ is well defined and belongs to $H^{2l,l+\frac1{2\kappa}}_w$ for every $l\ge0$. Indeed, by the support condition, $u$ is integrated only over $\Omega^{[<h]}$, where $w\ge c2^{-2h}$, so $u\in L^1(\Omega^{[<h]})$. The good regularization kernels $K^h$ constructed in \cite[Lemma~2.9]{ifrim2023compressible} are $C^\infty$ with bounded derivatives of every order for fixed $h$, so $\Psi_hu\in C^\infty(\bar\Omega)$ by differentiation under the integral sign, and all its weighted norms are finite since $\Omega$ is bounded. Thus, we have
    \begin{equation*}
        \Psi_hu \in \bigcap_{l\geq0}H^{2l,l+\frac{1}{2\kappa}}_w.
    \end{equation*}

    For integer indices $0\leq\alpha,\beta\leq m_*$, the proof of Proposition~2.10 of \cite{ifrim2023compressible} uses the kernels in two ways. It subtracts from the input a polynomial of degree at most $2\alpha-1$ (none when $\alpha=0$), which uses the moment conditions up to order $2\alpha-1$. It differentiates the output at most $2\beta$ times, which uses the derivative bounds up to order $2\beta$. Since we have
    \begin{align*}
        2\alpha-1&\leq 2m_*-1<2N_*,\\
        2\beta&\leq 2m_*\leq 2N_*,
    \end{align*}
    the stated moment and derivative bounds suffice simultaneously for every such pair of integer indices with our choice of $N_*$.

    We apply \cite[Proposition~2.10]{ifrim2023compressible} to $(0,u)$ with both defining functions in the proposition equal to $w$. For $\beta\geq\alpha$, the regularization bound in \cite[Proposition~2.10]{ifrim2023compressible} with $k=\alpha$ and $j=\beta-\alpha$ reads
    \begin{equation*}
        \norm{\Psi_h u}_{H^{2\beta,\beta+\frac{1}{2\kappa}}_w}\leq C2^{2(\beta-\alpha)h}\norm{u}_{H^{2\alpha,\alpha+\frac{1}{2\kappa}}_w}.
    \end{equation*}
    For $\beta\leq\alpha$, the same estimate with $j=0$ and the embedding \eqref{eq:Decreasing-weighted-scale} gives
    \begin{equation}\label{eq:dyadic-regularization-bound}
        \norm{\Psi_hu}_{H^{2\beta,\beta+\frac{1}{2\kappa}}_w} \leq C2^{2(\beta-\alpha)_+h}\norm{u}_{H^{2\alpha,\alpha+\frac{1}{2\kappa}}_w},
    \end{equation}
    where $(\beta-\alpha)_+=\max{(\beta-\alpha,0)}$.

    For $0\leq\beta\leq\alpha\leq m_*$, the difference and error bounds in \cite[Proposition~2.10]{ifrim2023compressible}, with $k=\alpha$ and $j=\beta-\alpha$, give
    \begin{align}
        \label{eq:dyadic-difference-bound-low}
        \norm{(\Psi_{h+1}-\Psi_h)u}_{H^{2\beta,\beta+\frac{1}{2\kappa}}_w} &\leq C2^{2(\beta-\alpha)h} \norm{u}_{H^{2\alpha,\alpha+\frac{1}{2\kappa}}_w}, \\
        \label{eq:dyadic-error-bound}
        \norm{(I-\Psi_h)u}_{H^{2\beta,\beta+\frac{1}{2\kappa}}_w} &\leq C2^{2(\beta-\alpha)h} \norm{u}_{H^{2\alpha,\alpha+\frac{1}{2\kappa}}_w}.
    \end{align}

    For $\beta>\alpha$, \eqref{eq:dyadic-regularization-bound} implies
    \begin{align*}
        \norm{(\Psi_{h+1}-\Psi_h)u}_{H^{2\beta,\beta+\frac{1}{2\kappa}}_w} &\leq \norm{\Psi_{h+1}u}_{H^{2\beta,\beta+\frac{1}{2\kappa}}_w} + \norm{\Psi_hu}_{H^{2\beta,\beta+\frac{1}{2\kappa}}_w} \\
        &\leq C\left(2^{2(\beta-\alpha)(h+1)}+2^{2(\beta-\alpha)h}\right) \norm{u}_{H^{2\alpha,\alpha+\frac{1}{2\kappa}}_w} \\
        &\leq C2^{2(\beta-\alpha)h} \norm{u}_{H^{2\alpha,\alpha+\frac{1}{2\kappa}}_w},
    \end{align*}
where the constant in the last inequality is uniform because $0<\beta-\alpha\leq m_*$. Together with \eqref{eq:dyadic-difference-bound-low}, this proves
\begin{equation}\label{eq:dyadic-difference-bound}
    \norm{(\Psi_{h+1}-\Psi_h)u}_{H^{2\beta,\beta+\frac{1}{2\kappa}}_w} \leq C2^{2(\beta-\alpha)h} \norm{u}_{H^{2\alpha,\alpha+\frac{1}{2\kappa}}_w}
\end{equation}
for all $0\leq\alpha,\beta\leq m_*$.

We next define the continuous parameter family. Let
\begin{align*}
    q(s) &:= \begin{cases}
        0, & s\leq0, \\
        \exp\left(-\frac{1}{s}\right), & s>0,
    \end{cases} \\
    b(s) &:= \frac{q\left(s-\frac{1}{3}\right)}{q\left(s-\frac{1}{3}\right)+q\left(\frac{2}{3}-s\right)}.
\end{align*}
Then $b\in C^\infty(\R)$, $0\leq b\leq1$, and
\begin{align*}
    b(s)&=0, \qquad s\leq\frac{1}{3}, \\
    b(s)&=1, \qquad s\geq\frac{2}{3}.
\end{align*}
For $h\in\N\cup\{0\}$, set $\theta_h:=2^{2h}$. If $\theta\in[\theta_h,\theta_{h+1}]$, set
\begin{subequations}{\label{eq:continuous-smoothing-operator}}
\begin{align}
    \tau(\theta) &:= \frac{\theta-\theta_h}{\theta_{h+1}-\theta_h}, \\
    S_\theta &:= \left(1-b(\tau(\theta))\right)\Psi_h +b(\tau(\theta))\Psi_{h+1}.
\end{align}
\end{subequations}
Because $b$ is constant on neighborhoods of $0$ and $1$, the two formulas for $S_\theta$ on the intervals meeting at $\theta_h$ both equal $\Psi_h$ on a neighborhood of $\theta_h$. Their derivatives of every positive order with respect to $\theta$ vanish there. Thus, the definitions in \eqref{eq:continuous-smoothing-operator} give a $C^\infty$ family for $\theta\geq1$. Since $S_\theta u$ is a linear combination of $\Psi_hu$ and $\Psi_{h+1}u$, it belongs to $\bigcap_{l\ge0}H^{2l,l+\frac{1}{2\kappa}}_w$, which is the range asserted in the lemma. Moreover, we have
\begin{subequations}{\label{eq:continuous-dyadic-comparison}}
\begin{align}
    \theta_h &\leq\theta\leq4\theta_h, \\
    \theta_{h+1}-\theta_h &=3\theta_h.
\end{align}
\end{subequations}

It follows from \eqref{eq:continuous-smoothing-operator}, \eqref{eq:dyadic-regularization-bound}, and \eqref{eq:continuous-dyadic-comparison} that
\begin{align*}
    \norm{S_\theta u}_{H^{2\beta,\beta+\frac{1}{2\kappa}}_w}
    \leq C\theta_h^{(\beta-\alpha)_+}\norm{u}_{H^{2\alpha,\alpha+\frac{1}{2\kappa}}_w}
    \leq C\theta^{(\beta-\alpha)_+}\norm{u}_{H^{2\alpha,\alpha+\frac{1}{2\kappa}}_w},
\end{align*}
which is \eqref{eq:S1}. If $0\leq\beta\leq\alpha$, then \eqref{eq:continuous-smoothing-operator}, \eqref{eq:dyadic-error-bound}, and \eqref{eq:continuous-dyadic-comparison} give
\begin{align*}
    \norm{S_\theta u-u}_{H^{2\beta,\beta+\frac{1}{2\kappa}}_w}
    \leq C\theta_h^{\beta-\alpha}\norm{u}_{H^{2\alpha,\alpha+\frac{1}{2\kappa}}_w}
    \leq C\theta^{\beta-\alpha}\norm{u}_{H^{2\alpha,\alpha+\frac{1}{2\kappa}}_w},
\end{align*}
which is \eqref{eq:S2}. Differentiating \eqref{eq:continuous-smoothing-operator} gives
\begin{align}
\label{eq:continuous-smoothing-operator-derivative}
    \frac{d}{d\theta}S_\theta
    =
    \frac{b^\prime(\tau(\theta))}{3\theta_h}
    (\Psi_{h+1}-\Psi_h).
\end{align}
It follows from \eqref{eq:dyadic-difference-bound}, \eqref{eq:continuous-dyadic-comparison}, and \eqref{eq:continuous-smoothing-operator-derivative} that
\begin{align*}
    \norm{\frac{d}{d\theta}S_\theta u}_{H^{2\beta,\beta+\frac{1}{2\kappa}}_w}
    \leq C\theta_h^{\beta-\alpha-1}\norm{u}_{H^{2\alpha,\alpha+\frac{1}{2\kappa}}_w}
    \leq C\theta^{\beta-\alpha-1}\norm{u}_{H^{2\alpha,\alpha+\frac{1}{2\kappa}}_w},
\end{align*}
which is \eqref{eq:S3}.
\end{proof}

\newpage

\bibliographystyle{abbrv}
\bibliography{reference}
\end{document}